\documentclass[aos,preprint]{imsart}
\input{glyphtounicode}
\usepackage{amsfonts,amsmath,amsthm,mathtools,bbm,eucal}
\usepackage{amssymb,mathabx}
\usepackage[notrig]{physics}
\usepackage{graphicx,hyperref,enumitem,multirow}
\usepackage{placeins}
\graphicspath{{figures/}}
\usepackage[ruled,vlined,linesnumbered]{algorithm2e}
\usepackage[nameinlink,capitalize,noabbrev,nosort]{cleveref}
\usepackage{crossreftools}
\usepackage[authoryear]{natbib}
\usepackage{bibunits}
\defaultbibliographystyle{plainnat}
\defaultbibliography{references}
\usepackage[title,toc]{appendix}
\usepackage[dvipsnames,x11names]{xcolor}
\usepackage{booktabs,makecell}
\allowdisplaybreaks
\hypersetup{colorlinks=True, urlcolor=Blue, citecolor=MidnightBlue, linkcolor=RedViolet}
\makeatletter
\def\nocontentsline{%
    \let\@@addcontentsline\addcontentsline
    \def\addcontentsline##1##2##3{\let\addcontentsline\@@addcontentsline}%
}
\makeatother
\theoremstyle{plain}
\newtheorem{theorem}{Theorem}[section]
\makeatletter\@ifundefined{iflatexml}{\let\iflatexml\iffalse}{}\makeatother
\iflatexml
	\newcommand{\sharedtheorem}[2]{\newtheorem{#1}[theorem]{#2}}
\else
	\usepackage{aliascnt}
	\newcommand{\sharedtheorem}[2]{\newaliascnt{#1}{theorem}\newtheorem{#1}[#1]{#2}\aliascntresetthe{#1}}
\fi
\sharedtheorem{proposition}{Proposition}
\sharedtheorem{lemma}{Lemma}
\sharedtheorem{corollary}{Corollary}
\theoremstyle{definition}
\sharedtheorem{assumption}{Assumption}
\sharedtheorem{example}{Example}
\sharedtheorem{remark}{Remark}
\newtheorem{conditioninner}{Condition}
\newenvironment{condition}[1]{\renewcommand\theconditioninner{#1}\conditioninner}{\endconditioninner}

\crefname{assumption}{Assumption}{Assumptions}
\crefname{conditioninner}{Condition}{Conditions}
\numberwithin{equation}{section}

\newcommand{\ind}{\mathbbm 1}
\newcommand{\cE}{\mathcal E}
\newcommand{\cP}{\mathcal P}
\newcommand{\cS}{\mathcal S}
\newcommand{\cX}{\mathcal X}
\newcommand{\cY}{\mathcal Y}
\newcommand{\E}{\mathbb E}
\newcommand{\bP}{\mathbb P}
\newcommand{\R}{\mathbb R}
\renewcommand{\d}{\mathrm d}
\newcommand{\e}{\mathrm e}
\newcommand{\TV}{\operatorname{TV}}
\newcommand{\Leb}{\lambda}
\newcommand{\eps}{\epsilon}
\newcommand{\cov}{\hat\gamma}
\newcommand{\Cstar}{\mathsf C^{\star}}
\newcommand{\hp}{\widehat{\mathsf p}}
\newcommand{\hC}{\widehat{\mathsf C}}
\newcommand{\hP}{\widehat{\mathsf P}}

\newcommand{\heta}{\widehat{\eta}}

\newcommand{\genericdel}[4]{%
  \ifcase#3\relax
  \ifx#1.\else#1\fi#4\ifx#2.\else#2\fi\or
  \bigl#1#4\bigr#2\or
  \Bigl#1#4\Bigr#2\or
  \biggl#1#4\biggr#2\or
  \Biggl#1#4\Biggr#2\else
  \left#1#4\right#2\fi
}
\newcommand{\del}[2][-1]{\genericdel(){#1}{#2}}
\newcommand{\cbr}[2][-1]{\genericdel\{\}{#1}{#2}}
\newcommand{\sbr}[2][-1]{\genericdel[]{#1}{#2}}
 \let\abs\Abs
 \let\norm\Norm

\begin{document}

\begin{frontmatter}
\title{Coverage and shape of conformal highest-density sets: a minimax theory}
\runtitle{Conformal highest-density sets}

\begin{aug}
\author[A]{\fnms{Ilsang} \snm{Ohn}\ead[label=e1]{ilsang.ohn@inha.ac.kr}}
 \address[A]{Department of Statistics, Inha University, \printead[presep={\ }]{e1}}
\end{aug}
\runauthor{I. Ohn}

\begin{abstract}
We study conditional coverage and geometric efficiency of conformal highest-density sets. The excess volume of such a set over the conditional highest-density set, the smallest set with the required conditional coverage, splits into a coverage error and an excess mass. We show that, among procedures whose volume is accurate up to a constant, the coverage error has the minimax rate of nonparametric regression on the covariate, and that the excess mass has the squared minimax rate of conditional density estimation. The usual split conformal construction calibrates these sets with a single threshold, which corrects the bias of the density estimate only on average over the covariate. We prove that its coverage error can then be of the order of the slower density rate. We therefore propose two-stage calibration, which estimates the coverage of the super-level sets of a density estimate by a regression on the covariate and uses this estimated coverage as the conformity score. Since these sets are fixed at each covariate value, their coverage is as smooth as the conditional density, and, under mild regularity conditions, the coverage error attains the regression rate even when the density estimate is inconsistent. The shape of the set still depends on the accuracy of the density estimate: sup-norm accuracy controls it at the minimax excess-mass rate, whereas accuracy in an integrated norm, even of minimax order, does not. With a local polynomial density estimate, two-stage calibration attains both minimax rates up to logarithmic factors. Simulations and an application to bike-sharing demand illustrate the results.
\end{abstract}



\end{frontmatter}
\begin{bibunit}
\addtocontents{toc}{\protect\setcounter{tocdepth}{-1}}

\section{Introduction}\label{sec:intro}
A predictive set for a multivariate response should contain the response with a prescribed probability given the covariate, and it should be small. Given a covariate $x$, the smallest set with conditional coverage $1-\alpha$ is the conditional highest-density set, called the highest predictive density (HPD) set by \citet{izbicki2022cd},
	\begin{align}
		\Cstar(x)=\{y:\,{\mathsf p}(y\mid x)\ge t_\star(x)\},
	\end{align}
where ${\mathsf p}(\cdot\mid x)$ is the conditional density and the level $t_\star(x)$ gives $\Cstar(x)$ coverage $1-\alpha$. HPD sets can be disconnected and adapt their scale to the covariate. They are therefore the natural target of conformal procedures built on conditional density estimators and generative models \citep{izbicki2020flexible,izbicki2022cd,sampson2024flexible,plassier2024conditionally}.

Split conformal calibration of the super-level sets of a density estimate guarantees finite-sample marginal coverage \citep{vovk2005algorithmic,lei2018distribution}. Exact conditional coverage, in contrast, is impossible without distributional assumptions \citep{vovk2012conditional,foygelbarber2021limits}. Asymptotic consistency and conditional validity of the resulting sets have been established \citep{izbicki2020flexible,izbicki2022cd}, but a minimax theory of their conditional coverage and geometric efficiency is missing. In particular, it is not known which property of the density estimate governs the conditional coverage, which governs the geometry of the set, and what the minimax rates are. This paper answers these questions.

\paragraph{Volume decomposition} To separate coverage from geometry, we show that the excess volume of any predictive set $\hC(x)$ splits into two terms (\cref{lem:identity}). Let ${\mathsf P}_x$ be the conditional distribution of $Y$ given $X=x$, and let $\cov(x):={\mathsf P}_x(\hC(x))$ be the conditional coverage of $\hC(x)$. Then
	\begin{align}
		\Leb\bigl(\hC(x)\bigr)-\Leb\bigl(\Cstar(x)\bigr)=\frac{\cov(x)-(1-\alpha)}{t_\star(x)}+\cE(\hC(x)),
	\end{align}
where $\Leb$ is Lebesgue measure. The first term is due to the coverage error. The second, the \emph{excess mass} $\cE(\hC(x))\ge0$, is due to points on the wrong side of the oracle level. The two terms have different minimax rates.

\paragraph{Minimax lower bounds} Suppose that the conditional density is $\beta$-H\"older jointly in the covariate and the response, which have dimensions $d_X$ and $d_Y$. Let
	\begin{align}
		\zeta_n:=(\log n/n)^{\beta/(2\beta+d_X+d_Y)},\qquad\rho_n:=(\log n/n)^{\beta/(2\beta+d_X)}
	\end{align}
be the density rate and the regression rate. We prove lower bounds of order $\rho_n$ for the uniform conditional coverage error and, up to a logarithmic factor, of order $\zeta_n^2$ for the excess mass (\cref{thm:lower}). Coverage alone does not pin down a procedure, since data-free sets can have negligible coverage error and a volume that exceeds the oracle volume by a constant (\cref{sec:coverage}). Coverage and efficiency must therefore be studied together. We prove the coverage lower bounds among procedures whose volume exceeds the oracle volume by at most a fixed constant, a class that contains every uniformly consistent procedure. The excess-mass lower bound holds for all procedures.

\paragraph{Two-stage calibration} We propose two-stage calibration, which corrects the conditional coverage at the regression rate, although the density itself can be estimated only at the slower density rate. A density estimate fitted on one part of the sample fixes, at each covariate value, a nested family of candidate predictive sets, its super-level sets. The coverage $x'\mapsto{\mathsf P}_{x'}(A)$ of a fixed predictive set $A$ is a regression function of the covariate that is as smooth as the conditional density, whatever the density estimate. Estimating it is the regression problem behind the coverage lower bound. Two-stage calibration therefore regresses, on the covariates of a second part of the sample, the indicators that its responses fall in the candidate predictive sets of the target covariate, and it uses this estimated coverage as the conformity score. Hence the coverage error attains the regression rate, and the accuracy of the density estimate affects only the shape of the sets.

We show that two-stage calibration has conditional coverage error of order $\rho_n$ uniformly in $x$, even for an inconsistent density estimate (\cref{thm:local-cov}). This holds for neural network estimates, with no condition on their training or accuracy (\cref{prop:nn-conditions}). If the density estimate has sup-norm error $r_n$, the excess mass is of order $r_n^2+\rho_n^2$, up to a logarithmic factor at covariates at which two connected components of the oracle set merge (\cref{thm:local}). The result covers a kernel estimate that is piecewise constant in the covariate, as tree-based estimates are (\cref{prop:partition}). For a local polynomial density estimate, $r_n$ is of order $\zeta_n$, and two-stage calibration attains both lower bounds in every dimension, exactly for the uniform coverage error and up to logarithmic factors for the excess mass (\cref{prop:lp-conditions}). The procedure evaluates the density estimate only at the regression points near each covariate and computes no integral over the response space.

\paragraph{One-stage calibration} Density-based conformal sets in common use, such as HPD-split \citep{izbicki2022cd} and conformal highest conditional density sets (CHCDS) \citep{sampson2024flexible}, instead calibrate a single conformal threshold, an approach that we call one-stage calibration. For HPD-split, the coverage error is, to first order, the variation in $x$ of the density error integrated over the oracle set, because a single threshold removes only the average of this error over $x$ (\cref{thm:onestage}). The integral averages out the noise of a smoothing density estimate but not its bias. For the local polynomial plug-in at the density-optimal bandwidth, one-stage calibration therefore has conditional coverage error of order $\zeta_n$, and, in general, neither this bandwidth nor a larger one attains the regression rate (\cref{thm:onestage-lower} and~\cref{cor:lp-onestage}). Smaller bandwidths reduce the bias but raise the noise, and the guarantee for the shape then degrades (\cref{rem:undersmoothing}).

\paragraph{Shape and sup-norm accuracy} We measure the shape of a set by its \emph{shape loss} $\cS(\hC(x))$, the extra volume of $\hC(x)$ over the smallest set with the same coverage. For procedures with accurate coverage, it has the minimax order of the excess mass. A density estimate orders the responses, and calibration selects an initial segment of this order. Points are therefore misplaced only where the density estimate reverses the order of responses across the level.

We show that a sup-norm error $r$ gives shape loss at most of order $r^2$, up to the same logarithmic factor (\cref{prop:supshape}). A density estimate attaining the density rate therefore has shape loss of order $\zeta_n^2$. In contrast, an $L^q$ error $r^{(q)}$ near the level set gives shape loss only of order $(r^{(q)})^{2q/(q+1)}$ (\cref{thm:Lq}), and a perturbation of the true density shows that this exponent is sharp for every $q\in[1,\infty]$ (\cref{thm:Lq-lower} in the supplement). Even within H\"older classes, an error of minimax order in $L^2$ or total variation does not guarantee shape loss of order $\zeta_n^2$ (\cref{thm:sep-smooth}).

\paragraph{Related work} Several conformal methods build predictive sets from estimates of the conditional distribution or density \citep{chernozhukov2021distributional,izbicki2022cd,sampson2024flexible,plassier2024conditionally}. Others let the calibration depend on the covariate through quantile regression over a function class \citep{gibbs2023conformal,duchi2025sample}, a rectifying quantile estimate \citep{plassier2025rectifying}, random localization \citep{hore2025conformal}, a localized probability integral transform \citep{zhou2025conformal}, a recalibration of HPD levels without a conformal step \citep[Appendix~D.2]{dey2025ladar} or minimax optimization over a class of sets \citep{bao2026shape}.

Excess-mass estimation of level sets goes back to \citet{hartigan1987estimation,polonik1995measuring,tsybakov1997nonparametric} (see also \citealp{mammen1995asymptotical,cadre2006kernel,scott2006learning,singh2009adaptive,rinaldo2010generalized,chen2017density}), and conditional minimum-volume regions to \citet{polonik2000conditional}. For plug-in level sets and classification at a \emph{fixed} level, pointwise concentration and a margin condition give excess-risk and excess-mass rates, with matching lower bounds \citep{tsybakov1997nonparametric,mammen1999smooth,audibert2007fast,rigollet2009optimal}. Here the level is set by calibration, and the comparison set depends on the achieved coverage.

For highest-density regions without covariates, \citet{mason2009asymptotic} establish asymptotic normality of the symmetric-difference error of plug-in level sets. \citet{samworth2010asymptotics,doss2018bandwidth} show that the bandwidth that is optimal for the density is not optimal for the region. \Cref{cor:lp-onestage} is a conditional and calibrated analogue for the coverage.

The shape loss is the excess of a mass-volume curve, a criterion introduced for anomaly ranking \citep{clemencon2013scoring,clemencon2018mass}. Our density estimates fit a network or a local polynomial to a contrastive regression function (\cref{sec:cov-nn,sec:cov-twostage-alg}). This device goes back to \citet{qin1998inferences} and to the noise-contrastive estimation of \citet{gutmann2012noise}, and \citet{dai2025cindes} use it with deep networks. Other conditional density estimators include those of \citet{efromovich2010dimension,izbicki2017converting}.

\paragraph{Notation} The marginal distribution of $X$ is $\Pi$, and the conditional distribution of $Y$ given $X=x$ is ${\mathsf P}_x$, with Lebesgue density ${\mathsf p}(\cdot\mid x)$. The estimates of ${\mathsf P}_x$ and ${\mathsf p}(\cdot\mid x)$ are $\hP_x$ and $\hp(\cdot\mid x)$. We write $\triangle$ for symmetric difference, $\ind\{\cdot\}$ for indicators, $|\cdot|$ for the Euclidean norm and $|\cdot|_\infty$ for the maximum norm. For a distribution ${\mathsf P}$ of $(X,Y)$, $\E_{\mathsf P}$ and $\bP_{\mathsf P}$ are the expectation and the probability under the product distribution of the i.i.d.\ sample from ${\mathsf P}$ and of the auxiliary randomness. We write $\E$ and $\bP$ when the distribution is fixed, and $\sigma(Z)$ for the $\sigma$-field generated by a random element $Z$. For positive sequences, $a_n\lesssim b_n$ means $a_n\le Cb_n$ for a constant $C>0$ that depends neither on $n$ nor on $x$ or the distribution, where these vary. Further, $a_n\gtrsim b_n$ means $b_n\lesssim a_n$, $a_n\asymp b_n$ means both, and $a_n\ll b_n$ means $a_n/b_n\to0$. Unless stated otherwise, $C$ and $c$ denote positive constants that depend only on $\alpha$, $\beta$, $d_X$, $d_Y$, the sets $\cY$ and $\cX_0$ of \cref{sec:setup}, the constants of the assumptions and the tuning parameters that do not depend on $n$. Their values may change from one occurrence to the next. The implied constants in $\lesssim$, $\gtrsim$, $\asymp$ and $O(\cdot)$ are of the same kind, and ``for all sufficiently large $n$'' means $n\ge n_0$ for such a constant $n_0$, which may also depend on the deterministic sequences in the statement, such as a bandwidth.

\paragraph{Organization} \Cref{sec:setup} defines HPD-split and states the assumptions, and \cref{sec:losses} introduces the losses and the volume decomposition. \Cref{sec:coverage} gives the minimax lower bounds. \Cref{sec:cov-twostage} develops two-stage calibration under conditions on the density estimate and verifies them for neural network estimates, a local polynomial estimate and an estimate that is piecewise constant in the covariate. \Cref{sec:cov-onestage} analyzes one-stage calibration and compares our guarantees with existing ones. \Cref{sec:shape} treats the shape loss under $L^q$ accuracy, and \cref{sec:experiments} reports the experiments. In the supplement, \cref{app:literature} states the existing guarantees in our notation, and the other appendices contain all proofs and the details of the experiments.

\section{Setting and assumptions}\label{sec:setup}
This section defines the oracle sets and HPD-split and states the conditions on the conditional density and the covariate distribution. Since HPD-split selects the level from the data, we define the oracle sets at every mass and impose the conditions at all levels near the oracle level.

\subsection{Oracle sets and HPD-split}\label{sec:setup-oracle}
Let $(X,Y)\in\cX\times\cY$, where $\cX\subset\R^{d_X}$ and $\cY\subset\R^{d_Y}$ is a compact box, and let
	\begin{align}
		G_x(v):={\mathsf P}_x\{{\mathsf p}(Y\mid x)\ge v\},\qquad t_\gamma(x):=\sup\{v:G_x(v)\ge \gamma\}\quad(\gamma\in(0,1]),
	\end{align}
so that $G_x(v)$ is the coverage of the super-level set $\{y:{\mathsf p}(y\mid x)\ge v\}$, and $t_\gamma(x)$ is the largest level whose super-level set has coverage at least $\gamma$. For $\gamma=0$, we set $t_0(x):=\sup_{y\in\cY}{\mathsf p}(y\mid x)$. Then $\gamma\mapsto t_\gamma(x)$ is nonincreasing on $[0,1]$. The \emph{oracle set} is the HPD set of the true conditional density with coverage at least $\gamma$,
	\begin{align}
		\Cstar_\gamma(x):=\{y:\,{\mathsf p}(y\mid x)\ge t_\gamma(x)\},\qquad V_x(\gamma):=\Leb\bigl(\Cstar_\gamma(x)\bigr),
	\end{align}
where $V_x(\gamma)$ is its volume. We write $\gamma_\star:=1-\alpha\in(0,1)$ for the target coverage, and $\Cstar:=\Cstar_{\gamma_\star}$ and $t_\star:=t_{\gamma_\star}$ for the oracle set at this coverage and its level. By the Neyman--Pearson lemma, $\Cstar_\gamma(x)$ has the smallest Lebesgue measure among all sets with coverage at least $\gamma$ whenever $G_x(t_\gamma(x))=\gamma$ and $t_\gamma(x)>0$. Both conditions hold for $\gamma$ near $\gamma_\star$ under \cref{ass:margin} below.

\paragraph{HPD-split} The procedure of \citet{izbicki2022cd} conformalizes the HPD sets of a fitted density. Let $\hp(\cdot\mid x)$ be a conditional density estimate fitted on data independent of a calibration sample $(X_i,Y_i)_{i\le n_{\rm cal}}$, and let $\widehat{\mathsf P}_x$ be the distribution with density $\hp(\cdot\mid x)$. For $v\ge0$ and $\widetilde Y\sim\widehat{\mathsf P}_x$, the estimated masses strictly above and at or above the level $v$ are
	\begin{align}
		\Lambda_x(v):=\widehat{\mathsf P}_x\{\hp(\widetilde Y\mid x)>v\},\qquad\overline\Lambda_x(v):=\widehat{\mathsf P}_x\{\hp(\widetilde Y\mid x)\ge v\}.
	\end{align}
The nonconformity score of $y$, which we call its \emph{HPD value}, is the estimated mass above $v:=\hp(y\mid x)$, with ties broken at random,
	\begin{align}\label{eq:score}
		\hat{\mathsf s}^{\rm hpd}(x,y):=\Lambda_x(v)+U\bigl\{\overline\Lambda_x(v)-\Lambda_x(v)\bigr\},
	\end{align}
where $U\sim U[0,1]$ is a tie-breaking variable. The score lies in $[0,1]$ and is uniform on $[0,1]$ for $\widetilde Y\sim\widehat{\mathsf P}_x$. When the level sets of $\hp(\cdot\mid x)$ are Lebesgue-null, \eqref{eq:score} does not depend on $U$ and equals one minus the conformity score of \citet{izbicki2022cd}. Given $\hp$ and the covariates, $U$ is generated independently at each calibration and test point and is \emph{shared by every $y$ scored at that point}. Hence the score of the test point is exchangeable with the scores of the calibration points, the \emph{calibration scores}. In practice, the masses $\Lambda_x$ are computed by Monte Carlo or on a grid.

HPD-split takes $\hat\tau$ to be the $\lceil(1-\alpha)(n_{\rm cal}+1)\rceil$-th smallest calibration score, with $\hat\tau:=1$ if this index exceeds $n_{\rm cal}$, in which case the set is $\cY$. It returns $\hC(x)=\{y:\hat{\mathsf s}^{\rm hpd}(x,y)\le\hat\tau\}$. Since the score is nonincreasing in $\hp(y\mid x)$, this set lies between $\{y:\hp(y\mid x)>\hat t(x)\}$ and $\{y:\hp(y\mid x)\ge\hat t(x)\}$ for some level $\hat t(x)$, that is, it is a super-level set of $\hp(\cdot\mid x)$ up to ties. For a new pair $(X,Y)$ independent of the data, by exchangeability, we have $\bP(Y\in\hC(X))\ge1-\alpha$, and $\bP(Y\in\hC(X))\le1-\alpha+1/(n_{\rm cal}+1)$ when the calibration scores are almost surely distinct \citep{papadopoulos2002inductive,lei2018distribution,angelopoulos2023conformal}.

\begin{remark}[Choice of the HPD value as score]\label{rem:hpd-value}
The estimated density $\hp(y\mid x)$ itself also orders the responses at each $x$, but as a score it calibrates a single density level for all covariates. Since the coverage $G_x(t)$ of a super-level set varies with $x$, a common level gives a coverage that varies with $x$ even when $\hp={\mathsf p}$. The HPD value instead thresholds an estimated mass, which has the same meaning at every $x$. With $\hp={\mathsf p}$ and Lebesgue-null level sets, the set $\{y:\hat{\mathsf s}^{\rm hpd}(x,y)\le u\}$ has conditional coverage exactly $u$ at every $x$.
\end{remark}
\subsection{Regularity conditions}\label{sec:setup-assumptions}
We use two regularity conditions. For $\beta>0$ and $L>0$, let $\ell:=\lceil\beta\rceil-1$ and let $\mathcal H^\beta_{\mathcal D}(L)$ be the H\"older class of functions $g$ on a domain $\mathcal D\subset\R^{d'}$ whose partial derivatives up to order $\ell$ exist and satisfy
	\begin{equation}\label{eq:holder}
		\max_{|\boldsymbol k|\le\ell}\sup_{\mathcal D}|\partial^{\boldsymbol k}g|\le L,\qquad\max_{|\boldsymbol k|=\ell}\ |\partial^{\boldsymbol k}g(z)-\partial^{\boldsymbol k}g(z')|\le L|z-z'|^{\beta-\ell}\quad(z,z'\in\mathcal D).
	\end{equation}
Here $\boldsymbol k=(k_1,\dots,k_{d'})\in\{0,1,\dots\}^{d'}$ is a multi-index, with $|\boldsymbol k|:=k_1+\dots+k_{d'}$, $\boldsymbol k!:=k_1!\cdots k_{d'}!$, $\partial^{\boldsymbol k}:=\partial_1^{k_1}\cdots\partial_{d'}^{k_{d'}}$ and $s^{\boldsymbol k}:=s_1^{k_1}\cdots s_{d'}^{k_{d'}}$ for $s\in\R^{d'}$. We write $\mathcal H^\beta_{\mathcal D}$ for the union over $L$, and drop the subscript $\mathcal D$ when the domain is clear. The first condition imposes smoothness.

\begin{assumption}[Smoothness]\label{ass:smooth}
The map $(x,y)\mapsto {\mathsf p}(y\mid x)$ belongs to the H\"older class $\mathcal H^\beta_{\cX\times\cY}(L)$ for some $\beta\ge1$, and $\sup_{(x,y)\in\cX\times\cY}{\mathsf p}(y\mid x)\le {\mathsf p}_{\max}$.
\end{assumption}

Throughout, $d_X,d_Y\ge1$ and $d:=d_X+d_Y$. Joint smoothness in $(x,y)$ gives the conditional density rate with exponent $\beta/(2\beta+d)$. In contrast, \citet[Assumption~1]{lei2014distribution} assume smoothness of order $\beta$ in $y$ and only Lipschitz continuity in $x$, which leads to a slower rate. The restriction $\beta\ge1$ makes $t_\star$ Lipschitz in $x$.

The second is a two-sided margin condition with exponent one at all levels near $t_\star(x)$, up to a logarithmic factor in its upper bound. Let $\mathfrak L$ be either $\mathfrak L(\eps):=1$ or $\mathfrak L(\eps):=1+\log_+(1/\eps)$ on $(0,\infty)$. For $k\ge1$, we set $u^k\mathfrak L(u):=0$ and $u^k/\mathfrak L(u):=0$ at $u=0$.

\begin{assumption}[Margin]\label{ass:margin}
There exist constants $\varepsilon_0,c_m,C_m>0$ and $a\in(0,1)$ such that, for all $x\in\cX$, all $t$ with $|t-t_\star(x)|\le\varepsilon_0$ and all $\eps\in(0,\varepsilon_0]$,
	\begin{align}
		&\Leb\{y:|{\mathsf p}(y\mid x)-t|\le\eps\}\le C_m\eps\,\mathfrak L(\eps),\label{eq:margin-upper}\\
		&\min\cbr[1]{{\mathsf P}_x\{t+a\eps\le{\mathsf p}(Y\mid x)<t+\eps\},\ {\mathsf P}_x\{t-\eps<{\mathsf p}(Y\mid x)\le t-a\eps\}}\ge c_m\eps.\label{eq:margin-lower}
	\end{align}
\end{assumption}

In words, ${\mathsf p}(\cdot\mid x)$ is neither flat nor too steep near the level \citep{lei2013distribution}. The upper bound \eqref{eq:margin-upper} rules out atoms of the distribution of ${\mathsf p}(Y\mid x)$ near $t_\star$, and the two disjoint one-sided bands of \eqref{eq:margin-lower} give ${\mathsf P}_x\{|{\mathsf p}(Y\mid x)-t|\le\eps\}\ge2c_m\eps$. Hence $G_x$ is continuous and strictly decreasing near $t_\star$. Consequently, with $c_0:=c_m\varepsilon_0/8$, we have $|t_\gamma(x)-t_\star(x)|\le\varepsilon_0/4$ for $|\gamma-\gamma_\star|\le c_0$, the map $\gamma\mapsto t_\gamma(x)$ is continuous and strictly decreasing there, and ${\mathsf P}_x(\Cstar_\gamma(x))=\gamma$. Since the lower band at the level $t_\star(x)-\varepsilon_0$ has positive probability, we have $t_\star(x)>\varepsilon_0$, and hence $t_\gamma(x)>\varepsilon_0/2$ for $|\gamma-\gamma_\star|\le c_0$. \Cref{lem:t-window} in the supplement collects these facts, together with a Lipschitz bound on $\gamma\mapsto t_\gamma(x)$. \citet[Assumption~2]{lei2014distribution} assume this lower bound on the levels. It bounds the oracle volume by $2/\varepsilon_0$ and the factor $1/t_\star(x)$ in \cref{lem:identity}.

The logarithmic factor admits covariates at which the topology of $\Cstar(x)$ changes, for instance where two of its connected components merge. The level set then passes through a critical point of ${\mathsf p}(\cdot\mid x)$. At a nondegenerate saddle point, the band $\{y:|{\mathsf p}(y\mid x)-t|\le\eps\}$ has volume of order $\eps\log(1/\eps)$ when $d_Y=2$ and of order $\eps$ when $d_Y\ge3$. For $d_Y=1$, the topology changes at local extrema, where the band has volume of order $\eps^{1/2}$ and neither choice of $\mathfrak L$ applies.

By the coarea formula, \eqref{eq:margin-upper} and \eqref{eq:margin-lower} hold with $\mathfrak L\equiv1$ when, uniformly in $x$, ${\mathsf p}(\cdot\mid x)$ is continuously differentiable on the band $\mathcal J_x(2\varepsilon_0)$, where
	\begin{align}\label{eq:Jx}
		\mathcal J_x(\eps):=\{y:|{\mathsf p}(y\mid x)-t_\star(x)|\le\eps\},\qquad\mathcal J_x:=\mathcal J_x(\varepsilon_0),
	\end{align}
this band lies in the interior of $\cY$, and the gradient is uniformly continuous and bounded above and away from zero on it. These are the gradient conditions of \citet{doss2018bandwidth} and \citet{qiao2020asymptotics}. A margin condition at all levels near $t_\star$ must have exponent one, since a fixed range of levels is the union of about $1/\eps$ disjoint bands of width $2\eps$ (\cref{rem:margin-exponent} in the supplement).

Margin conditions go back to the $\gamma$-exponent condition of \citet{polonik1995measuring}, an upper bound on the band uniformly over a set of levels. They are standard in level-set estimation and classification \citep{tsybakov1997nonparametric,audibert2007fast,rigollet2009optimal}. Two-sided bounds at the oracle level are assumed by \citet{lei2013distribution} for an unconditional density and by \citet{lei2014distribution} uniformly in $x$. \citet{izbicki2022cd} assume that, for every $x$, the distribution of ${\mathsf p}(Y\mid x)$ has a density bounded above and, near the HPD level, away from zero. Their condition implies \cref{ass:margin} when the neighborhood and the constants do not depend on $x$.

\subsection{The analysis set}\label{sec:rates}
We state the guarantees on an \emph{analysis set} $\cX_0\subset\cX$, for example a set slightly inside the support of $X$, and evaluate the losses at $x\in\cX_0$. From now on, we compute $\hat\tau$ from the $n_{\rm cal}^\circ$ calibration points with $X_i\in\cX_0$ and replace $n_{\rm cal}$ by $n_{\rm cal}^\circ$ in its definition. The coverage guarantees of \cref{sec:setup-oracle} then hold conditionally on $X\in\cX_0$.

\begin{assumption}[Design]\label{ass:design}
The analysis set $\cX_0$ is compact and convex with nonempty interior. For $\delta>0$, let $\cX_0^{\delta}:=\{x:\ \inf_{x'\in\cX_0}|x-x'|_\infty\le\delta\}$ be the $\delta$-neighborhood of $\cX_0$. There exist constants $\delta_0,f_{\min}>0$ such that $\Pi$ has a density $f_X$ on $\cX_0^{2\delta_0}$ with $f_X\in\mathcal H^1_{\cX_0^{2\delta_0}}(L)$ and $f_X\ge f_{\min}$.
\end{assumption}

Convexity of $\cX_0$ makes the sets $\cX_0^{\delta}\times\cY$ convex, and hence the H\"older bounds of the proofs apply without extension arguments. The condition on $\Pi$ is the density condition of \citet[Assumption~1(a)]{lei2014distribution} and the strong density assumption of \citet{audibert2007fast}, required here only near $\cX_0$. The condition $f_X\in\mathcal H^1_{\cX_0^{2\delta_0}}(L)$ means that $f_X$ is bounded by $L$ and $L$-Lipschitz.
\section{Losses and the volume decomposition}\label{sec:losses}
This section defines the losses and relates them to the volume. A predictive set can exceed the oracle volume by covering the wrong mass or by having the wrong shape. The coverage error, the excess mass and the shape loss separate the two effects.

\subsection{Excess mass and the volume decomposition}\label{sec:setup-losses}
A predictive set $\hC$ is a function of the data and of auxiliary randomness whose values $\hC(x)$ are Borel subsets of $\cY$, with $\ind\{y\in\hC(x)\}$ jointly measurable in $x$, $y$, the data and the auxiliary randomness. Let $\cov(x):={\mathsf P}_x(\hC(x))$ denote its conditional coverage given $X=x$, the data and all auxiliary randomness. Bounds on $\cov(x)$ at each $x$ with high probability over the data combine the object-conditional and training-conditional validity of \citet{vovk2012conditional}. Besides the coverage error $|\cov(x)-\gamma_\star|$ and the symmetric difference $\Leb(\hC(x)\triangle\Cstar(x))$, we use two geometric losses. The \emph{excess mass at the oracle level} is
	\begin{align}
		\cE(\hC(x)):=\frac1{t_\star(x)}\int_{\hC(x)\triangle\Cstar(x)}\bigl|{\mathsf p}(y\mid x)-t_\star(x)\bigr|\,\d y .
	\end{align}
Up to the factor $1/t_\star(x)$, this is the loss in the excess-mass functional of \citet{hartigan1987estimation} and \citet{polonik1995measuring} at the level $t_\star(x)$ relative to the oracle set. The coverage error and the excess mass together determine the excess volume of any set.

\begin{lemma}[Volume decomposition]\label{lem:identity}
Suppose that \cref{ass:margin} holds. Then, for every $x$ and every measurable set $\hC(x)$,
	\begin{align}
		\Leb\bigl(\hC(x)\bigr)-\Leb\bigl(\Cstar(x)\bigr)=\frac{\cov(x)-\gamma_\star}{t_\star(x)}+\cE(\hC(x)).
	\end{align}
\end{lemma}

\paragraph{Excess mass and symmetric difference} Under \cref{ass:smooth,ass:margin}, there exist constants $c,C>0$ such that, for any $x$ and any measurable set $\hC(x)$,
	\begin{equation}\label{eq:E-symdiff}
		\frac{c\,\Leb\bigl(\hC(x)\triangle\Cstar(x)\bigr)^2}{\mathfrak L\bigl(\Leb(\hC(x)\triangle\Cstar(x))\bigr)}\le\cE(\hC(x))\le C\,\Leb\bigl(\hC(x)\triangle\Cstar(x)\bigr).
	\end{equation}
The lower bound holds because at most $C_m\eps\,\mathfrak L(\eps)$ of the volume of $\hC(x)\triangle\Cstar(x)$ lies within $\eps$ of the level. With $\mathfrak L\equiv1$, it is sharp when the symmetric difference is a thin layer around the boundary of $\Cstar(x)$. Accordingly, the minimax lower bound for the excess mass in \cref{thm:lower} is the square of that for the symmetric difference.

\subsection{The shape loss and sup-norm accuracy}\label{sec:setup-shape}\label{sec:setup-supnorm}
The excess mass still depends on the coverage. A super-level set $A$ of ${\mathsf p}(\cdot\mid x)$ whose coverage differs from $\gamma_\star$ has $\cE(A)>0$, although no set of the same coverage is smaller. The \emph{shape loss} removes this dependence by comparing $\hC(x)$ with the best set of the same coverage,
	\begin{align}
		\cS(\hC(x)):=\Leb\bigl(\hC(x)\bigr)-V_x\bigl(\cov(x)\bigr)=\frac1{t_{\cov(x)}(x)}\int_{\hC(x)\triangle\Cstar_{\cov(x)}(x)}\bigl|{\mathsf p}(y\mid x)-t_{\cov(x)}(x)\bigr|\,\d y,
	\end{align}
defined when the coverage of $\hC(x)$ is \emph{in the window} $|\cov(x)-\gamma_\star|\le c_0$, where the two expressions are equal. Outside the window, we set $\cS(\hC(x)):=\cE(\hC(x))$.

We denote by $\cS^\circ(\hC(x))$ the integral expression at the achieved coverage, defined also outside the window, with $0/0:=0$ when $t_{\cov(x)}(x)=0$.

\begin{lemma}[Shape loss bound]\label{lem:shape-split}
Suppose that \cref{ass:margin} holds. Then there exists a constant $C>0$ such that, for any $x$ and any measurable set $\hC(x)$,
	\begin{align}
		\cS(\hC(x))\le\cE(\hC(x))\le\cS(\hC(x))+C\,\bigl|\cov(x)-\gamma_\star\bigr|^2 .
	\end{align}
\end{lemma}

Within the window,
	\begin{align}
		\cE(\hC(x))-\cS(\hC(x))=V_x(\cov(x))-V_x(\gamma_\star)-\frac{\cov(x)-\gamma_\star}{t_\star(x)} .
	\end{align}
This is the gap between the mass-volume function $V_x$, which is convex on $[\gamma_\star-c_0,\gamma_\star+c_0]$ with $V_x'(\gamma)=1/t_\gamma(x)$, and its tangent at $\gamma_\star$. The excess volume is therefore the sum of the linear coverage term of \cref{lem:identity}, a term of the order of the squared coverage error, and the shape loss.

\paragraph{Interpretation of the shape loss} At a fixed coverage in the window, the shape loss measures the extra volume caused by misranking responses relative to the true conditional density. By the integral expression, each point of $\hC(x)\triangle\Cstar_{\cov(x)}(x)$ contributes its distance in density from the level $t_{\cov(x)}(x)$, in units of that level. Hence $\cS(\hC(x))=0$ if and only if $\hC(x)$ is a super-level set of ${\mathsf p}(\cdot\mid x)$ up to a null set. For a set $\hC(x)=\{y:\hat a(y\mid x)\ge\hat t\}$, the shape loss thus depends only on the coverage and on the order in which the \emph{ranking function} $\hat a(\cdot\mid x)$ sorts the points of $\cY$.

For example, let ${\mathsf p}(\cdot\mid x)$ be a symmetric bimodal density on the line whose minimum between the modes lies well below $t_\star$. An interval centered between the modes with coverage exactly $\gamma_\star$ has no coverage error but positive shape loss, since it contains the low-density region that separates them. A super-level set of ${\mathsf p}(\cdot\mid x)$ with coverage in the window, here two intervals around the modes, has zero shape loss. Its excess volume $V_x(\cov(x))-V_x(\gamma_\star)$ is due entirely to the coverage error, which a different calibration can remove.

\paragraph{Sup-norm control} The shape loss of a super-level set of a ranking function is at most of the order of its squared sup-norm error near the level set.

\begin{proposition}[Shape loss under sup-norm accuracy]\label{prop:supshape}
Suppose that \cref{ass:smooth,ass:margin} hold. Then there exists a constant $C>0$ such that the following holds. For a fixed $x$, let $\hat a$ be a function on $\cY$ with $\sup_\cY|\hat a-{\mathsf p}(\cdot\mid x)|\le\varepsilon_0/8$ and $\sup_{\mathcal J_x}|\hat a-{\mathsf p}(\cdot\mid x)|\le r\le\varepsilon_0/8$, with $\mathcal J_x$ the band \eqref{eq:Jx}. Let $\hC(x)$ satisfy $\{y:\hat a(y)>\hat t\}\subset\hC(x)\subset\{y:\hat a(y)\ge\hat t\}$ for some level $\hat t$, and suppose that its coverage is in the window. Then
	\begin{align}
		\cS(\hC(x))\le Cr^2\mathfrak L(r),\qquad\cE(\hC(x))\le C\bigl(r^2\mathfrak L(r)+|\cov(x)-\gamma_\star|^2\bigr).
	\end{align}
\end{proposition}

\section{Lower bounds for coverage and excess mass}\label{sec:coverage}\label{sec:cov-lower}
By the volume decomposition, the coverage error and the excess mass determine the excess volume. We bound these two quantities and the symmetric difference from below, and we bound the coverage error both integrated over the covariate and uniformly in it. Up to logarithmic factors, the bounds are given by the conditional density rate $\zeta_n$ \citep{fan1996estimation,efromovich2007conditional,efromovich2010dimension} and the regression rate $\rho_n$ \citep{stone1982optimal} of \cref{sec:intro}.

The proofs of the lower bounds perturb a conditional density $f_0(y)$ that does not depend on $x$, namely a Gaussian density restricted to $\cY$ and renormalized (see \cref{app:lower} of the supplement for the construction). Part (i) of the following assumption ensures that the class below contains $f_0$ and its perturbations, and part (ii) fixes a uniform design.

\begin{assumption}[Regularity for the lower bounds]\label{ass:constants}
The constants and the design satisfy the following.
\begin{enumerate}[label=(\roman*)]
\item (Constants) Let $L_0,{\mathsf p}_\sharp,\bar C,\kappa_0,c_\natural>0$ be the thresholds of \eqref{eq:lb-thresholds} of the supplement, which depend only on $\beta$, $d_Y$, $\cY$ and $\alpha$. The constants of \cref{ass:smooth,ass:margin,ass:design} satisfy $L\ge L_0$, ${\mathsf p}_{\max}\ge{\mathsf p}_\sharp$, $C_m\kappa_0\ge\bar C$, $\varepsilon_0\le\kappa_0$, $c_m\le\min\{1-a,\frac12\}c_\natural$ and $f_{\min}\le1$.
\item (Design) The covariate distribution $\Pi$ is uniform on $[0,1]^{d_X}$, so that $\cX=[0,1]^{d_X}$, and $\cX_0$ is a fixed axis-parallel cube at distance at least $2\delta_0$ from the boundary of $[0,1]^{d_X}$.
\end{enumerate}
\end{assumption}

The conditions on $C_m$ and $c_m$ are stated relative to $\kappa_0$ and $a$, because the margin constants of $f_0$ and of its perturbations scale with these. We prove the lower bounds over the class
	\begin{align}
		\cP_\beta:=\bigl\{{\mathsf P}:\ {\mathsf P}\text{ satisfies \cref{ass:smooth,ass:margin,ass:design,ass:constants}}\bigr\}.
	\end{align}
Here $\beta$, $\mathfrak L$, $\cX_0$ and the other constants of these assumptions are fixed. The distributions in the proofs satisfy \cref{ass:margin} with $\mathfrak L\equiv1$, and the bounds therefore hold for either choice of $\mathfrak L$.

\begin{theorem}[Minimax lower bounds]\label{thm:lower}
Let $\mathfrak C_n$ be the class of predictive sets that are functions of a sample of size $n$ and of auxiliary randomness, and, for $c_1,c_2>0$, let
	\begin{equation}
		\mathfrak V_n(c_1,c_2):=\Biggl\{\hC\in\mathfrak C_n:\ \sup_{{\mathsf P}\in\cP_\beta}\E_{\mathsf P}\sbr{\int_{\cX_0}\ind\bigl\{\Leb\bigl(\hC(x)\bigr)-\Leb\bigl(\Cstar(x)\bigr)>c_1\bigr\}\,\d x}\le c_2\Biggr\}.
	\end{equation}
\begin{enumerate}[label=(\roman*)]
\item (Excess mass and symmetric difference) There exists a constant $c>0$ such that, for all $n\ge1$,
	\begin{align}
		\inf_{\hC\in\mathfrak C_n}\sup_{{\mathsf P}\in\cP_\beta}\E_{\mathsf P}\sbr{\int_{\cX_0}\cE(\hC(x))\,\d x}&\ge c\,n^{-2\beta/(2\beta+d)},\\
		\inf_{\hC\in\mathfrak C_n}\sup_{{\mathsf P}\in\cP_\beta}\E_{\mathsf P}\sbr{\int_{\cX_0}\Leb\bigl(\hC(x)\triangle\Cstar(x)\bigr)\,\d x}&\ge c\,n^{-\beta/(2\beta+d)} .
	\end{align}
\item (Coverage, under coarse volume accuracy) There exist constants $c_1,c>0$ and $c_2\in(0,\Leb(\cX_0))$ such that, for all $n\ge1$,
	\begin{equation}
		\inf_{\hC\in\mathfrak V_n(c_1,c_2)}\sup_{{\mathsf P}\in\cP_\beta}\E_{\mathsf P}\sbr{\int_{\cX_0}|\cov(x)-\gamma_\star|\,\d x}\ge c\,n^{-\beta/(2\beta+d_X)} .
	\end{equation}
\item (Uniform coverage) With the constants $c_1,c_2$ of (ii), there exists a constant $c\in(0,1)$ such that, for all sufficiently large $n$,
	\begin{equation}
		\inf_{\hC\in\mathfrak V_n(c_1,c_2)}\sup_{{\mathsf P}\in\cP_\beta}\bP_{\mathsf P}\cbr[3]{\sup_{x\in\cX_0}|\cov(x)-\gamma_\star|\ge c\Bigl(\frac{\log n}n\Bigr)^{\beta/(2\beta+d_X)}}\ge c ,
	\end{equation}
where $\bP_{\mathsf P}$ is an outer probability if the event is not measurable. Consequently, with $\E^*_{\mathsf P}$ the outer expectation,
	\begin{align}
		\inf_{\hC\in\mathfrak V_n(c_1,c_2)}\sup_{{\mathsf P}\in\cP_\beta}\E^*_{\mathsf P}\sbr{\sup_{x\in\cX_0}|\cov(x)-\gamma_\star|}\ge c^2\Bigl(\frac{\log n}n\Bigr)^{\beta/(2\beta+d_X)} .
	\end{align}
\end{enumerate}
\end{theorem}

The class $\mathfrak V_n(c_1,c_2)$ excludes only procedures whose sets exceed the oracle volume by more than $c_1$ on a set of covariates of expected measure more than $c_2$. It imposes no restriction on the shape of the sets or on sets that are too small, whose coverage falls short of $\gamma_\star$ by at least $t_\star(x)$ times the volume deficit (\cref{lem:identity}). Since $\Leb(\hC(x))-\Leb(\Cstar(x))\le\Leb(\hC(x)\triangle\Cstar(x))$, the class contains, for all sufficiently large $n$, every procedure that is consistent for $\Cstar$ uniformly over $\cP_\beta$, such as two-stage calibration with the local polynomial estimate of \cref{sec:cov-twostage-alg}.

Some restriction is necessary. To see this, we partition $\cY$ into boxes with sides of order $1/n$ and choose in each box a sub-box that occupies the fraction $\gamma_\star$ of its volume. The union of these sub-boxes has coverage error $O(1/n)\ll\rho_n$ at every $x$ and uses no data. Under the distributions constructed in \cref{app:lower} of the supplement, however, its symmetric difference with $\Cstar(x)$ has volume of constant order. By \cref{lem:identity} and \eqref{eq:E-symdiff}, its volume then exceeds $\Leb(\Cstar(x))$ by a constant, and the class excludes it for sufficiently small $c_1$.

Without covariates, (i) reduces to the known lower bounds for density level sets under a margin condition \citep{tsybakov1997nonparametric,rigollet2009optimal}, and (ii) and (iii) parallel the lower bounds of \citet{stone1982optimal} in $L^1$ and in sup-norm for the regression function $x\mapsto{\mathsf P}_x(A)$ of a fixed set $A$. The proof of (i) perturbs $f_0$ by small bumps near the boundary of its oracle set, as \citet[Theorem~2]{lei2014distribution} do for $d_Y=1$. It exchanges density values between antipodal bumps, and all hypotheses therefore share the oracle level.

\begin{remark}[Lower bound for the shape loss]\label{rem:shape-lower}
Since $\cS(\hC(x))\le\cE(\hC(x))$, every upper bound on the excess mass also bounds the shape loss. Conversely, $\cS(\hC(x))\ge\cE(\hC(x))-C|\cov(x)-\gamma_\star|^2$ by \cref{lem:shape-split}. Let $\hC$ be a procedure whose expected integrated squared coverage error is at most a small constant times $n^{-2\beta/(2\beta+d)}$ under every ${\mathsf P}\in\cP_\beta$. Averaging over the hypotheses in the proof of \cref{thm:lower}(i), we then obtain
	\begin{align}
		\sup_{{\mathsf P}\in\cP_\beta}\E_{\mathsf P}\sbr{\int_{\cX_0}\cS(\hC(x))\,\d x}\ge c\,n^{-2\beta/(2\beta+d)} .
	\end{align}
Two-stage calibration with the local polynomial estimate of \cref{sec:cov-twostage-alg} satisfies this condition, since its squared coverage error is $O(\rho_n^2)$ by \cref{prop:lp-conditions}(ii). For such procedures the shape loss and the excess mass have the same minimax order.
\end{remark}

\section{Minimax optimal two-stage calibration}\label{sec:cov-twostage}
The coverage lower bound is the regression rate $\rho_n$, which is faster than the conditional density rate $\zeta_n$. We therefore let a density estimate fix the candidate predictive sets at each covariate value, and we estimate their coverage by a regression on the covariate.

\subsection{The procedure}\label{sec:cov-procedure}
We split the sample into disjoint parts $\mathcal D_{\rm fit}$, $\mathcal D_{\rm reg}$ and $\mathcal D_{\rm cal}$, used to fit a ranking function, to regress the coverage of its super-level sets on the covariate and to calibrate the threshold. Let $\hat a(y\mid x)$ be a ranking function fitted on $\mathcal D_{\rm fit}$, for instance a conditional density estimate, not necessarily normalized; it may be any function for which $(x,y,\mathcal D_{\rm fit})\mapsto\hat a(y\mid x)$ is Borel. At a covariate $x$, its super-level sets
	\begin{align}
		\mathcal U_x(t):=\{y\in\cY:\hat a(y\mid x)\ge t\}\qquad(t\in\R)
	\end{align}
are fixed given $\mathcal D_{\rm fit}$, and their coverage $Q_x(t):={\mathsf P}_x(\mathcal U_x(t))$ is nonincreasing in $t$. The ranking function thus fixes the order of the responses at each covariate, and the candidate predictive sets are the initial segments of this order. Given this order, conditional coverage requires only a level at each $x$ at which the candidate predictive set has coverage $1-\alpha$, whether or not the order is accurate. This is a conditional analogue of split conformal prediction, which attains marginal coverage with any score: the level is now chosen at each $x$, from an estimate of $Q_x$. HPD-split replaces $Q_x$ by the mass that the density estimate assigns to the candidate predictive sets, and this mass is only as accurate as the estimate.

\begin{figure}[t]
\centering
\includegraphics[width=\textwidth]{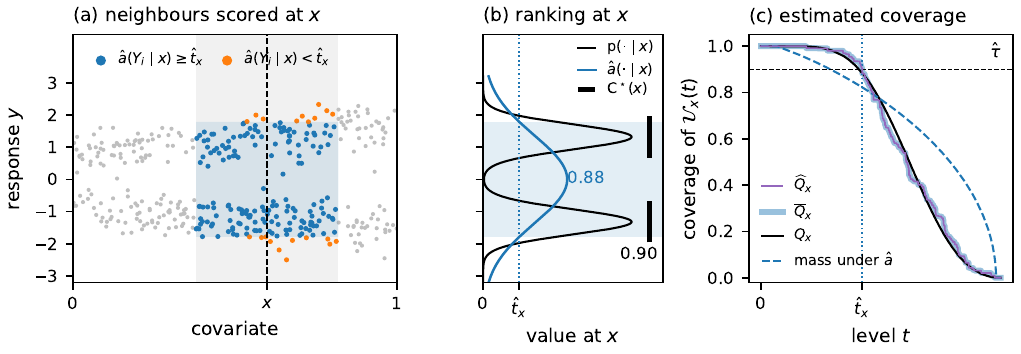}
\caption{Two-stage calibration at a covariate $x$, for a toy model with $d_X=d_Y=1$, a bimodal conditional distribution and, as the ranking function $\hat a$, the Gaussian density with the conditional mean and variance. (a) Regression points; the points in the window around $x$ are scored by $\hat a(\cdot\mid x)$, and the shaded candidate predictive set $\mathcal U_x(\hat t_x)$, fixed at $x$, is applied to all of them. (b) The ranking function $\hat a(\cdot\mid x)$ and the conditional density ${\mathsf p}(\cdot\mid x)$; $\mathcal U_x(\hat t_x)$ is the set where $\hat a(\cdot\mid x)\ge\hat t_x$, the bars show the oracle set $\Cstar(x)$, and the numbers are the true conditional coverages of the two sets. (c) The weighted counts $\widehat Q_x$, their envelope $\overline Q_x$, the true coverage $Q_x$ of the candidate predictive sets and, dashed, the mass that the normalized $\hat a(\cdot\mid x)$ assigns to them, which one-stage calibration uses in place of $Q_x$; $\overline Q_x$ falls to the threshold $\hat\tau=0.9$ at $\hat t_x$. The level at which the dashed curve falls to $0.9$ gives a set with true coverage $0.98$.}\label{fig:alg}
\end{figure}

The proposed two-stage calibration, \cref{alg:local}, instead estimates $Q_x(t)$ by counting how many of the responses observed near $x$ fall in $\mathcal U_x(t)$ (\cref{fig:alg}). It scores each response $Y_i$ of $\mathcal D_{\rm reg}$ by the ranking function at the target covariate, $\hat a(Y_i\mid x)$, and estimates $Q_x(t)$ for all $t$ at once by regressing the indicators $\ind\{\hat a(Y_i\mid x)\ge t\}$ on the covariates with a linear smoother, such as local polynomial regression. Since the sets $\mathcal U_x(t)$ are fixed given $\mathcal D_{\rm fit}$, the regression target is the coverage $x'\mapsto{\mathsf P}_{x'}(\mathcal U_x(t))$ of a fixed predictive set, evaluated at $x'=x$. This map is in $\mathcal H^\beta$ with norm at most $L\Leb(\cY)$, and the estimate therefore attains the rate $\rho_n$ for all candidate predictive sets, whatever the dependence of $\hat a$ on $x$. \Cref{alg:local} then scores $y$ by the estimated coverage of the smallest candidate predictive set that contains $y$ and selects the level at $x$ by applying split conformal prediction to this score. An inaccurate ranking function therefore changes the shape of the set but not the accuracy of its coverage.

\paragraph{Regression and score} The regression of \cref{alg:local} is a linear smoother on the covariates of $\mathcal D_{\rm reg}$, with weights $w_i(x)$ that are Borel functions of $x$ and of these covariates and do not depend on the responses, such as the local polynomial weights of \cref{app:second-stage} of the supplement. At a covariate $x$, the estimated coverage of the candidate predictive sets, its nonincreasing envelope and the score of $y$ are
	\begin{align}
		\widehat Q_x(t)&:=\sum_{i\in\mathcal D_{\rm reg}}w_i(x)\ind\{\hat a(Y_i\mid x)\ge t\},\label{eq:Qhat}\\
		\overline Q_x(t)&:=\min\Bigl\{1,\max\Bigl\{0,\sup_{v\ge t}\widehat Q_x(v)\Bigr\}\Bigr\},\label{eq:Qbar}\\
		\hat{\mathsf s}^{\rm two}(x,y)&:=\overline Q_x(\hat a(y\mid x))\qquad(t\in\R,\ y\in\cY).\label{eq:rank-score}
	\end{align}
The envelope $\overline Q_x$ is the smallest nonincreasing function above $\widehat Q_x$, clipped to $[0,1]$. The conformal threshold $\hat\tau$ is the $\lceil(1-\alpha)(n_{\rm cal}^\circ+1)\rceil$-th smallest of the scores $\hat{\mathsf s}^{\rm two}(X_i,Y_i)$ of the $n_{\rm cal}^\circ$ calibration points with $X_i\in\cX_0$, and $\hat\tau:=1$ if this index exceeds $n_{\rm cal}^\circ$, in which case the output is $\cY$, as for HPD-split. With $\hat t_x:=\inf\{t:\overline Q_x(t)\le\hat\tau\}$, the output $\hC^{\rm two}(x)=\{y:\hat{\mathsf s}^{\rm two}(x,y)\le\hat\tau\}$ of \cref{alg:local} lies between $\{y:\hat a(y\mid x)>\hat t_x\}$ and $\{y:\hat a(y\mid x)\ge\hat t_x\}$, since $\overline Q_x$ is nonincreasing. \Cref{sec:cov-theory} states the properties of the weights that the theory uses.

\begin{algorithm}[ht]
\caption{Two-stage calibration}\label{alg:local}
\KwIn{disjoint parts $\mathcal D_{\rm fit}$, $\mathcal D_{\rm reg}$, $\mathcal D_{\rm cal}$ of the sample, of sizes $n_{\rm fit}$, $n_{\rm reg}$ and $n_{\rm cal}$; miscoverage level $\alpha$; weights $w_i(\cdot)$ of a linear smoother on the covariates of $\mathcal D_{\rm reg}$.}
On $\mathcal D_{\rm fit}$: fit a ranking function $\hat a(y\mid x)$\;
At a covariate $x$: compute the weights $w_i(x)$, the values $\hat a(Y_i\mid x)$ \emph{at $x$} of the points of $\mathcal D_{\rm reg}$ with $w_i(x)\ne0$,
and the score $\hat{\mathsf s}^{\rm two}(x,\cdot)$ of \eqref{eq:rank-score}\;
On $\mathcal D_{\rm cal}$: compute the conformal threshold $\hat\tau$ from the scores $\hat{\mathsf s}^{\rm two}(X_i,Y_i)$ of its points with $X_i\in\cX_0$\;
\KwOut{$\hC^{\rm two}(x)=\{y:\,\hat{\mathsf s}^{\rm two}(x,y)\le\hat\tau\}$.}
\end{algorithm}

\paragraph{Computation} At a covariate $x$, $\overline Q_x$ is the running maximum, starting from zero, of the cumulative weights of the regression points with $w_i(x)\ne0$, sorted by decreasing $\hat a(Y_i\mid x)$ with tied values grouped. At every calibration and test covariate, the procedure evaluates $\hat a$ at these points, and it computes no integral over $\cY$. For local polynomial weights with bandwidth $h_{\rm reg}$, their number is $O(nh_{\rm reg}^{d_X})$. Computational details are given in \cref{app:numerics} of the supplement.

\paragraph{Relation to other calibrations} Localized calibrations of a conformity score, such as localized conformal prediction \citep{guan2023localized}, randomly localized conformal prediction (RLCP) \citep{hore2025conformal}, rectified conformal prediction (RCP) \citep{plassier2025rectifying} and the conformal sets of \citet{zhou2025conformal}, weight or regress the score of each sample point at its own covariate. Likewise, the recalibration of HPD levels in \citet[Appendix~D.2]{dey2025ladar}, which has no conformal step, estimates $Q_x$ from the HPD values of the sample points at their own covariates. Their regression targets are then only as smooth in the covariate as the fitted score. \Cref{alg:local} instead scores all regression points at the target covariate, so that the sets whose coverage it regresses do not move with the covariates of the points. For a univariate response, distributional conformal prediction \citep{chernozhukov2021distributional} uses the probability integral transform of an estimated conditional distribution function as the score. The score of \cref{alg:local} is the transform of the ranking function under the true conditional distribution. The regression estimates this transform at the target covariate, and no consistent estimate of the conditional distribution is therefore needed.

\subsection{Theory of \texorpdfstring{\cref{alg:local}}{Algorithm 1}}\label{sec:cov-theory}\label{sec:cov-local}
The theory assumes a sample split and uses a condition on the regression weights and three conditions on $\hat a$, all stated for the total sample size $n$.

\begin{assumption}[Sample split]\label{ass:split}
The observations are i.i.d.\ copies of $(X,Y)$. A procedure uses $n$ of them, split into disjoint parts: $\mathcal D_{\rm fit}$, $\mathcal D_{\rm reg}$ and $\mathcal D_{\rm cal}$ for two-stage calibration (\cref{alg:local}), and $\mathcal D_{\rm fit}$ and $\mathcal D_{\rm cal}$ for one-stage calibration. Their sizes $n_{\rm fit}$, $n_{\rm reg}$ and $n_{\rm cal}$ are each at least a fixed fraction of $n$.
\end{assumption}

The weights must give a regression with stochastic error and bias of order $\rho_n$ for $\beta$-H\"older regression functions and vary stably with $x$.

\begin{condition}{(W)}[Regression weights]\label{cond:W}
There exists a constant $C_W>0$ such that, for any $c>0$, there exists $C>0$ for which, for all sufficiently large $n$, there is an event $\Omega_X$ in the $\sigma$-field of the covariates of $\mathcal D_{\rm reg}$, of probability at least $1-Cn^{-c}$, on which the following hold for all $x,x'\in\cX_0$:
\begin{enumerate}[label=(\roman*)]
\item $w_i(x)=0$ unless $X_i\in\cX_0^{\delta_0}$, and $\sum_{i\in\mathcal D_{\rm reg}}w_i(x)^2\le C\rho_n^2/\log n$;
\item $\sum_{i\in\mathcal D_{\rm reg}}|w_i(x)-w_i(x')|\le n^{C_W}|x-x'|_\infty$;
\item $\bigl|\sum_{i\in\mathcal D_{\rm reg}}w_i(x)g(X_i)-g(x)\bigr|\le C\rho_n$ for all $g\in\mathcal H^\beta_{\cX_0^{\delta_0}}(1)$.
\end{enumerate}
\end{condition}

Local polynomial regression of degree $\ell$ in $x$ with a smooth kernel and a bandwidth of order $(\log n/n)^{1/(2\beta+d_X)}$ satisfies \cref{cond:W} (\cref{lem:loc-weights} in the supplement), and the experiments use local linear weights. The three conditions on $\hat a$ follow. The first two give the coverage bound, and the third, accuracy, gives the shape bounds.

\begin{condition}{(N)}[Null level sets]\label{cond:N}
For any $c>0$ there exists $C>0$ such that, for all sufficiently large $n$, there is an event in $\sigma(\mathcal D_{\rm fit})$ of probability at least $1-Cn^{-c}$ on which the level sets $\{y:\hat a(y\mid x)=v\}$ are Lebesgue-null for all $x\in\cX_0$ and $v\in\R$.
\end{condition}

Under \cref{cond:N}, the distribution of $\hat a(Y\mid x)$ with $Y\sim{\mathsf P}_x$ is atomless. Hence $Q_x$ is continuous, every $u\in(0,1)$ is the coverage of some candidate predictive set, and, by the probability integral transform, $Q_x(\hat a(Y\mid x))$ is uniform on $[0,1]$. \Cref{cond:N} cannot be dropped. If $\hat a(\cdot\mid x)$ is constant on a set of positive probability, $Q_x$ jumps at that value, and no candidate predictive set need have coverage close to $\gamma_\star$. For a constant $\hat a$, every candidate predictive set is $\emptyset$ or $\cY$.

The second condition bounds the complexity of the family of candidate predictive sets
	\begin{align}\label{eq:cand-family}
		\mathcal A_n:=\bigl\{\{y\in\cY:\hat a(y\mid x)\ge v\}:\ x\in\cX_0,\ v\in\R\bigr\}
	\end{align}
through its growth function $\mathfrak S_{\mathcal A_n}(m)$, the largest number of distinct sets $\{i\le m:y_i\in A\}$, $A\in\mathcal A_n$, over points $y_1,\dots,y_m\in\cY$.

\begin{condition}{(V)}[Complexity of the candidate predictive sets]\label{cond:V}
There exists a constant $C_{\rm gr}>0$ such that, for any $c>0$, there exists $C>0$ for which, for all sufficiently large $n$, there is an event in $\sigma(\mathcal D_{\rm fit})$ of probability at least $1-Cn^{-c}$ on which $\mathfrak S_{\mathcal A_n}(m)\le(nm)^{C_{\rm gr}}$ for all $m\ge1$.
\end{condition}

The candidate predictive sets vary with $x$, and \cref{cond:V} lets a single event control the regression at every covariate. It holds for the estimates of \cref{sec:cov-estimates}: networks with the usual smooth activation functions, a local polynomial estimate, whose super-level sets are given by sign conditions on polynomials in $(x,y)$, and estimates that take only polynomially many distinct forms as $x$ varies. Without \cref{cond:V}, the bounds of \cref{thm:local-cov,thm:local} hold at each $x$ outside a set of small $\Pi$-mass, because the event on which the estimate and the regression at $x$ are accurate then depends on $x$ (\cref{rem:local-each} in the supplement).

\begin{theorem}[Conditional coverage]\label{thm:local-cov}
Suppose that \cref{ass:smooth,ass:design,ass:split} hold. Let $\hC^{\rm two}$ be the output of \cref{alg:local} with a ranking function $\hat a$ and weights $w_i$, and let $\cov^{\rm two}(x)$ be its conditional coverage. Then the following hold.
\begin{enumerate}[label=(\roman*)]
\item If the weights satisfy \cref{cond:W} and $\hat a$ satisfies \cref{cond:N,cond:V}, then for every $c>0$ there exists a constant $C>0$ such that, for all sufficiently large $n$, with probability at least $1-Cn^{-c}$,
	\begin{align}
		\sup_{x\in\cX_0}\abs[0]{\cov^{\rm two}(x)-\gamma_\star}\le C\rho_n .
	\end{align}
\item For any $n$, whether or not these conditions hold, a new pair $(X,Y)$ independent of the data satisfies $\bP\{Y\in\hC^{\rm two}(X)\mid X\in\cX_0\}\ge1-\alpha$.
\end{enumerate}
\end{theorem}

\Cref{thm:local-cov} requires neither \cref{ass:margin} nor any accuracy of $\hat a$. Indeed, let $e_n:=\sup_{x\in\cX_0,\,t\in\R}|\widehat Q_x(t)-Q_x(t)|$, which is $O(\rho_n)$ with high probability under \cref{cond:W,cond:V} by the uniform regression bound of \cref{lem:loc-acc}(ii) in the supplement. The envelope preserves this accuracy, $|\hat{\mathsf s}^{\rm two}(x,y)-Q_x(\hat a(y\mid x))|\le e_n$, and, by \cref{cond:N}, $Q_x(\hat a(Y\mid x))$ is uniform on $[0,1]$ for $Y\sim{\mathsf P}_x$. Hence $|{\mathsf P}_x\{\hat{\mathsf s}^{\rm two}(x,Y)\le u\}-u|\le e_n$ for all $x\in\cX_0$ and $u\in[0,1]$, and $|\cov^{\rm two}(x)-\gamma_\star|\le2e_n$ up to the calibration error of $\hat\tau$, which is of order $(\log n/n)^{1/2}$. For instance, a Gaussian density with covariate-dependent mean and covariance satisfies \cref{cond:N,cond:V}, since its super-level sets are ellipsoids. If the conditional distribution is bimodal, these sets have the wrong shape, but their coverage error is still of order $\rho_n$, as in the example of \cref{fig:alg}. In contrast, the coverage error of one-stage calibration is governed by the variation in $x$ of the coverage residual, the density error integrated over the HPD set (\cref{sec:onestage-residual}).

The third condition requires accuracy near the level set at a deterministic rate $r_n\to0$, up to a strictly increasing transformation.

\begin{condition}{(A)}[Accuracy]\label{cond:A}
For any $c>0$ there exists $C>0$ such that, for all sufficiently large $n$, there is an event in $\sigma(\mathcal D_{\rm fit})$ of probability at least $1-Cn^{-c}$ on which, for every $x\in\cX_0$, some strictly increasing function $\phi_x:\R\to\R$ satisfies
	\begin{align*}
		\sup_{y\in\cY}\abs[0]{\phi_x(\hat a(y\mid x))-{\mathsf p}(y\mid x)}\le\frac{\varepsilon_0}8,\qquad\sup_{y\in\mathcal J_x}\abs[0]{\phi_x(\hat a(y\mid x))-{\mathsf p}(y\mid x)}\le Cr_n .
	\end{align*}
\end{condition}

\Cref{cond:A} allows a transformation $\phi_x$, such as the normalization of a density estimate, because \cref{alg:local} depends on $\hat a(\cdot\mid x)$ only through the order of its values. Moreover, \cref{cond:A} requires the rate $r_n$ only on the band $\mathcal J_x$ of \eqref{eq:Jx} around the oracle level. Under \cref{cond:A}, the super-level sets of $\hat a(\cdot\mid x)$ approximate those of ${\mathsf p}(\cdot\mid x)$, and \cref{prop:supshape} bounds the shape loss.

\begin{theorem}[Shape and volume]\label{thm:local}
Suppose that \cref{ass:smooth,ass:margin,ass:design,ass:split} hold, and let $\hC^{\rm two}$ be the output of \cref{alg:local} with weights that satisfy \cref{cond:W} and a ranking function $\hat a$ that satisfies \cref{cond:N,cond:V,cond:A}, the last with a deterministic sequence $r_n\to0$. Then for every $c>0$ there exists a constant $C>0$ such that, for all sufficiently large $n$, with probability at least $1-Cn^{-c}$,
	\begin{equation}\label{eq:local-bounds-gen}
		\begin{aligned}
		\sup_{x\in\cX_0}\cS(\hC^{\rm two}(x))&\le Cr_n^2\mathfrak L(r_n),\\
		\sup_{x\in\cX_0}\cE(\hC^{\rm two}(x))&\le C\bigl(r_n^2\mathfrak L(r_n)+\rho_n^2\bigr),\\
		\sup_{x\in\cX_0}\abs[0]{\Leb(\hC^{\rm two}(x))-\Leb(\Cstar(x))}&\le C\bigl(\rho_n+r_n^2\mathfrak L(r_n)\bigr).
		\end{aligned}
	\end{equation}
\end{theorem}

By \cref{lem:identity}, the volume bound follows from the excess-mass bound and the coverage bound of \cref{thm:local-cov}.

\subsection{Estimates that satisfy the conditions}\label{sec:cov-estimates}
We verify \cref{cond:N,cond:V} for neural network estimates, and all three conditions on the ranking function for a local polynomial estimate and for an estimate that is piecewise constant in the covariate. With the local polynomial estimate, two-stage calibration attains both lower bounds of \cref{sec:coverage}.

\subsubsection{Neural network estimates}\label{sec:cov-nn}
Both the network estimate and the local polynomial estimate of \cref{sec:cov-twostage-alg} recover ${\mathsf p}$ from a bounded regression function of $(x,y)$. Let ${\mathsf q}_0$ be the uniform density on $\cY$, let $K\ge1$ be a fixed integer, and set
	\begin{align}
		\eta(x,y):=\frac{{\mathsf p}(y\mid x)}{{\mathsf p}(y\mid x)+K{\mathsf q}_0(y)},\qquad {\mathsf p}=K{\mathsf q}_0\,\psi(\eta),\quad\psi(r):=\frac r{1-r}.
	\end{align}
For each observation $(X_i,Y_i)$ we draw $K$ reference points $Y'_{ij}\sim {\mathsf q}_0$. We label the row $(X_i,Y_i)$ by $1$ and the rows $(X_i,Y'_{ij})$ by $0$, so that $\eta$ is the probability of label $1$ given a row. This is the contrastive device of \citet{qin1998inferences} and \citet{gutmann2012noise}. The blocks $(X_i,Y_i,Y'_{i1},\dots,Y'_{iK})$ are i.i.d., and we include the reference points of the fitting observations in $\mathcal D_{\rm fit}$.

A network estimate fits the log-odds $\log\{\eta/(1-\eta)\}=\log\{{\mathsf p}/(K{\mathsf q}_0)\}$ by a feedforward network $f_\theta(x,y)$ on the labeled rows of $\mathcal D_{\rm fit}$, for instance by minimizing the cross-entropy. The density estimate $K{\mathsf q}_0\exp(f_{\hat\theta})$ is $K{\mathsf q}_0\psi(\heta)$ for the fitted probability $\heta:=1/(1+e^{-f_{\hat\theta}})$. The fitted parameter $\hat\theta$ may be any Borel function of $\mathcal D_{\rm fit}$, so that no property of the training is used. The architecture, that is, the numbers of layers and units, is fixed, and the activation function is Pfaffian: a polynomial in its argument $t$ and in functions $g_1,\dots,g_q$ of $t$ such that each derivative $g_i'$ is a polynomial in $t,g_1,\dots,g_i$. The logistic function $\sigma(t):=1/(1+e^{-t})$, with $\sigma'=\sigma-\sigma^2$, is Pfaffian, and so are $\tanh(t)=2\sigma(2t)-1$, the SiLU function $t\sigma(t)$ and the softplus function $\log(1+e^t)$, whose derivative is $\sigma$.

\Cref{cond:N} fails at $x$ if $f_{\hat\theta}(x,\cdot)$ is constant, and no property of the training excludes this. We therefore break ties by a small polynomial. Let $z(y)\in[-1,1]^{d_Y}$ be the image of $y$ under the affine map of $\cY$ onto $[-1,1]^{d_Y}$, and, for an integer $J\ge1$, let $\mathrm{pol}_J(y):=\sum_{k=1}^{d_Y}z_k(y)^J$. The network ranking function is
	\begin{align*}
		\hat a^{\rm nn}(y\mid x):=K{\mathsf q}_0\exp\{f_{\hat\theta}(x,y)\}+n^{-2}\mathrm{pol}_J(y) .
	\end{align*}
The perturbation is a deterministic tie-breaking rule: where the density estimate is flat, $\hat a^{\rm nn}$ orders the responses by $\mathrm{pol}_J$. Since $|\mathrm{pol}_J|\le d_Y$ on $\cY$, the perturbation changes the density estimate by at most $d_Yn^{-2}$.

\begin{proposition}[Network estimates]\label{prop:nn-conditions}
Let $J\ge1$. Then the following hold.
\begin{enumerate}[label=(\roman*)]
\item The ranking function $\hat a^{\rm nn}$ satisfies \cref{cond:N,cond:V}.
\item If \cref{ass:smooth,ass:design,ass:split} hold, then, with weights that satisfy \cref{cond:W}, the output of \cref{alg:local} with the ranking function $\hat a^{\rm nn}$ satisfies the bounds of \cref{thm:local-cov}.
\end{enumerate}
\end{proposition}

By (ii), the coverage error of two-stage calibration with a network estimate is $O(\rho_n)$ uniformly in $x$.

\subsubsection{The local polynomial estimate}\label{sec:cov-twostage-alg}
The shape bounds of \cref{thm:local} require \cref{cond:A}. We obtain this condition at the minimax rate by fitting the regression function $\eta$ of \cref{sec:cov-nn} with local polynomials in $(x,y)$. Local polynomial estimates of a conditional density go back to \citet{fan1996estimation}, who apply local polynomial regression in $x$ to a kernel of $Y_i-y$. Our fit is instead a linear smoother of the binary labels of the contrastive rows, with an explicit sup-norm rate and bias.

For a bandwidth $b>0$ and each $z=(x,y)$, we fit by weighted least squares a polynomial of degree $\ell$ in $(z'-z)/b$, where $z'$ ranges over the rows, with weights $\mathcal K((z'-z)/b)$ \citep{fan1996local}. Here $\mathcal K(s):=\prod_{j=1}^d\mathcal K_1(s_j)$ is the product kernel of a univariate kernel $\mathcal K_1$ that is smooth in the sense of \cref{app:tools} of the supplement. We clip the fitted intercept to $[0,\bar\eta']$ for a fixed $\bar\eta'\in({\mathsf p}_{\max}/({\mathsf p}_{\max}+K{\mathsf q}_0),1)$, denote the result by $\heta$, and, with $\mathrm{pol}_J$ as in \cref{sec:cov-nn}, set
	\begin{align*}
		\hat a^{\rm lp}_b(y\mid x):=K{\mathsf q}_0\psi(\heta(x,y))+n^{-2}\mathrm{pol}_J(y) .
	\end{align*}
We write $\tilde a^{\rm lp}_b:=K{\mathsf q}_0\psi(\heta)$ for the estimate before the perturbation. In this section, we use the density-optimal bandwidth $h\asymp(\log n/n)^{1/(2\beta+d)}$. The remaining conventions are given in \cref{app:lp-estimator} of the supplement.

The level sets of $\tilde a^{\rm lp}_h(\cdot\mid x)$ may have positive measure, for instance where the fitted intercept is clipped, and, as for the network, the perturbation breaks these ties. If $\mathcal K_1$ is a polynomial on $[-1,1]$, such as the triweight kernel $(1-s^2)^3$, then $\tilde a^{\rm lp}_h$ is piecewise a ratio of polynomials in $(x,y)$, which gives \cref{cond:V}. The perturbation gives \cref{cond:N} when $J$ is an even integer larger than a constant $D_0$ that depends only on $\beta$, $d$ and the degree of $\mathcal K_1$. The perturbation, of size at most $d_Yn^{-2}$, is far below the sup-norm error of $\tilde a^{\rm lp}_h$, which is of order $\zeta_n$, and therefore does not change the rates.

\begin{proposition}[The local polynomial estimate]\label{prop:lp-conditions}
Suppose that \cref{ass:smooth,ass:design,ass:split} hold, that $\mathcal K_1$ is a polynomial on $[-1,1]$ and that $J>D_0$ is an even integer. Then the following hold.
\begin{enumerate}[label=(\roman*)]
\item The ranking function $\hat a^{\rm lp}_h$ satisfies \cref{cond:N,cond:V} and, if \cref{ass:margin} holds, \cref{cond:A} with $r_n=\zeta_n$.
\item With weights that satisfy \cref{cond:W}, the output of \cref{alg:local} with the ranking function $\hat a^{\rm lp}_h$ satisfies the bounds of \cref{thm:local-cov} and, if in addition \cref{ass:margin} holds, those of \cref{thm:local} with $r_n=\zeta_n$.
\end{enumerate}
\end{proposition}

By \cref{prop:lp-conditions}(ii), two-stage calibration with $\hat a^{\rm lp}_h$ attains the lower bounds of \cref{thm:lower} in every dimension, exactly for the uniform coverage error and up to logarithmic factors otherwise.

\subsubsection{Estimates that are piecewise constant in the covariate}\label{sec:cov-partition}
Let $h_X=h_{X,n}$ and $h_Y=h_{Y,n}$ be deterministic sequences. For $x$ in a cube $\mathcal Q$ of a regular grid of side $h_X$, let $\hp^{\rm part}(\cdot\mid x)$ be the Gaussian kernel estimate with bandwidth $h_Y$ of the responses $Y_i$ of $\mathcal D_{\rm fit}$ with $X_i\in \mathcal Q$, divided at each $y$ by the mass that the kernel centered at $y$ assigns to $\cY$, and normalized on $\cY$. If no $X_i$ lies in $\mathcal Q$, let $\hp^{\rm part}(\cdot\mid x)$ be the uniform density on $\cY$. Density estimates that partition the covariate space by a tree, a forest or nearest neighbors share this step structure in $x$. Let
	\begin{align}
		r_n:=h_X+h_Y+\Bigl(\frac{\log n}{nh_X^{d_X}h_Y^{d_Y}}\Bigr)^{1/2}.
	\end{align}

\begin{proposition}[Estimates that are discontinuous in the covariate]\label{prop:partition}
Suppose that \cref{ass:smooth,ass:design,ass:split} hold and that $r_n\to0$. Then the following hold.
\begin{enumerate}[label=(\roman*)]
\item The estimate $\hp^{\rm part}$ satisfies \cref{cond:N,cond:V} and, if \cref{ass:margin} holds, \cref{cond:A} with the sequence $r_n$.
\item With weights that satisfy \cref{cond:W}, the output of \cref{alg:local} with the ranking function $\hp^{\rm part}$ satisfies the bounds of \cref{thm:local-cov} and, if in addition \cref{ass:margin} holds, those of \cref{thm:local} with this sequence $r_n$.
\end{enumerate}
\end{proposition}

By (ii), the coverage error of two-stage calibration with $\hp^{\rm part}$ is $O(\rho_n)$ uniformly in $x$, although the estimate jumps in $x$. For $\beta=1$ and $h_X\asymp h_Y\asymp(\log n/n)^{1/(2+d)}$, $r_n$ is of order $\zeta_n$, and two-stage calibration with this estimate attains both lower bounds of \cref{thm:lower}, up to logarithmic factors for the excess mass. For $\beta>1$, the density error $r_n$ of this estimate exceeds $\zeta_n$, and the excess mass therefore need not attain the rate $\zeta_n^2$. The coverage error of one-stage calibration, in contrast, can jump where the estimate jumps (\cref{prop:partition-jump}).

\section{One-stage calibration and existing guarantees}\label{sec:cov-onestage}
A single conformal threshold corrects the coverage only on average over $x$. We show that the coverage error of one-stage calibration is therefore, to first order, the variation in $x$ of the density error integrated over the oracle set, and that, for the local polynomial plug-in, this variation is of the order of the bias. For estimates that are piecewise constant in the covariate, the coverage can jump where the estimate jumps. We also compare our guarantees with existing ones.

\subsection{The coverage residual}\label{sec:onestage-residual}
\begin{figure}[t]
\centering
\includegraphics[width=\textwidth]{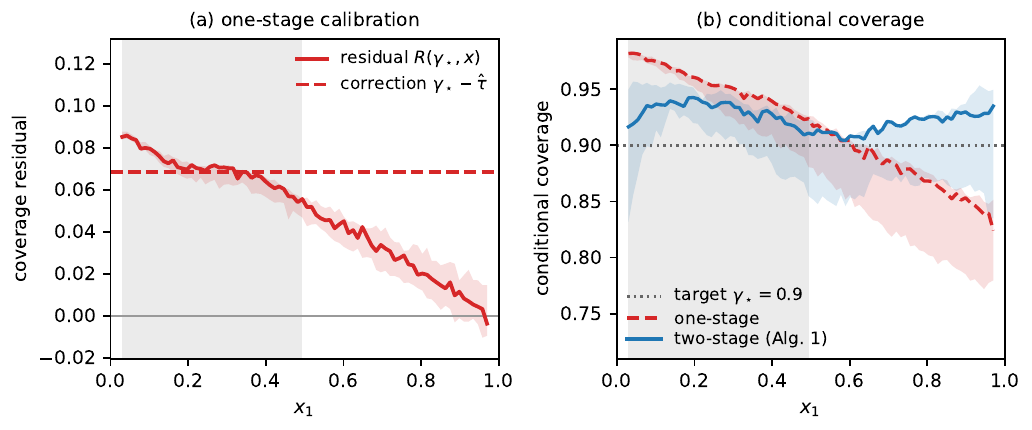}
\caption{Coverage residual and conditional coverage along the segment $\{(x_1,0.5,0.5):0<x_1<1\}$ for the local polynomial estimate in the multimodal model with $d_Y=2$, in the setting of \cref{sec:simulations}. (a) True coverage residual $R(\gamma_\star,x)$ of the one-stage fit (solid) and the correction $\gamma_\star-\hat\tau$ of its single threshold (dashed). (b) Conditional coverage of one-stage calibration (dashed) and of \cref{alg:local} (solid). Curves show one replication, and bands the pointwise range over replications. Gray: hard region.}\label{fig:residual}
\end{figure}

One-stage calibration applies the conformal threshold directly to the HPD values of an estimate fitted on $\mathcal D_{\rm fit}$. For fixed fitting data and $u\in[0,1]$, let
	\begin{align}\label{eq:F-R}
		F(u,x):={\mathsf P}_x\{\hat{\mathsf s}^{\rm hpd}(x,Y)\le u\},\qquad R(u,x):=F(u,x)-u,
	\end{align}
where the probability is over $Y\sim{\mathsf P}_x$ and the tie-breaking variable at $x$. On average over the tie-breaking variable, the HPD set $\{y:\hat{\mathsf s}^{\rm hpd}(x,y)\le u\}$ has estimated mass $u$ and coverage $F(u,x)$, and the \emph{coverage residual} $R(u,x)$ is the error of this estimated mass. One-stage calibration returns at every $x$ the HPD set with estimated mass $\hat\tau$, whose coverage, on average over the tie-breaking variable, is $\hat\tau+R(\hat\tau,x)$. Since $\hat\tau$ does not depend on $x$, the coverage error of one-stage calibration is governed by the variation of $R$ in $x$. When $\hp$ is close to ${\mathsf p}$ and its level sets near $t_\star(x)$ are Lebesgue-null, \cref{lem:hpd-sets} in the supplement shows that, for $|u-\gamma_\star|\le c_0$, the residual is the density error integrated over the super-level set $A_x(u)$ of $\hp(\cdot\mid x)$ with $\hP_x$-mass $u$,
	\begin{align}\label{eq:R-int}
		R(u,x)=\int_{A_x(u)}\bigl({\mathsf p}-\hp\bigr)(y\mid x)\,\d y .
	\end{align}
These definitions apply to every estimate, and the variation of $R$ in $x$ depends on the estimate. For a smoothing estimate at bandwidth $b$, such as the local polynomial estimate, the density error is the sum of a bias term and a noise term. The noise at a point averages over about $nb^d$ observations, whereas its integral over a set of fixed volume averages over about $nb^{d_X}$ observations near $x$. As for smooth integral functionals of a density \citep{birge1995estimation,laurent1996efficient}, the integral therefore averages out the noise as in a regression on $x$ alone, but not the bias, whose integral over the set varies with $x$ (\cref{fig:residual} and~\cref{cor:lp-onestage}).

\subsection{Coverage error and the integrated bias}\label{sec:cov-onestage-thm}
Let $\hC^{\rm one}$ be the output of one-stage calibration with an estimate $\hp$ fitted on $\mathcal D_{\rm fit}$ and exact HPD values, and let $\cov^{\rm one}(x)$ be its conditional coverage. Let
	\begin{align*}
		R^\star(x):=\int_{\Cstar(x)}({\mathsf p}-\hp)(y\mid x)\,\d y,\qquad\bar R^\star:=\frac1{\Pi(\cX_0)}\int_{\cX_0}R^\star\,\d\Pi
	\end{align*}
be the coverage residual \eqref{eq:R-int} at $u=\gamma_\star$ with the HPD set $A_x(\gamma_\star)$ of the estimate replaced by the oracle set, and its average over the analysis set. Since ${\mathsf P}_x(\Cstar(x))=\gamma_\star$, we have $R^\star(x)=\gamma_\star-\hP_x(\Cstar(x))$. \Cref{sec:onestage-residual} shows that the coverage error of one-stage calibration is governed by the variation of the residual in $x$. The next theorem makes this precise: to first order, the coverage error is the variation of $R^\star$ in $x$.

\begin{theorem}[One-stage coverage]\label{thm:onestage}
Suppose that \cref{ass:smooth,ass:design,ass:split} hold and that \cref{ass:margin} holds with $\mathfrak L\equiv1$. Let $r_n\to0$ be a deterministic sequence with $r_n\ge(\log n/n)^{1/2}$, and let $\Omega_r\in\sigma(\mathcal D_{\rm fit})$ be an event on which $\sup_{\cX_0\times\cY}|\hp-{\mathsf p}|\le r_n$ and, for all $x\in\cX_0$ and all $v$ with $|v-t_\star(x)|\le\varepsilon_0/2$, the level set $\{y:\hp(y\mid x)=v\}$ is Lebesgue-null. Then for every $c>0$ there exists a constant $C>0$ such that, for all sufficiently large $n$, with probability at least $\bP(\Omega_r)-Cn^{-c}$,
	\begin{align*}
		&\sup_{x\in\cX_0}\abs[1]{\cov^{\rm one}(x)-\gamma_\star-\{R^\star(x)-\bar R^\star\}}\le C\Bigl\{r_n^2+\Bigl(\frac{\log n}n\Bigr)^{1/2}\Bigr\},\\
		&\sup_{x\in\cX_0}\abs[0]{\cov^{\rm one}(x)-\gamma_\star}\le Cr_n,\qquad\sup_{x\in\cX_0}\bigl\{\cS(\hC^{\rm one}(x))+\cE(\hC^{\rm one}(x))\bigr\}\le Cr_n^2 .
	\end{align*}
\end{theorem}

The first bound shows that a single threshold removes the average $\bar R^\star$ of the residual but not its variation in $x$. By the other two bounds, one-stage calibration has coverage error at most of the order of the density error, and its sets have the shape loss and the excess mass of a super-level set of an estimate with this error. The integrated coverage error attains this order in expectation when the mean of $R^\star$ varies in $x$ at this order. Besides the conditions of \cref{thm:onestage}, the following lower bound requires only a condition on this mean.

\begin{theorem}[Coverage lower bound from the mean error]\label{thm:onestage-lower}
Suppose that the conditions of \cref{thm:onestage} hold with $\bP(\Omega_r)\ge1-\delta_n$ for a deterministic sequence $\delta_n$. Suppose also that there exist a deterministic sequence $a_n\to0$ with $r_n^2+(\log n/n)^{1/2}+\delta_n=o(a_n)$ and a bounded function $\mathfrak g$ on $\cX_0$, not depending on $n$, that is not almost everywhere constant, such that
	\begin{align*}
		\E\sbr[0]{R^\star(x)}=a_n\mathfrak g(x)+o(a_n)\qquad\text{uniformly in }x\in\cX_0 .
	\end{align*}
Then there exists a constant $c>0$, depending also on $\mathfrak g$, such that, for all sufficiently large $n$,
	\begin{align*}
		\E\sbr{\int_{\cX_0}\abs[0]{\cov^{\rm one}(x)-\gamma_\star}\,\d x}\ge c\,a_n .
	\end{align*}
\end{theorem}

For a smoothing estimate, $\E\sbr[0]{R^\star(x)}$ is, up to second-order terms, minus the bias integrated over the oracle set. For a deterministic bandwidth $b$, let $\hp^{\rm lp}_b$ be the estimate $\hat a^{\rm lp}_b$ of \cref{sec:cov-twostage-alg} normalized to a density, and let $\hC^{\rm one}_b$ and $\cov^{\rm one}_b(x)$ be the output and the conditional coverage of one-stage calibration with $\hp^{\rm lp}_b$. Under the conditions of the following corollary and for $b\ge h$, the mean residual of $\hp^{\rm lp}_b$ is $-b^\beta\int_{\Cstar(x)}B_\beta(x,y)\,\d y+o(b^\beta)$, where $B_\beta$ is the leading term of the bias, defined in \eqref{eq:leading-bias} of the supplement.

\begin{corollary}[The local polynomial estimate]\label{cor:lp-onestage}
Let $\beta\ge2$ be an integer. Suppose that \cref{ass:smooth,ass:design,ass:split} hold, that \cref{ass:margin} holds with $\mathfrak L\equiv1$, and that the regularity conditions of \cref{app:fixed-law} of the supplement hold. If $x\mapsto\int_{\Cstar(x)}B_\beta(x,y)\,\d y$ is not almost everywhere constant on $\cX_0$, then there exists a constant $c>0$ such that, for every deterministic bandwidth $b\to0$ with $b\ge h$ and all sufficiently large $n$,
	\begin{align*}
		\E\sbr{\int_{\cX_0}|\cov^{\rm one}_b(x)-\gamma_\star|\,\d x}\ge c\,b^\beta .
	\end{align*}
\end{corollary}

Here, as in \cref{thm:onestage-lower} and unlike the other constants, $c$ and the minimal sample size depend on the distribution. At $b=h$, the proof of \cref{cor:lp-onestage} shows that $\hp^{\rm lp}_h$ satisfies the conditions of \cref{thm:onestage} with $r_n$ a multiple of $\zeta_n$, so that one-stage calibration has, with high probability, coverage error $O(\zeta_n)$ and excess mass $O(\zeta_n^2)$. Since $h^\beta\asymp\zeta_n$, its coverage error at the density-optimal bandwidth is of order $\zeta_n$: at most $C\zeta_n$ uniformly in $x$ with high probability, and at least $c\zeta_n$ in expected integrated absolute value. This order exceeds $\rho_n$ by a polynomial factor, and no larger bandwidth improves it. With a symmetric kernel, the leading bias vanishes for odd $\beta$, and the condition on $B_\beta$ then requires $\beta$ to be even. For $\beta=2$, \cref{lem:gb-curv} in the supplement reduces it to a condition on curvatures in $x$, which a heteroscedastic Gaussian distribution satisfies (see \cref{app:gauss-bias} of the supplement for the computation).

\begin{remark}[Undersmoothing]\label{rem:undersmoothing}
A smaller bandwidth reduces the bias that yields the lower bound of \cref{cor:lp-onestage}, but it raises the noise of the estimate. For a deterministic bandwidth $b\to0$ with $nb^d/\log n\to\infty$, the estimate $\hp^{\rm lp}_b$ has sup-norm error of order $b^\beta+\{\log n/(nb^d)\}^{1/2}$ by \cref{lem:supnorm-b} in the supplement, and \cref{thm:onestage} bounds the excess mass of one-stage calibration with this estimate by the square of this error. For $b^\beta\le\rho_n$, this bound exceeds $\zeta_n^2$ by a polynomial factor. Heuristically, the noise misranks responses near the level set at the scale of its standard deviation $(nb^d)^{-1/2}$. For $d_Y\ge2$, the noise near two points of the level set farther apart than $2b$ is nearly independent, and a single threshold therefore cannot correct the misranking at both. Two-stage calibration instead keeps the density-optimal bandwidth for the shape and corrects the coverage by the regression.
\end{remark}

\subsection{Estimates that are piecewise constant in the covariate}\label{sec:onestage-partition}
For the estimate of \cref{sec:cov-partition}, which is piecewise constant in $x$, the residual $R^\star$ can jump where the estimate jumps. A single threshold cannot remove such jumps, and, when they exceed the remainder in the first bound of \cref{thm:onestage}, they appear in the coverage of one-stage calibration. The following proposition gives conditions under which they are of order $h_X$.

\begin{proposition}[One-stage calibration of estimates that are discontinuous in the covariate]\label{prop:partition-jump}
Let $\hp^{\rm part}$ and $r_n$ be as in \cref{sec:cov-partition}, and let $\hC^{\rm one}$ and $\cov^{\rm one}(x)$ be the output and the conditional coverage of one-stage calibration with $\hp^{\rm part}$ and exact HPD values. Suppose that \cref{ass:smooth,ass:design,ass:split} hold, that \cref{ass:margin} holds with $\mathfrak L\equiv1$, that $r_n\to0$ and $h_Y+\{\log n/(nh_X^{d_X}h_Y^{d_Y})\}^{1/2}=o(h_X)$, and that there exist constants $c_V>0$, $\varsigma\in\{-1,1\}$ and a closed cube $\mathcal B\subset\cX_0$, not depending on $n$, such that $\varsigma V_x(\gamma_\star)-c_Vx_1$ is nondecreasing in $x_1$ on $\mathcal B$. Then $h_X/\rho_n\to\infty$, and there exists a constant $c>0$ such that, with probability tending to one,
	\begin{align}
		\sup_{x\in\cX_0}\abs[0]{\cov^{\rm one}(x)-\gamma_\star}\ge c\,h_X .
	\end{align}
\end{proposition}

At a face shared by two cubes, the coverage of one-stage calibration jumps by $t_\star$ times the difference of the averages of $V_x(\gamma_\star)$ over the two cubes, up to $O(h_Y+\{\log n/(nh_X^{d_X}h_Y^{d_Y})\}^{1/2}+h_X^2)$. Under the conditions of \cref{prop:partition-jump}, this jump is of order $h_X$, which exceeds $\rho_n$. Two-stage calibration, in contrast, has coverage error $O(\rho_n)$ for the same estimate (\cref{prop:partition}). The bandwidth condition holds, for instance, for $h_Y\asymp(\log n/n)^{1/(2+d)}$ and $h_X\asymp h_Y\log n$.

\subsection{Existing guarantees}\label{sec:cov-existing}
Existing guarantees differ from ours mainly in what they average over. \Cref{tab:literature} in the supplement states them in our notation. Our bounds concern $\cov(x)$, conditional on all data, at each $x$. The guarantees of CP$^2$ \citep{plassier2024conditionally} and RCP average over the calibration sample, and those of HPD-split \citep{izbicki2022cd} and \citet{gao2025volume} over all data. Those of CP$^2$, HPD-split and \citet{gao2025volume} hold at a random test covariate or outside a set of small $\Pi$-mass. Our bounds imply bounds of the same order on these averages. Among density-based sets with one-stage calibration, HPD-split has consistency results, and the conditional coverage of distributional conformal prediction tends to $1-\alpha$ when the estimate of the conditional distribution function is consistent. The bounds for CHCDS \citep{sampson2024flexible} and CP$^2$ depend on the sup-norm and the total-variation accuracy of the density estimate and therefore hold at the density rate.

For $d_Y=1$, the local conformal sets of \citet{lei2014distribution} are obtained by running conformal prediction separately on cubes in $x$. They have coverage at least $1-\alpha$ on each cube, whereas our finite-sample guarantee is marginal over $\cX_0$. Their symmetric-difference rate, slower than $\zeta_n$ for $\beta>1$, is the only one among conformal level-set methods with a matching lower bound. \Cref{thm:lower}(i) is the counterpart of this lower bound in every dimension. CD-split \citep{izbicki2022cd} calibrates on cells built from the estimated density, with asymptotic conditional validity and no rate.

The guarantees of the localized calibrations of \cref{sec:cov-procedure} need assumptions on the fitted score. Those of localized conformal prediction, RCP and \citet{zhou2025conformal} assume that the conditional distribution of the score is smooth in the covariate. For RLCP and a univariate response, \citet{conrad2026beyond} give two bounds that hold simultaneously over a ball around a random localization center. For a fixed score whose conditional distribution function is H\"older continuous in $x$, with exponent at most one, the coverage error is of the order of the regression rate. For a learned score close in sup-norm to a fixed score whose conditional $(1-\alpha)$-quantile does not depend on $x$, as for conformalized quantile regression \citep{romano2019conformalized}, the coverage error is at most of the order of a calibration term plus this sup-norm distance. For the HPD value of a density estimate, the assumption of the first bound depends on the estimate and fails, for instance, for estimates that are discontinuous in the covariate. The second bound depends on the accuracy of the density estimate, as does the bound of \cref{thm:onestage} for one-stage calibration. \Cref{thm:local-cov} requires neither smoothness in the covariate nor accuracy of the ranking function, only \cref{cond:N,cond:V}, and its rate $\rho_n$ does not depend on $d_Y$.

Conditional calibration \citep{gibbs2023conformal} guarantees coverage at least $1-\alpha$ on average over the data under every nonnegative reweighting of the covariate distribution in a finite-dimensional class. \citet{duchi2025sample} bound the coverage given the data on groups of covariates, and minimax optimization predictive inference (MOPI) \citep{bao2026shape} bounds the squared coverage error integrated over the covariates, relative to the best set of a given class. None of them gives a rate at each $x$ under smoothness.

\section{Shape loss under $L^q$ accuracy}\label{sec:shape}
Since a sup-norm error $r$ gives shape loss $O(r^2)$, an estimate attaining the density rate has shape loss of the minimax excess-mass order $\zeta_n^2$. We show that accuracy in $L^q$, even of minimax order, cannot replace sup-norm accuracy.

Throughout this section, \cref{ass:margin} holds with $\mathfrak L\equiv1$. We fix $x$ and let $\hC(x)=\{y:\hat a(y\mid x)\ge\hat t\}$ for a ranking function $\hat a(\cdot\mid x)$ on $\cY$ and a level $\hat t$. For HPD-split, if the level sets of $\hp(\cdot\mid x)$ are Lebesgue-null, we may take $\hat a(\cdot\mid x)=\hp(\cdot\mid x)$. Since a strictly increasing transformation $\phi$ of $\hat a$ leaves its super-level sets unchanged, we measure accuracy through $\phi(\hat a(\cdot\mid x))$.

\begin{theorem}[From $L^q$ accuracy to shape]\label{thm:Lq}
Suppose that \cref{ass:smooth} holds. Then there exists a constant $C>0$, not depending on $x$, $\hat a$ or $\hat t$, such that the following holds for any $q\in[1,\infty]$. Let $\hC(x)$ have coverage in the window, let $\phi:\R\to\R$ be any strictly increasing function with $\sup_{y\in\cY}|\phi(\hat a(y\mid x))-{\mathsf p}(y\mid x)|\le\varepsilon_0/8$, and let
	\begin{align}
		r^{(q)}:=\norm[0]{\phi(\hat a(\cdot\mid x))-{\mathsf p}(\cdot\mid x)}_{L^q(\mathcal J_x)}
	\end{align}
be the $L^q$ error of $\phi(\hat a(\cdot\mid x))$ on the band $\mathcal J_x$ of \eqref{eq:Jx}. Then
	\begin{align}
		\cS(\hC(x))\le C\,(r^{(q)})^{2q/(q+1)} .
	\end{align}
\end{theorem}

Since accuracy matters only near the level, the theorem measures $L^q$ accuracy on the band $\mathcal J_x$ alone and requires merely the coarse sup-norm bound $\varepsilon_0/8$ elsewhere. In the bimodal example of \cref{sec:setup-shape}, let an estimate equal ${\mathsf p}(\cdot\mid x)$ on $\mathcal J_x$ and, elsewhere, be arbitrary but on the same side of the band as ${\mathsf p}(\cdot\mid x)$. Its shape loss is zero however large its $L^2$ error is.

For $q=1,2,\infty$, the exponents are $1$, $4/3$ and $2$. Hence, for an error of order $\zeta_n$ on $\mathcal J_x$, the shape loss is at most of order $\zeta_n$ under total-variation accuracy and of order $\zeta_n^{4/3}$ under $L^2$ accuracy. The exponent $4/3$ is that of plug-in classification under a margin condition with exponent one \citep{audibert2007fast}. Both rates are slower than the rate $\zeta_n^2$ under sup-norm accuracy. Neither exponent can be improved, even at the level that calibration selects: the estimate of \cref{thm:Lq-lower} in the supplement raises and lowers ${\mathsf p}(\cdot\mid x)$ by $2\delta$ on two sets of measure of order $\delta$ just below and just above the level, respectively, and has $L^q$ error of order $\delta^{1+1/q}$ and shape loss of order $\delta^2$ at every level within $\delta$ of $t_\star(x)$.

\begin{figure}[t]
\centering
\includegraphics[width=\textwidth]{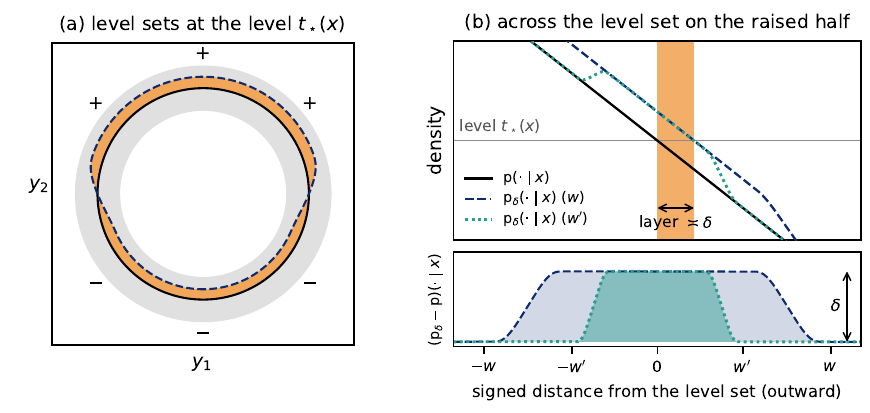}
\caption{The conditional density ${\mathsf p}_\delta$ of \cref{thm:sep-smooth}, drawn for a bivariate radial density. (a) Level sets of the density (solid) and of ${\mathsf p}_\delta(\cdot\mid x)$ (dashed) at the oracle level $t_\star(x)$. A smooth perturbation of height $\delta$, supported in the gray shell of half-width $w$ around the level set, raises the density on the upper half ($+$) and lowers it on the lower half ($-$). Shaded: the misplaced points, between the two level sets. (b) Top: the density across the level set on the raised half, perturbed on shells of half-widths $w$ and $w'<w$. Bottom: the two errors ${\mathsf p}_\delta(\cdot\mid x)-{\mathsf p}(\cdot\mid x)$. Both perturbations move the boundary by the same distance of order $\delta$. Hence the misplaced layer (shaded) and the shape loss are the same, whereas the $L^q$ norms of the error shrink with the width of the shell.}\label{fig:sep-smooth}
\end{figure}

The gap persists for smooth estimates (\cref{fig:sep-smooth}). We raise the density smoothly by $\delta$ on a shell of half-width $w$ around one half of the level set and lower it by $\delta$ on the other half. The result is still a density, and the coverage of its super-level sets is nearly unchanged. On the raised half, the super-level set moves out by a distance of order $\delta$, and on the lowered half it moves in by as much. The misplaced layer has volume of order $\delta$, and its points have density at distance of order $\delta$ from the level. Hence the shape loss is of order $\delta^2$, as for a sup-norm error $\delta$. In contrast, since the shell has volume of order $w$, the $L^2$ and total-variation norms of the error are of order $\delta w^{1/2}$ and $\delta w$. Within the $\beta$-H\"older class, we can choose $w\asymp\delta^{1/\beta}$. The shape loss $\delta^2$ is then much larger than the squares $\delta^2w$ and $\delta^2w^2$ of these norms.

\begin{theorem}[Separation within the H\"older class]\label{thm:sep-smooth}
Let $\beta>1$, let \cref{ass:constants} hold, and let ${\mathsf p}(\cdot\mid x):=f_0$ for every $x$, where $f_0$ is the base density of the lower bounds, defined in \eqref{eq:lb-f0} of the supplement. Then, for all sufficiently small $\delta>0$, there exists a conditional density ${\mathsf p}_\delta$, used in place of an estimate, with the following properties. Here ${\mathsf P}_{\delta,x}$ denotes the distribution with density ${\mathsf p}_\delta(\cdot\mid x)$, and the constants in $\asymp$ do not depend on $\delta$, $x$ or $s$.
\begin{enumerate}[label=(\roman*)]
\item The density ${\mathsf p}_\delta(\cdot\mid x)$ is the same for every $x$, and $(x,y)\mapsto{\mathsf p}_\delta(y\mid x)$ belongs to $\mathcal H^\beta(L)$.
\item The $L^2$ and total-variation distances satisfy
	\begin{align}
		\norm{{\mathsf p}_\delta(\cdot\mid x)-{\mathsf p}(\cdot\mid x)}_{L^2}\asymp\delta^{1+1/(2\beta)},\qquad\TV({\mathsf P}_{\delta,x},{\mathsf P}_x)\asymp\delta^{1+1/\beta}.
	\end{align}
\item For any level $s$ with $|s-t_\star(x)|\le\delta/4$,
	\begin{align}
		\cS\bigl(\{y:{\mathsf p}_\delta(y\mid x)\ge s\}\bigr)\asymp\cE\bigl(\{y:{\mathsf p}_\delta(y\mid x)\ge s\}\bigr)\asymp\delta^2 .
	\end{align}
\end{enumerate}
\end{theorem}

Consequently, with $\eps:=\delta^{1+1/(2\beta)}$, an estimate whose $L^2$ and total-variation errors are at most of order $\eps$ can have shape loss of order $\eps^{4\beta/(2\beta+1)}\gg\eps^2$ at every level near $t_\star(x)$.

\section{Numerical experiments}\label{sec:experiments}\label{sec:simulations}
We compare two-stage calibration with one-stage calibration and other conformal methods in conditional coverage and in shape, on three simulated models.

\paragraph{Models} In all models, $d_X=3$, $d_Y\in\{2,5\}$, $\cY=[-2.5,2.5]^{d_Y}$, $\alpha=0.1$, and $X\sim U[0,1]^3$. We write $\mathtt N(\mathrm{ctr},\mathrm{sd}^2)$ for the normal distribution on $\R^{d_Y}$ with mean $\mathrm{ctr}$ and covariance $\mathrm{sd}^2I_{d_Y}$. Each model has a hard region of covariate values, and its complement is the easy region. \Cref{app:numerics} of the supplement gives the exact specifications.
\begin{itemize}
\item \emph{Gaussian}: $\mathtt N(\mathrm{ctr}(x),\mathrm{sd}(x)^2)$, where the first coordinate of $\mathrm{ctr}(x)$ and the scale $\mathrm{sd}(x)\in[0.3,0.6]$ vary smoothly with $x$. The hard region is $\{x\in[0,1]^3:x_1>0.5\}$.
\item \emph{multimodal}: a mixture of four components $\mathtt N(\mathrm{ctr}_k,\mathrm{sd}_k^2)$ with fixed centers and scales $\mathrm{sd}_k\in\{0.02,0.059,0.17,0.5\}$ and with weights that depend on $x$. The hard region is the set of covariates at which the two components with the smallest scales carry more than half of the mass.
\item \emph{manifold}: $\{1-\mathrm{wt}(x)\}\,\mathtt N(0,0.5^2)+\mathrm{wt}(x)\,\mathsf T$, where $\mathsf T$ concentrates within about $0.02$ of a circle of radius $0.8$ and the weight $\mathrm{wt}(x)$ increases with $x_1$. The hard region is $\{x\in[0,1]^3:x_1>0.72\}$, where $\mathsf T$ carries most of the mass.
\end{itemize}
In the multimodal and manifold models, part of the mass lies at scale $0.02$. Hence the conditional density is hard to estimate, especially for $d_Y=5$. All three models satisfy \cref{ass:smooth}, and the Gaussian model also satisfies \cref{ass:margin}. \Cref{ass:margin} fails in the manifold model with $d_Y=2$: at $x_1\approx0.54$, a circle of critical points of ${\mathsf p}(\cdot\mid x)$ lies at the level $t_\star(x)$. We use $\cX_0=[0,1]^3$, which reaches the boundary of the support of $X$, whereas \cref{ass:design} requires a design density bounded below on a neighborhood of $\cX_0$.

\paragraph{Setup} Each replication has $1000$ calibration points and $4000$ points for fitting and for the auxiliary step. Methods with an auxiliary step fit the density on $3000$ of these $4000$ points. We compute the conditional coverage and the shape loss $\cS^\circ$ by Monte Carlo at ten covariate values in the hard region and ten in the easy region, and average the coverage errors over the two regions with equal weight. Results are means over twenty replications, with all methods run on the same data in each replication. The code and the results of all runs are available at \url{https://github.com/ilsangohn/conformal-hpd}.

\paragraph{Methods} We compare \cref{alg:local} with one-stage calibration, that is, HPD-split \citep{izbicki2022cd}, for three density estimates: a contrastive multilayer perceptron (NN), which, with a single uniform reference, is the estimator CINDES of \citet{dai2025cindes}; the contrastive local linear estimate (LP); and the partition estimate of \cref{prop:partition}, which is piecewise constant in $x$. Two variants of two-stage calibration use the same estimate and auxiliary points as \cref{alg:local}: two-stage (own covariate) regresses the HPD values of the auxiliary points at their own covariates, as in \citet[Appendix~D.2]{dey2025ladar}, and two-stage (quantile) applies RCP to the HPD score. The baselines are CHCDS \citep{sampson2024flexible}, a local ellipsoid, probabilistic conformal prediction (PCP) \citep{wang2022probabilistic}, RCP \citep{plassier2025rectifying} applied to the ellipsoid and PCP scores, and optimal-transport conformal prediction (OT-CP) \citep{thurin2025optimal,klein2025multivariate}, the last only for $d_Y=2$. One-stage calibration with the true density serves as a benchmark.

\begin{figure}[!t]
\centering
\includegraphics[width=\textwidth]{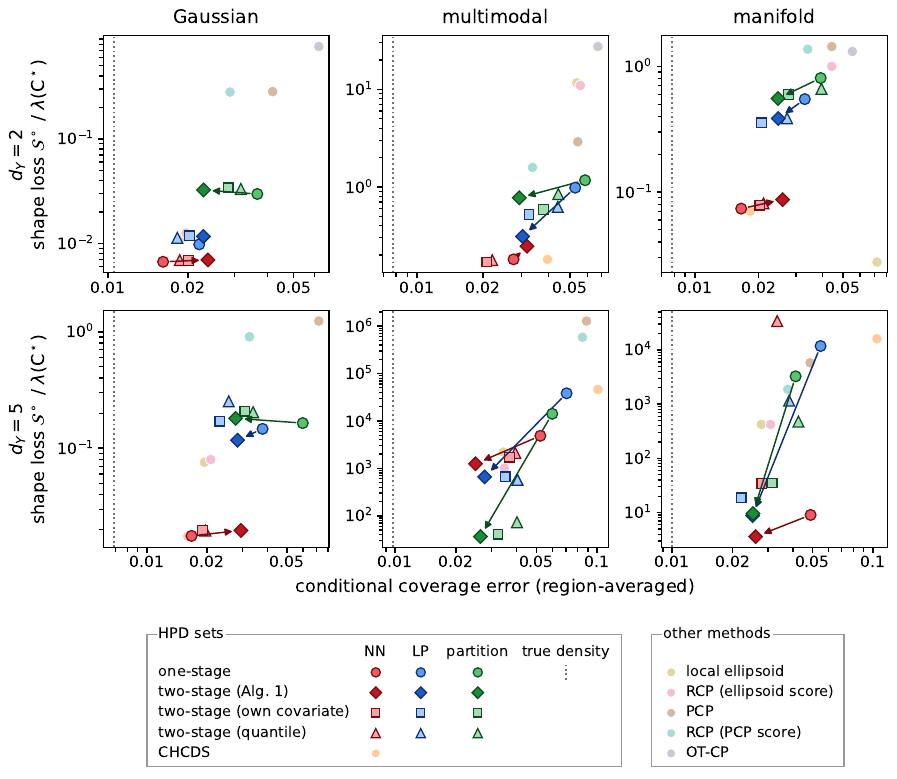}
\caption{Region-averaged conditional coverage error and shape loss relative to $\Leb(\Cstar)$ on logarithmic axes (means over twenty replications). Arrows connect the one-stage and \cref{alg:local} versions of each density estimate, the latter fitted on fewer points. Dotted line: one-stage calibration with the true density, whose shape loss is zero. Standard errors are reported in \cref{app:numerics} of the supplement.}\label{fig:cov-shape}
\end{figure}

\begin{figure}[!t]
\centering
\includegraphics[width=\textwidth]{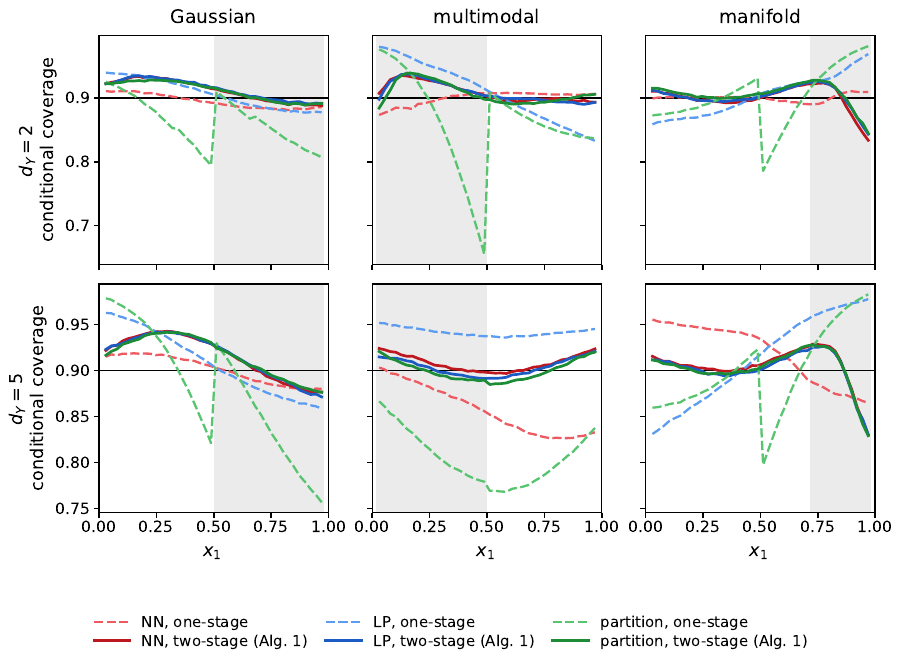}
\caption{Mean conditional coverage at the covariates $(x_1,0.5,0.5)$, $0.03\le x_1\le0.97$, of one-stage calibration (dashed) and \cref{alg:local} (solid) with the NN (red), the LP (blue) and the partition estimate (green). Gray: hard region; horizontal line: $1-\alpha=0.9$.}\label{fig:cov-curves}
\end{figure}

\paragraph{Results} \Cref{fig:cov-shape} reports the coverage error and the shape loss, \cref{fig:cov-curves} the conditional coverage along the segment $\{(x_1,0.5,0.5):0<x_1<1\}$, and \cref{app:numerics} of the supplement all values with standard errors. Marginal coverage is close to $0.9$ for every method.

The coverage error of \cref{alg:local} is about the same in every model, for both values of $d_Y$ and for all three density estimates. Its guaranteed rate does not depend on the accuracy of the estimate (\cref{thm:local-cov}), and its coverage curves with the three estimates nearly coincide. They stay close to $0.9$ except near the face $x_1=1$ in the manifold model, where the covariates reach the boundary of their support and cross-validation selects the largest regression bandwidth of the grid. The coverage error of one-stage calibration, in contrast, follows the accuracy of the estimate. It is small for the network in the Gaussian model and in the manifold model with $d_Y=2$, and large in the multimodal and manifold models with $d_Y=5$, where the estimates resolve the components of scale $0.02$ poorly. Two-stage calibration therefore lowers the coverage error where the estimate is poor and raises it where the estimate is accurate.

In the multimodal and manifold models with $d_Y=5$, the one-stage sets have shape losses many times $\Leb(\Cstar)$, and \cref{alg:local} reduces them by one to three orders of magnitude with the same density fit, as \cref{tab:samefit} in the supplement shows. These ratios are large because, at covariates where the components of scale $0.02$ carry most of the mass, $\Leb(\Cstar(x))$ is tiny, and a set of moderate volume is many times larger. Summed over random covariates, where covariates with large oracle sets dominate, the volumes of the sets of one-stage calibration and of \cref{alg:local} are at most about $2.3$ times the oracle volume (\cref{tab:main-gaussian,tab:main-multimodal,tab:main-manifold} in the supplement). In these two settings, the estimates lie outside the range of \cref{prop:supshape}, and the shape loss of a ranking depends strongly on the level at which it is cut. One-stage calibration cuts it where the estimated mass reaches the calibrated threshold, and the resulting set keeps large regions of small density over which the estimates spread their mass. \Cref{alg:local} instead cuts the ranking where the coverage estimated from the auxiliary points reaches $1-\alpha$. In the other cases, the shape losses of the two calibrations are of the same order. Although \cref{cond:A} is verified only for the LP and the partition estimate, \cref{alg:local} has its smallest shape loss with the network in every setting except the multimodal model with $d_Y=5$: \cref{thm:local} requires accuracy of the ranking function only near the level set and up to a strictly increasing transformation, and the minimax rate of the LP bounds its worst case over a H\"older class rather than its accuracy at these sample sizes.

Two-stage (own covariate) and two-stage (quantile) regress HPD values at the auxiliary covariates instead of fixing the candidate predictive sets at the target covariate. Their regression targets therefore vary with $x$ as the estimate does, whereas that of \cref{alg:local} is as smooth as the conditional density. Accordingly, their coverage errors are larger than that of \cref{alg:local} for the partition estimate, whose HPD values jump at the cube faces, and in the multimodal model with $d_Y=5$, but smaller where the estimate is accurate.

Among the baselines, PCP and OT-CP have shape losses well above those of \cref{alg:local} in every model. CHCDS, which uses the network density, has the largest coverage errors of all methods in the multimodal and manifold models with $d_Y=5$. The local ellipsoid is competitive in the Gaussian model, where its shape is correct.

\paragraph{Real data} \Cref{app:bike-mv} of the supplement applies the methods to the hourly bike-sharing demand of the Capital Bikeshare system \citep{fanaee2014event,fanaee2013bike}, with the hour, the temperature and the humidity as covariates and the log counts of casual and of registered rentals as the response. There, the sets of one-stage calibration with the network density are the smallest and follow the skewed shape of the demand. Two-stage calibration brings the coverage on groups of test points closer to $0.9$ for the LP and the partition estimate but not for the network.

\section{Discussion}\label{sec:discussion}
The excess volume of a predictive set over the oracle set splits into a coverage term and an excess mass, with minimax rates $\rho_n$ and $\zeta_n^2$. The coverage of a fixed predictive set is a regression function of the covariate. Two-stage calibration estimates it at the regression rate, uniformly in the covariate and whatever the accuracy of the density estimate. The shape loss is at most of the order of the squared sup-norm error of the estimate, a bound that accuracy in $L^q$ cannot give. With a local polynomial estimate, two-stage calibration attains both rates in every dimension. The coverage error of one-stage calibration is the variation in the covariate of the density error integrated over the oracle set. The integration averages out the noise of the estimate but not its bias. Three problems remain open.

\paragraph{Adaptation} In two-stage calibration, the regression bandwidth depends on $\beta$, and the excess-mass rate requires the density-optimal bandwidth. Adaptation to $\beta$, for instance by cross-validation in either stage, is open.

\paragraph{A shape lower bound under calibrated coverage} The lower bound for the excess mass holds for all procedures, but that for the shape loss holds only for procedures with accurate coverage (\cref{rem:shape-lower}). A shape lower bound among procedures whose coverage lies in a fixed window around $\gamma_\star$ is open. Since such a procedure chooses its level from the data, perturbations localized at the oracle level do not constrain it. With perturbations at all levels of the window, the loss depends on the perturbation through the data-dependent level, and the coordinatewise decomposition behind Assouad's lemma is not available.

\paragraph{Using the sample more efficiently} Two-stage calibration uses an additional split of the sample, which increases the density error. Cross-fitting \citep{chernozhukov2018double}, with an independent calibration sample kept for finite-sample validity, could reduce this loss of accuracy. Its analysis requires control of the dependence between the fold-specific density estimates and the regression targets.

\begin{funding}
This work was supported by the National Research Foundation of Korea (NRF) grant funded by the Korea government (MSIT) (RS-2026-25486448) and by an INHA UNIVERSITY Research Grant.
\end{funding}

\begin{acks}[Statement on AI use]
The author used the AI assistants Claude Code (Anthropic) and Codex (OpenAI) to edit the text, to check the proofs and the literature, and to write and run the code for the numerical experiments. The author reviewed all AI-assisted content and takes full responsibility for the paper.
\end{acks}

\begin{supplement}
\stitle{Supplement to ``Coverage and shape of conformal highest-density sets: a minimax theory''}
\sdescription{Guarantees of related procedures, proofs of all results, and details and further results of the numerical experiments, including the bike-sharing application.}
\end{supplement}

\begingroup
\renewcommand{\doi}[1]{\href{https://doi.org/#1}{doi:\nolinkurl{#1}}}
\putbib
\endgroup
\end{bibunit}

\newpage

\setcounter{page}{1}
\renewcommand{\thepage}{S-\arabic{page}}

\begin{appendices}
\begin{bibunit}
\crefalias{section}{appendix}
\crefalias{subsection}{appendix}
\crefalias{subsubsection}{appendix}

	\begin{center}
		{\large \textsc{Supplement to ``Coverage and shape of conformal highest-density sets: a minimax theory''}} \\
		\medskip
		{Ilsang Ohn}
		\medskip
	\end{center}

Section, equation and result numbers without a letter prefix refer to the main paper.

\setcounter{tocdepth}{2}
\begingroup
\makeatletter
\def\@pnumwidth{2.5em}
\makeatother
\tableofcontents
\endgroup
\addtocontents{toc}{\protect\setcounter{tocdepth}{2}}

\section{Guarantees of related procedures}\label{app:literature}
\Cref{tab:literature} states the guarantees compared in \cref{sec:cov-existing} in the notation of this paper.
\begin{table}[!t]
\centering\footnotesize
\caption{Guarantees for conformal and plug-in predictive sets, in the notation of this paper ($\mathsf C$ is the oracle set of each work). Here $\cov(x)$ is the conditional coverage given all data; $\E_{\rm cal}\cov(x)$ averages it over the calibration sample and $\E\cov(x)$ over all data; $X$ denotes a random test covariate. Level-set rates are for a margin condition with exponent one; constants are those of the cited results.}\label{tab:literature}
\setlength{\tabcolsep}{3pt}\renewcommand{\arraystretch}{1.03}
\begin{tabular}{@{}>{\raggedright\arraybackslash}p{0.16\textwidth}>{\raggedright\arraybackslash}p{0.035\textwidth}>{\raggedright\arraybackslash}p{0.14\textwidth}>{\raggedright\arraybackslash}p{0.62\textwidth}@{}}
\toprule
work & $d_Y$ & calibration & guarantee\\
\midrule
\citet{lei2013distribution} & any & no covariate & $\Leb(\hC\triangle\mathsf C)=O_P\bigl((\log n/n)^{\beta/(2\beta+d_Y)}+(\log n/n)^{1/2}\bigr)$\\
\citet{lei2014distribution} & $1$ & local, on cells in $x$ & $\sup_x\Leb\bigl(\hC(x)\triangle\mathsf C(x)\bigr)=O_P\bigl((\log n/n)^{\beta/(\beta(d_X+2)+1)}\bigr)$, with a matching lower bound\\
HPD-split \citep{izbicki2022cd} & any & one-stage & ${\mathsf P}\{Y\in\hC(X)\triangle\mathsf C(X)\}=o(1)$ and $|\E\cov(x)-\gamma_\star|=o(1)$ outside a set of vanishing mass\\
CD-split \citep{izbicki2022cd} & any & local, on cells built from $\hp$ & ${\mathsf P}\{Y\in\hC(X)\mid X\in A\}\ge\gamma_\star$ on each cell $A$; asymptotic conditional validity, without a rate, as the cells shrink in the metric defined by $\hp$\\
\citet{gibbs2023conformal} & any & quantile regression over a class $\mathcal F$ & $\E\sbr{f(X)(\ind\{Y\in\hC(X)\}-\gamma_\star)}\ge0$ for nonnegative $f$ in a finite-dimensional $\mathcal F$, on average over all data; up to a regularization term for infinite-dimensional $\mathcal F$\\
CP$^2$ \citep{plassier2024conditionally} & any & one-stage & $\E_{\rm cal}\cov(x)\ge\gamma_\star-\TV(\hP_x,{\mathsf P}_x)-\eps_n(x)$, with $\eps_n(x)$ a calibration error term; $|\E_{\rm cal}\cov(X)-\gamma_\star|=O_P\bigl((\log n/n)^{1/2}+\overline{\TV}_n\bigr)$, with $\overline{\TV}_n$ the average total-variation error\\
RCP \citep{plassier2025rectifying} & any & quantile regression of the score & $|\E_{\rm cal}\cov(x)-\gamma_\star|$ bounded by the square root of the excess pinball loss of the kernel quantile estimate, up to a calibration term\\
RLCP \citep{hore2025conformal}, analyzed by \citet{conrad2026beyond} & $1$ & random localization at the test covariate, with bandwidth $h$ & for a fixed score whose conditional distribution function is $\beta$-H\"older in $x$ with $\beta\le1$: given the localization center, with probability at least $1-\delta$ over the calibration sample, $|\cov(x)-\gamma_\star|\lesssim h^\beta+\{\log(4/\delta)/(nh^{d_X})\}^{1/2}+\log(4/\delta)/(nh^{d_X})$ simultaneously for $x$ within $h$ of the center; the rate $n^{-\beta/(2\beta+d_X)}$ up to logarithmic factors. For a learned score close in sup-norm to a fixed score whose conditional $\gamma_\star$-quantile does not depend on $x$: the same calibration term plus this sup-norm distance, without $h^\beta$\\
\citet{duchi2025sample} & any & quantile regression on $m$ features & given the data, coverage on each of $m$ groups $\mathcal G$ at least $\gamma_\star-c\{m\log n/(n\Pi(\mathcal G))\}^{1/2}$; minimax up to $\log n$\\
MOPI \citep{bao2026shape} & any & minimax over a class of sets & given the data, with high probability: $\int(\cov-\gamma_\star)^2\,\d\Pi$ at most the best value in the class plus approximation and complexity terms, e.g.\ $O(\log^{d_X+1}n/n)$ for a Vapnik--Chervonenkis (VC) class of sets that contains the target and for weight functions in a Gaussian-kernel reproducing kernel Hilbert space\\
\citet{gao2025volume} & $1$ & unions of $k$ intervals & if the c.d.f.\ estimate is $\delta$-accurate over unions of $k$ intervals: $\cov(X)\ge\gamma_\star-3\delta$ with probability at least $1-\delta$, and, with probability at least $1-2\delta$, volume at most that of the best $k$ intervals at coverage $\gamma_\star+1/n+4\delta+\gamma_{\rm grid}$, where $\gamma_{\rm grid}$ is a discretization slack\\
\midrule
this paper (\cref{thm:onestage} and~\cref{cor:lp-onestage}) & any & one-stage & for an estimate with sup-norm error $r_n\ge(\log n/n)^{1/2}$ and null level sets near the level, with high probability: $|\cov^{\rm one}(x)-\gamma_\star|\lesssim r_n$ and $\cE(\hC^{\rm one}(x))\lesssim r_n^2$; for integer $\beta\ge2$ and the local polynomial plug-in at bandwidths $b\ge h$, under (R1)--(R3) of \cref{app:fixed-law}: $\E\int_{\cX_0}|\cov^{\rm one}_b(x)-\gamma_\star|\,\d x\gtrsim b^\beta$ at distributions whose integrated leading bias varies with $x$\\
(\cref{thm:local-cov,thm:local} and~\cref{prop:lp-conditions}) & any & two-stage & for $\beta\ge1$ and the perturbed local polynomial estimate, with high probability, $|\cov(x)-\gamma_\star|\lesssim\rho_n$ and $\cE(\hC(x))\lesssim\zeta_n^2\mathfrak L(\zeta_n)$, uniformly in $x\in\cX_0$; minimax orders, exactly for uniform coverage and up to logarithmic factors otherwise (\cref{thm:lower})\\
\bottomrule
\end{tabular}
\end{table}

\FloatBarrier
\section{Proofs for Sections~\ref{sec:setup} and~\ref{sec:losses}}\label{app:setup}
\emph{H\"older norms.} Let $\mathcal O$ be a convex set with nonempty interior. For a function $g$ on $\mathcal O$ and $\beta'>0$, let $\norm[0]{g}_{\mathcal H^{\beta'}_{\mathcal O}}$ denote the smallest $L$ in \eqref{eq:holder} with $\beta'$ in place of $\beta$. Since $\mathcal O$ need not be open, the derivatives in \eqref{eq:holder} at its boundary points are defined as the continuous extensions of those in the interior. These extensions exist because, by \eqref{eq:holder}, the derivatives of order at most $\lceil\beta'\rceil-1$ are uniformly continuous on the interior. For $\beta'>1$, Taylor's formula with integral remainder then holds along every segment in $\mathcal O$, by approximation with segments in the interior.

For a fixed $x$ we write
	\begin{align*}
		&t_\gamma:=t_\gamma(x),\quad t_\star:=t_\star(x),\quad\Cstar_\gamma:=\Cstar_\gamma(x),\quad\Cstar:=\Cstar(x),\quad G:=G_x,\quad V:=V_x,\\
		&\hC:=\hC(x),\quad\cov:=\cov(x),\quad\cE:=\cE(\hC(x)),\quad\cS:=\cS(\hC(x)),
	\end{align*}
where the second line applies when $\hC$ is also fixed.

We first show that $G$ is continuous and strictly decreasing on $[t_\star-\varepsilon_0,t_\star+\varepsilon_0]$. Let $s$ lie in this interval. By \eqref{eq:margin-upper} and the fact that $\eps\mathfrak L(\eps)\to0$ as $\eps\to0$, we have
	\begin{align*}
		\Leb\{y:{\mathsf p}(y\mid x)=s\}\le\inf_{0<\eps\le\varepsilon_0}C_m\eps\,\mathfrak L(\eps)=0,
	\end{align*}
and hence ${\mathsf P}_x\{{\mathsf p}(Y\mid x)=s\}=0$. Since $G$ is left-continuous and its jump at $s$ is ${\mathsf P}_x\{{\mathsf p}(Y\mid x)=s\}$, it is continuous at $s$. Now let $s_1<s_2$ lie in $[t_\star-\varepsilon_0,t_\star+\varepsilon_0]$. The lower one-sided band of \eqref{eq:margin-lower} at the level $(s_1+s_2)/2$, with $\eps=(s_2-s_1)/2$, lies in $\{y:s_1<{\mathsf p}(y\mid x)<s_2\}$. By the probability bound of \cref{ass:margin} for this band, we obtain
	\begin{align}\label{eq:G-lower}
		G(s_1)-G(s_2)\ge{\mathsf P}_x\Bigl\{s_1<{\mathsf p}(Y\mid x)\le\frac{s_1+s_2}2-a\,\frac{s_2-s_1}2\Bigr\}\ge\frac{c_m}2(s_2-s_1).
	\end{align}
Therefore $G$ is strictly decreasing on $[t_\star-\varepsilon_0,t_\star+\varepsilon_0]$. Let
	\begin{align*}
		\mathcal W:=[t_\star-\varepsilon_0/4,t_\star+\varepsilon_0/4].
	\end{align*}

\begin{lemma}[Levels near the oracle level]\label{lem:t-window}
Suppose that \cref{ass:margin} holds, fix $x$, and let $t_{\min}:=\varepsilon_0/2$. Then $t_\star>(1+a)\varepsilon_0$. Moreover, every $\gamma$ with $|\gamma-\gamma_\star|\le c_0$ satisfies $t_\gamma\in\mathcal W$, $G(t_\gamma)=\gamma$ and $t_\gamma\ge t_{\min}$, and
	\begin{align}\label{eq:t-lip}
		|t_\gamma-t_{\gamma'}|\le\frac2{c_m}|\gamma-\gamma'|\qquad(|\gamma'-\gamma_\star|\le c_0).
	\end{align}
\end{lemma}

\begin{proof}
We first locate $t_\star$, $t_0$ and $t_1$. By \eqref{eq:margin-lower} with $\eps=\varepsilon_0$, the upper one-sided band at the level $t_\star+\varepsilon_0$ and the lower one-sided band at the level $t_\star-\varepsilon_0$ each have probability at least $c_m\varepsilon_0>0$. In particular, each of them contains a point at which ${\mathsf p}(\cdot\mid x)>0$. Since ${\mathsf p}(\cdot\mid x)\ge t_\star+(1+a)\varepsilon_0$ on the upper band, we have $t_0=\sup_\cY{\mathsf p}(\cdot\mid x)\ge t_\star+(1+a)\varepsilon_0$. On the lower band, ${\mathsf p}(\cdot\mid x)\le t_\star-(1+a)\varepsilon_0$. Thus $t_\star-(1+a)\varepsilon_0>0$, that is, $t_\star>(1+a)\varepsilon_0$, and ${\mathsf P}_x\{{\mathsf p}(Y\mid x)\le t_\star-(1+a)\varepsilon_0\}>0$. The last bound implies that $G(v)<1$ for $v>t_\star-(1+a)\varepsilon_0$, and hence $t_1\le t_\star-(1+a)\varepsilon_0$. Therefore $t_0,t_1\notin\mathcal W$. Since $G$ is continuous at $t_\star=t_{\gamma_\star}$, we have $G(t_\star)=\gamma_\star$. Hence, by \eqref{eq:G-lower} and $c_0=c_m\varepsilon_0/8$, we have $G(t_\star-\varepsilon_0/4)\ge\gamma_\star+c_0$ and $G(v)<\gamma_\star-c_0$ for $v>t_\star+\varepsilon_0/4$. If $|\gamma-\gamma_\star|\le c_0$, it follows that $G(t_\star-\varepsilon_0/4)\ge\gamma$ and $G(v)<\gamma$ for $v>t_\star+\varepsilon_0/4$, and hence $t_\gamma\in\mathcal W$. Consequently, the interval $[\gamma_\star-c_0,\gamma_\star+c_0]$ contains neither $0$ nor $1$, and it therefore lies in $(0,1)$.

Now let $|\gamma-\gamma_\star|\le c_0$. By the definition of $t_\gamma$ and the monotonicity of $G$, we have $G(v)\ge\gamma$ for $v<t_\gamma$ and $G(v)<\gamma$ for $v>t_\gamma$. Since $t_\gamma\in\mathcal W$, the function $G$ is continuous at $t_\gamma$, and hence $G(t_\gamma)=\gamma$. If also $|\gamma'-\gamma_\star|\le c_0$, then $G(t_{\gamma'})=\gamma'$, and, by \eqref{eq:G-lower} with $\{s_1,s_2\}=\{t_\gamma,t_{\gamma'}\}\subset\mathcal W$, we obtain $|\gamma-\gamma'|\ge(c_m/2)|t_\gamma-t_{\gamma'}|$, which is \eqref{eq:t-lip}. Finally, since $t_\gamma\in\mathcal W$ and $t_\star>(1+a)\varepsilon_0$, we have $t_\gamma\ge t_\star-\varepsilon_0/4>(3/4+a)\varepsilon_0>t_{\min}$.
\end{proof}

\begin{remark}[The exponent in the margin condition]\label{rem:margin-exponent}
Suppose that the band probabilities are bounded above and below by multiples of $\eps^\upsilon$, uniformly over a nontrivial interval of levels. The set of points at which the density lies in a fixed range of levels inside that interval is the union of about $1/\eps$ disjoint bands of width $2\eps$. Its probability, which does not depend on $\eps$, is therefore of order $\eps^{\upsilon-1}$ for all sufficiently small $\eps$. Hence $\upsilon=1$.
\end{remark}

\subsection{Proof of \texorpdfstring{\cref{lem:identity}}{Lemma \ref{lem:identity}}}
\begin{proof}
Let $A$ be a measurable set, and let $\gamma$ satisfy ${\mathsf P}_x(\Cstar_\gamma)=\gamma$ and $t_\gamma>0$. By \cref{lem:t-window}, every $\gamma$ with $|\gamma-\gamma_\star|\le c_0$ satisfies both conditions under \cref{ass:margin}. Writing volumes and probabilities as integrals of indicators and using ${\mathsf P}_x(\Cstar_\gamma)=\gamma$, we obtain
	\begin{align*}
		\Leb(A)-V(\gamma)-\frac{{\mathsf P}_x(A)-\gamma}{t_\gamma}
		=\int\bigl(\ind_A-\ind_{\Cstar_\gamma}\bigr)\Bigl(1-\frac{{\mathsf p}(\cdot\mid x)}{t_\gamma}\Bigr)\d y .
	\end{align*}
Since ${\mathsf p}(\cdot\mid x)<t_\gamma$ on $A\setminus\Cstar_\gamma$ and ${\mathsf p}(\cdot\mid x)\ge t_\gamma$ on $\Cstar_\gamma\setminus A$, the integrand equals $\ind_{A\triangle\Cstar_\gamma}|{\mathsf p}(\cdot\mid x)-t_\gamma|/t_\gamma$. Therefore,
	\begin{align}\label{eq:basic-identity}
		\Leb(A)-V(\gamma)=\frac{{\mathsf P}_x(A)-\gamma}{t_\gamma}+\frac1{t_\gamma}\int_{A\triangle\Cstar_\gamma}|{\mathsf p}(y\mid x)-t_\gamma|\,\d y .
	\end{align}
With $A=\hC$ and $\gamma=\gamma_\star$, this is \cref{lem:identity}, since $V(\gamma_\star)=\Leb(\Cstar)$ and the last term is $\cE$.
\end{proof}

\subsection{Proof of \texorpdfstring{\eqref{eq:E-symdiff}}{the comparison of excess mass and symmetric difference}}\label{app:E-symdiff}
\begin{proof}
We write $V_\triangle:=\Leb(\hC\triangle\Cstar)$. If $V_\triangle=0$, then $\cE=0$. By \cref{ass:smooth} and \cref{lem:t-window}, we have $t_{\min}\le t_\star\le{\mathsf p}_{\max}$ and $|{\mathsf p}(\cdot\mid x)-t_\star|\le{\mathsf p}_{\max}$. Therefore $\cE\le({\mathsf p}_{\max}/t_{\min})\,V_\triangle$, which is the upper bound.

For the lower bound, let $\eps\in(0,\varepsilon_0]$. By \eqref{eq:margin-upper} at the level $t_\star$ and by $t_\star\le{\mathsf p}_{\max}$,
	\begin{align*}
		\Leb\bigl((\hC\triangle\Cstar)\cap\mathcal J_x(\eps)\bigr)\le C_m\eps\,\mathfrak L(\eps),\qquad\frac{|{\mathsf p}(\cdot\mid x)-t_\star|}{t_\star}\ge\frac\eps{{\mathsf p}_{\max}}\ \text{on }\cY\setminus\mathcal J_x(\eps).
	\end{align*}
Hence
	\begin{align*}
		\cE\ge\frac\eps{{\mathsf p}_{\max}}\,\bigl(V_\triangle-C_m\eps\,\mathfrak L(\eps)\bigr)_+.
	\end{align*}
Both choices of $\mathfrak L$ are nonincreasing, $u\mapsto u\mathfrak L(u)$ is nondecreasing, and, for $u>0$ and $s\ge1$,
	\begin{align}\label{eq:L-props}
		1\le\mathfrak L(u),\qquad\mathfrak L(u/s)\le(1+\log s)\,\mathfrak L(u),\qquad\mathfrak L\bigl(u/\mathfrak L(u)\bigr)\le2\mathfrak L(u).
	\end{align}
We fix $c_\triangle\in(0,1]$ with $c_\triangle\Leb(\cY)\le\varepsilon_0$ and $2C_mc_\triangle\{1+\log(1/c_\triangle)\}\le\frac12$, and we choose $\eps:=c_\triangle V_\triangle/\mathfrak L(V_\triangle)$. Since $\mathfrak L\ge1$ and $V_\triangle\le\Leb(\cY)$, we have $\eps\le\varepsilon_0$. By \eqref{eq:L-props} with $s=1/c_\triangle$ and $u=V_\triangle/\mathfrak L(V_\triangle)$,
	\begin{align*}
		\mathfrak L(\eps)\le2\bigl\{1+\log(1/c_\triangle)\bigr\}\mathfrak L(V_\triangle),\qquad C_m\eps\,\mathfrak L(\eps)\le\frac{V_\triangle}2.
	\end{align*}
By the lower bound on $\cE$ above, we obtain $\cE\ge c_\triangle V_\triangle^2/\{2{\mathsf p}_{\max}\mathfrak L(V_\triangle)\}$. This proves the lower bound with $c:=c_\triangle/(2{\mathsf p}_{\max})$.
\end{proof}

\subsection{Proof of \texorpdfstring{\cref{lem:shape-split}}{Lemma \ref{lem:shape-split}}}
\begin{proof}
We first prove the properties of $V=V_x$ stated after the lemma. Let $\gamma>\gamma'$ lie in $[\gamma_\star-c_0,\gamma_\star+c_0]$. The oracle sets are nested, and $\Cstar_\gamma\setminus\Cstar_{\gamma'}=\{y:t_\gamma\le{\mathsf p}(y\mid x)<t_{\gamma'}\}$. By \cref{lem:t-window}, ${\mathsf P}_x(\Cstar_\gamma\setminus\Cstar_{\gamma'})=\gamma-\gamma'$ and $\Leb(\Cstar_\gamma\setminus\Cstar_{\gamma'})=V(\gamma)-V(\gamma')$. Hence
	\begin{align*}
		\frac1{t_{\gamma'}}\le\frac{V(\gamma)-V(\gamma')}{\gamma-\gamma'}\le\frac1{t_\gamma}.
	\end{align*}
By \eqref{eq:t-lip} and $t_\gamma\ge t_{\min}$, the function $\gamma\mapsto1/t_\gamma$ is Lipschitz with constant $2/(c_mt_{\min}^2)$. Thus $V$ is differentiable with $V'(\gamma)=1/t_\gamma$, and $V'$ is Lipschitz with the same constant. Since $\gamma\mapsto t_\gamma$ is nonincreasing, $V'$ is nondecreasing, and therefore $V$ is convex.

If $|\cov-\gamma_\star|>c_0$, then $\cS=\cE$ by definition. Otherwise, let $\gamma:=\cov$. By \cref{lem:t-window}, ${\mathsf P}_x(\Cstar_\gamma)=\gamma={\mathsf P}_x(\hC)$. Hence \eqref{eq:basic-identity} at coverage $\gamma$ with $A=\hC$ reads
	\begin{align*}
		\Leb(\hC)-V(\gamma)=\frac1{t_\gamma}\int_{\hC\triangle\Cstar_\gamma}|{\mathsf p}(y\mid x)-t_\gamma|\,\d y,
	\end{align*}
that is, the two expressions for $\cS$ in \cref{sec:setup-shape} agree. By the definition of $\cS$ and by \cref{lem:identity}, we have $\Leb(\hC)=V(\gamma)+\cS$ and $\Leb(\hC)=V(\gamma_\star)+(\gamma-\gamma_\star)/t_\star+\cE$. Subtracting the first identity from the second and using $V'(\gamma_\star)=1/t_\star$, we obtain
	\begin{align*}
		\cE-\cS=V(\gamma)-V(\gamma_\star)-(\gamma-\gamma_\star)V'(\gamma_\star)=\int_{\gamma_\star}^{\gamma}\bigl\{V'(u)-V'(\gamma_\star)\bigr\}\,\d u.
	\end{align*}
Since $V'$ is nondecreasing, this integral is nonnegative, and, since $V'$ is Lipschitz with constant $2/(c_mt_{\min}^2)$, it is at most $|\cov-\gamma_\star|^2/(c_mt_{\min}^2)$.
\end{proof}

\subsection{Two deterministic comparisons}
Let $\hat a$ be a function on $\cY$, $\hat t$ a level, and $\hC$ a set with $\{y:\hat a(y)>\hat t\}\subset\hC\subset\{y:\hat a(y)\ge\hat t\}$. We also use the following inclusion. Let $g,g'$ be functions on $\cY$, let $s,s'$ be levels, and let $A,A'$ be sets with $\{g>s\}\subset A\subset\{g\ge s\}$ and $\{g'>s'\}\subset A'\subset\{g'\ge s'\}$. If $y\in A\setminus A'$, then $s'\ge g'(y)\ge g(y)-\sup_\cY|g-g'|\ge s-\sup_\cY|g-g'|$, and symmetrically if $y\in A'\setminus A$. Hence
	\begin{align}\label{eq:sandwich}
		A\triangle A'\subset\Bigl\{y:|g'(y)-s'|\le\sup_\cY|g-g'|+|s-s'|\Bigr\}.
	\end{align}

\begin{lemma}[Level comparison]\label{lem:threshold}
Under \cref{ass:margin}, let $\hC$ be as above with $|\cov-\gamma_\star|\le c_0$, and suppose that $\sup_\cY|\hat a-{\mathsf p}(\cdot\mid x)|\le r$. Then $|\hat t-t_{\cov}|\le r$.
\end{lemma}

\begin{proof}
Since $|\hat a-{\mathsf p}(\cdot\mid x)|\le r$ everywhere,
	\begin{align*}
		\{y:{\mathsf p}(y\mid x)>\hat t+r\}\subset\{y:\hat a(y)>\hat t\}\subset\hC\subset\{y:\hat a(y)\ge\hat t\}\subset\{y:{\mathsf p}(y\mid x)\ge\hat t-r\},
	\end{align*}
and hence
	\begin{align*}
		{\mathsf P}_x\{{\mathsf p}(Y\mid x)>\hat t+r\}\le\cov\le G(\hat t-r).
	\end{align*}
Since $t_{\cov}=\sup\{v:G(v)\ge\cov\}$, the upper bound gives $\hat t-r\le t_{\cov}$. By \cref{lem:t-window}, we have $t_{\cov}\in\mathcal W$ and $G(t_{\cov})=\cov$. Suppose, to reach a contradiction, that $\hat t+r<t_{\cov}$, and let $s\in(\hat t+r,t_{\cov})$ satisfy $s\ge t_\star-\varepsilon_0$. Since $G$ is strictly decreasing on $[t_\star-\varepsilon_0,t_\star+\varepsilon_0]$ by \eqref{eq:G-lower}, the lower bound gives
	\begin{align*}
		\cov\ge{\mathsf P}_x\{{\mathsf p}(Y\mid x)>\hat t+r\}\ge G(s)>G(t_{\cov})=\cov,
	\end{align*}
a contradiction. Therefore $|\hat t-t_{\cov}|\le r$.
\end{proof}

\begin{lemma}[Shape loss and the error on the misplaced set]\label{lem:shape-error}
Under \cref{ass:smooth,ass:margin}, there exists a constant $c>0$ such that the following holds. Let $\hC$ be as above with $|\cov-\gamma_\star|\le c_0$ and $\sup_\cY|\hat a-{\mathsf p}(\cdot\mid x)|\le\varepsilon_0/8$, and let $\mathfrak e:=\hat a-{\mathsf p}(\cdot\mid x)$ and $\mathfrak M:=\hC\triangle\Cstar_{\cov}$. Then $\mathfrak M\subset\mathcal J_x$, the band \eqref{eq:Jx}, and
	\begin{align*}
		\frac{7t_{\min}}8\,\cS(\hC)\le\int_{\mathfrak M}|\mathfrak e|\,\d y,\qquad\frac{c\,\Leb(\mathfrak M)^2}{\mathfrak L\bigl(\Leb(\mathfrak M)\bigr)}\le\cS(\hC).
	\end{align*}
\end{lemma}

\begin{proof}
By \cref{lem:threshold}, we have $|\hat t-t_{\cov}|\le\varepsilon_0/8$. By \eqref{eq:sandwich} with $(g,s,A)=(\hat a,\hat t,\hC)$ and $(g',s',A')=({\mathsf p}(\cdot\mid x),t_{\cov},\Cstar_{\cov})$, every $y\in\mathfrak M$ satisfies
	\begin{align*}
		|{\mathsf p}(y\mid x)-t_{\cov}|\le\sup_\cY|\mathfrak e|+|\hat t-t_{\cov}|\le\varepsilon_0/4.
	\end{align*}
Since $t_{\cov}\in\mathcal W$ by \cref{lem:t-window}, it follows that $|{\mathsf p}(y\mid x)-t_\star|\le\varepsilon_0/2$, that is, $y\in\mathcal J_x(\varepsilon_0/2)$. Therefore $\mathfrak M\subset\mathcal J_x$.

Since $\hat a\ge\hat t$ on $\hC$ and $\hat a\le\hat t$ off $\hC$, we have $(\ind_{\hC}-\ind_{\Cstar_{\cov}})(\hat a-\hat t)\ge0$ on $\cY$. Recall that ${\mathsf P}_x(\Cstar_{\cov})=\cov={\mathsf P}_x(\hC)$ by \cref{lem:t-window} and that $\cS(\hC)=\Leb(\hC)-\Leb(\Cstar_{\cov})$ by definition. Integrating the last inequality and writing $\hat a={\mathsf p}(\cdot\mid x)+\mathfrak e$, we therefore obtain
	\begin{align*}
		0&\le\int(\ind_{\hC}-\ind_{\Cstar_{\cov}})({\mathsf p}(\cdot\mid x)+\mathfrak e-\hat t)\,\d y\\
		&=-\hat t\,\cS(\hC)+\int(\ind_{\hC}-\ind_{\Cstar_{\cov}})\,\mathfrak e\,\d y\le-\hat t\,\cS(\hC)+\int_{\mathfrak M}|\mathfrak e|\,\d y .
	\end{align*}
Moreover, by \cref{lem:t-window}, we have $t_{\cov}\ge t_\star-\varepsilon_0/4>3\varepsilon_0/4$, and hence
	\begin{align*}
		\hat t\ge t_{\cov}-\varepsilon_0/8>5\varepsilon_0/8>7t_{\min}/8.
	\end{align*}
Since $\cS(\hC)\ge0$, the last two displays yield the first bound of the lemma.

For the second bound, we use the integral expression $\cS(\hC)=t_{\cov}^{-1}\int_{\mathfrak M}|{\mathsf p}(y\mid x)-t_{\cov}|\,\d y$ of \cref{sec:setup-shape}. Since \eqref{eq:margin-upper} holds at the level $t_{\cov}\in\mathcal W$, the proof of \eqref{eq:E-symdiff} applies with $t_{\cov}$ in place of $t_\star$ and $\cS$ in place of $\cE$, and it yields the second bound.
\end{proof}

\subsection{Proof of \texorpdfstring{\cref{prop:supshape}}{Proposition \ref{prop:supshape}}}
\begin{proof}
By \cref{lem:threshold} with $\varepsilon_0/8$ in place of $r$ and by \cref{lem:t-window}, we have $|\hat t-t_\star|\le|\hat t-t_{\cov}|+|t_{\cov}-t_\star|\le3\varepsilon_0/8$. Off $\mathcal J_x$, we have $|{\mathsf p}(\cdot\mid x)-t_\star|>\varepsilon_0$, and hence $|{\mathsf p}(\cdot\mid x)-\hat t|>5\varepsilon_0/8>|\hat a-{\mathsf p}(\cdot\mid x)|$. Thus $\hat a-\hat t$ and ${\mathsf p}(\cdot\mid x)-\hat t$ have the same sign. On $\mathcal J_x$, we have $|\hat a-{\mathsf p}(\cdot\mid x)|\le r$. Therefore
	\begin{align*}
		\{y:{\mathsf p}(y\mid x)>\hat t+r\}\subset\{y:\hat a(y)>\hat t\}\subset\hC\subset\{y:\hat a(y)\ge\hat t\}\subset\{y:{\mathsf p}(y\mid x)\ge\hat t-r\}.
	\end{align*}
The proof of \cref{lem:threshold} uses only these inclusions, and it therefore gives $|\hat t-t_{\cov}|\le r$. Let $\mathfrak M:=\hC\triangle\Cstar_{\cov}$. By the inclusions, $\hat t-r\le{\mathsf p}(\cdot\mid x)\le t_{\cov}$ on $\hC\setminus\Cstar_{\cov}$ and $t_{\cov}\le{\mathsf p}(\cdot\mid x)\le\hat t+r$ on $\Cstar_{\cov}\setminus\hC$, and hence $|{\mathsf p}(\cdot\mid x)-t_{\cov}|\le2r$ on $\mathfrak M$. Since $|t_{\cov}-t_\star|\le\varepsilon_0/4$ and $2r\le\varepsilon_0$, \eqref{eq:margin-upper} at the level $t_{\cov}$ and the monotonicity of $\mathfrak L$ give, for $r>0$,
	\begin{align*}
		\Leb(\mathfrak M)&\le2C_mr\,\mathfrak L(2r)\le2C_mr\,\mathfrak L(r),\\
		\cS&=\frac1{t_{\cov}}\int_{\mathfrak M}|{\mathsf p}(y\mid x)-t_{\cov}|\,\d y\le\frac{2r}{t_{\min}}\Leb(\mathfrak M)\le Cr^2\mathfrak L(r).
	\end{align*}
For $r=0$, the same integral vanishes. The second bound follows from the first and the bound $\cE\le\cS+C|\cov-\gamma_\star|^2$ of \cref{lem:shape-split}.
\end{proof}

\section{Proof of \texorpdfstring{\cref{thm:lower}}{Theorem \ref{thm:lower}}}\label{app:lower}
\emph{Conventions.} Since the assumptions are invariant under translation, we translate the $y$-coordinates so that the box $\cY$ is centered at the origin and $\cY=-\cY$. Let $w_\cY$ be the smallest side length of $\cY$. Let $\mathbb S^{d_Y-1}$ be the unit sphere of $\R^{d_Y}$ and $S_{d_Y-1}$ its surface measure ($\mathbb S^0=\{-1,1\}$ and $S_0$ the counting measure if $d_Y=1$), so that $\d y=r^{d_Y-1}\,\d r\,\d S_{d_Y-1}(\theta)$ for $y=r\theta$. We write $|\mathbb S^{d_Y-1}|:=S_{d_Y-1}(\mathbb S^{d_Y-1})$ and let $\mathcal B(y,r)$ and $\bar{\mathcal B}(y,r)$ be the open and the closed Euclidean balls with center $y$ and radius $r$. We write $\norm[0] g_{\mathcal H^\beta}$ for the H\"older norm $\norm[0] g_{\mathcal H^\beta_{\mathcal O}}$ of \cref{app:setup} with $\mathcal O=[0,1]^{d_X}\times\cY$, which is subadditive. By \cref{ass:constants}(ii), $\cX_0$ is an axis-parallel cube of side length $w_0>0$.

\subsection{The base density and its perturbations}\label{app:lower-lemmas}
\emph{Base density.} Let $Z\sim \mathtt N(0,I_{d_Y})$, let $q$ be the $(1+\gamma_\star)/2$-quantile of $|Z|$, and let $\sigma_0:=w_\cY/(4q)$. The base density is the density of $\mathtt N(0,\sigma_0^2I_{d_Y})$ restricted to $\cY$ and renormalized,
	\begin{align}\label{eq:lb-f0}
		f_0(y):=\bar f(|y|),\qquad\bar f(r):=\frac{\e^{-r^2/(2\sigma_0^2)}}{\int_\cY\e^{-|y'|^2/(2\sigma_0^2)}\,\d y'}\qquad(y\in\cY,\ r\ge0).
	\end{align}
Since $\cY$ is centered at the origin and has smallest side length $w_\cY$, we have $\bar{\mathcal B}(0,r)\subset\cY$ for $r\le w_\cY/2$. For these $r$, let
	\begin{align*}
		\mathfrak m(r):=\int_{\mathcal B(0,r)}f_0\,\d y=\frac{\bP\{|\sigma_0 Z|\le r\}}{\bP\{\sigma_0 Z\in\cY\}}.
	\end{align*}
The function $\mathfrak m$ is continuous and strictly increasing on $[0,w_\cY/2]$, with $\mathfrak m(0)=0$. Since $w_\cY/(4\sigma_0)=q$ and $\bP\{\sigma_0Z\in\cY\}\le1$, we also have $\mathfrak m(w_\cY/4)\ge\bP\{|Z|\le q\}=(1+\gamma_\star)/2>\gamma_\star$. Hence there is a unique $r_\star\in(0,w_\cY/4)$ with $\mathfrak m(r_\star)=\gamma_\star$. Let
	\begin{align}\label{eq:lb-radii}
		r_-:=\frac{r_\star}2,\qquad r_+:=\frac{w_\cY}4,\qquad\bar t:=\bar f(r_\star),\qquad C_{\bar f'}:=\max_{[r_-,r_+]}|\bar f'|,\qquad c_{\bar f'}:=\min_{[r_-,r_+]}|\bar f'| .
	\end{align}
The thresholds of \cref{ass:constants} are
	\begin{equation}\label{eq:lb-thresholds}
		\begin{gathered}
		L_0:=\max\{1,2\norm{f_0}_{\mathcal H^\beta}\},\qquad{\mathsf p}_\sharp:=2\bar f(0),\qquad\bar C:=2|\mathbb S^{d_Y-1}|r_+^{d_Y-1},\qquad{\mathsf p}_\flat:=\frac{\min_\cY f_0}2,\\
		c_\natural:=\frac{{\mathsf p}_\flat|\mathbb S^{d_Y-1}|r_-^{d_Y-1}}{2C_{\bar f'}},\\
		\kappa_0:=\min\cbr[3]{{\mathsf p}_\flat,\ \frac{\bar t}2,\ \frac{\bar f(r_-)-\bar t}4,\ \frac{\bar t-\bar f(r_+)}4,\ \frac{c_{\bar f'}}2}.
		\end{gathered}
	\end{equation}
All the quantities in \eqref{eq:lb-f0}--\eqref{eq:lb-thresholds} depend only on $d_Y$, $\cY$ and $\alpha$, except $L_0$, which also depends on $\beta$. Under \cref{ass:constants}(i), $f_0$ has the following properties, which are the conditions used in the rest of the proof and are verified in \cref{lem:constants-suffice} below.
\begin{enumerate}[label=(L\arabic*)]
\item The function $f_0$, viewed as a function of $(x,y)$ that does not depend on $x$, satisfies $\norm[0]{f_0}_{\mathcal H^\beta}\le L/2$ and $2{\mathsf p}_\flat\le f_0\le{\mathsf p}_{\max}/2$.
\item The level $\bar t$ is the $\gamma_\star$-level of $f_0$. Moreover, $f_0>\bar t+4\varepsilon_0$ on $\mathcal B(0,r_-)$, $f_0<\bar t-4\varepsilon_0$ on $\cY\setminus\bar{\mathcal B}(0,r_+)$, and $-C_{\bar f'}\le \bar f'\le-2\kappa_0$ on $[r_-,r_+]$.
\end{enumerate}
By the choice of $r_\star$ and \eqref{eq:lb-radii}, we have $r_-<r_\star<r_+$. Let
	\begin{align}\label{eq:lb-annulus}
		\mathcal A:=\{y:r_-\le|y|\le r_+\}
	\end{align}
be the closed annulus. By (L1), (L2), \eqref{eq:lb-thresholds} and \cref{ass:constants}(i), we have
	\begin{align*}
		2\kappa_0\le C_{\bar f'},\qquad\varepsilon_0\le\kappa_0\le{\mathsf p}_\flat\le\frac{{\mathsf p}_{\max}}4,\qquad\bar t\ge2\kappa_0\ge2\varepsilon_0 .
	\end{align*}

\begin{lemma}[The thresholds suffice]\label{lem:constants-suffice}
The norm $\norm[0]{f_0}_{\mathcal H^\beta}$ is finite and the thresholds \eqref{eq:lb-thresholds} are positive. Under \cref{ass:constants}(i), $f_0$ satisfies (L1) and (L2).
\end{lemma}

\begin{proof}
\emph{Thresholds.} The density $f_0$ does not depend on $x$ and is a multiple of $\e^{-|y|^2/(2\sigma_0^2)}$ on the compact box $\cY$. Thus its derivatives of order at most $\ell+1$ are bounded on $\cY$. Since $\cY$ is convex, the mean value theorem then yields $\norm[0]{f_0}_{\mathcal H^\beta}<\infty$. Next, $\bar f'(r)=-(r/\sigma_0^2)\bar f(r)<0$ for $r>0$. Consequently, $\bar f$ is strictly decreasing on $[0,\infty)$, and $c_{\bar f'}\ge r_-\bar f(r_+)/\sigma_0^2>0$. Since $r_-<r_\star<r_+$, we also have $\bar f(r_-)>\bar t>\bar f(r_+)$. Hence $\kappa_0>0$. The other thresholds are positive by their definitions.

\emph{Level.} Let
	\begin{align}\label{eq:lb-G0-def}
		G_0(v):=\int_{\{y:f_0(y)\ge v\}}f_0\,\d y .
	\end{align}
For $r\le w_\cY/2$, we have $\{y:f_0(y)\ge\bar f(r)\}=\bar{\mathcal B}(0,r)$, and hence $G_0(\bar f(r))=\mathfrak m(r)$. In particular, $G_0(\bar t)=\gamma_\star$. For $v>\bar t$, either $v>\bar f(0)$ and $G_0(v)=0$, or $v=\bar f(r)$ with $r<r_\star$ and $G_0(v)=\mathfrak m(r)<\gamma_\star$. Hence $\bar t=\sup\{v:G_0(v)\ge\gamma_\star\}$ is the $\gamma_\star$-level of $f_0$.

\emph{(L1) and (L2).} Since $\bar f$ is strictly decreasing, \cref{ass:constants}(i) and \eqref{eq:lb-thresholds} imply
	\begin{align*}
		&\norm{f_0}_{\mathcal H^\beta}\le\frac{L_0}2\le\frac L2,\qquad2{\mathsf p}_\flat=\min_\cY f_0,\qquad\max_\cY f_0=\frac{{\mathsf p}_\sharp}2\le\frac{{\mathsf p}_{\max}}2,\\
		&\bar f(r_-)\ge\bar t+4\varepsilon_0,\qquad\bar f(r_+)\le\bar t-4\varepsilon_0,\\
		&-C_{\bar f'}\le\bar f'\le-c_{\bar f'}\le-2\kappa_0\ \text{on }[r_-,r_+],
	\end{align*}
which establish (L1) and the remaining claims of (L2).
\end{proof}

\emph{Small perturbations.} We call a function $\varphi$ on $\cY$ a \emph{small perturbation} if $\varphi\in C^1(\cY)$, $\int_\cY\varphi\,\d y=0$ and
	\begin{align}\label{eq:lb-admissible}
		\sup_\cY|\varphi|\le\varepsilon_0/2,\qquad\sup_\cY|\nabla\varphi|\le \kappa_0 .
	\end{align}
For a small perturbation $\varphi$, let ${\mathsf P}_\varphi$ be the distribution with density ${\mathsf p}_\varphi$, where
	\begin{align}\label{eq:lb-pphi}
		{\mathsf p}_\varphi:=f_0+\varphi,\qquad G_\varphi(v):={\mathsf P}_\varphi\{{\mathsf p}_\varphi(Y)\ge v\},\qquad t_\varphi:=\sup\{v:G_\varphi(v)\ge\gamma_\star\}.
	\end{align}
By (L1), \eqref{eq:lb-admissible} and $\varepsilon_0\le{\mathsf p}_\flat$, we have ${\mathsf p}_\varphi\ge2{\mathsf p}_\flat-\varepsilon_0/2\ge{\mathsf p}_\flat$. For $\varphi=0$, the function $G_\varphi$ and the level $t_\varphi$ are the function $G_0$ of \eqref{eq:lb-G0-def} and the level $\bar t$ of (L2). Let
	\begin{align}\label{eq:lb-Delta}
		\Delta_\varphi:=\sup_{|s-\bar t|\le2\varepsilon_0}|G_\varphi(s)-G_0(s)|,\qquad\delta_\varphi:=2\Delta_\varphi/c_\natural .
	\end{align}

\begin{lemma}[Small perturbations of the base density]\label{lem:lb-perturb}
Suppose that \cref{ass:constants}(i) holds, and let $\varphi$ be a small perturbation. Then the following hold.
\begin{enumerate}[label=(\roman*)]
\item For $\bar t-3\varepsilon_0\le v_1<v_2\le\bar t+3\varepsilon_0$,
	\begin{align}
		\Leb\{y:v_1\le{\mathsf p}_\varphi(y)\le v_2\}&\le\frac{\bar C}{2\kappa_0}(v_2-v_1),\notag\\
		{\mathsf P}_\varphi\{y:v_1<{\mathsf p}_\varphi(y)<v_2\}&\ge c_\natural(v_2-v_1).\label{eq:lb-band}
	\end{align}
\item (Level) If $\delta_\varphi<2\varepsilon_0$, then $|t_\varphi-\bar t|\le\delta_\varphi$.
\item (Membership) Suppose, in addition, that \cref{ass:constants}(ii) holds. Let $(\varphi_x)_{x\in[0,1]^{d_X}}$ be small perturbations with $\Delta_{\varphi_x}\le c_0/2$ for every $x$, and suppose that $\norm[0]{(x,y)\mapsto\varphi_x(y)}_{\mathcal H^\beta}\le L/2$. Then the distribution ${\mathsf P}$ with $\Pi$ uniform on $[0,1]^{d_X}$ and ${\mathsf p}(\cdot\mid x)=f_0+\varphi_x$ belongs to $\cP_\beta$.
\end{enumerate}
\end{lemma}

\begin{proof}
Part (i) follows from bounds along the rays $r\mapsto r\theta$, integrated in polar coordinates. Part (ii) compares $G_\varphi$ with $G_0$ of \eqref{eq:lb-G0-def}, which decreases at rate at least $c_\natural$ near $\bar t$. Part (iii) checks the assumptions defining $\cP_\beta$ at each $x$, using (i) and (ii).

\emph{(i)} We fix $\theta\in\mathbb S^{d_Y-1}$ and let $\bar f_\theta(r):={\mathsf p}_\varphi(r\theta)$. By (L2), \eqref{eq:lb-admissible} and $2\kappa_0\le C_{\bar f'}$, we have, on $[r_-,r_+]$,
	\begin{align*}
		\bar f_\theta'(r)=\bar f'(r)+\theta\cdot\nabla\varphi(r\theta)\in[-C_{\bar f'}-\kappa_0,\,-2\kappa_0+\kappa_0]\subset[-2C_{\bar f'},\,-\kappa_0].
	\end{align*}
Since $|\varphi|\le\varepsilon_0/2$, it follows from (L2) and continuity that
	\begin{align*}
		\bar f_\theta(r_-)\ge\bar t+\frac{7\varepsilon_0}2,\qquad\bar f_\theta(r_+)\le\bar t-\frac{7\varepsilon_0}2,\qquad|{\mathsf p}_\varphi-\bar t|>\frac{7\varepsilon_0}2\ \text{on }\cY\setminus\mathcal A.
	\end{align*}
Let $\bar t-3\varepsilon_0\le v_1<v_2\le\bar t+3\varepsilon_0$ and $I_\theta:=\bar f_\theta^{-1}([v_1,v_2])$. By the two displays above, $\bar f_\theta$ is strictly decreasing on $[r_-,r_+]$ and takes no value in $[v_1,v_2]$ outside $(r_-,r_+)$. Therefore $I_\theta$ is an interval in $(r_-,r_+)$ with $(v_2-v_1)/(2C_{\bar f'})\le\Leb(I_\theta)\le(v_2-v_1)/\kappa_0$, and
	\begin{align*}
		\{y:v_1\le{\mathsf p}_\varphi(y)\le v_2\}=\{r\theta:\theta\in\mathbb S^{d_Y-1},\ r\in I_\theta\}.
	\end{align*}
Integrating $\d y=r^{d_Y-1}\,\d r\,\d S_{d_Y-1}(\theta)$ over these intervals, with $r_-^{d_Y-1}\le r^{d_Y-1}\le r_+^{d_Y-1}$, and using ${\mathsf p}_\varphi\ge{\mathsf p}_\flat$, we obtain
	\begin{align*}
		\Leb\{y:v_1\le{\mathsf p}_\varphi(y)\le v_2\}&\le|\mathbb S^{d_Y-1}|\,r_+^{d_Y-1}\,\frac{v_2-v_1}{\kappa_0}=\frac{\bar C}{2\kappa_0}(v_2-v_1),\\
		{\mathsf P}_\varphi\{y:v_1<{\mathsf p}_\varphi(y)<v_2\}&\ge{\mathsf p}_\flat|\mathbb S^{d_Y-1}|\,r_-^{d_Y-1}\,\frac{v_2-v_1}{2C_{\bar f'}}=c_\natural(v_2-v_1),
	\end{align*}
where the equalities follow from the definitions of $\bar C$ and $c_\natural$ in \eqref{eq:lb-thresholds}.

\emph{(ii)} Recall from the proof of \cref{lem:constants-suffice} (\emph{Level}) that $G_0(\bar t)=\gamma_\star$. By \eqref{eq:lb-band} with $\varphi=0$, we have, for $0\le\delta'\le2\varepsilon_0$,
	\begin{align}
		G_0(\bar t-\delta')-\gamma_\star&={\mathsf P}_0\{y:\bar t-\delta'\le f_0(y)<\bar t\}\ge c_\natural\delta',\notag\\
		\gamma_\star-G_0(\bar t+\delta')&={\mathsf P}_0\{y:\bar t\le f_0(y)<\bar t+\delta'\}\ge c_\natural\delta' .\label{eq:lb-G0}
	\end{align}
Since $\delta_\varphi<2\varepsilon_0$, we can choose $\delta'\in(\delta_\varphi,2\varepsilon_0]$. By \eqref{eq:lb-G0}, the definition \eqref{eq:lb-Delta} of $\Delta_\varphi$ and $c_\natural\delta_\varphi=2\Delta_\varphi$, we have
	\begin{align*}
		G_\varphi(\bar t-\delta_\varphi)&\ge\gamma_\star+c_\natural\delta_\varphi-\Delta_\varphi\ge\gamma_\star,\\
		G_\varphi(\bar t+\delta')&\le\gamma_\star-c_\natural\delta'+\Delta_\varphi<\gamma_\star,
	\end{align*}
and hence $t_\varphi\ge\bar t-\delta_\varphi$. Since $G_\varphi$ is nonincreasing, the second bound implies $t_\varphi\le\bar t+\delta'$. Letting $\delta'\downarrow\delta_\varphi$, we obtain $t_\varphi\le\bar t+\delta_\varphi$.

\emph{(iii)} We check the conditions at a fixed $x$, with $\varphi:=\varphi_x$.
\begin{itemize}
\item \cref{ass:smooth}. By (L1), the assumed bound on $\norm[0]{(x,y)\mapsto\varphi_x(y)}_{\mathcal H^\beta}$ and subadditivity, we have $\norm[0]{\mathsf p}_{\mathcal H^\beta}\le L$. Since $\varepsilon_0\le{\mathsf p}_{\max}/4$, we also have ${\mathsf p}_\varphi\le{\mathsf p}_{\max}/2+\varepsilon_0/2\le{\mathsf p}_{\max}$.
\item \cref{ass:margin}. Since $\Delta_\varphi\le c_0/2$ and $c_0=c_m\varepsilon_0/8\le c_\natural\varepsilon_0/16$ by \cref{ass:constants}(i), we have $\delta_\varphi\le c_0/c_\natural\le\varepsilon_0/16$, and hence $|t_\varphi-\bar t|\le\varepsilon_0/16$ by (ii). For the band conditions, let $|s-t_\varphi|\le\varepsilon_0$ and $\eps\in(0,\varepsilon_0]$. Then $[s-\eps,s+\eps]\subset[\bar t-3\varepsilon_0,\bar t+3\varepsilon_0]$. Hence, by (i) and \cref{ass:constants}(i), we have
	\begin{align*}
		\Leb\{y:|{\mathsf p}_\varphi(y)-s|\le\eps\}&\le\frac{\bar C}{\kappa_0}\eps\le C_m\eps\le C_m\eps\,\mathfrak L(\eps),\\
		{\mathsf P}_\varphi\{y:s+a\eps\le{\mathsf p}_\varphi(y)<s+\eps\}&\ge(1-a)c_\natural\eps\ge c_m\eps,\\
		{\mathsf P}_\varphi\{y:s-\eps<{\mathsf p}_\varphi(y)\le s-a\eps\}&\ge c_m\eps .
	\end{align*}
\item \cref{ass:design}. By \cref{ass:constants}(ii), $\cX_0$ is a cube and $\cX_0^{2\delta_0}\subset[0,1]^{d_X}$, and hence $f_X\equiv1$ on $\cX_0^{2\delta_0}$. Since $f_{\min}\le1$ by \cref{ass:constants}(i) and $L\ge L_0\ge1$, the uniform design satisfies \cref{ass:design}.
\item \cref{ass:constants}. Part (i) is assumed. For part (ii), $\Pi$ is uniform on $[0,1]^{d_X}$, and $\cX_0$ is the cube of the conventions.\qedhere
\end{itemize}
\end{proof}

\begin{lemma}[Mass of a small perturbation]\label{lem:lb-mass}
Suppose that \cref{ass:constants}(i) holds, and let $\varphi$ be a small perturbation. Then
	\begin{align}\label{eq:lb-Delta-bound}
		\Delta_\varphi\le\norm\varphi_1+\frac{3\bar t\bar C}{\kappa_0}\norm\varphi_\infty .
	\end{align}
\end{lemma}

\begin{proof}
Let $|s-\bar t|\le2\varepsilon_0$. By \eqref{eq:sandwich} with $g={\mathsf p}_\varphi$, $g'=f_0$ and $s'=s$,
	\begin{align*}
		\{y:{\mathsf p}_\varphi(y)\ge s\}\triangle\{y:f_0(y)\ge s\}\subset\{y:|f_0(y)-s|\le\norm\varphi_\infty\}.
	\end{align*}
Since $\norm[0]\varphi_\infty\le\varepsilon_0/2$, it follows from \cref{lem:lb-perturb}(i) with $\varphi=0$ that $\Leb\{y:|f_0(y)-s|\le\norm[0]\varphi_\infty\}\le\bar C\norm[0]\varphi_\infty/\kappa_0$. On this set, $f_0\le s+\norm[0]\varphi_\infty\le\bar t+3\varepsilon_0\le3\bar t$, since $2\varepsilon_0\le\bar t$. We decompose
	\begin{align*}
		G_\varphi(s)-G_0(s)=\int\bigl(\ind\{{\mathsf p}_\varphi\ge s\}-\ind\{f_0\ge s\}\bigr)f_0\,\d y+\int_{\{y:{\mathsf p}_\varphi(y)\ge s\}}\varphi\,\d y.
	\end{align*}
By the first display of this proof, the first integrand vanishes outside $\{y:|f_0(y)-s|\le\norm[0]\varphi_\infty\}$. Therefore
	\begin{align*}
		|G_\varphi(s)-G_0(s)|\le3\bar t\,\Leb\{y:|f_0(y)-s|\le\norm\varphi_\infty\}+\norm\varphi_1\le\frac{3\bar t\bar C}{\kappa_0}\norm\varphi_\infty+\norm\varphi_1 .
	\end{align*}
Taking the supremum over $|s-\bar t|\le2\varepsilon_0$, we obtain \eqref{eq:lb-Delta-bound}.
\end{proof}

\subsection{Tools and covariate cells}\label{app:lower-tools}
The proofs of parts (i) and (ii) of \cref{thm:lower} use the following form of Assouad's lemma \citep{yu1997assouad,tsybakov2009introduction}. In this lemma and below, $\hC$ is a predictive set as in \cref{thm:lower}, and we write $\bP_\omega:=\bP_{{\mathsf P}_\omega}$ and $\E_\omega:=\E_{{\mathsf P}_\omega}$.

\begin{lemma}[Assouad's lemma]\label{lem:lb-assouad}
Let $J\ge1$, let $\delta_{\rm sep}>0$, and let ${\mathsf P}_\omega$, $\omega\in\{0,1\}^{J}$, be distributions. For $i\le J$ and $b\in\{0,1\}$, let $\mathrm{Err}_i(\hC,b)\ge0$ be functions of the predictive set such that $\mathrm{Err}_i(\hC,0)+\mathrm{Err}_i(\hC,1)\ge\delta_{\rm sep}$ for every value of $\hC$. Suppose that $\TV({\mathsf P}_\omega^{\otimes n},{\mathsf P}_{\omega'}^{\otimes n})\le1/4$ for all $\omega,\omega'$ that differ in exactly one coordinate. Then
	\begin{align}\label{eq:lb-assouad}
		2^{-J}\sum_{\omega\in\{0,1\}^J}\E_\omega\sbr{\sum_{i\le J}\mathrm{Err}_i(\hC,\omega_i)}\ge\frac38J\delta_{\rm sep} .
	\end{align}
The same holds for a family indexed by $\{-1,1\}^J$, with $b\in\{-1,1\}$ in place of $b\in\{0,1\}$ and the sum over $\omega\in\{-1,1\}^J$.
\end{lemma}

\begin{proof}
We fix $i\le J$ and $\omega$ with $\omega_i=1$, and let $\omega'$ be $\omega$ with $\omega_i$ replaced by $0$. Since the functions $\mathrm{Err}_i(\hC,b)$ are nonnegative, integrating the inequality $\mathrm{Err}_i(\hC,0)+\mathrm{Err}_i(\hC,1)\ge\delta_{\rm sep}$ against the minimum of the densities of ${\mathsf P}_\omega^{\otimes n}$ and ${\mathsf P}_{\omega'}^{\otimes n}$, we obtain
	\begin{align*}
		\E_\omega\sbr{\mathrm{Err}_i(\hC,1)}+\E_{\omega'}\sbr{\mathrm{Err}_i(\hC,0)}\ge\cbr{1-\TV({\mathsf P}_\omega^{\otimes n},{\mathsf P}_{\omega'}^{\otimes n})}\delta_{\rm sep}\ge\frac34\delta_{\rm sep}.
	\end{align*}
If $\hC$ also depends on auxiliary randomness, the same bound holds with the joint distributions of the data and this randomness, which have the same total-variation distance. Summing the last display over the $2^{J-1}$ such pairs $(\omega,\omega')$ for each $i$, then over $i\le J$, and dividing by $2^J$, we obtain \eqref{eq:lb-assouad}. The relabeling $b\mapsto(1+b)/2$ yields the version for $\{-1,1\}^J$.
\end{proof}

\begin{lemma}[Kullback--Leibler divergence]\label{lem:lb-kl}
Suppose that \cref{ass:constants}(i) holds. Let ${\mathsf P}$ and ${\mathsf P}'$ be distributions with $\Pi$ uniform on $[0,1]^{d_X}$ and conditional densities ${\mathsf p}(\cdot\mid x)=f_0+\varphi_x$ and ${\mathsf p}'(\cdot\mid x)=f_0+\varphi_x'$, where $\varphi_x$ and $\varphi_x'$ are small perturbations for every $x$. Then
	\begin{align*}
		{\rm KL}\bigl({\mathsf P}^{\otimes n}\,\big\|\,{\mathsf P}'^{\otimes n}\bigr)&\le n\int_{[0,1]^{d_X}}\int_\cY\frac{\{({\mathsf p}-{\mathsf p}')(y\mid x)\}^2}{{\mathsf p}'(y\mid x)}\,\d y\,\d x\\
		&\le\frac n{{\mathsf p}_\flat}\int_{[0,1]^{d_X}}\int_\cY(\varphi_x-\varphi_x')^2\,\d y\,\d x .
	\end{align*}
If, moreover, ${\rm KL}({\mathsf P}^{\otimes n}\,\|\,{\mathsf P}'^{\otimes n})\le1/8$, then $\TV({\mathsf P}^{\otimes n},{\mathsf P}'^{\otimes n})\le1/4$.
\end{lemma}

\begin{proof}
Since $\Pi$ has density one on $[0,1]^{d_X}$,
	\begin{align*}
		{\rm KL}\bigl({\mathsf P}^{\otimes n}\,\big\|\,{\mathsf P}'^{\otimes n}\bigr)=n\int_{[0,1]^{d_X}}\int_\cY{\mathsf p}(y\mid x)\log\frac{{\mathsf p}(y\mid x)}{{\mathsf p}'(y\mid x)}\,\d y\,\d x .
	\end{align*}
Since $\log u\le u-1$ and $\int_\cY({\mathsf p}-{\mathsf p}')(y\mid x)\,\d y=0$, the inner integral is at most $\int_\cY\{({\mathsf p}-{\mathsf p}')(y\mid x)\}^2/{\mathsf p}'(y\mid x)\,\d y$, which gives the first inequality. The second follows from ${\mathsf p}(\cdot\mid x)-{\mathsf p}'(\cdot\mid x)=\varphi_x-\varphi_x'$ and ${\mathsf p}'(\cdot\mid x)\ge{\mathsf p}_\flat$, which holds since $\varphi_x'$ is a small perturbation. Finally, by Pinsker's inequality, $\TV({\mathsf P}^{\otimes n},{\mathsf P}'^{\otimes n})\le({\rm KL}/2)^{1/2}\le1/4$ when ${\rm KL}\le1/8$.
\end{proof}

The perturbations in the hypotheses of parts (i) and (ii) are sums of rescaled bumps with disjoint supports. The next lemma bounds their H\"older norms.

\begin{lemma}[H\"older norms of bumps]\label{lem:lb-holder}
The following hold.
\begin{enumerate}[label=(\roman*)]
\item Let $g_1,\dots,g_N\in C_c^\infty(\R^d)$ have pairwise disjoint supports. Then $\norm[0]{\sum_{i\le N}g_i}_{\mathcal H^\beta}\le2\max_{i\le N}\norm[0]{g_i}_{\mathcal H^\beta}$.
\item Let $\chi_{\rm gen}\in C_c^\infty(\R^d)$, $z_0\in\R^d$, $\mu>0$ and $m\ge1$. Then the function $z\mapsto\mu\chi_{\rm gen}(m(z-z_0))$ has $\norm[0]\cdot_{\mathcal H^\beta}$-norm at most $\mu m^\beta\norm[0]{\chi_{\rm gen}}_{\mathcal H^\beta_{\R^d}}$.
\item Let $\chi_1\in C_c^\infty(\R^{d_X})$ and $\chi_2\in C_c^\infty(\R^{d_Y})$. Then there exists a constant $C_H'>0$, depending only on $\chi_1$, $\chi_2$ and $\beta$, such that, for all $x_0\in\R^{d_X}$, $\mu>0$ and $m\ge1$, the function $(x,y)\mapsto\mu\chi_1(m(x-x_0))\chi_2(y)$ has $\norm[0]\cdot_{\mathcal H^\beta}$-norm at most $C_H'\mu m^\beta$.
\end{enumerate}
\end{lemma}

\begin{proof}
\emph{(i)} Let $g:=\sum_{i\le N}g_i$. Since the supports are disjoint, the bounds on the derivatives of $g$ of order at most $\ell$ are the maxima of those of the summands. Let $|\boldsymbol k|=\ell$ and $z,z'\in[0,1]^{d_X}\times\cY$. At most two summands have $z$ or $z'$ in their supports, and the other terms of
	\begin{align*}
		\partial^{\boldsymbol k}g(z)-\partial^{\boldsymbol k}g(z')=\sum_{i\le N}\bigl\{\partial^{\boldsymbol k}g_i(z)-\partial^{\boldsymbol k}g_i(z')\bigr\}
	\end{align*}
vanish. Hence the H\"older seminorm of $g$ is at most twice the largest H\"older seminorm of the summands. This proves (i).

\emph{(ii)} By the chain rule, the derivative $\partial^{\boldsymbol k}$ of $z\mapsto\mu\chi_{\rm gen}(m(z-z_0))$ is $\mu m^{|\boldsymbol k|}(\partial^{\boldsymbol k}\chi_{\rm gen})(m(z-z_0))$. Since $m\ge1$ and $|\boldsymbol k|\le\ell<\beta$, it is bounded by $\mu m^\beta\norm[0]{\partial^{\boldsymbol k}\chi_{\rm gen}}_\infty$. For $|\boldsymbol k|=\ell$, its H\"older quotient of order $\beta-\ell$ is at most $\mu m^\ell m^{\beta-\ell}=\mu m^\beta$ times that of $\partial^{\boldsymbol k}\chi_{\rm gen}$. This proves (ii).

\emph{(iii)} The function $(x,y)\mapsto\chi_1(x)\chi_2(y)$ is in $C_c^\infty(\R^d)$, and the map $(x,y)\mapsto(m(x-x_0),y)$ multiplies each derivative in $x$ by $m$ and distances by at most $m$. The argument of (ii) therefore applies, and (iii) holds with $C_H'$ the norm $\norm[0]\cdot_{\mathcal H^\beta_{\R^d}}$ of $(x,y)\mapsto\chi_1(x)\chi_2(y)$.
\end{proof}

\emph{Covariate cells.} Throughout the proof of \cref{thm:lower}, we fix $\chi_X$ such that
	\begin{align}\label{eq:lb-chiX}
		\chi_X\in C_c^\infty\Bigl(\Bigl(-\frac12,\frac12\Bigr)^{d_X}\Bigr),\qquad0\le\chi_X\le1,\qquad\chi_X=1\ \text{on }\Bigl[-\frac14,\frac14\Bigr]^{d_X}.
	\end{align}
Let $m\ge1/w_0$ be an integer. Since $\cX_0$ is a cube of side length $w_0$, it contains $J_X\ge c_Jm^{d_X}$ cubes $Q_j$ with disjoint interiors, where $c_J:=(w_0/2)^{d_X}$ and
	\begin{align}\label{eq:lb-cells}
		Q_j:=x_j+\Bigl[-\frac1{2m},\frac1{2m}\Bigr]^{d_X},\qquad Q_j':=x_j+\Bigl[-\frac1{4m},\frac1{4m}\Bigr]^{d_X}\qquad(j\le J_X).
	\end{align}
We call the cubes $Q_j$ the cells and the cubes $Q_j'$ their inner parts, with $\Leb(Q_j')=(2m)^{-d_X}$. The function $x\mapsto\chi_X(m(x-x_j))$ has compact support in the interior of $Q_j$ and equals one on $Q_j'$.

\subsection{Proof of \texorpdfstring{\cref{thm:lower}(i)}{Theorem \ref{thm:lower}(i)}}\label{app:lower-proof}
For (i), each hypothesis perturbs $f_0$ by a sum of localized functions of amplitude $\mu\asymp m^{-\beta}$, with $m\asymp n^{1/(2\beta+d)}$. They are localized in $x$ on the cells and in $y$ near pairs of antipodal points $\pm y_k$ on the level sphere $\{y:|y|=r_\star\}$ of $f_0$, with opposite signs at $y_k$ and $-y_k$. The reflection $y\mapsto-y$ near one pair reverses these signs and leaves the level unchanged. Hence all hypotheses share one level on the inner parts of the cells (\cref{lem:lb-i-family}). The sign attached to a cell and a pair determines whether a set of volume of order $\mu m^{-(d_Y-1)}$ near $y_k$ lies in the oracle set. On this set, the density differs from the level by at least $\mu/4$ (\cref{lem:lb-i-detect}).

\emph{Bumps and points.} We fix a radial $\chi_Y$ such that
	\begin{align}\label{eq:lb-chiY}
		\chi_Y\in C_c^\infty\Bigl(\mathcal B\Bigl(0,\frac12\Bigr)\Bigr),\qquad0\le\chi_Y\le\chi_Y(0)=1.
	\end{align}
Let $r_\chi\in(0,\frac12)$ be such that $\chi_Y\ge\frac12$ on $\bar{\mathcal B}(0,r_\chi)$, and let $C_\chi:=\norm[0]{\nabla\chi_Y}_\infty$.

For an integer $m\ge1/w_0$, let $y_1,\dots,y_{k_m}$ be a maximal family of points, with pairwise distances at least $1/m$, in the spherical cap
	\begin{align*}
		\Bigl\{y:|y|=r_\star,\ \langle y,\mathbf e_1\rangle\ge\frac1{2m}\Bigr\},
	\end{align*}
where $\mathbf e_1$ is the first canonical basis vector. Since $|y_k+y_l|\ge\langle y_k+y_l,\mathbf e_1\rangle\ge1/m$, the $2k_m$ points $\pm y_k$ are pairwise at distance at least $1/m$, and the balls $\mathcal B(\pm y_k,\frac1{2m})$ are disjoint. By maximality, the balls $\mathcal B(y_k,1/m)$ cover this cap. Each of them meets the sphere $\{y:|y|=r_\star\}$ in a set of surface measure at most $Cm^{-(d_Y-1)}$. For all sufficiently large $m$, the cap has surface measure at least half that of the hemisphere $\{y:|y|=r_\star,\langle y,\mathbf e_1\rangle\ge0\}$. Consequently, there exists a constant $c_{\mathbb S}>0$ such that $k_m\ge c_{\mathbb S}m^{d_Y-1}$ for all sufficiently large $m$ ($k_m=1$ if $d_Y=1$).

\emph{Constants.} Let $(\chi_X\otimes\chi_Y)(x,y):=\chi_X(x)\chi_Y(y)$, and let
	\begin{align*}
		C_H:=\norm{\chi_X\otimes\chi_Y}_{\mathcal H^\beta_{\R^d}},\qquad C_\Delta:=\frac{3\bar t\bar C}{\kappa_0}+\Leb(\mathcal A),\qquad C_{\rm KL}:=\frac{8\norm{\chi_X}_2^2\norm{\chi_Y}_2^2}{{\mathsf p}_\flat}.
	\end{align*}
We fix $c_\mu\in(0,1]$ with
	\begin{align}\label{eq:lb-ca}
		4C_Hc_\mu\le1,\qquad C_\chi c_\mu L\le \kappa_0,\qquad(c_\natural+8C_\Delta)c_\mu L\le4c_\natural r_\chi \kappa_0,\qquad8C_{\rm KL}c_\mu^2L^2\le1 .
	\end{align}
We then fix an integer $m_0\ge1/w_0$ such that, for all $m\ge m_0$,
	\begin{align*}
		&c_\mu L/m\le\varepsilon_0/2,\qquad1/(2m)\le\min\{r_\star-r_-,\,r_+-r_\star\},\\
		&k_m\ge c_{\mathbb S}m^{d_Y-1},\qquad C_\Delta c_\mu L/m\le c_0/2.
	\end{align*}

\emph{The hypotheses.} For an integer $m\ge m_0$, let $\mu:=c_\mu Lm^{-\beta}$ and, for $k\le k_m$,
	\begin{align*}
		\mathrm{dip}_k(y):=\chi_Y\bigl(m(y-y_k)\bigr)-\chi_Y\bigl(m(y+y_k)\bigr).
	\end{align*}
For $\omega\in\{-1,1\}^{J_X\times k_m}$, let ${\mathsf P}_\omega$ be the distribution with $\Pi$ uniform on $[0,1]^{d_X}$ and
	\begin{align*}
		{\mathsf p}_\omega(y\mid x):=f_0(y)+\varphi_{\omega,x}(y),\qquad\varphi_{\omega,x}(y):=\mu\sum_{j\le J_X,\,k\le k_m}\omega_{jk}\chi_X\bigl(m(x-x_j)\bigr)\mathrm{dip}_k(y).
	\end{align*}
We write $t^{(\omega)}_\star(x)$ and ${\mathsf C}^{\star,(\omega)}(x)$ for the level and the oracle set of ${\mathsf p}_\omega(\cdot\mid x)$. Let
	\begin{align}\label{eq:lb-phim}
		\varphi_m:=\mu\sum_{k\le k_m}\mathrm{dip}_k,\qquad t_m:=t_{\varphi_m},
	\end{align}
with $t_\varphi$ as in \eqref{eq:lb-pphi}. Since $\beta\ge1$,
	\begin{align}\label{eq:lb-am}
		\mu m=c_\mu Lm^{1-\beta}\le c_\mu L,\qquad\mu\le c_\mu L/m .
	\end{align}
For $m\ge m_0$, since $1/(2m)\le\min\{r_\star-r_-,\,r_+-r_\star\}$ by the choice of $m_0$, the disjoint balls $\mathcal B(\pm y_k,\frac1{2m})$ lie in the annulus $\mathcal A\subset\cY$ of \eqref{eq:lb-annulus}.

\begin{lemma}[Hypotheses for part (i)]\label{lem:lb-i-family}
Suppose that \cref{ass:constants} holds. Let $m\ge m_0$ be an integer, and let ${\mathsf P}_\omega$, $\omega\in\{-1,1\}^{J_X\times k_m}$, be the distributions defined above. Then the following hold.
\begin{enumerate}[label=(\roman*)]
\item For every $\omega$ and $x$, the function $\varphi_{\omega,x}$ is a small perturbation, and
	\begin{align}\label{eq:lb-firstorder}
		\Delta_{\varphi_{\omega,x}}\le C_\Delta\mu\le c_0/2 .
	\end{align}
\item For every $\omega$, we have $\norm[0]{(x,y)\mapsto\varphi_{\omega,x}(y)}_{\mathcal H^\beta}\le L/2$ and ${\mathsf P}_\omega\in\cP_\beta$.
\item For every $\omega$, $j\le J_X$ and $x\in Q_j'$,
	\begin{align}\label{eq:lb-common}
		t^{(\omega)}_\star(x)=t_m .
	\end{align}
Moreover,
	\begin{align}\label{eq:lb-tm}
		|t_m-\bar t|\le\frac{2C_\Delta}{c_\natural}\mu\le\frac{\varepsilon_0}8 .
	\end{align}
\end{enumerate}
\end{lemma}

\begin{proof}
We first show that each $\varphi_{\omega,x}$ is a small perturbation, then bound the H\"older norm of $(x,y)\mapsto\varphi_{\omega,x}(y)$, and finally use the reflection symmetry of $f_0$ to show that the level does not depend on $\omega$ for $x$ in the inner parts of the cells.

\emph{Step 1: small perturbations.} We fix $\omega$ and $x$, and let $\varphi:=\varphi_{\omega,x}$. If $x$ lies in the interior of no cell, then $\varphi=0$. Otherwise $x$ lies in the interior of exactly one cell $Q_j$, and $\varphi=\mu\sum_{k\le k_m}b_k\,\mathrm{dip}_k$ with $b_k:=\omega_{jk}\chi_X(m(x-x_j))\in[-1,1]$. Since the balls $\mathcal B(\pm y_k,\frac1{2m})$ are disjoint subsets of $\mathcal A\subset\cY$, we have $\int\mathrm{dip}_k\,\d y=0$. Since the functions $\chi_Y(m(\cdot\mp y_k))$ have disjoint supports and $|b_k|\le1$, we also have $|\varphi|\le\mu$ and $|\nabla\varphi|\le C_\chi\mu m$. By \eqref{eq:lb-am}, \eqref{eq:lb-ca} and the choice of $m_0$, we have
	\begin{align*}
		\mu\le c_\mu L/m\le\varepsilon_0/2,\qquad C_\chi\mu m\le C_\chi c_\mu L\le \kappa_0.
	\end{align*}
Hence $\varphi$ is a small perturbation. Since $\varphi$ is supported in $\mathcal A$ and $|\varphi|\le\mu$, we have $\norm[0]\varphi_1\le\mu\Leb(\mathcal A)$ and $\norm[0]\varphi_\infty\le\mu$. By \cref{lem:lb-mass}, we obtain
	\begin{align*}
		\Delta_\varphi\le\mu\Leb(\mathcal A)+\frac{3\bar t\bar C}{\kappa_0}\mu=C_\Delta\mu\le\frac{C_\Delta c_\mu L}m\le\frac{c_0}2,
	\end{align*}
where the last two inequalities follow from \eqref{eq:lb-am} and the choice of $m_0$. This proves (i).

\emph{Step 2: membership.} By \cref{lem:lb-holder}(ii), each of the $2J_Xk_m$ functions $(x,y)\mapsto\mu\chi_X(m(x-x_j))\chi_Y(m(y\mp y_k))$, $j\le J_X$, $k\le k_m$, has H\"older norm at most $C_H\mu m^\beta$. These functions have disjoint compact supports, since the cells have disjoint interiors and the balls $\mathcal B(\pm y_k,\frac1{2m})$ are disjoint. The function $(x,y)\mapsto\varphi_{\omega,x}(y)$ is a sum of these functions multiplied by $\pm1$. Hence, by \cref{lem:lb-holder}(i), we have
	\begin{align}\label{eq:lb-holder}
		\norm{(x,y)\mapsto\varphi_{\omega,x}(y)}_{\mathcal H^\beta}\le2C_H\mu m^\beta=2C_Hc_\mu L\le\frac L2,
	\end{align}
where the last inequality follows from $4C_Hc_\mu\le1$. By (i) and \eqref{eq:lb-holder}, the conditions of \cref{lem:lb-perturb}(iii) hold. Therefore ${\mathsf P}_\omega\in\cP_\beta$, which proves (ii).

\emph{Step 3: a common level.} We fix $\omega$ and $x$ in the interior of a cell $Q_j$, with $\varphi$ and $b_k$ as in Step 1. For $k\le k_m$, let $R_k:\cY\to\cY$ be the map $y\mapsto-y$ on $\mathcal B(y_k,\frac1{2m})\cup\mathcal B(-y_k,\frac1{2m})$ and the identity elsewhere. Since $f_0$ is radial and $\chi_Y$ is even, and since $\mathrm{dip}_{k'}$, $k'\ne k$, vanishes on these balls, we have $f_0\circ R_k=f_0$, $\mathrm{dip}_k\circ R_k=-\mathrm{dip}_k$ and $\mathrm{dip}_{k'}\circ R_k=\mathrm{dip}_{k'}$ for $k'\ne k$. Hence, with ${\mathsf p}_\varphi$ and $G_\varphi$ as in \eqref{eq:lb-pphi}, we have ${\mathsf p}_\varphi\circ R_k={\mathsf p}_{\varphi'}$, where $\varphi'$ is $\varphi$ with $b_k$ replaced by $-b_k$. Since $R_k$ is a bijection that preserves Lebesgue measure, it follows that, for every $s$,
	\begin{align*}
		G_{\varphi'}(s)=\int_\cY{\mathsf p}_\varphi(R_ky)\ind\{{\mathsf p}_\varphi(R_ky)\ge s\}\,\d y=\int_\cY{\mathsf p}_\varphi\ind\{{\mathsf p}_\varphi\ge s\}\,\d y=G_\varphi(s).
	\end{align*}
Replacing each negative $b_k$ by $|b_k|$ in turn shows that $G_\varphi$, and hence the level $t_\varphi$, are those of $\mu\chi_X(m(x-x_j))\sum_{k\le k_m}\mathrm{dip}_k$, which do not depend on $\omega$. For $x\in Q_j'$, where $\chi_X(m(x-x_j))=1$, this function is $\varphi_m$ in \eqref{eq:lb-phim}. Therefore $t^{(\omega)}_\star(x)=t_{\varphi_{\omega,x}}=t_m$, which proves \eqref{eq:lb-common}.

Next, $\varphi_m=\varphi_{\omega,x}$ for $x\in Q_j'$ and any $\omega$ with $\omega_{jk}=1$ for all $k$. Thus $\Delta_{\varphi_m}\le C_\Delta\mu\le c_0/2$ by \eqref{eq:lb-firstorder}, and hence, by \cref{ass:constants}(i),
	\begin{align*}
		\delta_{\varphi_m}=2\Delta_{\varphi_m}/c_\natural\le2C_\Delta\mu/c_\natural\le c_0/c_\natural\le\varepsilon_0/8.
	\end{align*}
By \cref{lem:lb-perturb}(ii), we therefore have $|t_m-\bar t|\le\delta_{\varphi_m}$. Together with the last display, this proves \eqref{eq:lb-tm}.
\end{proof}

\emph{Detection sets.} For an integer $m\ge m_0$ and $k\le k_m$, let $\theta_k:=y_k/r_\star$ and
	\begin{align*}
		\Theta_k&:=\Bigl\{\theta\in\mathbb S^{d_Y-1}:|\theta-\theta_k|\le\frac{r_\chi}{2mr_\star}\Bigr\},\qquad I_m:=\Bigl\{r\in[r_-,r_+]:|\bar f(r)-t_m|\le\frac\mu4\Bigr\},\\
		T_k&:=\{r\theta:\theta\in\Theta_k,\ r\in I_m\}.
	\end{align*}
For $A\subset\cY$, a predictive set $\hC$, $j\le J_X$, $k\le k_m$ and $b\in\{-1,1\}$, let
	\begin{align*}
		&\mathrm{Err}_k(A,1):=\Leb(T_k\setminus A),\qquad\mathrm{Err}_k(A,-1):=\Leb(T_k\cap A),\\
		&\mathrm{Err}_{jk}(\hC,b):=\int_{Q_j'}\mathrm{Err}_k\bigl(\hC(x),b\bigr)\,\d x .
	\end{align*}
Let $\cE_\omega(\hC(x))$ be the excess mass of $\hC(x)$ under ${\mathsf P}_\omega$.

\begin{lemma}[Detection sets for part (i)]\label{lem:lb-i-detect}
Suppose that \cref{ass:constants} holds. Let $m\ge m_0$ be an integer. Then the following hold.
\begin{enumerate}[label=(\roman*)]
\item The sets $T_k$, $k\le k_m$, are disjoint subsets of $\mathcal A$, with $T_k\subset\bar{\mathcal B}(y_k,r_\chi/m)\subset\mathcal B(y_k,\frac1{2m})$ and
	\begin{align}\label{eq:lb-Kvol}
		\Leb(T_k)\ge c_T\mu\,m^{-(d_Y-1)},
	\end{align}
where $c_T>0$ depends only on $d_Y$, $r_\chi$, $r_\star$, $r_-$ and $C_{\bar f'}$.
\item For every $\omega$, $j\le J_X$, $k\le k_m$ and $x\in Q_j'$,
	\begin{align}\label{eq:lb-sep}
		T_k\cap{\mathsf C}^{\star,(\omega)}(x)=\begin{cases}T_k,&\omega_{jk}=1,\\\emptyset,&\omega_{jk}=-1,\end{cases}\qquad|{\mathsf p}_\omega(\cdot\mid x)-t_m|\ge\frac\mu4\ \text{on }T_k .
	\end{align}
\item For every $\omega$ and every value of $\hC$,
	\begin{align}
		\int_{\cX_0}\Leb\bigl(\hC(x)\triangle{\mathsf C}^{\star,(\omega)}(x)\bigr)\d x&\ge\sum_{j\le J_X,\,k\le k_m}\mathrm{Err}_{jk}(\hC,\omega_{jk}),\notag\\
		\int_{\cX_0}\cE_\omega(\hC(x))\,\d x&\ge\frac\mu{8\bar t}\sum_{j\le J_X,\,k\le k_m}\mathrm{Err}_{jk}(\hC,\omega_{jk}).\label{eq:lb-loss}
	\end{align}
Moreover, for every $j\le J_X$, $k\le k_m$ and every value of $\hC$,
	\begin{align}\label{eq:lb-percoord}
		\mathrm{Err}_{jk}(\hC,-1)+\mathrm{Err}_{jk}(\hC,1)=\Leb(Q_j')\Leb(T_k)\ge2^{-d_X}c_T\mu\,m^{-(d-1)}.
	\end{align}
\end{enumerate}
\end{lemma}

\begin{proof}
We first locate the sets $T_k$ near the points $y_k$, then compare ${\mathsf p}_\omega(\cdot\mid x)$ with $t_m$ on them, and finally sum the resulting bounds over the sets $T_k$.

\emph{Step 1: geometry of the detection sets.} By \eqref{eq:lb-tm} and $\mu\le\varepsilon_0/2$, which follows from \eqref{eq:lb-am} and the choice of $m_0$, we have
	\begin{align*}
		[t_m-\mu/4,t_m+\mu/4]\subset[\bar t-\varepsilon_0/4,\bar t+\varepsilon_0/4]\subset[\bar f(r_+),\bar f(r_-)],
	\end{align*}
where the last inclusion follows from (L2). Since $\bar f$ decreases on $[r_-,r_+]$ with $2\kappa_0\le|\bar f'|\le C_{\bar f'}$, the set $I_m$ is an interval with $\Leb(I_m)\ge\mu/(2C_{\bar f'})$. Since $I_m\subset[r_-,r_+]$, we also have $T_k\subset\mathcal A$. By the definition of $I_m$, \eqref{eq:lb-tm} and $\bar f(r_\star)=\bar t$, we have $|\bar f(r)-\bar f(r_\star)|\le\mu/4+2C_\Delta\mu/c_\natural$ for $r\in I_m$. Since $|\bar f'|\ge2\kappa_0$ on $[r_-,r_+]$, this bound implies that
	\begin{align*}
		|r-r_\star|\le\frac{(c_\natural+8C_\Delta)\mu}{8c_\natural \kappa_0}\le\frac{(c_\natural+8C_\Delta)c_\mu L}{8c_\natural \kappa_0\,m}\le\frac{r_\chi}{2m} ,
	\end{align*}
where the last two inequalities follow from \eqref{eq:lb-am} and the third condition in \eqref{eq:lb-ca}. For $r\theta\in T_k$, it follows that $|r\theta-y_k|\le|r-r_\star|+r_\star|\theta-\theta_k|\le r_\chi/m$. Therefore $T_k\subset\bar{\mathcal B}(y_k,r_\chi/m)\subset\mathcal B(y_k,\frac1{2m})$. Since the balls $\mathcal B(y_k,\frac1{2m})$ are disjoint, the sets $T_k$ are disjoint as well.

A spherical cap of Euclidean radius $r'\le1$ has $S_{d_Y-1}$-measure at least $c\,r'^{d_Y-1}$, where $c>0$ depends only on $d_Y$. The set $\Theta_k$ is such a cap with $r'=r_\chi/(2mr_\star)$, and $r'\le1$ since $r_\chi<\frac12$ and $1/(2m)\le r_\star-r_-=r_\star/2$ by the choice of $m_0$. Integrating in polar coordinates over $\Theta_k\times I_m$, we therefore obtain
	\begin{align*}
		\Leb(T_k)\ge S_{d_Y-1}(\Theta_k)\,r_-^{d_Y-1}\frac\mu{2C_{\bar f'}}\ge c_T\mu\,m^{-(d_Y-1)},
	\end{align*}
which proves (i).

\emph{Step 2: separation.} We fix $(j,k)$, $\omega$ and $x\in Q_j'$. Since $\chi_X(m(x-x_j))=1$, and since the functions $\mathrm{dip}_{k'}$, $k'\ne k$, and $\chi_Y(m(\cdot+y_k))$ vanish on $\mathcal B(y_k,\frac1{2m})$, we have ${\mathsf p}_\omega(y\mid x)=f_0(y)+\mu\omega_{jk}\chi_Y(m(y-y_k))$ there. On $T_k$, we have $|f_0-t_m|\le\frac\mu4$ by the definition of $I_m$, and $\chi_Y(m(y-y_k))\ge\frac12$ by (i) and the choice of $r_\chi$. Consequently,
	\begin{align*}
		\omega_{jk}=1:&\qquad{\mathsf p}_\omega(y\mid x)-t_m\ge\frac\mu2-\frac\mu4=\frac\mu4,\\
		\omega_{jk}=-1:&\qquad t_m-{\mathsf p}_\omega(y\mid x)\ge\frac\mu2-\frac\mu4=\frac\mu4.
	\end{align*}
Since $t^{(\omega)}_\star(x)=t_m$ by \eqref{eq:lb-common}, these bounds prove \eqref{eq:lb-sep}.

\emph{Step 3: decomposition of the loss.} We fix $\omega$, $j$ and $x\in Q_j'$. By \eqref{eq:lb-sep}, we have, for every $k$,
	\begin{align*}
		\Leb\bigl(T_k\cap(\hC(x)\triangle{\mathsf C}^{\star,(\omega)}(x))\bigr)=\mathrm{Err}_k\bigl(\hC(x),\omega_{jk}\bigr).
	\end{align*}
The sets $T_k$ are disjoint by (i), and $|{\mathsf p}_\omega(\cdot\mid x)-t_m|\ge\mu/4$ on each of them by \eqref{eq:lb-sep}. Moreover, $t^{(\omega)}_\star(x)=t_m\le2\bar t$ by \eqref{eq:lb-common}, \eqref{eq:lb-tm} and $\varepsilon_0\le\bar t$. Summing the last display over $k$ and using the definition of the excess mass, we therefore obtain
	\begin{align*}
		\Leb\bigl(\hC(x)\triangle{\mathsf C}^{\star,(\omega)}(x)\bigr)&\ge\sum_{k\le k_m}\mathrm{Err}_k\bigl(\hC(x),\omega_{jk}\bigr),\\
		\cE_\omega(\hC(x))&\ge\frac\mu{4t_m}\sum_{k\le k_m}\mathrm{Err}_k\bigl(\hC(x),\omega_{jk}\bigr)\ge\frac\mu{8\bar t}\sum_{k\le k_m}\mathrm{Err}_k\bigl(\hC(x),\omega_{jk}\bigr).
	\end{align*}
Integrating the last display over $x$ in the disjoint inner parts $Q_j'\subset\cX_0$, we obtain \eqref{eq:lb-loss}. Since $\mathrm{Err}_k(A,-1)+\mathrm{Err}_k(A,1)=\Leb(T_k)$ for every $A\subset\cY$ and $\Leb(Q_j')=(2m)^{-d_X}$, \eqref{eq:lb-percoord} follows from \eqref{eq:lb-Kvol}.
\end{proof}

\begin{proof}[Proof of \cref{thm:lower}(i)]
We apply Assouad's lemma to the hypotheses of \cref{lem:lb-i-family}, with $m$ of order $n^{1/(2\beta+d)}$ and the losses of \cref{lem:lb-i-detect}. Let
	\begin{align}\label{eq:lb-m}
		m:=\max\bigl\{\lceil n^{1/(2\beta+d)}\rceil,\ m_0\bigr\},
	\end{align}
and let ${\mathsf P}_\omega$, $\omega\in\{-1,1\}^{J_X\times k_m}$, be the distributions defined before \cref{lem:lb-i-family}. By \cref{lem:lb-i-family}(ii), each ${\mathsf P}_\omega$ lies in $\cP_\beta$.

\emph{Step 1: information.} Let $\omega,\omega'$ differ only at $(j,k)$. Then ${\mathsf p}_\omega-{\mathsf p}_{\omega'}=\pm2\mu\chi_X(m(x-x_j))\mathrm{dip}_k(y)$, where the two terms of $\mathrm{dip}_k$ have disjoint supports. By \cref{lem:lb-i-family}(i), ${\mathsf p}_\omega(\cdot\mid x)$ and ${\mathsf p}_{\omega'}(\cdot\mid x)$ are $f_0$ plus small perturbations for every $x$. Hence, by \cref{lem:lb-kl}, we have
	\begin{align}
		{\rm KL}\bigl({\mathsf P}_\omega^{\otimes n}\,\big\|\,{\mathsf P}_{\omega'}^{\otimes n}\bigr)&\le n\int_{[0,1]^{d_X}}\int_\cY\frac{({\mathsf p}_\omega-{\mathsf p}_{\omega'})^2}{{\mathsf p}_{\omega'}}\,\d y\,\d x\notag\\
		&\le\frac{4n\mu^2}{{\mathsf p}_\flat}\,m^{-d_X}\norm{\chi_X}_2^2\cdot2m^{-d_Y}\norm{\chi_Y}_2^2\notag\\
		&=C_{\rm KL}c_\mu^2L^2\,n\,m^{-(2\beta+d)}\le\frac18 ,\label{eq:lb-KL}
	\end{align}
where the last inequality follows from \eqref{eq:lb-m}, which gives $m\ge n^{1/(2\beta+d)}$, and the last condition in \eqref{eq:lb-ca}. By \cref{lem:lb-kl}, it follows that $\TV({\mathsf P}_\omega^{\otimes n},{\mathsf P}_{\omega'}^{\otimes n})\le1/4$.

\emph{Step 2: Assouad's lemma.} We apply \cref{lem:lb-assouad}, in its form for $\{-1,1\}^J$, to the $J:=J_Xk_m\ge c_Jc_{\mathbb S}m^{d-1}$ coordinates $(j,k)$, with the functions $\mathrm{Err}_{jk}$ of \cref{lem:lb-i-detect} and $\delta_{\rm sep}:=2^{-d_X}c_T\mu m^{-(d-1)}$. The required separation holds by \eqref{eq:lb-percoord}, and the total-variation distances are at most $1/4$ by Step 1. Combining \eqref{eq:lb-assouad} with \eqref{eq:lb-loss}, we obtain, for a constant $c>0$,
	\begin{align*}
		2^{-J}\sum_{\omega\in\{-1,1\}^{J_X\times k_m}}\E_\omega\sbr{\int_{\cX_0}\Leb\bigl(\hC(x)\triangle{\mathsf C}^{\star,(\omega)}(x)\bigr)\,\d x}\ge c\mu,\\
		2^{-J}\sum_{\omega\in\{-1,1\}^{J_X\times k_m}}\E_\omega\sbr{\int_{\cX_0}\cE_\omega(\hC(x))\,\d x}\ge c\mu^2 .
	\end{align*}
The supremum over $\cP_\beta$ is at least the average over $\omega$.

\emph{Step 3: the rate.} Since $\lceil n^{1/(2\beta+d)}\rceil\le2n^{1/(2\beta+d)}$, we have, for every $n\ge1$,
	\begin{align*}
		m\le2m_0n^{1/(2\beta+d)},\qquad\mu=c_\mu Lm^{-\beta}\ge c_\mu L(2m_0)^{-\beta}n^{-\beta/(2\beta+d)}.
	\end{align*}
Together with the bounds of Step 2, this proves (i).
\end{proof}

\subsection{Proof of \texorpdfstring{\cref{thm:lower}(ii)}{Theorem \ref{thm:lower}(ii)}}\label{app:lower-cov}
For (ii), each hypothesis adds to $f_0$, in each cell, either nothing or a fixed perturbation of amplitude $\mu\asymp m^{-\beta}$, now with $m\asymp n^{1/(2\beta+d_X)}$. This perturbation changes the oracle set only by a set of volume $O(\mu)$ (\cref{lem:lb-ii-family}), but it changes the mass of every set close to the oracle set of $f_0$ by an amount of order $\mu$. By \cref{lem:identity} and \eqref{eq:E-symdiff}, a set with small coverage error and small excess volume is close to the oracle set of $f_0$. Hence no predictive set has small coverage error under both hypotheses of a cell unless its volume exceeds the oracle volume by a constant (\cref{lem:lb-ii-twopoint}). The class $\mathfrak V_n(c_1,c_2)$ bounds the expected measure of the covariates at which this happens.

\emph{The hypotheses.} Let ${\mathsf C}_0:=\{y:f_0(y)\ge \bar t\}$ be the oracle set of $f_0$, and let $\mathfrak u\in C_c^\infty(\R^{d_Y})$ be supported in the interior of $\cY$, with
	\begin{align*}
		\int\mathfrak u\,\d y=0,\qquad c_{\mathfrak u}:=\int_{{\mathsf C}_0}\mathfrak u\,\d y>0 .
	\end{align*}
For example, $\mathfrak u$ can be positive on a ball inside ${\mathsf C}_0$ and negative on an annulus outside it. For an integer $m\ge m_0'$, with the cells of \eqref{eq:lb-cells}, let $\mu:=c_\mu' Lm^{-\beta}$, where the constants $c_\mu'$ and $m_0'\ge1/w_0$ are fixed below. For $\omega\in\{0,1\}^{J_X}$, let ${\mathsf P}_\omega$ be the distribution with $\Pi$ uniform on $[0,1]^{d_X}$ and
	\begin{align}\label{eq:lb-pomega-ii}
		{\mathsf p}_\omega(y\mid x):=f_0(y)+\mu\sum_{j\le J_X}\omega_j\chi_X\bigl(m(x-x_j)\bigr)\mathfrak u(y).
	\end{align}
As in \cref{app:lower-proof}, we write $t^{(\omega)}_\star(x)$ and ${\mathsf C}^{\star,(\omega)}(x)$ for the level and the oracle set of ${\mathsf p}_\omega(\cdot\mid x)$, and ${\mathsf P}_{\omega,x}$ for the distribution with density ${\mathsf p}_\omega(\cdot\mid x)$. For a predictive set $\hC$, let
	\begin{align}\label{eq:lb-cov-omega}
		\cov_\omega(x):={\mathsf P}_{\omega,x}(\hC(x)).
	\end{align}
For $j\le J_X$, $x\in Q_j'$ and $b\in\{0,1\}$, let ${\mathsf P}_{b,x}$ be the distribution with density $f_0+b\,\mu\,\mathfrak u$. Since $\chi_X(m(\cdot-x_j))=1$ on $Q_j'$ and the other summands of \eqref{eq:lb-pomega-ii} vanish there, we have ${\mathsf P}_{\omega,x}={\mathsf P}_{\omega_j,x}$ for every $\omega$, every $j\le J_X$ and every $x\in Q_j'$.

\emph{Constants.} Let $C_H'$ be the constant of \cref{lem:lb-holder}(iii) for $\chi_1=\chi_X$ and $\chi_2=\mathfrak u$, and let
	\begin{align*}
		C_{\mathfrak u}:=\norm{\mathfrak u}_1+\frac{3\bar t\bar C}{\kappa_0}\norm{\mathfrak u}_\infty,\qquad C_{\rm sd}:=\frac{\bar C(\norm{\mathfrak u}_\infty+2C_{\mathfrak u}/c_\natural)}{\kappa_0},\qquad C_{\rm KL}':=\frac{\norm{\chi_X}_2^2\norm{\mathfrak u}_2^2}{{\mathsf p}_\flat}.
	\end{align*}
Let $c_E$ be the constant $c$ of the lower bound in \eqref{eq:E-symdiff} with $\mathfrak L\equiv1$, and let
	\begin{align*}
		\delta_{\rm sd}:=\frac{c_{\mathfrak u}}{4\norm{\mathfrak u}_\infty},\qquad c_1:=\frac{c_E\delta_{\rm sd}^2}4,\qquad c_2:=2^{-d_X-9}c_J .
	\end{align*}
Note that $0<c_2<c_J\le\Leb(\cX_0)$. We choose $c_\mu'\in(0,1]$ with
	\begin{align*}
		4C_H'c_\mu'\le1,\qquad8C_{\rm KL}'c_\mu'^2L^2\le1,\qquad2(2\beta+d_X)C_{\rm KL}'c_\mu'^2L^2\le d_X,
	\end{align*}
and then an integer $m_0'\ge1/w_0$ such that, for all $m\ge m_0'$,
	\begin{align*}
		&\mu\norm{\mathfrak u}_\infty\le\varepsilon_0/2,\qquad\mu\norm{\nabla \mathfrak u}_\infty\le \kappa_0,\qquad C_{\mathfrak u}\mu\le c_0/2,\\
		&C_{\rm sd}\mu\le c_1,\qquad c_{\mathfrak u}\mu\le4\bar tc_1 .
	\end{align*}

\begin{lemma}[Hypotheses for part (ii)]\label{lem:lb-ii-family}
Suppose that \cref{ass:constants} holds. Let $m\ge m_0'$ be an integer, and let ${\mathsf P}_\omega$, $\omega\in\{0,1\}^{J_X}$, be the distributions of \eqref{eq:lb-pomega-ii}. Then the following hold.
\begin{enumerate}[label=(\roman*)]
\item For every $\omega$ and $x$, the function $\varphi:={\mathsf p}_\omega(\cdot\mid x)-f_0$ is a small perturbation with $\Delta_\varphi\le C_{\mathfrak u}\mu\le c_0/2$. Moreover, ${\mathsf P}_\omega\in\cP_\beta$ for every $\omega$.
\item For every $\omega$ and $x$, we have $|t^{(\omega)}_\star(x)-\bar t|\le2C_{\mathfrak u}\mu/c_\natural$ and
	\begin{align}\label{eq:lb-symdiff0}
		\Leb\bigl({\mathsf C}^{\star,(\omega)}(x)\triangle{\mathsf C}_0\bigr)\le C_{\rm sd}\mu .
	\end{align}
\end{enumerate}
\end{lemma}

\begin{proof}
We check the conditions of \cref{lem:lb-perturb} for ${\mathsf p}_\omega(\cdot\mid x)-f_0$ and then compare the oracle sets with ${\mathsf C}_0$.

\emph{Step 1: membership.} We fix $\omega$ and $x$. Since the supports of the functions $\chi_X(m(\cdot-x_j))$ lie in the interiors of disjoint cells and $0\le\chi_X\le1$, the function $\varphi$ equals $\mu\,b\,\mathfrak u$ with $b\in[0,1]$. Thus $\int\varphi\,\d y=0$ and
	\begin{align*}
		|\varphi|\le\mu\norm{\mathfrak u}_\infty\le\varepsilon_0/2,\qquad|\nabla\varphi|\le\mu\norm{\nabla\mathfrak u}_\infty\le \kappa_0,
	\end{align*}
where the bounds by $\varepsilon_0/2$ and $\kappa_0$ follow from the choice of $m_0'$. Hence $\varphi$ is a small perturbation. Since $\norm[0]\varphi_1\le\mu\norm[0]{\mathfrak u}_1$ and $\norm[0]\varphi_\infty\le\mu\norm[0]{\mathfrak u}_\infty$, it follows from \cref{lem:lb-mass} and the choice of $m_0'$ that $\Delta_\varphi\le C_{\mathfrak u}\mu\le c_0/2$.

Next, by \cref{lem:lb-holder}(iii), each summand $(x,y)\mapsto\mu\chi_X(m(x-x_j))\mathfrak u(y)$ of ${\mathsf p}_\omega-f_0$ has H\"older norm at most $C_H'\mu m^\beta$. These summands have disjoint compact supports. Hence, by \cref{lem:lb-holder}(i), we have
	\begin{align*}
		\norm{{\mathsf p}_\omega-f_0}_{\mathcal H^\beta}\le2C_H'\mu m^\beta=2C_H'c_\mu' L\le\frac L2,
	\end{align*}
where the last inequality follows from $4C_H'c_\mu'\le1$. Therefore the conditions of \cref{lem:lb-perturb}(iii) hold, and ${\mathsf P}_\omega\in\cP_\beta$. This proves (i).

\emph{Step 2: level and oracle sets.} We fix $\omega$ and $x$. By (i) and \cref{ass:constants}(i), we have $\delta_\varphi=2\Delta_\varphi/c_\natural\le c_0/c_\natural<2\varepsilon_0$. Hence, by \cref{lem:lb-perturb}(ii), we have $|t^{(\omega)}_\star(x)-\bar t|\le\delta_\varphi\le2C_{\mathfrak u}\mu/c_\natural$. By \eqref{eq:sandwich} with $g=f_0+\varphi$, $s=t^{(\omega)}_\star(x)$, $g'=f_0$ and $s'=\bar t$, every $y\in{\mathsf C}^{\star,(\omega)}(x)\triangle{\mathsf C}_0$ satisfies
	\begin{align*}
		|f_0(y)-\bar t|\le\norm{\varphi}_\infty+|t^{(\omega)}_\star(x)-\bar t|\le\bigl(\norm{\mathfrak u}_\infty+2C_{\mathfrak u}/c_\natural\bigr)\mu.
	\end{align*}
By the choice of $m_0'$ and \cref{ass:constants}(i), the right-hand side is at most $5\varepsilon_0/8$. Hence, by the first bound of \cref{lem:lb-perturb}(i) with $\varphi=0$, we obtain \eqref{eq:lb-symdiff0}.
\end{proof}

\begin{lemma}[Two-point bound for part (ii)]\label{lem:lb-ii-twopoint}
Suppose that \cref{ass:constants} holds. Let $m\ge m_0'$ be an integer, and let ${\mathsf P}_\omega$, $\omega\in\{0,1\}^{J_X}$, be the distributions of \eqref{eq:lb-pomega-ii}. For $x\in\cX_0$, let
	\begin{align*}
		\Omega_{{\rm vol},x}:=\cbr[1]{\Leb\bigl(\hC(x)\bigr)-\Leb\bigl({\mathsf C}_0\bigr)>2c_1}.
	\end{align*}
Then the following hold for every value of $\hC$.
\begin{enumerate}[label=(\roman*)]
\item Let $j\le J_X$ and $x\in Q_j'$. If $\Omega_{{\rm vol},x}$ fails, then
	\begin{align}\label{eq:lb-twopoint}
		\bigl|{\mathsf P}_{1,x}(\hC(x))-\gamma_\star\bigr|+\bigl|{\mathsf P}_{0,x}(\hC(x))-\gamma_\star\bigr|\ge\frac{\mu c_{\mathfrak u}}2 .
	\end{align}
\item For every $x\in\cX_0$ and every $\omega$, if $\Omega_{{\rm vol},x}$ holds, then
	\begin{align}\label{eq:lb-penalty}
		\Leb\bigl(\hC(x)\bigr)-\Leb\bigl({\mathsf C}^{\star,(\omega)}(x)\bigr)>c_1 .
	\end{align}
\end{enumerate}
\end{lemma}

\begin{proof}
Part (i) uses the volume decomposition under ${\mathsf P}_{0,x}$ to show that $\hC(x)$ is close to ${\mathsf C}_0$ when its coverage error is small, and part (ii) follows from \cref{lem:lb-ii-family}(ii).

\emph{Step 1: the two-point bound.} Suppose that $\Omega_{{\rm vol},x}$ fails. If $|{\mathsf P}_{0,x}(\hC(x))-\gamma_\star|>2\bar tc_1$, then \eqref{eq:lb-twopoint} holds since $c_{\mathfrak u}\mu\le4\bar tc_1$. Suppose now that $|{\mathsf P}_{0,x}(\hC(x))-\gamma_\star|\le2\bar tc_1$, and let $\cE_0$ be the excess mass under ${\mathsf P}_{0,x}$, whose density $f_0$ has level $\bar t$ and oracle set ${\mathsf C}_0$. Since ${\mathsf P}_\omega\in\cP_\beta$ for $\omega=0$ by \cref{lem:lb-ii-family}(i), we have, by \cref{lem:identity} and the definition of $\Omega_{{\rm vol},x}$,
	\begin{align*}
		\cE_0(\hC(x))=\Leb\bigl(\hC(x)\bigr)-\Leb\bigl({\mathsf C}_0\bigr)-\frac{{\mathsf P}_{0,x}(\hC(x))-\gamma_\star}{\bar t}\le2c_1+2c_1=c_E\delta_{\rm sd}^2 .
	\end{align*}
The proof of \cref{lem:lb-perturb}(iii) shows that $f_0$ satisfies \cref{ass:margin} with $\mathfrak L\equiv1$. Hence, by \eqref{eq:E-symdiff} with $\mathfrak L\equiv1$, we have $\Leb(\hC(x)\triangle{\mathsf C}_0)\le\delta_{\rm sd}$. Therefore
	\begin{align*}
		\int_{\hC(x)}\mathfrak u\,\d y=\int_{{\mathsf C}_0}\mathfrak u\,\d y+\int\bigl(\ind_{\hC(x)}-\ind_{{\mathsf C}_0}\bigr)\mathfrak u\,\d y\ge c_{\mathfrak u}-\norm{\mathfrak u}_\infty\delta_{\rm sd}=\frac{3c_{\mathfrak u}}4 ,
	\end{align*}
and \eqref{eq:lb-twopoint} follows from
	\begin{align*}
		\bigl|{\mathsf P}_{1,x}(\hC(x))-\gamma_\star\bigr|+\bigl|{\mathsf P}_{0,x}(\hC(x))-\gamma_\star\bigr|&\ge\bigl|{\mathsf P}_{1,x}(\hC(x))-{\mathsf P}_{0,x}(\hC(x))\bigr|\\
		&=\mu\Bigl|\int_{\hC(x)}\mathfrak u\,\d y\Bigr|\ge\frac{3\mu c_{\mathfrak u}}4 .
	\end{align*}

\emph{Step 2: the penalty.} We fix $\omega$ and $x\in\cX_0$. By \eqref{eq:lb-symdiff0} and the choice of $m_0'$, we have
	\begin{align*}
		\bigl|\Leb\bigl({\mathsf C}^{\star,(\omega)}(x)\bigr)-\Leb\bigl({\mathsf C}_0\bigr)\bigr|\le C_{\rm sd}\mu\le c_1 .
	\end{align*}
Hence, on $\Omega_{{\rm vol},x}$, we have
	\begin{align*}
		\Leb\bigl(\hC(x)\bigr)-\Leb\bigl({\mathsf C}^{\star,(\omega)}(x)\bigr)\ge\Leb\bigl(\hC(x)\bigr)-\Leb\bigl({\mathsf C}_0\bigr)-c_1>c_1,
	\end{align*}
which proves \eqref{eq:lb-penalty}.
\end{proof}

\begin{proof}[Proof of \cref{thm:lower}(ii)]
We apply Assouad's lemma to the integrated coverage error plus a penalty on the events $\Omega_{{\rm vol},x}$ of \cref{lem:lb-ii-twopoint}, and then bound the penalty by the definition of $\mathfrak V_n(c_1,c_2)$. Let
	\begin{align*}
		m:=\max\bigl\{\lceil n^{1/(2\beta+d_X)}\rceil,\ m_0'\bigr\},
	\end{align*}
let ${\mathsf P}_\omega$, $\omega\in\{0,1\}^{J_X}$, be the distributions of \eqref{eq:lb-pomega-ii}, and let $\hC\in\mathfrak V_n(c_1,c_2)$. By \cref{lem:lb-ii-family}(i), each ${\mathsf P}_\omega$ lies in $\cP_\beta$.

\emph{Step 1: a penalized loss.} We consider the loss
	\begin{align*}
		\mathrm{Err}(\hC,\omega):=\int_{\cX_0}|\cov_\omega(x)-\gamma_\star|\,\d x+\frac{\mu c_{\mathfrak u}}2\int_{\cX_0}\ind_{\Omega_{{\rm vol},x}}\,\d x .
	\end{align*}
The penalty term does not depend on $\omega$. With
	\begin{align*}
		\mathrm{Err}_j(\hC,b):=\int_{Q_j'}\Bigl(\bigl|{\mathsf P}_{b,x}(\hC(x))-\gamma_\star\bigr|+\frac{\mu c_{\mathfrak u}}2\ind_{\Omega_{{\rm vol},x}}\Bigr)\d x ,
	\end{align*}
we have $\mathrm{Err}(\hC,\omega)\ge\sum_{j\le J_X}\mathrm{Err}_j(\hC,\omega_j)$, since $\cov_\omega(x)={\mathsf P}_{\omega_j,x}(\hC(x))$ for $x\in Q_j'$ and the sets $Q_j'\subset\cX_0$ are disjoint. At each $x\in Q_j'$, either $\Omega_{{\rm vol},x}$ holds and the penalty contributes $\mu c_{\mathfrak u}/2$ to each of $\mathrm{Err}_j(\hC,0)$ and $\mathrm{Err}_j(\hC,1)$, or it fails and \eqref{eq:lb-twopoint} applies by \cref{lem:lb-ii-twopoint}(i). Therefore, for every value of $\hC$,
	\begin{align*}
		\mathrm{Err}_j(\hC,0)+\mathrm{Err}_j(\hC,1)\ge\frac{\mu c_{\mathfrak u}}2(2m)^{-d_X} .
	\end{align*}

\emph{Step 2: Assouad's lemma.} Let $\omega,\omega'$ differ only at $j$. Then ${\mathsf p}_\omega-{\mathsf p}_{\omega'}=\pm\mu\chi_X(m(x-x_j))\mathfrak u(y)$. By \cref{lem:lb-ii-family}(i), ${\mathsf p}_\omega(\cdot\mid x)$ and ${\mathsf p}_{\omega'}(\cdot\mid x)$ are $f_0$ plus small perturbations for every $x$. Hence, by \cref{lem:lb-kl}, we have
	\begin{align*}
		{\rm KL}\bigl({\mathsf P}_\omega^{\otimes n}\,\big\|\,{\mathsf P}_{\omega'}^{\otimes n}\bigr)\le\frac{n\mu^2}{{\mathsf p}_\flat}m^{-d_X}\norm{\chi_X}_2^2\norm{\mathfrak u}_2^2=C_{\rm KL}'c_\mu'^2L^2\,n\,m^{-(2\beta+d_X)}\le\frac18,
	\end{align*}
where the last inequality follows from $m\ge n^{1/(2\beta+d_X)}$ and $8C_{\rm KL}'c_\mu'^2L^2\le1$. By \cref{lem:lb-kl}, it follows that $\TV({\mathsf P}_\omega^{\otimes n},{\mathsf P}_{\omega'}^{\otimes n})\le1/4$. Hence, by \cref{lem:lb-assouad} with $J:=J_X\ge c_Jm^{d_X}$ and Step 1, we obtain
	\begin{align*}
		2^{-J}\sum_{\omega\in\{0,1\}^J}\E_\omega\sbr{\mathrm{Err}(\hC,\omega)}\ge\frac38J\cdot\frac{\mu c_{\mathfrak u}}2(2m)^{-d_X}\ge\frac{3c_Jc_{\mathfrak u}}{2^{d_X+4}}\mu .
	\end{align*}

\emph{Step 3: the penalty and the rate.} Since ${\mathsf P}_\omega\in\cP_\beta$ and $\hC\in\mathfrak V_n(c_1,c_2)$, it follows from \eqref{eq:lb-penalty}, the definition of $\mathfrak V_n(c_1,c_2)$ and $c_2=2^{-d_X-9}c_J$ that
	\begin{align*}
		\E_\omega\sbr{\frac{\mu c_{\mathfrak u}}2\int_{\cX_0}\ind_{\Omega_{{\rm vol},x}}\,\d x}\le\frac{\mu c_{\mathfrak u}}2c_2=2^{-d_X-10}c_Jc_{\mathfrak u}\mu .
	\end{align*}
Subtracting this bound from the bound of Step 2, we obtain
	\begin{align*}
		2^{-J}\sum_{\omega\in\{0,1\}^J}\E_\omega\sbr{\int_{\cX_0}|\cov_\omega(x)-\gamma_\star|\,\d x}\ge\frac{3c_Jc_{\mathfrak u}}{2^{d_X+4}}\mu-\frac{c_Jc_{\mathfrak u}}{2^{d_X+10}}\mu\ge\frac{c_Jc_{\mathfrak u}}{2^{d_X+3}}\mu .
	\end{align*}
The supremum over $\cP_\beta$ is at least this average. Finally, as in Step 3 of the proof of (i), we have 
  \begin{align*}
    \mu\ge c_\mu' L(2m_0')^{-\beta}n^{-\beta/(2\beta+d_X)},
  \end{align*}
which proves (ii).
\end{proof}

\subsection{Proof of \texorpdfstring{\cref{thm:lower}(iii)}{Theorem \ref{thm:lower}(iii)}}\label{app:lower-sup}
The proof of part (iii) uses the family \eqref{eq:lb-pomega-ii} with cells of side of order $(\log n/n)^{1/(2\beta+d_X)}$, a uniform prior on it, and the constants of \cref{app:lower-cov}. Given the data, the coordinates of $\omega$ are independent (\cref{lem:lb-posterior}). A predictive set with small uniform coverage error determines the coordinate of every cell that the class $\mathfrak V_n(c_1,c_2)$ does not exclude (\cref{lem:lb-cells}). At this cell size, the Kullback--Leibler divergence between two hypotheses that differ in one cell is a small multiple of the logarithm of the number $J$ of cells. Hence, on a positive fraction of the cells, the posterior probability of each value of the coordinate is at least of order $1/J$, while $\mathfrak V_n(c_1,c_2)$ excludes only a small fraction of the cells (\cref{lem:lb-uncertain}). With probability bounded away from zero, the predictive set therefore determines at least one coordinate incorrectly, and its coverage error on that cell is of order $\mu$.

\emph{The prior.} Let $n$ be large enough that $(n/\log n)^{1/(2\beta+d_X)}\ge m_0'$, with $c_\mu'$ and $m_0'$ of \cref{app:lower-cov}, and let
	\begin{align*}
		m:=\bigl\lceil(n/\log n)^{1/(2\beta+d_X)}\bigr\rceil,\qquad\mu:=c_\mu' Lm^{-\beta},\qquad J:=J_X\ge c_Jm^{d_X}.
	\end{align*}
Let $\{{\mathsf P}_\omega:\omega\in\{0,1\}^J\}$ be the distributions of \eqref{eq:lb-pomega-ii} with these $m$ and $\mu$. Since $m\ge m_0'$, \cref{lem:lb-ii-family,lem:lb-ii-twopoint} apply to this family. Hence each ${\mathsf P}_\omega$ lies in $\cP_\beta$, and \eqref{eq:lb-twopoint} and \eqref{eq:lb-penalty} hold. Let $\hC\in\mathfrak V_n(c_1,c_2)$.

Let $\bar\bP$ and $\bar\E$ denote the probability and the expectation with respect to the joint distribution under which $\omega$ is uniform on $\{0,1\}^J$, the sample $(X_i,Y_i)_{i\le n}$ is i.i.d.\ from ${\mathsf P}_\omega$ given $\omega$, and the auxiliary randomness $\mathfrak X$ is independent of $\omega$ and the sample. Let $\mathcal D:=((X_i,Y_i)_{i\le n},\mathfrak X)$, which determines $\hC$.

\emph{Likelihood ratios and cell averages.} Let $\chi_j(x):=\chi_X(m(x-x_j))$, $z_j(x,y):=\mu\chi_j(x)\mathfrak u(y)/f_0(y)$ and
	\begin{align*}
		\Lambda^{\rm lr}_j:=\prod_{i\le n}\bigl\{1+z_j(X_i,Y_i)\bigr\}\qquad(j\le J).
	\end{align*}
By (L1) and the choice of $m_0'$, we have $f_0\ge2{\mathsf p}_\flat$ and $\mu\norm[0]{\mathfrak u}_\infty\le\varepsilon_0/2\le{\mathsf p}_\flat$. Hence $|z_j|\le\frac12$ and $\Lambda^{\rm lr}_j>0$. For $b\in\{0,1\}$ and $j\le J$, let
	\begin{align*}
		A_{b,j}:=(2m)^{d_X}\int_{Q_j'}\bigl|{\mathsf P}_{b,x}(\hC(x))-\gamma_\star\bigr|\,\d x,\qquad F_j:=(2m)^{d_X}\int_{Q_j'}\ind_{\Omega_{{\rm vol},x}}\,\d x,
	\end{align*}
with ${\mathsf P}_{b,x}$ of \cref{app:lower-cov} and $\Omega_{{\rm vol},x}$ of \cref{lem:lb-ii-twopoint}. By the joint measurability of $\hC$ and Tonelli's theorem, $A_{b,j}$ and $F_j$ are measurable functions of $\mathcal D$. Finally, let
	\begin{align*}
		\mathcal G:=\{j:F_j\le\tfrac12\},\qquad s_n:=\mu c_{\mathfrak u}/8,\qquad\mathcal R:=\{j:|\log\Lambda^{\rm lr}_j|\le\log(J/32)\}.
	\end{align*}

\begin{lemma}[The posterior]\label{lem:lb-posterior}
Suppose that \cref{ass:constants} holds, and let $n$, $m$, $\mu$, $J$ and $\bar\bP$ be as above. Then, given $\mathcal D$, the coordinates $\omega_j$, $j\le J$, are independent under $\bar\bP$, with $\bar\bP(\omega_j=1\mid\mathcal D)=\varpi_j:=\Lambda^{\rm lr}_j/(1+\Lambda^{\rm lr}_j)$. Moreover, $\underline\varpi_j:=\min\{\varpi_j,1-\varpi_j\}$ satisfies
	\begin{align}\label{eq:lb-post}
		\underline\varpi_j=\frac{\min\{\Lambda^{\rm lr}_j,1\}}{1+\Lambda^{\rm lr}_j}\ge\frac12\e^{-|\log\Lambda^{\rm lr}_j|}.
	\end{align}
\end{lemma}

\begin{proof}
Since the supports of the $\chi_j$ lie in the interiors of disjoint cells, for each $x$ at most one $z_j(x,\cdot)$ is nonzero. Thus ${\mathsf p}_\omega(y\mid x)=f_0(y)\prod_j\{1+z_j(x,y)\}^{\omega_j}$. Since $\Pi$ has density one on $[0,1]^{d_X}$, the joint density of $\omega$ and the sample with respect to the product of the counting measure and Lebesgue measure is
	\begin{align*}
		2^{-J}\prod_{i\le n}f_0(Y_i)\prod_{j\le J}(\Lambda^{\rm lr}_j)^{\omega_j}.
	\end{align*}
As a function of $\omega$, this density is a product of functions of the single coordinates $\omega_j$. Since $\mathfrak X$ is independent of $\omega$ and the sample, it follows that, given $\mathcal D$, the coordinates $\omega_j$ are independent, with $\bar\bP(\omega_j=1\mid\mathcal D)=\varpi_j$. The equality in \eqref{eq:lb-post} follows from the definition of $\varpi_j$, and the inequality follows from $1+\Lambda^{\rm lr}_j\le2\max\{\Lambda^{\rm lr}_j,1\}$ and $\min\{\Lambda,1\}/\max\{\Lambda,1\}=\e^{-|\log\Lambda|}$ for $\Lambda>0$.
\end{proof}

\begin{lemma}[Reduction to the cells]\label{lem:lb-cells}
Suppose that \cref{ass:constants} holds, and let $n$, $m$, $\mu$, $J$, $\bar\bP$ and $\hC$ be as above. Then, for every $\omega$ and every value of $\mathcal D$,
	\begin{align}\label{eq:lb-sup-cells}
		\sup_{x\in\cX_0}|\cov_\omega(x)-\gamma_\star|\ge\max_{j\le J}A_{\omega_j,j},
	\end{align}
with $\cov_\omega$ as in \eqref{eq:lb-cov-omega}. Moreover, with $\omega$ random under $\bar\bP$,
	\begin{align}\label{eq:lb-prod}
		\bar\bP\Bigl(\max_{j\in\mathcal G}A_{\omega_j,j}<s_n\Bigm|\mathcal D\Bigr)\le\exp\Bigl(-\sum_{j\in\mathcal G}\underline\varpi_j\Bigr).
	\end{align}
\end{lemma}

\begin{proof}
Since ${\mathsf P}_{\omega,x}={\mathsf P}_{\omega_j,x}$ for $x\in Q_j'$ and $\Leb(Q_j')=(2m)^{-d_X}$, the quantity $A_{\omega_j,j}$ is the average of $|\cov_\omega(x)-\gamma_\star|$ over $Q_j'$. Since $Q_j'\subset\cX_0$, \eqref{eq:lb-sup-cells} follows.

Next, let $j\in\mathcal G$. Averaging \eqref{eq:lb-twopoint} over the points $x\in Q_j'$ at which $\Omega_{{\rm vol},x}$ fails, we obtain $A_{0,j}+A_{1,j}\ge(1-F_j)\mu c_{\mathfrak u}/2\ge2s_n$. Consequently, at most one $b\in\{0,1\}$ has $A_{b,j}<s_n$. Let $\hat b_j$ be this $b$, or $\hat b_j:=0$ if there is none. Then $\{A_{\omega_j,j}<s_n\}\subset\{\omega_j=\hat b_j\}$ for $j\in\mathcal G$. Since $\mathcal G$ and the $\hat b_j$ are functions of $\mathcal D$, \cref{lem:lb-posterior} and \eqref{eq:lb-post} imply that
	\begin{align*}
		\bar\bP\Bigl(\max_{j\in\mathcal G}A_{\omega_j,j}<s_n\Bigm|\mathcal D\Bigr)\le\prod_{j\in\mathcal G}\bar\bP(\omega_j=\hat b_j\mid\mathcal D)\le\prod_{j\in\mathcal G}(1-\underline\varpi_j)\le\exp\Bigl(-\sum_{j\in\mathcal G}\underline\varpi_j\Bigr),
	\end{align*}
which proves \eqref{eq:lb-prod}.
\end{proof}

\begin{lemma}[Uncertain cells]\label{lem:lb-uncertain}
Suppose that \cref{ass:constants} holds. Then, for all sufficiently large $n$ and every $\hC\in\mathfrak V_n(c_1,c_2)$, with $m$, $\mu$, $J$ and $\bar\bP$ as above, the event $\{\#\mathcal R\ge J/8,\ \#(\mathcal G^c)\le J/16\}$ has $\bar\bP$-probability at least $1/16$, and $\sum_{j\in\mathcal G}\underline\varpi_j\ge1$ on this event.
\end{lemma}

\begin{proof}
We show that each cell lies in $\mathcal R$ with $\bar\bP$-probability at least $1/4$, bound the expected number of cells outside $\mathcal G$ by the definition of $\mathfrak V_n(c_1,c_2)$, and combine the two bounds.

\emph{Step 1: moments of the log-likelihood ratios.} We fix $j$ and $\omega$ with $\omega_j=0$. Where $z_j(x,\cdot)\ne0$, we have $\chi_{j'}(x)=0$ for $j'\ne j$, and hence ${\mathsf p}_\omega(\cdot\mid x)=f_0$ there. Consequently, under ${\mathsf P}_\omega$, the variables $z_j(X_i,Y_i)$ are i.i.d.\ with the distribution of $z_j(X,Y)$ for independent $X\sim\Pi$ and $Y\sim f_0$. We write $\E_0$ for the expectation under this distribution. Since $z-z^2\le\log(1+z)\le z$ and $|\log(1+z)|\le2|z|$ for $|z|\le\frac12$, and $\E_0z_j=\mu\int\chi_j\,\d x\int\mathfrak u\,\d y=0$, we have
	\begin{align*}
		&-\E_0z_j^2\le\E_0\log(1+z_j)\le0,\qquad\E_0\log^2(1+z_j)\le4\E_0z_j^2,\\
		&\E_0z_j^2\le\frac{\mu^2\norm{\chi_X}_2^2\norm{\mathfrak u}_2^2}{2{\mathsf p}_\flat m^{d_X}}=\frac{C_{\rm KL}'\mu^2}{2m^{d_X}} .
	\end{align*}
Let $S:=n\E_0z_j^2$. By the choice of $m$, we have $nm^{-(2\beta+d_X)}\le\log n$. Thus $S\le\frac12C_{\rm KL}'c_\mu'^2L^2\log n$. Since $\log\Lambda^{\rm lr}_j$ is a sum of $n$ i.i.d.\ copies of $\log(1+z_j)$, we have $-S\le\E_\omega\log\Lambda^{\rm lr}_j\le0$ and ${\rm Var}_\omega(\log\Lambda^{\rm lr}_j)\le4S$. By Chebyshev's inequality, we therefore have
	\begin{align*}
		\bP_\omega\bigl\{|\log\Lambda^{\rm lr}_j|\le S+(8S)^{1/2}\bigr\}\ge\frac12 .
	\end{align*}

\emph{Step 2: cells with uncertain coordinates.} Since $J\ge c_Jm^{d_X}$ and $m\ge(n/\log n)^{1/(2\beta+d_X)}$, we have, for all sufficiently large $n$,
	\begin{align*}
		\log J\ge\log c_J+\frac{d_X}{2\beta+d_X}(\log n-\log\log n)\ge\frac{d_X}{2(2\beta+d_X)}\log n.
	\end{align*}
Together with the last condition on $c_\mu'$ in \cref{app:lower-cov} and Step 1, this bound implies $S\le\frac12\log J$. For all sufficiently large $n$, it follows that
	\begin{align*}
		S+(8S)^{1/2}\le\frac12\log J+2(\log J)^{1/2}\le\log(J/32).
	\end{align*}
By Step 1, we therefore have $\bP_\omega(j\in\mathcal R)\ge\frac12$ for every $\omega$ with $\omega_j=0$. Averaging over $\omega$ yields $\bar\bP(j\in\mathcal R)\ge\frac12\bar\bP(\omega_j=0)=\frac14$. Summing over $j$, we obtain $\bar\E\,\#\mathcal R\ge J/4$. Since $\#\mathcal R\le J$, it follows that $\bar\bP(\#\mathcal R\ge J/8)\ge1/8$.

\emph{Step 3: cells excluded by the volume condition.} Since $F_j>\frac12$ for $j\notin\mathcal G$, we have
	\begin{align*}
		\#(\mathcal G^c)\le2\sum_{j\le J}F_j\le2(2m)^{d_X}\int_{\cX_0}\ind_{\Omega_{{\rm vol},x}}\,\d x.
	\end{align*}
By \eqref{eq:lb-penalty} and $\hC\in\mathfrak V_n(c_1,c_2)$, we have $\E_\omega\sbr[0]{\int_{\cX_0}\ind_{\Omega_{{\rm vol},x}}\,\d x}\le c_2$ for every $\omega$. Hence
	\begin{align*}
		\bar\E\,\#(\mathcal G^c)\le2^{d_X+1}m^{d_X}c_2\le\frac{2^{d_X+1}c_2}{c_J}J=\frac J{256},
	\end{align*}
and, by Markov's inequality, $\bar\bP(\#(\mathcal G^c)>J/16)\le1/16$. Together with Step 2, this bound gives
	\begin{align*}
		\bar\bP\del[1]{\#\mathcal R\ge J/8,\ \#(\mathcal G^c)\le J/16}\ge\bar\bP\del[1]{\#\mathcal R\ge J/8}-\bar\bP\del[1]{\#(\mathcal G^c)>J/16}\ge\frac18-\frac1{16}=\frac1{16}.
	\end{align*}
On this event, $\#(\mathcal G\cap\mathcal R)\ge J/16$. Since $\underline\varpi_j\ge16/J$ for $j\in\mathcal R$ by \eqref{eq:lb-post}, it follows that, on this event,
	\begin{align*}
		\sum_{j\in\mathcal G}\underline\varpi_j\ge\frac{16}J\,\#(\mathcal G\cap\mathcal R)\ge1 .
	\end{align*}
\end{proof}

\begin{proof}[Proof of \cref{thm:lower}(iii)]
We combine \cref{lem:lb-cells,lem:lb-uncertain} under the prior and then pass from the prior to a single distribution of the family.

\emph{Step 1: a bound under the prior.} Taking expectations in \eqref{eq:lb-prod}, and using that $\bar\bP(\sum_{j\in\mathcal G}\underline\varpi_j\ge1)\ge1/16$ by \cref{lem:lb-uncertain}, we obtain
	\begin{align*}
		\bar\bP\Bigl(\max_{j\le J}A_{\omega_j,j}\ge s_n\Bigr)\ge\bar\E\sbr{1-\exp\Bigl(-\sum_{j\in\mathcal G}\underline\varpi_j\Bigr)}\ge\frac{1-\e^{-1}}{16}.
	\end{align*}

\emph{Step 2: a single distribution and the rate.} The left-hand side of the last display is the average over $\omega$ of $\bP_\omega\{\max_jA_{\omega_j,j}\ge s_n\}$. Therefore $\bP_\omega\{\max_jA_{\omega_j,j}\ge s_n\}\ge(1-\e^{-1})/16$ for some $\omega$. By \eqref{eq:lb-sup-cells}, this ${\mathsf P}_\omega\in\cP_\beta$ satisfies
	\begin{align*}
		\bP_\omega\cbr[2]{\sup_{x\in\cX_0}|\cov_\omega(x)-\gamma_\star|\ge s_n}\ge\frac{1-\e^{-1}}{16},
	\end{align*}
as an outer probability if the event is not measurable. Finally, since $s_n=c_{\mathfrak u}c_\mu' Lm^{-\beta}/8$ and $m\le2(n/\log n)^{1/(2\beta+d_X)}$, we have $s_n\ge c(\log n/n)^{\beta/(2\beta+d_X)}$ for a constant $c\in(0,(1-\e^{-1})/16]$, and (iii) holds with this $c$. For $Z:=\sup_{x\in\cX_0}|\cov(x)-\gamma_\star|$, Markov's inequality gives
	\begin{align*}
		\E^*_{\mathsf P}[Z]\ge c\Bigl(\frac{\log n}n\Bigr)^{\beta/(2\beta+d_X)}\bP_{\mathsf P}\cbr[2]{Z\ge c\Bigl(\frac{\log n}n\Bigr)^{\beta/(2\beta+d_X)}},
	\end{align*}
which yields the bound on the outer expectation.
\end{proof}

\section{Proofs for \texorpdfstring{\cref{sec:cov-twostage}}{Section \ref{sec:cov-twostage}}}\label{app:coverage}
\subsection{Calibration}\label{app:calibration}
The procedures of \cref{sec:cov-twostage,sec:cov-onestage} calibrate a score, namely $\hat{\mathsf s}^{\rm two}$ or the HPD value. This score is a function of a point $(x,y)$ and of the tie-breaking variable at $x$, if any. This function is determined by the data used before calibration, which are $\mathcal D_{\rm fit}$ and, for two-stage calibration, $\mathcal D_{\rm reg}$.

\begin{lemma}[Calibration events]\label{lem:cal-events}
Suppose that \cref{ass:design,ass:split} hold, and let $\pi_0:=f_{\min}\Leb(\cX_0)>0$ and $\Omega^\circ_{\rm cal}:=\{n_{\rm cal}^\circ\ge\pi_0n_{\rm cal}/2\}$. Then the following hold.
\begin{enumerate}[label=(\roman*)]
\item $\bP(\Omega^\circ_{\rm cal})\ge1-\e^{-\pi_0^2n_{\rm cal}/2}$.
\item Conditionally on the data used before calibration and on $n_{\rm cal}^\circ$, if $n_{\rm cal}^\circ\ge1$, the scores of the $n_{\rm cal}^\circ$ calibration points in $\cX_0$ are i.i.d.\ with some distribution function $H$.
\item For every $c>0$ there is $C>0$ such that, under the conditioning of (ii), with conditional probability at least $1-Cn^{-c}$, the empirical distribution function $\widehat H$ of these scores satisfies
	\begin{align}\label{eq:ts-dkw}
		\sup_u|\widehat H(u)-H(u)|\le\eps_{{\rm cal},n}:=C\Bigl(\frac{\log n}{n_{\rm cal}^\circ}\Bigr)^{1/2}.
	\end{align}
\end{enumerate}
\end{lemma}

\begin{proof}
Since $f_X\ge f_{\min}$ on $\cX_0$, we have $\Pi(\cX_0)\ge\pi_0$. Hence $n_{\rm cal}^\circ$ is a binomial variable with $n_{\rm cal}$ trials and success probability at least $\pi_0$, and (i) follows from Hoeffding's inequality. Given the data used before calibration and $n_{\rm cal}^\circ$, the calibration points in $\cX_0$ are i.i.d.\ from the distribution of $(X,Y)$ given $X\in\cX_0$, and their tie-breaking variables are drawn afresh at each point as in \cref{sec:setup-oracle}. Consequently, their scores are i.i.d. By the Dvoretzky--Kiefer--Wolfowitz inequality with Massart's constant, the probability that \eqref{eq:ts-dkw} fails is at most $2\exp(-2C^2\log n)\le Cn^{-c}$ for $C\ge2$ and $2C^2\ge c$.
\end{proof}

\subsection{Regression weights}\label{app:second-stage}\label{app:tools}
The kernels of the local polynomial density estimate and of the second-stage regression are fixed product kernels $\mathcal K(s)=\prod_{j=1}^{d'}\mathcal K_1(s_j)$ on $\R^{d'}$, where $\mathcal K_1:\R\to[0,\infty)$ is continuously differentiable, vanishes outside $[-1,1]$, and is bounded below by a positive constant on $[-\frac12,\frac12]$. We call such kernels smooth. The lower bound on $[-\frac12,\frac12]$ keeps the local Gram matrices nondegenerate, even in windows truncated at the boundary of $\cY$.

Local polynomial regression of degree $\ell$ in $x$ on the covariates of $\mathcal D_{\rm reg}$ uses the following weights. Let $\mathcal K_X$ be a smooth product kernel on $\R^{d_X}$, let $D_X(x';x)$ be the vector of monomials of degree at most $\ell$ in $(x'-x)/h_{\rm reg}$, and let $W_X(x';x):=\mathcal K_X((x'-x)/h_{\rm reg})$. The fit at $x$ to responses attached to the points $i\in\mathcal D_{\rm reg}$ is their weighted sum with weights
	\begin{align}
		w_i(x)&:=\frac{\mathbf e_1^\top\widehat\Sigma_X(x)^{-1}D_X(X_i;x)W_X(X_i;x)}{n_{\rm reg}h_{\rm reg}^{d_X}},\label{eq:loc-weights}\\
		\widehat\Sigma_X(x)&:=\frac1{n_{\rm reg}h_{\rm reg}^{d_X}}\sum_{i\in \mathcal D_{\rm reg}}W_X(X_i;x)D_X(X_i;x)D_X(X_i;x)^\top,\notag
	\end{align}
where $\mathbf e_1:=(1,0,\dots,0)^\top$ selects the intercept, and $w_i(x):=0$ for all $i$ when $\widehat\Sigma_X(x)$ is singular.

\begin{lemma}[Local polynomial weights]\label{lem:loc-weights}
Suppose that \cref{ass:design,ass:split} hold and that $h_{\rm reg}\asymp(\log n/n)^{1/(2\beta+d_X)}$. Then the weights \eqref{eq:loc-weights} satisfy \cref{cond:W}.
\end{lemma}

\begin{proof}
Since $n_{\rm reg}\ge cn$ by \cref{ass:split}, it suffices to prove the bounds of \cref{cond:W} with $n_{\rm reg}$ in place of $n$.

\emph{(i)} For all sufficiently large $n_{\rm reg}$, we have $h_{\rm reg}\le\delta_0$ and $n_h:=n_{\rm reg}h_{\rm reg}^{d_X}\ge C\log n_{\rm reg}$. By \cref{ass:design}, the covariates of $\mathcal D_{\rm reg}$ have a density bounded above and below on $\cX_0^{2\delta_0}$. Let $\Sigma_X(x):=\E\sbr[0]{\widehat\Sigma_X(x)}$. By the step \emph{Gram matrix} of the proof of \cref{lem:vc}, with $x$ alone in place of $(x,y)$, there is a constant $c_{\Sigma,X}>0$ such that $\lambda_{\min}(\Sigma_X(x))\ge c_{\Sigma,X}$ for all $x\in\cX_0$. By the argument of the step \emph{Concentration} of that proof,
	\begin{align*}
		\bP(\Omega_X)\ge1-Cn_{\rm reg}^{-c},\qquad\Omega_X:=\cbr[3]{\sup_{x\in\cX_0}\norm[1]{\widehat\Sigma_X(x)-\Sigma_X(x)}\le\frac{c_{\Sigma,X}}2}.
	\end{align*}
Hence, by Weyl's inequality, $\norm[0]{\widehat\Sigma_X(x)^{-1}}\le2/c_{\Sigma,X}$ and $\norm[0]{\widehat\Sigma_X(x)}\le C$ for all $x\in\cX_0$ on $\Omega_X$. We work on $\Omega_X$ and write $\widehat\Sigma_X:=\widehat\Sigma_X(x)$, $W_i:=W_X(X_i;x)$ and $D_i:=D_X(X_i;x)$. The weights vanish unless $|X_i-x|_\infty\le h_{\rm reg}$, because $\mathcal K_X$ is supported in $[-1,1]^{d_X}$. We have 
    \begin{align*}
      w_i(x)=n_h^{-1}W_i\mathbf e_1^\top\widehat\Sigma_X^{-1}D_i, \qquad 
      \widehat\Sigma_X=n_h^{-1}\sum_{i\in\mathcal D_{\rm reg}}W_iD_iD_i^\top
    \end{align*}
and, since the first entry of $D_i$ is one, $n_h^{-1}\sum_{i\in\mathcal D_{\rm reg}}W_i=\mathbf e_1^\top\widehat\Sigma_X\mathbf e_1$. Hence, by the Cauchy--Schwarz inequality with the weights $0\le W_i\le\norm[0]{\mathcal K_X}_\infty$,
	\begin{align*}
		\sum_{i\in\mathcal D_{\rm reg}}|w_i(x)|&\le\bigl(\mathbf e_1^\top\widehat\Sigma_X\mathbf e_1\cdot\mathbf e_1^\top\widehat\Sigma_X^{-1}\mathbf e_1\bigr)^{1/2}\le C,\\
		\sum_{i\in\mathcal D_{\rm reg}}w_i(x)^2&\le\frac{\norm{\mathcal K_X}_\infty}{n_h}\,\mathbf e_1^\top\widehat\Sigma_X^{-1}\mathbf e_1\le\frac C{n_h} .
	\end{align*}

\emph{(ii)} Let $\Xi_i(x):=W_X(X_i;x)D_X(X_i;x)$. Since $D_X(X_i;x)$ is bounded on the window, we have
	\begin{align*}
		\sum_{i\in\mathcal D_{\rm reg}}|\Xi_i(x)|\le C\sum_{i\in\mathcal D_{\rm reg}}W_i=Cn_h\mathbf e_1^\top\widehat\Sigma_X\mathbf e_1\le Cn_h .
	\end{align*}
Since the kernel is smooth and compactly supported, the maps $x\mapsto\Xi_i(x)$ and $x\mapsto\Xi_i(x)D_X(X_i;x)^\top$ are Lipschitz with constant $C/h_{\rm reg}$. Summing over the $n_{\rm reg}$ indices, we obtain
	\begin{align*}
		\norm{\widehat\Sigma_X(x)-\widehat\Sigma_X(x')}&\le Ch_{\rm reg}^{-d_X-1}|x-x'|_\infty,\\
		\sum_{i\in\mathcal D_{\rm reg}}|\Xi_i(x)-\Xi_i(x')|&\le Cn_{\rm reg}h_{\rm reg}^{-1}|x-x'|_\infty.
	\end{align*}
By the identity
	\begin{align*}
		\widehat\Sigma_X(x)^{-1}-\widehat\Sigma_X(x')^{-1}=\widehat\Sigma_X(x)^{-1}\bigl\{\widehat\Sigma_X(x')-\widehat\Sigma_X(x)\bigr\}\widehat\Sigma_X(x')^{-1}
	\end{align*}
and the bounds above, it follows that
	\begin{align*}
		\sum_{i\in\mathcal D_{\rm reg}}|w_i(x)-w_i(x')|&\le C\norm{\widehat\Sigma_X(x)-\widehat\Sigma_X(x')}+\frac C{n_h}\sum_{i\in\mathcal D_{\rm reg}}|\Xi_i(x)-\Xi_i(x')|\\
		&\le Ch_{\rm reg}^{-d_X-1}|x-x'|_\infty .
	\end{align*}

\emph{(iii)} Let $g\in\mathcal H^\beta_{\cX_0^{\delta_0}}(1)$, and let $\mathcal T_x$ be its Taylor polynomial of degree $\ell$ at $x$. Since the local polynomial weights reproduce polynomials of degree at most $\ell$, we have $\sum_{i\in\mathcal D_{\rm reg}}w_i(x)\mathcal T_x(X_i)=\mathcal T_x(x)=g(x)$. Since $h_{\rm reg}\le\delta_0$, the window $\{x':|x'-x|_\infty\le h_{\rm reg}\}$ lies in $\cX_0^{\delta_0}$, and, by Taylor's theorem and \eqref{eq:holder}, $|g-\mathcal T_x|\le Ch_{\rm reg}^\beta$ on it. Hence, by (i), $|\sum_{i\in\mathcal D_{\rm reg}}w_i(x)g(X_i)-g(x)|\le Ch_{\rm reg}^\beta$.

\emph{The rates.} Since $h_{\rm reg}\le\delta_0$, the weights vanish unless $X_i\in\cX_0^{\delta_0}$. Since $h_{\rm reg}\asymp(\log n/n)^{1/(2\beta+d_X)}$, we have $h_{\rm reg}^\beta\le C\rho_n$, $1/(nh_{\rm reg}^{d_X})\le C\rho_n^2/\log n$ and, for all sufficiently large $n$, $Ch_{\rm reg}^{-d_X-1}\le n^{d_X+2}$. Hence (i)--(iii) give \cref{cond:W} with $C_W=d_X+2$.
\end{proof}

\subsection{Uniform regression over candidate predictive sets}\label{app:local}\label{app:condV}
Throughout this subsection, \cref{ass:smooth,ass:design,ass:split} hold, and the weights of \cref{alg:local} satisfy \cref{cond:W} on the event $\Omega_X$.

For $x\in\cX_0$ and a set $A$ that is fixed given $\mathcal D_{\rm fit}$, we write
	\begin{align*}
		\sum_{i\in\mathcal D_{\rm reg}}w_i(x)\ind\{Y_i\in A\}-{\mathsf P}_x(A)&=\mathfrak n_x(A)+\mathfrak b_x(A),
	\end{align*}
where
	\begin{align*}
    	  \mathfrak n_x(A)&:=\sum_{i\in\mathcal D_{\rm reg}}w_i(x)\bigl\{\ind\{Y_i\in A\}-{\mathsf P}_{X_i}(A)\bigr\},\\
    	\mathfrak b_x(A)&:=\sum_{i\in\mathcal D_{\rm reg}}w_i(x){\mathsf P}_{X_i}(A)-{\mathsf P}_x(A).
	\end{align*}
With $A=\mathcal U_x(t)$, the left-hand side is $\widehat Q_x(t)-Q_x(t)$. For a measurable $A\subset\cY$, let $g_A(x'):={\mathsf P}_{x'}(A)=\int_A{\mathsf p}(y\mid x')\,\d y$. Since $(x,y)\mapsto{\mathsf p}(y\mid x)$ is in $\mathcal H^\beta(L)$, the derivatives in $x'$ of order at most $\ell$ can be taken under the integral sign, and hence $g_A\in\mathcal H^\beta(L\Leb(\cY))$. Therefore, by \cref{cond:W}(i) and by (iii) applied to $g_A/\{L\Leb(\cY)\}$, we have, on $\Omega_X$, simultaneously for all $x\in\cX_0$ and all measurable $A\subset\cY$,
	\begin{align*}
		|\mathfrak b_x(A)|\le C\rho_n,\qquad\sum_{i\in\mathcal D_{\rm reg}}w_i(x)^2\le\frac{C\rho_n^2}{\log n} .
	\end{align*}
Given $\mathcal D_{\rm fit}$ and the covariates of $\mathcal D_{\rm reg}$, the $Y_i$ are independent with distributions ${\mathsf P}_{X_i}$, because $\mathcal D_{\rm reg}$ is independent of $\mathcal D_{\rm fit}$.

\begin{lemma}[A weighted empirical bound]\label{lem:weighted-sets}
Let $Y_1,\ldots,Y_N$ be independent, with possibly different distributions $P_1,\ldots,P_N$, let $w_1,\ldots,w_N\in\R$ be fixed weights, and let $\mathcal A$ be a Borel-parametrized family of measurable sets with growth function $\mathfrak S_{\mathcal A}$. Let $s_w:=\sum_{i\le N} w_i^2$, and let
	\begin{align*}
		Z_w:=\sup_{A\in\mathcal A}\abs{\sum_{i\le N} w_i\{\ind_A(Y_i)-P_i(A)\}}.
	\end{align*}
Then, for every $\xi>0$, with probability at least $1-\e^{-\xi}$,
	\begin{align}\label{eq:weighted-growth}
		Z_w\le 2\sqrt{2s_w\log\{2\mathfrak S_{\mathcal A}(N)\}}+\sqrt{s_w\xi/2}.
	\end{align}
The same assertion holds conditionally on a $\sigma$-field when, given it, the $Y_i$ are independent with distributions $P_i$ and the weights and the family are fixed.
\end{lemma}

\begin{proof}
If $s_w=0$, the assertion is immediate. We therefore assume that $s_w>0$. Let $Y_i'$ be independent copies of the respective $Y_i$, and let $\epsilon_i$ be independent Rademacher variables, independent of both samples. By Jensen's inequality, and since exchanging $Y_i$ and $Y_i'$ does not change the joint distribution, we have
	\begin{align*}
		\E Z_w&\le \E\sup_{A\in\mathcal A}\abs{\sum_{i\le N} w_i\{\ind_A(Y_i)-\ind_A(Y_i')\}}\\
		&=\E\sup_{A\in\mathcal A}\abs{\sum_{i\le N}\epsilon_i w_i\{\ind_A(Y_i)-\ind_A(Y_i')\}}\\
		&\le 2\E\sup_{A\in\mathcal A}\abs{\sum_{i\le N}\epsilon_i w_i\ind_A(Y_i)}.
	\end{align*}
Conditionally on $(Y_i)$, the vector $(w_i\ind_A(Y_i))_{i\le N}$ takes at most $\mathfrak S_{\mathcal A}(N)$ values as $A$ ranges over $\mathcal A$, and each value has Euclidean norm at most $s_w^{1/2}$. By Lemma~5.2 of \citet{massart2000some}, applied to these vectors and their negatives, at most $2\mathfrak S_{\mathcal A}(N)$ in number, we have
	\begin{align*}
		\E\sbr[3]{\sup_{A\in\mathcal A}\abs[2]{\sum_{i\le N}\epsilon_i w_i\ind_A(Y_i)}\,\biggm|\,(Y_i)_{i\le N}}\le\sqrt{2s_w\log\{2\mathfrak S_{\mathcal A}(N)\}}.
	\end{align*}
Therefore, by the first display,
	\begin{align*}
		\E Z_w\le2\sqrt{2s_w\log\{2\mathfrak S_{\mathcal A}(N)\}}.
	\end{align*}
Finally, $Z_w$ does not depend on the $Y_i$ with $w_i=0$, and changing any other $Y_i$ alone changes $Z_w$ by at most $|w_i|$. Consequently, by the bounded-differences inequality \citep[Theorem~D.8]{mohri2018foundations}, $\bP\{Z_w>\E Z_w+\sqrt{s_w\xi/2}\}\le\e^{-\xi}$, which proves \eqref{eq:weighted-growth}. All suprema are understood on completed probability spaces. For the Borel-parametrized family \eqref{eq:cand-family}, the suprema are universally measurable because projections of Borel sets are universally measurable. The same argument applies to the independent copies and signs used above.
\end{proof}

\begin{lemma}[Accuracy of the second-stage regression]\label{lem:loc-acc}
For every $c>0$ and $C_{\rm gr}>0$ there is $C>0$ such that, for all sufficiently large $n$ and every $\mathcal D_{\rm fit}$, the following hold.
\begin{enumerate}[label=(\roman*)]
\item For each fixed $x\in\cX_0$, with conditional probability at least $1-Cn^{-c}$ given $\mathcal D_{\rm fit}$,
	\begin{align*}
		\sup_{t\in\R}|\widehat Q_x(t)-Q_x(t)|\le C\rho_n.
	\end{align*}
\item If $\mathfrak S_{\mathcal A_n}(m)\le(nm)^{C_{\rm gr}}$ for every $m\ge1$, then, with conditional probability at least $1-Cn^{-c}$ given $\mathcal D_{\rm fit}$, the bound of (i) holds simultaneously for all $x\in\cX_0$.
\end{enumerate}
\end{lemma}

\begin{proof}
\emph{(i)} We condition on the fitting data and on regression covariates for which $\Omega_X$ occurs. The sets $\mathcal U_x(t)$, $t\in\R$, are nested, and hence their growth function is at most $m+1$ at $m$ points. We apply \cref{lem:weighted-sets} with $N=n_{\rm reg}$, the weights $w_i(x)$ and $\xi$ a sufficiently large constant times $\log n$. Hence, by \eqref{eq:weighted-growth} and the bound on $\sum_{i\in\mathcal D_{\rm reg}}w_i(x)^2$, we have, outside an event of conditional probability at most $n^{-c}$,
	\begin{align*}
		\sup_{t\in\R}|\mathfrak n_x(\mathcal U_x(t))|\le C\rho_n .
	\end{align*}
Together with the bias bound and $\bP(\Omega_X^c)\le Cn^{-c}$, this proves (i).

\emph{(ii)} We fix fitting data satisfying the stated growth bound and regression covariates for which $\Omega_X$ occurs. Let $\mathcal G_n\subset\cX_0$ be a deterministic net of mesh $n^{-(C_W+2)}$ and cardinality at most $Cn^{d_X(C_W+2)}$. At each $x_g\in\mathcal G_n$, we apply \cref{lem:weighted-sets} with $N=n_{\rm reg}$, the weights $w_i(x_g)$ and $\mathcal A=\mathcal A_n$, for which $\log\{2\mathfrak S_{\mathcal A_n}(n_{\rm reg})\}\le C\log n$. We take $\xi$ to be a sufficiently large constant times $\log n$. By a union bound over $\mathcal G_n$, with conditional probability at least $1-Cn^{-c}$, we have
	\begin{align}\label{eq:loc-grid-bound}
		\max_{x_g\in\mathcal G_n}\sup_{A\in\mathcal A_n}|\mathfrak n_{x_g}(A)|\le C\rho_n.
	\end{align}
For any $x\in\cX_0$, we choose $x_g\in\mathcal G_n$ with $|x-x_g|_\infty\le n^{-(C_W+2)}$. For every $A\in\mathcal A_n$, since $|\ind\{Y_i\in A\}-{\mathsf P}_{X_i}(A)|\le1$, \cref{cond:W}(ii) gives
	\begin{align*}
		|\mathfrak n_x(A)-\mathfrak n_{x_g}(A)|\le\sum_{i\in\mathcal D_{\rm reg}}|w_i(x)-w_i(x_g)|\le n^{C_W}n^{-(C_W+2)}=n^{-2}.
	\end{align*}
Since each $\mathcal U_x(t)$ belongs to $\mathcal A_n$, by \eqref{eq:loc-grid-bound}, the last display and the uniform bias bound, we obtain
	\begin{align*}
		\sup_{x\in\cX_0}\sup_{t\in\R}|\widehat Q_x(t)-Q_x(t)|\le C(\rho_n+n^{-2})\le C\rho_n.
	\end{align*}
Since $\bP(\Omega_X^c)\le Cn^{-c}$, this proves (ii).
\end{proof}

\subsection{Proof of \texorpdfstring{\cref{thm:local-cov,thm:local} and~\cref{rem:local-each}}{Theorems \ref{thm:local-cov} and \ref{thm:local} and Remark \ref{rem:local-each}}}\label{app:local-proof}
For $x\in\cX_0$ and $Y\sim{\mathsf P}_x$, let
	\begin{align*}
		T_x:=Q_x(\hat a(Y\mid x)),\qquad H_x(u):={\mathsf P}_x\{\hat{\mathsf s}^{\rm two}(x,Y)\le u\}\qquad(u\in\R).
	\end{align*}
The function $H_x$ depends on the fitting and regression data, and $\cov^{\rm two}(x)=H_x(\hat\tau)$.

\begin{remark}[Guarantees at each covariate]\label{rem:local-each}
Suppose that \cref{ass:smooth,ass:margin,ass:design,ass:split} hold and that the weights of \cref{alg:local} satisfy \cref{cond:W}. Let $r_n\to0$ be a deterministic sequence, and suppose that $\hat a$, which need not satisfy \cref{cond:V}, satisfies \cref{cond:N,cond:A} at each $x$: for every $c>0$, there exist $C>0$ and, for all sufficiently large $n$, events $\Omega_A(x)\in\sigma(\mathcal D_{\rm fit})$, $x\in\cX_0$, of probability at least $1-Cn^{-c}$ such that the set of pairs $(x,\mathcal D_{\rm fit})$ with $\mathcal D_{\rm fit}\in\Omega_A(x)$ is measurable and, on $\Omega_A(x)$, the level sets of $\hat a(\cdot\mid x)$ are Lebesgue-null and some strictly increasing $\phi_x$ satisfies the two bounds of \cref{cond:A} at $x$. Then for every $c>0$ there exists $C>0$ such that, for all sufficiently large $n$, with probability at least $1-Cn^{-c}$, the bounds of \cref{thm:local-cov}(i) and \eqref{eq:local-bounds-gen} hold with the suprema taken over $\cX_0\setminus\mathcal X_{\rm bad}$, for a random set $\mathcal X_{\rm bad}\subset\cX_0$ with $\Pi(\mathcal X_{\rm bad})\le n^{-c}$.
\end{remark}

\begin{proof}[Proof of \cref{thm:local-cov,thm:local} and~\cref{rem:local-each}]
We fix $c>0$ and use the events of \cref{cond:N,cond:V,cond:A}, the event of \cref{lem:loc-acc}(ii), the calibration events of \cref{lem:cal-events} and, for \cref{rem:local-each}, the events $\Omega_A(x)$, each with $2c+2$ in place of $c$. Let $C_{\rm acc}$ be the constant of \cref{lem:loc-acc}, with $C_{\rm gr}$ of \cref{cond:V}, or $C_{\rm gr}:=1$ for \cref{rem:local-each}. Let $e_n$ be as in \cref{sec:cov-theory}.

\emph{Step 1: the envelope.} Let $x\in\cX_0$ and $t\in\R$. For $v\ge t$, we have $Q_x(v)\le Q_x(t)$, and hence $\widehat Q_x(v)\le Q_x(t)+e_n$. Therefore
	\begin{align*}
		Q_x(t)-e_n\le\widehat Q_x(t)\le\sup_{v\ge t}\widehat Q_x(v)\le Q_x(t)+e_n .
	\end{align*}
Since $Q_x(t)\in[0,1]$, clipping to $[0,1]$ does not increase the distance to $Q_x(t)$. Consequently, for all $x\in\cX_0$ and $y\in\cY$,
	\begin{align}\label{eq:env}
		\abs[1]{\hat{\mathsf s}^{\rm two}(x,y)-Q_x(\hat a(y\mid x))}\le\sup_{t\in\R}|\overline Q_x(t)-Q_x(t)|\le e_n .
	\end{align}

\emph{Step 2: uniform scores.} On the event of \cref{cond:N}, since ${\mathsf P}_x$ has a Lebesgue density, we have ${\mathsf P}_x\{\hat a(Y\mid x)=t\}=0$ for all $x\in\cX_0$ and $t\in\R$. Hence $Q_x(t)={\mathsf P}_x\{\hat a(Y\mid x)\ge t\}$ is one minus the continuous distribution function of $\hat a(Y\mid x)$, and, by the probability integral transform,
	\begin{align}\label{eq:pit-unif}
		{\mathsf P}_x(T_x\le u)=u\qquad(x\in\cX_0,\ u\in[0,1]).
	\end{align}

\emph{Step 3: calibration.} By \eqref{eq:env},
	\begin{align*}
		\{T_x\le u-e_n\}\subset\{\hat{\mathsf s}^{\rm two}(x,Y)\le u\}\subset\{T_x\le u+e_n\},
	\end{align*}
and hence, by \eqref{eq:pit-unif},
	\begin{align}\label{eq:Hx}
		|H_x(u)-u|\le e_n\qquad(x\in\cX_0,\ u\in[0,1]).
	\end{align}
Let $m:=n_{\rm cal}^\circ$ and $k:=\lceil\gamma_\star(m+1)\rceil$. By \cref{lem:cal-events}(ii), given the fitting and regression data and $m$, the scores of the calibration points in $\cX_0$ are i.i.d.\ with distribution function $H(u):=\Pi(\cX_0)^{-1}\int_{\cX_0}H_x(u)\,\Pi(\d x)$, which also satisfies \eqref{eq:Hx}. Let $\widehat H$ be their empirical distribution function, and suppose that $\sup_u|\widehat H(u)-H(u)|\le\eps_{{\rm cal},n}$ and that $\delta^{\rm thr}:=e_n+\eps_{{\rm cal},n}+2/m\le\min\{\gamma_\star,1-\gamma_\star\}$. Since $\gamma_\star\le k/m\le\gamma_\star+2/m$, we have
	\begin{align*}
		\widehat H(\gamma_\star-\delta^{\rm thr})\le\gamma_\star-\frac2m<\frac km,\qquad\widehat H(\gamma_\star+\delta^{\rm thr})\ge\gamma_\star+\frac2m\ge\frac km .
	\end{align*}
Hence $k\le m$ and $|\hat\tau-\gamma_\star|\le \delta^{\rm thr}$. Since $\cov^{\rm two}(x)=H_x(\hat\tau)$, by \eqref{eq:Hx}, we obtain
	\begin{align}\label{eq:two-cov}
		\sup_{x\in\cX_0}|\cov^{\rm two}(x)-\gamma_\star|\le e_n+\delta^{\rm thr}\le2e_n+\eps_{{\rm cal},n}+\frac2m .
	\end{align}

\emph{Step 4: rates.} On the events of \cref{lem:cal-events}, we have $m\ge\pi_0n_{\rm cal}/2\gtrsim n$ by \cref{ass:split} and $\eps_{{\rm cal},n}\le C\sqrt{\log n/n}\le C\rho_n$. On the event of \cref{cond:V}, it follows from \cref{lem:loc-acc}(ii) that $e_n\le C_{\rm acc}\rho_n$ with conditional probability at least $1-Cn^{-2c-2}$. Hence, for all sufficiently large $n$, the conditions of Step 3 hold, and \eqref{eq:two-cov} gives the coverage bound of \cref{thm:local-cov}(i). The union of the failure events has probability at most $Cn^{-c}$.

Under the conditions of \cref{thm:local}, with $\phi_x$ as in \cref{cond:A} and $\hat t_x$ as in \cref{sec:cov-procedure}, we have
	\begin{align*}
		\{y:\phi_x(\hat a(y\mid x))>\phi_x(\hat t_x)\}\subset\hC^{\rm two}(x)\subset\{y:\phi_x(\hat a(y\mid x))\ge\phi_x(\hat t_x)\}\qquad(x\in\cX_0),
	\end{align*}
and, by the coverage bound, the coverage of $\hC^{\rm two}(x)$ is in the window. By \cref{cond:A}, and since $Cr_n\le\varepsilon_0/8$ for all sufficiently large $n$, \cref{prop:supshape} with $r=Cr_n$ and \eqref{eq:L-props} give the shape and excess-mass bounds, and \cref{lem:identity,lem:t-window} give the volume bound.

\emph{The exceptional set for \cref{rem:local-each}.} Let
	\begin{align*}
		\mathcal X_{\rm bad}:=\bigl\{x\in\cX_0:\ \mathcal D_{\rm fit}\notin\Omega_A(x)\ \text{or}\ \sup_{t\in\R}|\widehat Q_x(t)-Q_x(t)|>C_{\rm acc}\rho_n\bigr\}.
	\end{align*}
Both $\widehat Q_x(t)$ and $Q_x(t)$ are left-continuous in $t$, and hence the supremum may be taken over rational $t$. With the measurability required in \cref{rem:local-each}, the indicator of $\{x\in\mathcal X_{\rm bad}\}$ is jointly measurable in $x$ and the data. By Fubini's theorem and \cref{lem:loc-acc}(i),
	\begin{align*}
		\E\sbr[1]{\Pi(\mathcal X_{\rm bad})}\le\int_{\cX_0}\Bigl\{\bP\bigl(\Omega_A(x)^c\bigr)+Cn^{-2c-2}\Bigr\}\,\Pi(\d x)\le Cn^{-2c-2}.
	\end{align*}
Hence, by Markov's inequality, $\bP\{\Pi(\mathcal X_{\rm bad})>n^{-c-1}\}\le Cn^{-c-1}$. On the complement of this event, we have $\Pi(\mathcal X_{\rm bad})\le n^{-c-1}\le n^{-c}$, Steps 1 and 2 apply at every $x\in\cX_0\setminus\mathcal X_{\rm bad}$ with $C_{\rm acc}\rho_n$ in place of $e_n$, and \eqref{eq:Hx} holds there. Since $|H_x(u)-u|\le1$ for every $x$ and $\bP(X\in\mathcal X_{\rm bad}\mid X\in\cX_0)\le n^{-c-1}/\pi_0\le C\rho_n$, we have $\sup_{u\in[0,1]}|H(u)-u|\le C\rho_n$. Step 3 then gives the coverage bound at every $x\in\cX_0\setminus\mathcal X_{\rm bad}$, and Step 4 gives the other bounds there.

\emph{Step 5: finite-sample validity.} For arbitrary fitting and regression data, $\hat{\mathsf s}^{\rm two}$ is a fixed Borel function of $(x,y)$. Conditionally on these data and on which calibration covariates lie in $\cX_0$, the retained calibration pairs and an independent test pair drawn from the distribution of $(X,Y)$ given $X\in\cX_0$ are independent and identically distributed, and hence their scores are exchangeable. Let $m$ and $k$ be as in Step 3. If $k\le m$, the probability that the test score does not exceed the $k$-th calibration order statistic is at least $1-\alpha$, even with ties, by the quantile argument for exchangeable scores \citep[Section~2]{lei2018distribution}. If $k>m$, the threshold is $1$, the output is $\cY$, and the probability is one. Averaging over the conditioning gives the stated marginal validity for every $n$ and every ranking function.
\end{proof}

\subsection{Proof of \texorpdfstring{\cref{prop:nn-conditions}}{Proposition \ref{prop:nn-conditions}}}\label{app:nn}
We fix the fitting sample and write $f:=f_{\hat\theta}$ and $\hat a:=\hat a^{\rm nn}$. Let $\operatorname{act}$ be the activation function, a polynomial in $t$ and $g_1(t),\dots,g_q(t)$ as in \cref{sec:cov-nn}. The functions $g_1,\dots,g_q$ solve a triangular system of ordinary differential equations with polynomial right-hand sides, and are therefore real-analytic. Hence $\operatorname{act}$ and $f$ are real-analytic.

\begin{proof}
\emph{\Cref{cond:N}.} Let $x\in\cX_0$ and $v\in\R$. The function
	\begin{align*}
		g(y):=K{\mathsf q}_0\exp\{f(x,y)\}+n^{-2}\mathrm{pol}_J(y)-v
	\end{align*}
is real-analytic on $\R^{d_Y}$. We fix $y_0\in\R^{d_Y}$ and let $\mathbf e_1$ be the first coordinate vector. As $s\to\infty$, the first coordinate of $z(y_0+s\mathbf e_1)$ tends to infinity, and the other coordinates stay fixed. Hence $\mathrm{pol}_J(y_0+s\mathbf e_1)\to\infty$, and, since the exponential term is positive, $g(y_0+s\mathbf e_1)\to\infty$. Thus $g$ is not identically zero, and its zero set, which contains the level set $\{y\in\cY:\hat a(y\mid x)=v\}$, is Lebesgue-null \citep{mityagin2020zero}. This holds for every $x\in\cX_0$, $v\in\R$ and fitting sample.

\emph{\Cref{cond:V}.} Let $u_1,\dots,u_U$ be the pre-activations of the units of $f$, ordered by layer, as functions of $w:=(y,x)\in\R^{d_Y+d_X}$. By induction on the layers, each $u_j$ is a polynomial in $w$ and in the values $g_i(u_k)$ at the units $k$ of earlier layers. By the chain rule and $g_i'(t)=P_i(t,g_1(t),\dots,g_i(t))$ for polynomials $P_i$, every partial derivative of $g_i(u_j)$ is a polynomial in $w$, in the $g_{i'}(u_k)$ for earlier units $k$, and in $g_1(u_j),\dots,g_i(u_j)$. Hence the functions $g_i(u_j)$, ordered by unit and then by $i$, followed by $\exp\{f(w)\}$, form a Pfaffian chain in $w$ of length $qU+1$, whose degree depends only on the architecture and $\operatorname{act}$. The function
	\begin{align*}
		\Psi(y;x,v):=\hat a(y\mid x)-v=K{\mathsf q}_0\exp\{f(x,y)\}+n^{-2}\mathrm{pol}_J(y)-v
	\end{align*}
is a polynomial in $(y,x,v)$ and in the members of this chain, of degree at most $J$. Its coefficients depend on the fitting sample and on $n$, but its degrees and the length and degree of the chain do not. Therefore, by \citet[Section~4.2]{karpinski1997polynomial}, the family of the sets $\{y:\Psi(y;x,v)>0\text{ or }\Psi(y;x,v)=0\}$, $(x,v)\in\R^{d_X+1}$, has VC dimension at most a constant $D_{\rm nn}$ that depends only on the architecture, $\operatorname{act}$, $J$, $d_X$ and $d_Y$. Every set of $\mathcal A_n$ is the intersection of such a set with $\cY$. By the lemma of \citet{sauer1972density}, for $n\ge2$ and all $m\ge1$,
	\begin{align*}
		\mathfrak S_{\mathcal A_n}(m)\le\sum_{i=0}^{D_{\rm nn}}\binom mi\le(m+1)^{D_{\rm nn}}\le(nm)^{D_{\rm nn}}.
	\end{align*}
This is \cref{cond:V} with $C_{\rm gr}=D_{\rm nn}$, for every fitting sample.

Part (ii) follows from (i) and \cref{thm:local-cov}.
\end{proof}

\subsection{Accuracy of the local polynomial estimate}\label{app:lp}\label{app:lp-estimator}
This subsection collects the properties of the local polynomial estimate used in \cref{sec:cov-twostage-alg,sec:cov-onestage}. In this subsection and in \cref{app:lp-mix}, we write $\tilde a$ and $\hat a$ for $\tilde a^{\rm lp}_b$ and $\hat a^{\rm lp}_b$ at the bandwidth $b$ under consideration, and we let $\nu_b:=\{\log n/(nb^d)\}^{1/2}$ and $\eps_b:=b^\beta+\nu_b$. The bias of the estimate is of order $b^\beta$ and its sup-norm noise of order $\nu_b$, and at $b=h$ both are of order $\zeta_n$. Let $\bar h:=\min\{\delta_0,w_\cY\}$, with $w_\cY$ as in \cref{app:lower}.

The fitted intercept $\heta^{\rm raw}(z)$ is clipped to $[0,\bar\eta']$, and the result is $\heta$. By the choice of $\bar\eta'$ in \cref{sec:cov-twostage-alg}, we have $\sup\eta<\bar\eta'$. Let $\widetilde Z:=\int\tilde a\,\d y$, $\hp:=\tilde a/\widetilde Z$ if $\widetilde Z\ge c_Z$, and $\hp:={\mathsf q}_0$ otherwise, for a fixed $c_Z\in(0,1)$. Let $\tilde a(\cdot\mid x):=\tilde a(x,\cdot)$.

\emph{Notation.} Let $f(x,y):=f_X(x)\{{\mathsf p}(y\mid x)+K{\mathsf q}_0(y)\}$, which is $K+1$ times the density of the pooled rows on $\cX_0^{2\delta_0}\times\cY$. By \cref{ass:smooth,ass:design}, $f$ is bounded and Lipschitz on $\cX_0^{2\delta_0}\times\cY$. Moreover, $\eta=\Psi({\mathsf p})$ with
	\begin{align}\label{eq:Psi}
		\Psi(u):=\frac u{u+K{\mathsf q}_0},
	\end{align}
which has bounded derivatives of all orders on $[0,\infty)$. For a row $(z',\iota)$, with label $\iota\in\{0,1\}$, and a point $z$, let $D(z';z)$ be the vector of monomials of degree at most $\ell$ in $(z'-z)/b$, and let $W(z';z):=\mathcal K((z'-z)/b)$. Let
	\begin{align*}
		\widehat\Sigma(z):=(n_{\rm fit}b^d)^{-1}\sum_{\rm rows}WDD^\top,\qquad\widehat{\rm rhs}(z):=(n_{\rm fit}b^d)^{-1}\sum_{\rm rows}WD\iota,
	\end{align*}
where $\sum_{\rm rows}$ runs over the $n_{\rm fit}(K+1)$ rows $(z',\iota)$ of $\mathcal D_{\rm fit}$, and $W$ and $D$ are evaluated at $(z';z)$. The raw estimate is $\heta^{\rm raw}=\mathbf e_1^\top\widehat\Sigma^{-1}\widehat{\rm rhs}$ when $\widehat\Sigma(z)$ is invertible, and $\heta^{\rm raw}(z)=0$ otherwise, so that the estimator is defined on every event. Let
	\begin{align*}
		\bar\eta_b:=\mathbf e_1^\top \Sigma_b^{-1}{\rm rhs}_b,\qquad\Sigma_b:=\E\sbr{\widehat\Sigma},\qquad{\rm rhs}_b:=\E\sbr{\widehat{\rm rhs}}.
	\end{align*}
By the step \emph{Gram matrix} of the proof of \cref{lem:vc}, there is a constant $c_\Sigma>0$ such that $\lambda_{\min}(\Sigma_b(z))\ge c_\Sigma$ for all $z\in\cX_0^{\delta_0}\times\cY$ and $b\le\bar h$.

\begin{lemma}[Uniform stochastic error]\label{lem:vc}
Suppose that \cref{ass:smooth,ass:design,ass:split} hold. Then for every $c>0$ there is $C>0$ such that the following holds for all sufficiently large $n$. Let $b$ be a deterministic bandwidth with $(C\log n/n)^{1/d}\le b\le\bar h$, and let $\Omega_{\rm vc}\in\sigma(\mathcal D_{\rm fit})$ be the event on which, for all $z\in\cX_0^{\delta_0}\times\cY$, we have $\lambda_{\min}(\widehat\Sigma(z))\ge c_\Sigma/2$ and
	\begin{align*}
		|\heta^{\rm raw}(z)-\bar\eta_b(z)|+\norm{\widehat\Sigma(z)-\Sigma_b(z)}+\norm{\widehat{\rm rhs}(z)-{\rm rhs}_b(z)}\le C\sqrt{\log n/(nb^d)} .
	\end{align*}
Then $\bP(\Omega_{\rm vc})\ge1-Cn^{-c}$.
\end{lemma}

\begin{proof}
Since the rows generated by one observation share its covariate, they are dependent. We therefore apply Bernstein's inequality to the i.i.d.\ blocks of rows rather than to single rows.

\emph{Block variables.} We group the rows into the i.i.d.\ blocks $O_i:=(X_i,Y_i,Y'_{i1},\dots,Y'_{iK})$, $i\le n_{\rm fit}$. With $s=(z'-z)/b$, for a function $\phi$ on $\R^d$ and $k_0\in\{0,1\}$ let
	\begin{align*}
		V_{\phi,z,b}(O_i):=b^{-d}\sum_{\text{rows }(z',\iota)\text{ of block }i}\phi(s)\,\iota^{k_0},\qquad\phi_{\boldsymbol k}(s):=\mathcal K(s)s^{\boldsymbol k} .
	\end{align*}
Every entry of $\widehat\Sigma(z)$ and $\widehat{\rm rhs}(z)$ is $n_{\rm fit}^{-1}\sum_{i\le n_{\rm fit}}V_{\phi_{\boldsymbol k},z,b}(O_i)$ for a multi-index $\boldsymbol k$ with $|\boldsymbol k|\le2\ell$ and some $k_0$, and the corresponding entry of $\Sigma_b(z)$ or ${\rm rhs}_b(z)$ is $\E V_{\phi_{\boldsymbol k},z,b}(O_1)$. We write $\phi_{\boldsymbol k}(s)=\prod_{j=1}^d\mathcal K_{1,k_j}(s_j)$, with $\mathcal K_{1,k}(s):=\mathcal K_1(s)s^k$ for integers $k\ge0$. Since $\mathcal K_1$ is continuously differentiable and vanishes outside $[-1,1]$, each $\phi_{\boldsymbol k}$ is continuously differentiable, and $\phi_{\boldsymbol k}$ and its first derivatives vanish outside $[-1,1]^d$. For a continuously differentiable $\phi$, we have $\partial_{z_j}\{\phi(s)\}=-b^{-1}(\partial_j\phi)(s)$, and hence
	\begin{align*}
		b\,\partial_{z_j}V_{\phi,z,b}=-V_{\partial_j\phi,z,b} .
	\end{align*}
For $\phi=\phi_{\boldsymbol k}$ and fixed $b$, the right-hand side $-V_{\partial_j\phi_{\boldsymbol k},z,b}$ is bounded.

Let $\Phi$ be the finite set of the functions $\phi_{\boldsymbol k}$, $|\boldsymbol k|\le2\ell$. We take $C\ge1$, so that $nb^d\ge\log n$ and, for all sufficiently large $n$, $b\ge n^{-1}$.

\emph{Gram matrix.} We show that there is a constant $c_\Sigma>0$ such that $\lambda_{\min}(\Sigma_b(z))\ge c_\Sigma$ for all $z\in\cX_0^{\delta_0}\times\cY$ and $b\le\bar h$. Let $S_z:=\{s\in[-\frac12,\frac12]^d:\ z+bs\in\cX_0^{2\delta_0}\times\cY\}$, a product of intervals because $\cY$ is a box. Since $b\le\delta_0$, the $x$-window lies in $\cX_0^{2\delta_0}$, and hence the intervals of $S_z$ in the $x$-coordinates equal $[-\frac12,\frac12]$. Since $b\le w_\cY$, in each $y$-coordinate the face of $\cY$ farther from $z$ is at distance at least $w_\cY/2\ge b/2$. Thus the interval of $S_z$ in this coordinate contains the half of $[-\frac12,\frac12]$ that points toward this face. Therefore $S_z$ contains one of the $2^{d_Y}$ boxes $[-\frac12,\frac12]^{d_X}\times\prod_{j\le d_Y}I_j$ with $I_j\in\{[-\frac12,0],[0,\frac12]\}$. On these boxes, $\mathcal K$ is bounded below, and the rescaled rows have a density bounded below, since $f\ge f_{\min}K{\mathsf q}_0$ on $\cX_0^{2\delta_0}\times\cY$. The monomials of degree at most $\ell$ are linearly independent on each of these finitely many boxes. It follows that $\lambda_{\min}(\Sigma_b(z))\ge c_\Sigma>0$.

\emph{Concentration.} Let $z\in\cX_0^{\delta_0}\times\cY$ and $\phi\in\Phi$. Since $b\le\delta_0$, the window $\{z':|z'-z|_\infty\le b\}$ meets $\R^{d_X}\times\cY$ only in $\cX_0^{2\delta_0}\times\cY$, where every row of a block has a bounded density. Therefore $\bP(|z'-z|_\infty\le b)\le Cb^d$ for each row $(z',\iota)$ of a block, and hence
	\begin{align*}
		\norm{V_{\phi,z,b}}_\infty\le Cb^{-d},\qquad{\rm Var}\,V_{\phi,z,b}(O_1)\le Cb^{-d}.
	\end{align*}
By Bernstein's inequality with $u=\{c+d(d+3)\}\log n$, together with $n_{\rm fit}\ge cn$ by \cref{ass:split} and $\log n\le nb^d$, we have
	\begin{align*}
		\bP\cbr[3]{\abs[3]{\frac1{n_{\rm fit}}\sum_{i\le n_{\rm fit}}V_{\phi,z,b}(O_i)-\E V_{\phi,z,b}(O_1)}>C\Bigl\{\frac{\log n}{nb^d}\Bigr\}^{1/2}}\le Cn^{-c-d(d+3)}.
	\end{align*}
By the identity $b\,\partial_{z_j}V_{\phi,z,b}=-V_{\partial_j\phi,z,b}$, the block variables and their expectations are Lipschitz in $z$ with constant $Cb^{-d-1}\le Cn^{d+1}$. Let $\mathcal G_n\subset\cX_0^{\delta_0}\times\cY$ be a set of at most $Cn^{d(d+3)}$ points such that every $z\in\cX_0^{\delta_0}\times\cY$ is within $|\cdot|_\infty$-distance $n^{-(d+3)}$ of a point of $\mathcal G_n$. Between $z$ and this point, every entry of $\widehat\Sigma-\Sigma_b$ and $\widehat{\rm rhs}-{\rm rhs}_b$ changes by at most $Cn^{-2}\le C\{\log n/(nb^d)\}^{1/2}$. By a union bound over $\mathcal G_n$ and the finitely many entries, the event
	\begin{align*}
		\Omega_{\rm grid}:=\cbr[3]{\sup_{z\in\cX_0^{\delta_0}\times\cY}\del[2]{\norm{\widehat\Sigma(z)-\Sigma_b(z)}+\norm{\widehat{\rm rhs}(z)-{\rm rhs}_b(z)}}\le C_1\Bigl\{\frac{\log n}{nb^d}\Bigr\}^{1/2}}
	\end{align*}
has probability at least $1-Cn^{-c}$, where $C_1$ does not depend on the constant $C$ in the range of $b$. Since $b\ge(C\log n/n)^{1/d}$, the right-hand side is at most $c_\Sigma/2$ for $C\ge(2C_1/c_\Sigma)^2$.

\emph{The raw fit.} On $\Omega_{\rm grid}$, by the step \emph{Gram matrix} and Weyl's inequality, we have $\lambda_{\min}(\widehat\Sigma)\ge c_\Sigma/2$. Therefore, $\widehat\Sigma$ is invertible and $\norm[0]{\widehat\Sigma^{-1}}\le2/c_\Sigma$, and we can write
	\begin{align*}
		\heta^{\rm raw}-\bar\eta_b&=\mathbf e_1^\top\bigl(\widehat\Sigma^{-1}\widehat{\rm rhs}-\Sigma_b^{-1}{\rm rhs}_b\bigr)=\mathbf e_1^\top\widehat\Sigma^{-1}\bigl[\widehat{\rm rhs}-\widehat\Sigma\Sigma_b^{-1}{\rm rhs}_b\bigr]\\
		&=\mathbf e_1^\top\widehat\Sigma^{-1}\bigl[(\widehat{\rm rhs}-{\rm rhs}_b)-(\widehat\Sigma-\Sigma_b)\Sigma_b^{-1}{\rm rhs}_b\bigr].
	\end{align*}
Since $|\iota|\le1$ and $\lambda_{\min}(\Sigma_b)\ge c_\Sigma$, the vector $\Sigma_b^{-1}{\rm rhs}_b$ is bounded, and hence
	\begin{align*}
		|\heta^{\rm raw}-\bar\eta_b|\le C\bigl(\norm{\widehat\Sigma-\Sigma_b}+\norm{\widehat{\rm rhs}-{\rm rhs}_b}\bigr).
	\end{align*}
With the bound of the step \emph{Concentration}, this gives the bound of the lemma on $\Omega_{\rm grid}$. Thus $\Omega_{\rm grid}\subset\Omega_{\rm vc}$ and $\bP(\Omega_{\rm vc})\ge1-Cn^{-c}$.
\end{proof}

\begin{lemma}[Bias]\label{lem:bias}
Suppose that \cref{ass:smooth,ass:design} hold. Then there exists a constant $C>0$ such that $|\bar\eta_b(z)-\eta(z)|\le Cb^\beta$ for all $z\in\cX_0^{\delta_0}\times\cY$ and $b\le\bar h$.
\end{lemma}

\begin{proof}
Let $\mathcal T_z$ be the Taylor polynomial of $\eta$ at $z$ of degree $\ell$. In each block, only the row $(X_i,Y_i)$ has label $1$, and $\eta({\mathsf p}+K{\mathsf q}_0)={\mathsf p}$. Therefore, for every bounded function $g$ on rows that vanishes outside $\cX_0^{2\delta_0}\times\cY$, we have
	\begin{align*}
		\E\sbr{g(X_i,Y_i)}=\int f_X{\mathsf p}\,g\,\d z'=\int f\,\eta\,g\,\d z'=\E\sbr[2]{\sum_{\text{rows }z'\text{ of block }i}g(z')\eta(z')} .
	\end{align*}
It follows that
	\begin{align*}
		{\rm rhs}_b(z)=\E\sbr[2]{(n_{\rm fit}b^d)^{-1}\sum_{\rm rows}W\,D\,\eta(z')}.
	\end{align*}
Since $\mathcal T_z(z')$ is a polynomial of degree $\ell$ in $(z'-z)/b$, it equals $D(z';z)^\top\theta^{\rm cf}_z$ for a coefficient vector $\theta^{\rm cf}_z$ with first entry $\mathcal T_z(z)=\eta(z)$. Thus
	\begin{align*}
		\mathbf e_1^\top\Sigma_b(z)^{-1}\E\sbr[2]{(n_{\rm fit}b^d)^{-1}\sum_{\rm rows}W\,D\,\mathcal T_z(z')}=\mathbf e_1^\top\Sigma_b(z)^{-1}\Sigma_b(z)\theta^{\rm cf}_z=\eta(z).
	\end{align*}
Subtracting this identity from $\bar\eta_b(z)=\mathbf e_1^\top\Sigma_b(z)^{-1}{\rm rhs}_b(z)$ and using the expression for ${\rm rhs}_b(z)$ above, we obtain
	\begin{align*}
		\bar\eta_b(z)-\eta(z)=\mathbf e_1^\top \Sigma_b(z)^{-1}\E\sbr[2]{(n_{\rm fit}b^d)^{-1}\sum_{\rm rows}W\,D\,\bigl\{\eta(z')-\mathcal T_z(z')\bigr\}}.
	\end{align*}
It remains to bound the expectation in the last display by $Cb^\beta$. Since $\eta=\Psi({\mathsf p})$, the chain rule shows that $\eta$ is in $\mathcal H^\beta$ on the convex set $\cX_0^{2\delta_0}\times\cY$. This set contains the segments from $z$ to the rows in the window of $z$. Since the derivatives of $\eta$ of order $\ell$ are H\"older of order $\beta-\ell\in(0,1]$, it follows from Taylor's theorem that $|\eta(z')-\mathcal T_z(z')|\le C|z'-z|^\beta\le Cb^\beta$ on the window. Moreover, the product $W\,D$ is bounded, and, by the steps \emph{Gram matrix} and \emph{Concentration} of the proof of \cref{lem:vc}, $\norm[0]{\Sigma_b(z)^{-1}}\le1/c_\Sigma$ and $\E\sbr[1]{\sum_{\rm rows}\ind\{|z'-z|_\infty\le b\}}\le Cn_{\rm fit}b^d$. Therefore, the norm of the expectation is at most $(n_{\rm fit}b^d)^{-1}\cdot Cn_{\rm fit}b^d\cdot Cb^\beta\le Cb^\beta$.
\end{proof}

Let
	\begin{align}\label{eq:a-bar}
		\bar a:=K{\mathsf q}_0\psi(\bar\eta_b),\qquad\bar Z(x):=\int\bar a(x,y)\,\d y,\qquad\bar{\mathsf p}:=\bar a/\bar Z.
	\end{align}
Unlike $\heta$, the function $\bar\eta_b$ is not clipped. Let $\bar\eta'':=(\sup\eta+\bar\eta')/2$. Under \cref{ass:smooth,ass:design}, by \cref{lem:bias}, we have $|\bar\eta_b-\eta|\le Cb^\beta$ on $\cX_0^{\delta_0}\times\cY$. Since $0\le\eta\le\sup\eta<\bar\eta''$, it follows that $\bar\eta_b$ takes values in $[-\frac12,\bar\eta'']$ for all sufficiently small $b$. Since $\psi$ is Lipschitz on this interval, $K{\mathsf q}_0\psi(\eta)={\mathsf p}$ and $\int{\mathsf p}(y\mid x)\,\d y=1$, we have
	\begin{align*}
		|\bar Z(x)-1|\le\Leb(\cY)\sup_y|\bar a-{\mathsf p}|\le Cb^\beta,
	\end{align*}
which is at most $\frac12(1-c_Z)$ for all sufficiently small $b$. Therefore, there is $h_0\in(0,\bar h]$ such that, for $b\le h_0$, the function $\bar\eta_b$ takes values in $[-\frac12,\bar\eta'']$ on $\cX_0^{\delta_0}\times\cY$ and $\bar Z\ge\max\{c_Z,\frac12\}$ on $\cX_0^{\delta_0}$. In particular, $\bar{\mathsf p}$ is well defined on $\cX_0^{\delta_0}\times\cY$ for $b\le h_0$.

We write $B:=\bar{\mathsf p}-{\mathsf p}$ for the bias and $N:=\hp-\bar{\mathsf p}$ for the noise, where $\bar{\mathsf p}$ need not equal $\E\hp$. For $b\le h_0$, since $\bar Z\ge\frac12$ and ${\mathsf p}\le{\mathsf p}_{\max}$, we have, on $\cX_0^{\delta_0}\times\cY$,
	\begin{align*}
		|B|\le\frac{|\bar a-{\mathsf p}|}{\bar Z}+\frac{{\mathsf p}|1-\bar Z|}{\bar Z}\le Cb^\beta .
	\end{align*}

\begin{lemma}[Sup-norm error]\label{lem:supnorm-b}
Suppose that \cref{ass:smooth,ass:design,ass:split} hold, and let $b\to0$ be a deterministic bandwidth with $nb^d/\log n\to\infty$. Then for every $c>0$ there are constants $C,C_{\rm clip}>0$ such that, for all sufficiently large $n$, we have $b\le h_0$, the bandwidth $b$ lies in the range of \cref{lem:vc} with this $c$, and the following hold on the event $\Omega_{\rm vc}$ of that lemma.
\begin{enumerate}[label=(\roman*)]
\item $\widetilde Z\ge\max\{c_Z,\frac12\}$ on $\cX_0^{\delta_0}$, and, on $\cX_0^{\delta_0}\times\cY$,
	\begin{align*}
		|\tilde a-{\mathsf p}|\le C\eps_b,\qquad|\hp-{\mathsf p}|\le C\eps_b,\qquad|N|\le C\eps_b .
	\end{align*}
\item The clipping is inactive at every $(x,y)$ with $x\in\cX_0^{\delta_0}$ and ${\mathsf p}(y\mid x)>C_{\rm clip}\eps_b$.
\end{enumerate}
\end{lemma}

\begin{proof}
Since $b\to0$ and $nb^d/\log n\to\infty$, for all sufficiently large $n$, we have $b\le h_0$ and $nb^d\ge C\log n$ with the constant $C$ of \cref{lem:vc}, and hence $b$ is in the range of that lemma. Moreover, $\eps_b\to0$. On $\Omega_{\rm vc}$, by \cref{lem:vc,lem:bias}, we have
	\begin{align*}
		|\heta^{\rm raw}-\eta|\le|\heta^{\rm raw}-\bar\eta_b|+|\bar\eta_b-\eta|\le C\nu_b+Cb^\beta\le C_1\eps_b
	\end{align*}
for a constant $C_1>0$. Since clipping to $[0,\bar\eta']$ is $1$-Lipschitz and fixes $\eta$, it follows that $|\heta-\eta|\le C_1\eps_b$. Since $b\le h_0$, the function $\bar\eta_b$ takes values in $[-\frac12,\bar\eta'']$, and $\psi$ is Lipschitz on $[-\frac12,\bar\eta']$. Hence $|\tilde a-{\mathsf p}|\le C\eps_b$ and $|\tilde a-\bar a|\le C\eps_b$. Therefore
	\begin{align*}
		|\widetilde Z-1|\le\Leb(\cY)C\eps_b,\qquad|\widetilde Z-\bar Z|\le\Leb(\cY)C\eps_b,
	\end{align*}
which implies that, for all sufficiently large $n$, we have $\widetilde Z\ge\max\{c_Z,\frac12\}$ and $\hp=\tilde a/\widetilde Z$. Moreover, since $b\le h_0$, we have $\bar Z\ge\max\{c_Z,\frac12\}$. It follows that
	\begin{align*}
		|\hp-{\mathsf p}|&=\Bigl|\frac{\tilde a-{\mathsf p}}{\widetilde Z}+{\mathsf p}\Bigl(\frac1{\widetilde Z}-1\Bigr)\Bigr|\le C\eps_b,\\
		|N|&=\Bigl|\frac{\tilde a-\bar a}{\widetilde Z}-\frac{\bar a(\widetilde Z-\bar Z)}{\widetilde Z\bar Z}\Bigr|\le C\eps_b.
	\end{align*}

Let $C_{\rm clip}:=2C_1({\mathsf p}_{\max}+K{\mathsf q}_0)$. Since $\eta={\mathsf p}/({\mathsf p}+K{\mathsf q}_0)$, if ${\mathsf p}(y\mid x)>C_{\rm clip}\eps_b$, then
	\begin{align*}
		\eta(x,y)\ge\frac{{\mathsf p}(y\mid x)}{{\mathsf p}_{\max}+K{\mathsf q}_0}>2C_1\eps_b,
	\end{align*}
and hence $\heta^{\rm raw}(x,y)\ge\eta(x,y)-C_1\eps_b>C_1\eps_b>0$. Moreover, $\heta^{\rm raw}\le\sup\eta+C_1\eps_b<\bar\eta'$ for all sufficiently large $n$. Therefore the clipping is inactive at $(x,y)$.
\end{proof}

\begin{lemma}[Linearization of the local polynomial fit]\label{lem:linear}
Suppose that \cref{ass:smooth,ass:design,ass:split} hold, let $c>0$, and let $b$ be a deterministic bandwidth in the range of \cref{lem:vc} with this $c$. Let $\bar\theta^{\rm cf}:=\Sigma_b^{-1}{\rm rhs}_b$, so that $\mathbf e_1^\top\bar\theta^{\rm cf}=\bar\eta_b$, and, for the rows $(z_{ij},\iota_{ij})$, $j=0,\dots,K$, of block $i$, let
	\begin{align*}
		g_{ij}(z):=W(z_{ij};z)\,\mathbf e_1^\top\Sigma_b(z)^{-1}D(z_{ij};z)\bigl(\iota_{ij}-D(z_{ij};z)^\top\bar\theta^{\rm cf}(z)\bigr),\qquad\xi_i:=\sum_{j=0}^Kg_{ij}.
	\end{align*}
Then there is a constant $C>0$, not depending on $b$, such that, for all sufficiently large $n$ and all $z\in\cX_0^{\delta_0}\times\cY$, the variables $\xi_i(z)$, $i\le n_{\rm fit}$, are i.i.d.\ with mean zero, $|g_{ij}(z)|\le C$, and $g_{ij}(z)=0$ unless $|z_{ij}-z|_\infty\le b$. Moreover, on the event $\Omega_{\rm vc}$ of \cref{lem:vc}, the linear term $\xi:=(n_{\rm fit}b^d)^{-1}\sum_{i\le n_{\rm fit}}\xi_i$ and the remainder $N_\Sigma:=\heta^{\rm raw}-\bar\eta_b-\xi$ satisfy $|N_\Sigma|\le C\nu_b^2$ on $\cX_0^{\delta_0}\times\cY$.
\end{lemma}

\begin{proof}
Since $\sum_{\rm rows}WD(\iota-D^\top\bar\theta^{\rm cf})=n_{\rm fit}b^d(\widehat{\rm rhs}-\widehat\Sigma\bar\theta^{\rm cf})$, we have
	\begin{align*}
		\xi=\mathbf e_1^\top\Sigma_b^{-1}\bigl(\widehat{\rm rhs}-\widehat\Sigma\bar\theta^{\rm cf}\bigr).
	\end{align*}
The variables $\xi_i(z)$ are i.i.d.\ because the blocks are, and $\E\xi_i(z)=0$ since ${\rm rhs}_b=\Sigma_b\bar\theta^{\rm cf}$. Since $\int_{\{z':|z'-z|_\infty\le b\}}f(z')\,\d z'\le Cb^d$, the vector ${\rm rhs}_b$ is bounded. Hence, since $\mathcal K$ vanishes outside $[-1,1]^d$ and $\norm[0]{\Sigma_b^{-1}}\le1/c_\Sigma$, the function $g_{ij}$ is bounded and vanishes unless $|z_{ij}-z|_\infty\le b$. On $\Omega_{\rm vc}$, the matrix $\widehat\Sigma$ is invertible, and
	\begin{align*}
		\heta^{\rm raw}-\bar\eta_b=\mathbf e_1^\top\widehat\Sigma^{-1}\bigl(\widehat{\rm rhs}-\widehat\Sigma\bar\theta^{\rm cf}\bigr).
	\end{align*}
Hence
	\begin{align*}
		N_\Sigma=\mathbf e_1^\top\widehat\Sigma^{-1}\bigl(\Sigma_b-\widehat\Sigma\bigr)\Sigma_b^{-1}\bigl\{(\widehat{\rm rhs}-{\rm rhs}_b)-(\widehat\Sigma-\Sigma_b)\bar\theta^{\rm cf}\bigr\},
	\end{align*}
and, by \cref{lem:vc}, $|N_\Sigma|\le C\nu_b^2$.
\end{proof}

\subsection{Proof of \texorpdfstring{\cref{prop:lp-conditions}}{Proposition \ref{prop:lp-conditions}}}\label{app:lp-mix}
Let $k_{\mathcal K}$ be the degree of the polynomial that agrees with $\mathcal K_1$ on $[-1,1]$, and let $s_\ell:=\binom{d+\ell}{\ell}$ be the number of monomials in the local polynomial fit. One valid choice of the constant $D_0$ of \cref{sec:cov-twostage-alg} is
	\begin{align*}
		D_0:=s_\ell(k_{\mathcal K}d+2\ell).
	\end{align*}
Since $J$ is even, $\mathrm{pol}_J\ge0$ and $\hat a\ge\tilde a\ge0$.

\begin{proof}
\emph{\Cref{cond:A}.} Since $h\asymp(\log n/n)^{1/(2\beta+d)}$, we have $h\to0$, $nh^d/\log n\to\infty$ and $\eps_h\asymp\zeta_n$. Hence, for all sufficiently large $n$, the bounds of \cref{lem:supnorm-b}(i) hold at $b=h$ on the event $\Omega_{\rm vc}$ of \cref{lem:vc}, whose probability is at least $1-Cn^{-c}$, and give $\widetilde Z\ge\frac12$ and $|\tilde a/\widetilde Z-{\mathsf p}|\le C\zeta_n$ on $\cX_0\times\cY$. Let $\widehat Z(x):=\int_\cY\hat a(y\mid x)\,\d y$. Since $|\hat a-\tilde a|\le d_Yn^{-2}$, we have $|\widehat Z-\widetilde Z|\le d_Y\Leb(\cY)n^{-2}$, and, since $0\le\tilde a\le K{\mathsf q}_0\psi(\bar\eta')$, also $|\hat a/\widehat Z-\tilde a/\widetilde Z|\le Cn^{-2}$ on $\cX_0\times\cY$ for all sufficiently large $n$. Since $n^{-2}=o(\zeta_n)$, the normalization $a\mapsto a/\widehat Z(x)$ satisfies both bounds of \cref{cond:A} with $r_n=\zeta_n$ for all sufficiently large $n$.

\emph{Polynomial representations.} We fix the fitting sample and the bandwidth $h$, and we write $c_K:=K\mathsf q_0$ and $a_{\max}:=c_K\bar\eta'/(1-\bar\eta')$. We divide the query space by the coordinate hyperplanes at distance $h$ from the coordinates of the pooled rows. On each resulting piece, the set of rows in the kernel window is fixed, and every kernel weight is a polynomial in $(x,y)$ of total degree at most $k_{\mathcal K}d$. Hence each entry of $\widehat\Sigma$ has degree at most $k_{\mathcal K}d+2\ell$, and each entry of $\widehat{\rm rhs}$ has degree at most $k_{\mathcal K}d+\ell$. Let
	\begin{align*}
		\mathfrak D:=\det\widehat\Sigma,\qquad\mathfrak N:=\mathbf e_1^\top\operatorname{adj}(\widehat\Sigma)\widehat{\rm rhs}.
	\end{align*}
Both are polynomials of degree at most $D_0$, because $\widehat\Sigma$ is an $s_\ell\times s_\ell$ matrix. Since $\widehat\Sigma$ is positive semidefinite, it is invertible exactly when $\mathfrak D>0$, and then $\heta^{\rm raw}=\mathfrak N/\mathfrak D$. When $\mathfrak D=0$, we have $\heta^{\rm raw}=0$ by convention, and hence $\tilde a=0$. For $\mathfrak D>0$, the unnormalized estimate takes the three forms
	\begin{align}\label{eq:mixture-branches}
		\tilde a=
		\begin{cases}
		0,&\mathfrak N\le0,\\
		c_K\mathfrak N/(\mathfrak D-\mathfrak N),&0<\mathfrak N<\bar\eta'\mathfrak D,\\
		a_{\max},&\mathfrak N\ge\bar\eta'\mathfrak D.
		\end{cases}
	\end{align}
In the middle branch, we have $\mathfrak D-\mathfrak N>0$, since $\bar\eta'<1$. At a clipping boundary the adjacent formulas agree. The formulas remain valid on kernel-window boundaries by the continuity of the kernel, whose value at $\pm1$ is zero.

\emph{Null level sets.} We fix $x$ and $v\in\R$. The response-coordinate hyperplanes of the step \emph{Polynomial representations} partition $\cY$, up to a Lebesgue-null set, into at most $C n^{d_Y}$ boxes. On each box and each branch of \eqref{eq:mixture-branches}, $\tilde a(\cdot\mid x)$ is a ratio $\mathfrak P(y)/\mathfrak Q(y)$ of polynomials of degree at most $D_0$, with $\mathfrak Q\ne0$ on the branch. The singular-matrix and lower-clipping branches have $\mathfrak P=0$ and $\mathfrak Q=1$, the upper-clipping branch has $\mathfrak P=a_{\max}$ and $\mathfrak Q=1$, and the middle branch has $\mathfrak P=c_K\mathfrak N$ and $\mathfrak Q=\mathfrak D-\mathfrak N$. On such a branch, the equation $\hat a(y\mid x)=v$ implies
	\begin{align*}
		\mathfrak H(y):=\mathfrak P(y)+\{n^{-2}\mathrm{pol}_J(y)-v\}\mathfrak Q(y)=0.
	\end{align*}
The polynomial $\mathfrak H$ is nonzero. Indeed, $\mathrm{pol}_J$ has degree $J$, and hence $\mathrm{pol}_J\mathfrak Q$ has degree $J+\deg \mathfrak Q>D_0$, whereas both $\mathfrak P$ and $v\mathfrak Q$ have degree at most $D_0$. Therefore the highest-degree term of $\mathrm{pol}_J\mathfrak Q$ cannot cancel. A nonzero polynomial has a Lebesgue-null zero set \citep{mityagin2020zero}. Since there are finitely many boxes and branches, every level set of $\hat a(\cdot\mid x)$ is Lebesgue-null, for every $x$ and every fitting sample.

\emph{A sign-counting bound.} We use the following elementary fact, in which zero signs are counted. If $r$ real polynomials of degree at most $D$ depend on $d'$ real parameters, the number of vectors of their signs in $\{-1,0,1\}^r$ is at most
	\begin{align}\label{eq:polynomial-sign-count}
		C_{d',D}(r+1)^{L^{\rm mon}_{d',D}},\qquad L^{\rm mon}_{d',D}:=\binom{d'+D}{D}.
	\end{align}
To see this, we lift the parameters to the vector of all monomials of degree at most $D$, so that each polynomial becomes a linear functional on $\R^{L^{\rm mon}_{d',D}}$. The number of faces of all dimensions in an arrangement of $r$ hyperplanes in this space is at most $\sum_{j=0}^{L^{\rm mon}_{d',D}}2^j\binom rj$. This bound follows by induction: adding a hyperplane increases the number of faces by at most twice the number of faces in the induced arrangement on that hyperplane. Each realizable sign vector corresponds to a face, which proves \eqref{eq:polynomial-sign-count}.

\emph{\Cref{cond:V}.} We fix arbitrary response points $y_1,\ldots,y_m$, and we divide the covariate space by the coordinate hyperplanes at distance $h$ from the coordinates of the covariates $X_i$. Counting the faces of this subdivision as well as its full-dimensional pieces gives at most $C n^{d_X}$ pieces. On each piece and for each $y_j$, the active rows in the $x$-coordinates are fixed, and their activity in the $y$-coordinates is already fixed by $y_j$. Thus the polynomials $\mathfrak D_j(x):=\mathfrak D(x,y_j)$ and $\mathfrak N_j(x):=\mathfrak N(x,y_j)$ have degree at most $D_0$ on this piece. Let $q_j:=n^{-2}\mathrm{pol}_J(y_j)$. Whether $\hat a(y_j\mid x)\ge v$ is determined by the signs of the following six polynomials in $(x,v)$:
	\begin{align*}
		&\mathfrak D_j,\quad\mathfrak N_j,\quad\mathfrak N_j-\bar\eta'\mathfrak D_j,\\
		&q_j-v,\quad q_j+a_{\max}-v,\quad(q_j-v)(\mathfrak D_j-\mathfrak N_j)+c_K\mathfrak N_j.
	\end{align*}
The fourth polynomial applies on the singular-matrix or lower-clipping branch, the fifth on the upper-clipping branch, and the sixth on the middle branch after multiplication by the positive denominator $\mathfrak D_j-\mathfrak N_j$. All these polynomials have total degree at most $D_0+1$. Their zero signs specify the weak inequalities and the branch boundaries. By \eqref{eq:polynomial-sign-count} applied to the $6m$ polynomials and the $d_X+1$ parameters, there are at most $C(m+1)^{L^{\rm mon}}$ patterns on each covariate piece, where $L^{\rm mon}:=L^{\rm mon}_{d_X+1,D_0+1}$ is fixed. Combining the covariate pieces gives, for every fitting sample and all $m\ge1$,
	\begin{align*}
		\mathfrak S_{\mathcal A_n}(m)\le Cn^{d_X}(m+1)^{L^{\rm mon}}\le(nm)^{C_{\rm gr}}
	\end{align*}
for all sufficiently large $n$ and a fixed $C_{\rm gr}$. This is \cref{cond:V}.

Part (ii) follows from (i) and \cref{thm:local-cov,thm:local}.
\end{proof}

\subsection{Proof of \texorpdfstring{\cref{prop:partition}}{Proposition \ref{prop:partition}}}\label{app:partition}
\begin{proof}
We write $\hp$ for $\hp^{\rm part}$. Let
	\begin{align*}
		\nu_{\rm part}:=\Bigl\{\frac{\log n}{nh_X^{d_X}h_Y^{d_Y}}\Bigr\}^{1/2},\qquad r_n':=h_Y+\nu_{\rm part},
	\end{align*}
so that $r_n=h_X+r_n'$. Since $r_n\to0$, we take $n$ sufficiently large, so that $h_X\le\min\{1,\delta_0/2\}$, $h_Y\le\min\{1,w_\cY\}$ and $\nu_{\rm part}\le1$. Hence $h_Y\ge n^{-1}$ and $h_X^{-d_X}\le n$.

Let $\phi$ be the standard Gaussian density on $\R^{d_Y}$, $\phi_{h_Y}:=h_Y^{-d_Y}\phi(\cdot/h_Y)$ and $Z_{h_Y}(y):=\int_\cY\phi_{h_Y}(y-u)\,\d u$. The cubes of the grid are half-open, so that every $x\in\cX$ lies in exactly one of them. For a cube $\mathcal Q$ of the grid with $n_{\mathcal Q}:=\#\{i\in \mathcal D_{\rm fit}:X_i\in \mathcal Q\}\ge1$, the estimate at $x\in \mathcal Q$ is $\hp_{\mathcal Q}:=\hat a_{\mathcal Q}/\int_\cY\hat a_{\mathcal Q}\,\d y$, where
	\begin{align*}
		\hat a_{\mathcal Q}(y):=\frac1{n_{\mathcal Q}\,Z_{h_Y}(y)}\sum_{i\in \mathcal D_{\rm fit}:\,X_i\in \mathcal Q}\phi_{h_Y}(y-Y_i)\qquad(y\in\cY),
	\end{align*}
and $\hp(\cdot\mid x):={\mathsf q}_0$ if $n_{\mathcal Q}=0$.

Below, $\mathcal Q$ denotes a cube that meets $\cX_0^{\delta_0}$. There are at most $Ch_X^{-d_X}\le Cn$ such cubes. Since $h_X\le\delta_0/2$, we have $\bar{\mathcal Q}\subset\cX_0^{2\delta_0}$, and hence, by \cref{ass:design}, $f_{\min}h_X^{d_X}\le\Pi(\mathcal Q)\le Lh_X^{d_X}$. Let ${\mathsf p}_{\mathcal Q}:=\Pi(\mathcal Q)^{-1}\int_{\mathcal Q}{\mathsf p}(\cdot\mid x)\,\Pi(\d x)$ be the density of $Y$ given $X\in \mathcal Q$.

\emph{Step 1: sup-norm accuracy.} We show that, for every $c>0$, on an event $\Omega_{\rm p}$ of probability at least $1-Cn^{-c}$, every $\mathcal Q$ has $n_{\mathcal Q}\ge1$ and
	\begin{align}\label{eq:part-sup}
		\sup_{y\in\cY}\abs[0]{\hp_{\mathcal Q}(y)-{\mathsf p}_{\mathcal Q}(y)}\le Cr_n',\qquad\sup_{y\in\cY}\abs[0]{\hp_{\mathcal Q}(y)-{\mathsf p}(y\mid x)}\le Cr_n\quad(x\in\bar{\mathcal Q}).
	\end{align}
By \eqref{eq:holder} with $\beta\ge1$, ${\mathsf p}$ is Lipschitz on the convex set $\cX_0^{2\delta_0}\times\cY$. Hence ${\mathsf p}_{\mathcal Q}$ is Lipschitz, and $\abs[0]{{\mathsf p}_{\mathcal Q}(y)-{\mathsf p}(y\mid x)}\le Ch_X$ for $x\in\bar{\mathcal Q}$. Consequently, the second bound in \eqref{eq:part-sup} follows from the first.

Let $y\in\cY$. Since $h_Y\le w_\cY$, for each $j$, one of the two parts into which $y_j$ splits the $j$-th side of $\cY$ has length at least $h_Y/2$. Hence $Z_{h_Y}(y)\ge c_\phi$ for a constant $c_\phi>0$. Let
	\begin{align*}
		\bar a_{\mathcal Q}(y):=\frac1{Z_{h_Y}(y)}\int_\cY\phi_{h_Y}(y-u){\mathsf p}_{\mathcal Q}(u)\,\d u.
	\end{align*}
Since ${\mathsf p}_{\mathcal Q}$ is Lipschitz and $Z_{h_Y}\ge c_\phi$, we have
	\begin{align*}
		\abs[0]{\bar a_{\mathcal Q}(y)-{\mathsf p}_{\mathcal Q}(y)}\le\frac1{c_\phi}\int_\cY\phi_{h_Y}(y-u)\abs[0]{{\mathsf p}_{\mathcal Q}(u)-{\mathsf p}_{\mathcal Q}(y)}\,\d u\le C\int_{\R^{d_Y}}\phi_{h_Y}(s)\abs[0]{s}\,\d s\le Ch_Y.
	\end{align*}

Let
	\begin{align*}
		V^{\rm ker}_i(y):=\ind\{X_i\in \mathcal Q\}\phi_{h_Y}(y-Y_i),\qquad T(y):=\frac1{n_{\rm fit}}\sum_{i\in \mathcal D_{\rm fit}}V^{\rm ker}_i(y),
	\end{align*}
so that $\hat a_{\mathcal Q}=n_{\rm fit}T/(n_{\mathcal Q}Z_{h_Y})$ and $\E T=\Pi(\mathcal Q)Z_{h_Y}\bar a_{\mathcal Q}\le{\mathsf p}_{\max}\Pi(\mathcal Q)Z_{h_Y}$. Let $\mathcal G_n$ be a grid in $\cY$ of mesh $n^{-d_Y-2}$. By the boundedness of $\phi$ and ${\mathsf p}$, the bound $\Pi(\mathcal Q)\le Lh_X^{d_X}$ and the choice of the mesh, we have
	\begin{align*}
		0\le V^{\rm ker}_i(y)\le Ch_Y^{-d_Y},\qquad\E V^{\rm ker}_i(y)^2\le Ch_X^{d_X}h_Y^{-d_Y},\qquad\#\mathcal G_n\le Cn^{d_Y(d_Y+2)},
	\end{align*}
and $T$ and $\E T$ are Lipschitz with constant $Ch_Y^{-d_Y-1}\le Cn^{d_Y+1}$. Moreover, $n_{\rm fit}\asymp n$ by \cref{ass:split}, and $1/n\le\log n/(nh_Y^{d_Y})=h_X^{d_X}\nu_{\rm part}^2\le h_X^{d_X}\nu_{\rm part}$. Hence, by Bernstein's inequality at the points of $\mathcal G_n$, a union bound over $\mathcal G_n$ and the at most $Cn$ cubes, and the Lipschitz bounds between grid points, we have, outside an event of probability at most $Cn^{-c}$, for all $\mathcal Q$,
	\begin{align*}
		\sup_{y\in\cY}\abs[0]{T(y)-\E T(y)}\le C\Bigl\{\Bigl(\frac{h_X^{d_X}\log n}{n_{\rm fit}h_Y^{d_Y}}\Bigr)^{1/2}+\frac{\log n}{n_{\rm fit}h_Y^{d_Y}}+\frac1n\Bigr\}\le Ch_X^{d_X}\nu_{\rm part}.
	\end{align*}

Likewise, by Bernstein's inequality for the indicators $\ind\{X_i\in \mathcal Q\}$, we have $\abs[0]{n_{\mathcal Q}/n_{\rm fit}-\Pi(\mathcal Q)}\le Ch_X^{d_X}\nu_{\rm part}$ for all $\mathcal Q$, outside an event of probability at most $Cn^{-c}$. Thus the event
	\begin{align*}
		\Omega_{\rm p}:=\bigcap_{\mathcal Q}\cbr[3]{\sup_{y\in\cY}\abs[0]{T(y)-\E T(y)}\le Ch_X^{d_X}\nu_{\rm part},\ \abs[3]{\frac{n_{\mathcal Q}}{n_{\rm fit}}-\Pi(\mathcal Q)}\le Ch_X^{d_X}\nu_{\rm part}}
	\end{align*}
has probability at least $1-Cn^{-c}$. Since $\Pi(\mathcal Q)\ge f_{\min}h_X^{d_X}$ and $\nu_{\rm part}\le r_n\to0$, we have $Ch_X^{d_X}\nu_{\rm part}\le\Pi(\mathcal Q)/2$ for all sufficiently large $n$. Hence, on $\Omega_{\rm p}$, $n_{\mathcal Q}\ge n_{\rm fit}\Pi(\mathcal Q)/2\ge1$ and
	\begin{align*}
		\abs[0]{\hat a_{\mathcal Q}-\bar a_{\mathcal Q}}&=\frac1{Z_{h_Y}}\abs[3]{\frac{T-\E T}{n_{\mathcal Q}/n_{\rm fit}}+\E T\,\frac{\Pi(\mathcal Q)-n_{\mathcal Q}/n_{\rm fit}}{\Pi(\mathcal Q)n_{\mathcal Q}/n_{\rm fit}}}\\
		&\le\frac{2\abs[0]{T-\E T}}{c_\phi\Pi(\mathcal Q)}+\frac{2{\mathsf p}_{\max}\abs[0]{n_{\mathcal Q}/n_{\rm fit}-\Pi(\mathcal Q)}}{\Pi(\mathcal Q)}\le C\nu_{\rm part} .
	\end{align*}
Together with $\abs[0]{\bar a_{\mathcal Q}-{\mathsf p}_{\mathcal Q}}\le Ch_Y$, this gives $\abs[0]{\hat a_{\mathcal Q}-{\mathsf p}_{\mathcal Q}}\le Cr_n'$. Since $\int_\cY{\mathsf p}_{\mathcal Q}\,\d y=1$, it follows that $\abs[0]{\int_\cY\hat a_{\mathcal Q}\,\d y-1}\le C\Leb(\cY)r_n'\le\frac12$ for all sufficiently large $n$. Hence
	\begin{align*}
		\abs[0]{\hp_{\mathcal Q}-{\mathsf p}_{\mathcal Q}}\le2\abs[0]{\hat a_{\mathcal Q}-{\mathsf p}_{\mathcal Q}}+2{\mathsf p}_{\max}\abs[3]{\int_\cY\hat a_{\mathcal Q}\,\d y-1}\le Cr_n',
	\end{align*}
which is the first bound in \eqref{eq:part-sup}.

\emph{Step 2: \cref{cond:N,cond:V,cond:A}.} By \eqref{eq:part-sup}, for all sufficiently large $n$, the bounds of \cref{cond:A} hold on $\Omega_{\rm p}$ at every $x\in\cX_0^{\delta_0}$, with the identity transformation. We fix $\mathcal Q$ with $n_{\mathcal Q}\ge1$ and $v\in\R$. By the definition of $\hp_{\mathcal Q}$, the level set $\{y\in\cY:\hp_{\mathcal Q}(y)=v\}$ lies in the zero set of $g_\mu:=\phi_{h_Y}*\mu$, where
	\begin{align*}
		\mu:=\sum_{i:X_i\in \mathcal Q}\delta_{Y_i}-\Bigl(vn_{\mathcal Q}\int_\cY\hat a_{\mathcal Q}\,\d y\Bigr)\Leb(\cdot\cap\cY).
	\end{align*}
Since $\mu$ is a finite signed measure with compact support, $g_\mu$ extends to an entire function on $\mathbb C^{d_Y}$, and the Fourier transform of $g_\mu$ is $\e^{-h_Y^2\abs[0]{\xi}^2/2}\hat\mu(\xi)$. Since $\mu$ has atoms, $\mu\ne0$. Hence $\hat\mu\ne0$, and thus $g_\mu\ne0$. Since $g_\mu$ is real-analytic on $\R^{d_Y}$ and not identically zero, its zero set is Lebesgue-null. Consequently, every level set of $\hp_{\mathcal Q}$ is Lebesgue-null, and $\hp$ satisfies \cref{cond:N} on $\Omega_{\rm p}$ and \cref{cond:A} with the sequence $r_n$.

For \cref{cond:V}, we observe that, within each cube, $\hp(\cdot\mid x)$ does not depend on $x$ and its super-level sets form a nested family. Such a family picks out at most $m+1$ different subsets of any $m$ response points, even when some of these points lie at the same level. Therefore
	\begin{align*}
		\mathfrak S_{\mathcal A_n}(m)\le Ch_X^{-d_X}(m+1)\le Cn(m+1)\le(nm)^{C_{\rm gr}}
	\end{align*}
for a fixed $C_{\rm gr}$ and all sufficiently large $n$.

Part (ii) follows from (i) and \cref{thm:local-cov,thm:local}.
\end{proof}

\section{Proofs for \texorpdfstring{\cref{sec:cov-onestage}}{Section \ref{sec:cov-onestage}}}\label{app:onestage-proofs}
\subsection{Regularity of the oracle level}\label{app:oracle-reg}
\begin{lemma}[Lipschitz continuity of the oracle level]\label{lem:tstar-lip}
Under \cref{ass:smooth,ass:margin,ass:design}, there exists a constant $C>0$ such that, for all $x,x'\in\cX_0^{\delta_0}$:
\begin{enumerate}[label=(\roman*)]
\item $|t_\star(x)-t_\star(x')|\le C|x-x'|_\infty$;
\item $\Leb\bigl(\Cstar(x)\triangle\Cstar(x')\bigr)\le C|x-x'|_\infty\,\mathfrak L(|x-x'|_\infty)$.
\end{enumerate}
\end{lemma}

\begin{proof}
\emph{(i)} Let $\eps':=\sup_y|{\mathsf p}(y\mid x)-{\mathsf p}(y\mid x')|$. Since $\cX_0^{\delta_0}$ is convex, the segment from $x$ to $x'$ lies in it. Since $\beta\ge1$, ${\mathsf p}$ is Lipschitz on $\cX_0^{\delta_0}\times\cY$ by \eqref{eq:holder}, and hence $\eps'\le C|x-x'|_\infty$.

Next, $|{\mathsf P}_x(A)-{\mathsf P}_{x'}(A)|\le\Leb(\cY)\eps'$ for every set $A$. If $\Leb(\cY)\eps'>c_0$, then $|x-x'|_\infty$ is bounded below by a positive constant, and (i) holds because $t_\star\le{\mathsf p}_{\max}$. Suppose therefore that $\Leb(\cY)\eps'\le c_0$. By \cref{lem:t-window}, we have ${\mathsf P}_{x'}(\Cstar(x'))=\gamma_\star$. Therefore $\gamma':={\mathsf P}_x(\Cstar(x'))$ satisfies $|\gamma'-\gamma_\star|\le\Leb(\cY)\eps'\le c_0$. Applying \cref{lem:threshold} at $x$ with $\hat a={\mathsf p}(\cdot\mid x')$, $\hC=\Cstar(x')$, $\hat t=t_\star(x')$ and $r=\eps'$, we obtain $|t_\star(x')-t_{\gamma'}(x)|\le\eps'$. Moreover, by \eqref{eq:t-lip}, we have $|t_{\gamma'}(x)-t_\star(x)|\le2\Leb(\cY)\eps'/c_m$. Hence
	\begin{align*}
		|t_\star(x)-t_\star(x')|&\le|t_\star(x')-t_{\gamma'}(x)|+|t_{\gamma'}(x)-t_\star(x)|\\
		&\le\Bigl(1+\frac{2\Leb(\cY)}{c_m}\Bigr)\eps'\le C|x-x'|_\infty.
	\end{align*}

\emph{(ii)} For $y\in\Cstar(x)\setminus\Cstar(x')$ and for $y\in\Cstar(x')\setminus\Cstar(x)$, we have, respectively,
	\begin{align*}
		t_\star(x)&\le{\mathsf p}(y\mid x)\le{\mathsf p}(y\mid x')+\eps'<t_\star(x')+\eps',\\
		t_\star(x)&>{\mathsf p}(y\mid x)\ge t_\star(x')-\eps'.
	\end{align*}
In both cases, by (i),
	\begin{align*}
		|{\mathsf p}(y\mid x)-t_\star(x)|\le\eps'+|t_\star(x)-t_\star(x')|\le C_1|x-x'|_\infty
	\end{align*}
for a constant $C_1>0$. Therefore $\Cstar(x)\triangle\Cstar(x')\subset\mathcal J_x(C_1|x-x'|_\infty)$. If $C_1|x-x'|_\infty\le\varepsilon_0$, then, by \eqref{eq:margin-upper} at $x$ and \eqref{eq:L-props}, $\Leb\bigl(\mathcal J_x(C_1|x-x'|_\infty)\bigr)\le C|x-x'|_\infty\mathfrak L(|x-x'|_\infty)$. Otherwise, $|x-x'|_\infty$ is bounded below by a positive constant, and the same bound holds with a larger constant, because $\Leb(\Cstar(x)\triangle\Cstar(x'))\le2\Leb(\cY)$.
\end{proof}

\subsection{HPD sets of an accurate estimate}\label{app:hpd-sets}
\begin{lemma}[HPD sets of an accurate estimate]\label{lem:hpd-sets}
Suppose that \cref{ass:margin} holds. For fixed $\mathcal D_{\rm fit}$ and $x$, let
	\begin{align*}
		\eps_{\rm sup}:=\sup_y|\hp(y\mid x)-{\mathsf p}(y\mid x)|,
	\end{align*}
and suppose that $\eps_{\rm sup}\le\varepsilon_0/8$, that $\Leb(\cY)\eps_{\rm sup}<\min\{c_m\varepsilon_0/32,c_0/2\}$ and that the level sets $\{y:\hp(y\mid x)=v\}$ with $|v-t_\star(x)|\le\varepsilon_0/2$ are Lebesgue-null. Let $u,u'\in[\gamma_\star-c_0,\gamma_\star+c_0]$, $v_x(u):=\inf\{v:\Lambda_x(v)\le u\}$ and $A_x(u):=\{y:\hp(y\mid x)\ge v_x(u)\}$. Then the following hold.
\begin{enumerate}[label=(\roman*)]
\item The level $v_x(u)$ lies in $(t_\star(x)-\varepsilon_0/2,t_\star(x)+\varepsilon_0/2)$.
\item For every value of the tie-breaking variable,
	\begin{align}\label{eq:A-sets}
		\{y:\hat{\mathsf s}^{\rm hpd}(x,y)\le u\}=A_x(u),\qquad\hP_x(A_x(u))=u.
	\end{align}
Moreover, $F(u,x)={\mathsf P}_x(A_x(u))$, and $R(u,x)$ satisfies \eqref{eq:R-int}. If, in addition, $|u-\gamma_\star|\le c_0/2$, then
	\begin{align*}
		|v_x(u)-t_\star(x)|\le\eps_{\rm sup}+\frac2{c_m}\bigl(|u-\gamma_\star|+\Leb(\cY)\eps_{\rm sup}\bigr).
	\end{align*}
\item If $u'\le u$, then $|R(u,x)-R(u',x)|\le(2\eps_{\rm sup}/\varepsilon_0)(u-u')$.
\end{enumerate}
\end{lemma}

\begin{proof}
\emph{(i)} Let $v_\pm:=t_\star(x)\pm\varepsilon_0/2$. Since $\eps_{\rm sup}\le\varepsilon_0/8$, we have
	\begin{align*}
		\{y:\hp(y\mid x)\ge v_-\}&\supset\{y:{\mathsf p}(y\mid x)\ge t_\star(x)-3\varepsilon_0/8\},\\
		\{y:\hp(y\mid x)\ge v_+\}&\subset\{y:{\mathsf p}(y\mid x)\ge t_\star(x)+3\varepsilon_0/8\}.
	\end{align*}
By \cref{lem:t-window}, we have $|t_{\gamma_\star\pm c_0}(x)-t_\star(x)|\le\varepsilon_0/4$ and $G_x(t_{\gamma_\star\pm c_0}(x))=\gamma_\star\pm c_0$. Hence, by \eqref{eq:G-lower},
	\begin{align*}
		G_x\bigl(t_\star(x)-3\varepsilon_0/8\bigr)\ge\gamma_\star+c_0+\frac{c_m\varepsilon_0}{16},\qquad G_x\bigl(t_\star(x)+3\varepsilon_0/8\bigr)\le\gamma_\star-c_0-\frac{c_m\varepsilon_0}{16}.
	\end{align*}
Since $|\hP_x(A)-{\mathsf P}_x(A)|\le\Leb(\cY)\eps_{\rm sup}<c_m\varepsilon_0/32$ for every set $A\subset\cY$, it follows that
	\begin{align*}
		\overline\Lambda_x(v_-)>\gamma_\star+c_0\ge u,\qquad\overline\Lambda_x(v_+)<\gamma_\star-c_0\le u.
	\end{align*}

By the continuity of $\hP_x$ along monotone sequences of sets, $\Lambda_x$ is nonincreasing and right-continuous, and
	\begin{align*}
		\lim_{v'\uparrow v}\Lambda_x(v')=\overline\Lambda_x(v),\qquad\overline\Lambda_x(v)-\Lambda_x(v)=\hP_x\{y:\hp(y\mid x)=v\}.
	\end{align*}
The right-hand side of the second identity vanishes when the level set is Lebesgue-null. Since the level set at $v_-$ is Lebesgue-null, we have $\Lambda_x(v_-)>u$, and, by right-continuity, $\Lambda_x(v')>u$ for some $v'>v_-$. Hence $v_x(u)\ge v'>v_-$. Since $\overline\Lambda_x(v_+)<u$ is the left limit of $\Lambda_x$ at $v_+$, we have $\Lambda_x(v')<u$ for some $v'<v_+$, and hence $v_x(u)\le v'<v_+$.

\emph{(ii)} We write $v_u:=v_x(u)$. Since $\Lambda_x$ is nonincreasing and right-continuous,
	\begin{align*}
		\Lambda_x(v)\le u\iff v\ge v_u .
	\end{align*}
By (i), the level set at $v_u$ is Lebesgue-null. Hence $\Lambda_x(v_u)=\overline\Lambda_x(v_u)$, and, since $\overline\Lambda_x(v_u)$ is the left limit of $\Lambda_x$ at $v_u$,
	\begin{align*}
		u\le\lim_{v\uparrow v_u}\Lambda_x(v)=\overline\Lambda_x(v_u)=\Lambda_x(v_u)\le u .
	\end{align*}
Therefore $\hP_x(A_x(u))=\overline\Lambda_x(v_u)=u$. Now let $y\in\cY$ and $v:=\hp(y\mid x)$. By \eqref{eq:score}, we have $\Lambda_x(v)\le\hat{\mathsf s}^{\rm hpd}(x,y)\le\overline\Lambda_x(v)$, and hence, by the equivalence above,
	\begin{align*}
		v\ge v_u&\implies\hat{\mathsf s}^{\rm hpd}(x,y)\le\overline\Lambda_x(v_u)=u,\\
		v<v_u&\implies\hat{\mathsf s}^{\rm hpd}(x,y)\ge\Lambda_x(v)>u .
	\end{align*}
This proves \eqref{eq:A-sets}. Consequently, $F(u,x)={\mathsf P}_x(A_x(u))$ and
	\begin{align*}
		R(u,x)={\mathsf P}_x(A_x(u))-\hP_x(A_x(u))=\int_{A_x(u)}({\mathsf p}-\hp)(y\mid x)\,\d y,
	\end{align*}
which is \eqref{eq:R-int}.

For the last claim, let $|u-\gamma_\star|\le c_0/2$ and $\gamma_A:={\mathsf P}_x(A_x(u))$. Since $\hP_x(A_x(u))=u$, we have $|\gamma_A-u|=|{\mathsf P}_x(A_x(u))-\hP_x(A_x(u))|\le\Leb(\cY)\eps_{\rm sup}<c_0/2$, and hence
	\begin{align*}
		|\gamma_A-\gamma_\star|\le|u-\gamma_\star|+\Leb(\cY)\eps_{\rm sup}\le c_0 .
	\end{align*}
By \cref{lem:threshold} with $\hat a=\hp(\cdot\mid x)$, $\hC=A_x(u)$, $\hat t=v_u$ and $r=\eps_{\rm sup}$, we have $|v_u-t_{\gamma_A}(x)|\le\eps_{\rm sup}$. Combining this bound with \eqref{eq:t-lip}, we obtain
	\begin{align*}
		|v_u-t_\star(x)|\le|v_u-t_{\gamma_A}(x)|+|t_{\gamma_A}(x)-t_\star(x)|\le\eps_{\rm sup}+\frac2{c_m}\bigl(|u-\gamma_\star|+\Leb(\cY)\eps_{\rm sup}\bigr).
	\end{align*}

\emph{(iii)} Since $v_x$ is nonincreasing, we have $A_x(u')\subset A_x(u)$. On $A_x(u)\setminus A_x(u')$, by (i) and \cref{lem:t-window}, $\hp(\cdot\mid x)\ge v_x(u)>t_\star(x)-\varepsilon_0/2>\varepsilon_0/2$, and hence $|({\mathsf p}-\hp)(\cdot\mid x)|\le\eps_{\rm sup}\le(2\eps_{\rm sup}/\varepsilon_0)\hp(\cdot\mid x)$. By (ii), we have $\int_{A_x(u)\setminus A_x(u')}\hp(y\mid x)\,\d y=u-u'$. Hence, by \eqref{eq:R-int} and the last bound,
	\begin{align*}
		|R(u,x)-R(u',x)|&=\Bigl|\int_{A_x(u)\setminus A_x(u')}({\mathsf p}-\hp)(y\mid x)\,\d y\Bigr|\\
		&\le\frac{2\eps_{\rm sup}}{\varepsilon_0}\int_{A_x(u)\setminus A_x(u')}\hp(y\mid x)\,\d y=\frac{2\eps_{\rm sup}}{\varepsilon_0}(u-u').
	\end{align*}
\end{proof}

\subsection{Coverage from total variation}\label{app:tv}
Step~1 of the proof of \cref{prop:tv} and the event of probability at least $1-\vartheta$ constructed there apply to a general $\hp$. Only the first paragraph of Step~2 uses the assumption that the level sets of $\hp(\cdot\mid x)$ are Lebesgue-null.

\begin{lemma}[Uniformity of the HPD value under the estimate]\label{lem:hpd-uniform}
Fix $x$ and $\hp$, and let $Y_0\sim\hP_x$ be independent of $U$ at $x$. Then $\hat{\mathsf s}^{\rm hpd}(x,Y_0)$ is uniform on $[0,1]$. Equivalently, for every $\tau\in[0,1]$, the inclusion probability
	\begin{align}\label{eq:pi-u}
		\pi_\tau(y):=\bP\{\hat{\mathsf s}^{\rm hpd}(x,y)\le\tau\},
	\end{align}
where the probability is taken over $U$ at $x$, satisfies $\int\pi_\tau\,\d\hP_x=\tau$.
\end{lemma}

\begin{proof}
Let $T_0:=\hp(Y_0\mid x)$, and let $F_0$ be its distribution function, so that $\Lambda_x(v)=1-F_0(v)$ and $\overline\Lambda_x(v)=1-F_0(v-)$. By \eqref{eq:score}, the score of $Y_0$ is
	\begin{align*}
		\hat{\mathsf s}^{\rm hpd}(x,Y_0)=1-\bigl\{F_0(T_0-)+(1-U)\bigl(F_0(T_0)-F_0(T_0-)\bigr)\bigr\}.
	\end{align*}
The term in braces is the distributional transform of $T_0$ with the variable $1-U$, which is uniform on $[0,1]$ and independent of $T_0$. By Proposition~2.1 of \citet{ruschendorf2009distributional}, this term is uniform on $[0,1]$, and hence so is the score. In particular, $\int\pi_\tau\,\d\hP_x=\bP\{\hat{\mathsf s}^{\rm hpd}(x,Y_0)\le\tau\}=\tau$ for $\tau\in[0,1]$.
\end{proof}

\begin{proposition}[Coverage from total variation]\label{prop:tv}
Let $\hC$ be the output of HPD-split of \cref{sec:setup-oracle} with an estimate $\hp$, and let $\cov(x)$ be its conditional coverage. Let $\vartheta\in(0,1)$, and let the calibration points be i.i.d.\ copies of $(X,Y)$, independent of the fitting data. Suppose that the fitting data and the value of $n_{\rm cal}^\circ$ are fixed and that $n_{\rm cal}^\circ\ge1$ (otherwise $\hat\tau=1$ and the output is $\cY$), and let
	\begin{align}
		&\TV_x:=\TV(\hP_x,{\mathsf P}_x),\qquad\overline{\TV}:=\Pi(\cX_0)^{-1}\int_{\cX_0}\TV_x\,\Pi(\d x),\notag\\
		&\eps_{\rm cal}:=\Bigl\{\frac{\log(2/\vartheta)}{2n_{\rm cal}^\circ}\Bigr\}^{1/2},\qquad\Delta:=\overline{\TV}+\eps_{\rm cal}+\frac2{n_{\rm cal}^\circ}.\label{eq:tv-consts}
	\end{align}
Let $x$ be a point at which the level sets of $\hp(\cdot\mid x)$ are Lebesgue-null. Then, with probability at least $1-\vartheta$ over the calibration sample,
	\begin{align*}
		|\cov(x)-\gamma_\star|\le\TV_x+\Delta .
	\end{align*}
\end{proposition}

\begin{proof}[Proof of \cref{prop:tv}]
\emph{Step 1: the calibrated threshold.} We write $m:=n_{\rm cal}^\circ$ and $k:=\lceil\gamma_\star(m+1)\rceil$. Let $H$ be the distribution function of a calibration score, conditionally on the fitting data and on $m$. The covariate of a calibration point has the distribution of $X$ given $X\in\cX_0$, and its response, given the covariate $x$, has distribution ${\mathsf P}_x$. Hence, for $\tau\in[0,1]$, with $\pi_\tau$ computed at $x$,
	\begin{align*}
		H(\tau)=\Pi(\cX_0)^{-1}\int_{\cX_0}\int\pi_\tau\,\d{\mathsf P}_x\,\Pi(\d x).
	\end{align*}
Since $0\le\pi_\tau\le1$ and $\int\pi_\tau\,\d\hP_x=\tau$ by \cref{lem:hpd-uniform}, the definition of $\TV_x$ gives $|\int\pi_\tau\,\d{\mathsf P}_x-\tau|\le\TV_x$ for every $x$. Averaging this bound, and then letting $\tau'\uparrow\tau$ in the resulting bound at $\tau'$, we obtain
	\begin{align*}
		|H(\tau)-\tau|\le\overline{\TV},\qquad|H(\tau-)-\tau|\le\overline{\TV}\qquad(\tau\in[0,1]).
	\end{align*}

Since each calibration score uses its own covariate and tie-breaking variable, the calibration scores are i.i.d.\ with distribution function $H$, conditionally on the fitting data and on $m$. The Dvoretzky--Kiefer--Wolfowitz inequality \citep{dvoretzky1956asymptotic} with Massart's constant also holds for discontinuous distribution functions \citep[Corollary~1 and Comment~2(iii)]{massart1990tight}. By this inequality, with probability at least $1-\vartheta$, the empirical distribution function $\widehat H$ of the calibration scores satisfies
	\begin{align*}
		\sup_\tau|\widehat H(\tau)-H(\tau)|\le\eps_{\rm cal},\qquad\sup_\tau|\widehat H(\tau-)-H(\tau-)|\le\eps_{\rm cal}.
	\end{align*}
We work on this event.

In the first case, where $k>m$, the output is $\cY$, and hence $\cov(x)=1$. Since $\gamma_\star(m+1)>m$ in this case, we have $|\cov(x)-\gamma_\star|=1-\gamma_\star<1/(m+1)\le\Delta$. In the second case, where $k\le m$, the threshold $\hat\tau$ is the $k$-th smallest calibration score. Since $\gamma_\star(m+1)\le k<\gamma_\star(m+1)+1$, we have
	\begin{align*}
		\widehat H(\hat\tau)\ge\frac km\ge\gamma_\star,\qquad\widehat H(\hat\tau-)\le\frac{k-1}m\le\gamma_\star+\frac1m .
	\end{align*}
Combining these bounds with the bounds on $H$ and $\widehat H$ above, we obtain
	\begin{align*}
		\hat\tau&\ge H(\hat\tau)-\overline{\TV}\ge\widehat H(\hat\tau)-\eps_{\rm cal}-\overline{\TV}\ge\gamma_\star-\overline{\TV}-\eps_{\rm cal},\\
		\hat\tau&\le H(\hat\tau-)+\overline{\TV}\le\widehat H(\hat\tau-)+\eps_{\rm cal}+\overline{\TV}\le\gamma_\star+\overline{\TV}+\eps_{\rm cal}+\frac1m,
	\end{align*}
and hence $|\hat\tau-\gamma_\star|\le\overline{\TV}+\eps_{\rm cal}+2/m=\Delta$.

\emph{Step 2: the coverage.} We work in the second case of Step~1 and write $\Lambda:=\Lambda_x$ and $\overline\Lambda:=\overline\Lambda_x$. Let $T:=\hp(\widetilde Y\mid x)$ with $\widetilde Y\sim\hP_x$, so that $\Lambda(v)=\hP_x(T>v)$ and $\overline\Lambda(v)=\hP_x(T\ge v)$. By \eqref{eq:score}, we have $\Lambda(\hp(y\mid x))\le\hat{\mathsf s}^{\rm hpd}(x,y)\le\overline\Lambda(\hp(y\mid x))$ for every $y$ and every value of $U$. Since $\hC(x)=\{y:\hat{\mathsf s}^{\rm hpd}(x,y)\le\hat\tau\}$, it follows that
	\begin{align*}
		\bigl\{y:\overline\Lambda(\hp(y\mid x))\le\hat\tau\bigr\}\subset\hC(x)\subset\bigl\{y:\Lambda(\hp(y\mid x))\le\hat\tau\bigr\}.
	\end{align*}
Since the level sets of $\hp(\cdot\mid x)$ are Lebesgue-null, we have, for every $v$,
	\begin{align*}
		\hP_x(T=v)=\int_{\{y:\hp(y\mid x)=v\}}\hp(y\mid x)\,\d y=v\,\Leb\{y:\hp(y\mid x)=v\}=0.
	\end{align*}
Hence $\Lambda=\overline\Lambda$, and, by \cref{lem:hpd-uniform}, the variable $\Lambda(T)$ is uniform on $[0,1]$. Thus both sets in the inclusions have $\hP_x$-mass $\hat\tau$, and $\hP_x(\hC(x))=\hat\tau$. By the definition of $\TV_x$, we have $|\cov(x)-\hP_x(\hC(x))|\le\TV_x$. Combining this bound with Step~1, we obtain $|\cov(x)-\gamma_\star|\le\TV_x+|\hat\tau-\gamma_\star|\le\TV_x+\Delta$.

For a general $\hp$, the variable $T$ can have atoms. Let $U'\sim U[0,1]$ be independent of $T$. By \cref{lem:hpd-uniform}, the variable $S:=\Lambda(T)+U'\{\overline\Lambda(T)-\Lambda(T)\}$ is uniform on $[0,1]$. Since $\Lambda(T)\le S\le\overline\Lambda(T)$, it follows that
	\begin{align*}
		\hP_x\cbr{\overline\Lambda(T)\le\hat\tau}\le\hat\tau\le\hP_x\cbr{\Lambda(T)\le\hat\tau},
	\end{align*}
and the outer terms differ by $g:=\hP_x\{\Lambda(T)\le\hat\tau<\overline\Lambda(T)\}$. Since $\overline\Lambda(v')\le\Lambda(v)$ for $v<v'$, at most one level $v$ satisfies $\Lambda(v)\le\hat\tau<\overline\Lambda(v)$, and $g=\overline\Lambda(v)-\Lambda(v)$ for this $v$, or $g=0$ if there is none. Since $|\hat\tau-\gamma_\star|\le\Delta$, such a level satisfies $\Lambda(v)\le\gamma_\star+\Delta$ and $\overline\Lambda(v)>\gamma_\star-\Delta$, and hence $g\le\eps^{\rm at}_x$, with $\eps^{\rm at}_x$ as in \cref{rem:tv-atoms}. By the inclusions above, we have $|\hP_x(\hC(x))-\hat\tau|\le g$, and the argument of the previous paragraph gives the bound of \cref{rem:tv-atoms}.
\end{proof}

\begin{remark}[Estimates with flat pieces]\label{rem:tv-atoms}
For estimates whose level sets may have positive Lebesgue measure, such as histograms and tree-based densities, let
	\begin{align*}
		\eps^{\rm at}_x:=\sup\cbr{\overline\Lambda_x(v)-\Lambda_x(v):\ \Lambda_x(v)\le\gamma_\star+\Delta,\ \overline\Lambda_x(v)\ge\gamma_\star-\Delta}
	\end{align*}
be the largest $\hP_x$-mass of a level set $\{y:\hp(y\mid x)=v\}$ for which the interval $[\Lambda_x(v),\overline\Lambda_x(v)]$ meets $[\gamma_\star-\Delta,\gamma_\star+\Delta]$. By Step~2 of the proof of \cref{prop:tv}, on the event of probability at least $1-\vartheta$ constructed in that proof, we have $|\cov(x)-\gamma_\star|\le\TV_x+\Delta+\eps^{\rm at}_x$.
\end{remark}

\subsection{One-stage calibration}\label{app:onestage-cal}
\begin{lemma}[One-stage calibration]\label{lem:onestage-cal}
Suppose that \cref{ass:smooth,ass:design,ass:split} hold and that \cref{ass:margin} holds with $\mathfrak L\equiv1$. Let $\hp$ be an estimate fitted on $\mathcal D_{\rm fit}$, let $\hC$ be the output of one-stage calibration with this $\hp$ and exact HPD values, calibrated on the $n_{\rm cal}^\circ$ points of $\mathcal D_{\rm cal}$ in $\cX_0$, let $\hat\tau$ be its threshold and $\cov(x)$ its conditional coverage. Let $r_n\to0$ be a deterministic sequence with $r_n\ge(\log n/n)^{1/2}$, and let $\Omega_r\in\sigma(\mathcal D_{\rm fit})$ be an event on which $\sup_{\cX_0\times\cY}|\hp-{\mathsf p}|\le r_n$. Then for every $c>0$ there exists a constant $C>0$ such that, for all sufficiently large $n$, there exists an event $\Omega_{\rm thr}\subset\Omega_r$ with $\bP(\Omega_{\rm thr})\ge\bP(\Omega_r)-Cn^{-c}$ on which the following hold.
\begin{enumerate}[label=(\roman*)]
\item $|\hat\tau-\gamma_\star|\le Cr_n$.
\item Suppose, in addition, that on $\Omega_r$, for all $x\in\cX_0$ and all $v$ with $|v-t_\star(x)|\le\varepsilon_0/2$, the level set $\{y:\hp(y\mid x)=v\}$ is Lebesgue-null, and let $\bar R(\gamma_\star):=\Pi(\cX_0)^{-1}\int_{\cX_0}R(\gamma_\star,x)\,\Pi(\d x)$. Then, for all $x\in\cX_0$, $\hC(x)=A_x(\hat\tau)$ for every value of the tie-breaking variable at $x$, and
	\begin{equation}\label{eq:onestage-centered}
		\cov(x)-\gamma_\star=R(\gamma_\star,x)-\bar R(\gamma_\star)+O\Bigl(\Bigl(\frac{\log n}n\Bigr)^{1/2}+r_n|\bar R(\gamma_\star)|\Bigr).
	\end{equation}
\end{enumerate}
\end{lemma}

\begin{proof}
We write $m:=n_{\rm cal}^\circ$ and $k:=\lceil\gamma_\star(m+1)\rceil$. Let $H$ be the distribution function of the HPD value of a calibration point in $\cX_0$, given $\mathcal D_{\rm fit}$ and $m$, and let $\widehat H$ be the empirical distribution function of the $m$ calibration scores, as in \cref{lem:cal-events}. By \cref{lem:cal-events} and \cref{ass:split},
	\begin{align*}
		\bP(\Omega_{\rm thr})\ge\bP(\Omega_r)-Cn^{-c},\qquad\Omega_{\rm thr}:=\Omega_r\cap\Omega^\circ_{\rm cal}\cap\cbr[2]{\sup_u|\widehat H(u)-H(u)|\le\eps_{{\rm cal},n}},
	\end{align*}
with $\Omega^\circ_{\rm cal}$ and $\eps_{{\rm cal},n}$ as in \cref{lem:cal-events}. We work on $\Omega_{\rm thr}$. Then $m\ge\pi_0n_{\rm cal}/2$, which is of order $n$, and hence $\eps_{{\rm cal},n}\le C(\log n/n)^{1/2}\le Cr_n$.

\emph{Step 1: the threshold.} Let $x\in\cX_0$, and let $\pi_u$ be the inclusion probability $\pi_\tau$ of \eqref{eq:pi-u} with $\tau=u$. By \cref{lem:hpd-uniform}, we have $\int\pi_u\,\d\hP_x=u$ for $u\in[0,1]$, and hence $|F(u,x)-u|=|\int\pi_u({\mathsf p}-\hp)(y\mid x)\,\d y|\le\Leb(\cY)r_n$. Averaging over $X$, we obtain $|H(u)-u|\le\Leb(\cY)r_n$ for $u\in[0,1]$. Let $\delta^{\rm thr}:=\Leb(\cY)r_n+\eps_{{\rm cal},n}+2/m$. For all sufficiently large $n$, we have $\delta^{\rm thr}\le\min\{\gamma_\star,1-\gamma_\star\}/2$, and hence $\gamma_\star\pm\delta^{\rm thr}\in[0,1]$. Since $\gamma_\star\le k/m\le\gamma_\star+2/m$, the bound on $H$ and \eqref{eq:ts-dkw} give
	\begin{align*}
		\widehat H(\gamma_\star-\delta^{\rm thr})\le\gamma_\star-\frac2m<\frac km,\qquad\widehat H(\gamma_\star+\delta^{\rm thr})\ge\gamma_\star+\frac2m\ge\frac km .
	\end{align*}
Hence $k\le m$ and $|\hat\tau-\gamma_\star|\le \delta^{\rm thr}\le Cr_n$.

\emph{Step 2: the expansion.} Suppose that the level sets are Lebesgue-null as in (ii). For all sufficiently large $n$, the conditions of \cref{lem:hpd-sets} hold with $\eps_{\rm sup}\le r_n$, and $\delta^{\rm thr}\le c_0/2$. Hence, by Step 1 and \cref{lem:hpd-sets}(ii), $\hC(x)=A_x(\hat\tau)$ for every value of the tie-breaking variable, and $\cov(x)=F(\hat\tau,x)$. By \cref{lem:hpd-sets}(iii), we have $|R(u,x)-R(u',x)|\le L_R|u-u'|$ for $u,u'\in[\gamma_\star-c_0,\gamma_\star+c_0]$, with $L_R:=2r_n/\varepsilon_0$. Averaging over $X$, we obtain $H(u)=u+\bar R(u)$, where $\bar R(u):=\Pi(\cX_0)^{-1}\int_{\cX_0}R(u,x)\,\Pi(\d x)$. Since $|\bar R(u)-\bar R(u')|\le L_R|u-u'|$ on the same interval, $H$ is continuous on $[\gamma_\star-c_0,\gamma_\star+c_0]$.

Since $\widehat H(\hat\tau)\ge k/m\ge\gamma_\star$, $\widehat H(\hat\tau-)\le(k-1)/m\le\gamma_\star+1/m$, and \eqref{eq:ts-dkw} also holds for the left limits, the continuity of $H$ at $\hat\tau$ gives $|H(\hat\tau)-\gamma_\star|\le\delta^{\rm cal}:=\eps_{{\rm cal},n}+1/m$. Since $H(\hat\tau)=\hat\tau+\bar R(\hat\tau)$, we obtain
	\begin{align*}
		|\hat\tau-\gamma_\star+\bar R(\gamma_\star)|\le|H(\hat\tau)-\gamma_\star|+|\bar R(\hat\tau)-\bar R(\gamma_\star)|\le\delta^{\rm cal}+L_R|\hat\tau-\gamma_\star|,
	\end{align*}
and therefore $(1-L_R)|\hat\tau-\gamma_\star|\le|\bar R(\gamma_\star)|+\delta^{\rm cal}$. Since $\cov(x)=\hat\tau+R(\hat\tau,x)$, it follows that
	\begin{align*}
		\bigl|\cov(x)-\gamma_\star-\{R(\gamma_\star,x)-\bar R(\gamma_\star)\}\bigr|&\le|\hat\tau-\gamma_\star+\bar R(\gamma_\star)|+|R(\hat\tau,x)-R(\gamma_\star,x)|\\
		&\le\delta^{\rm cal}+\frac{2L_R(|\bar R(\gamma_\star)|+\delta^{\rm cal})}{1-L_R}.
	\end{align*}
Since $\delta^{\rm cal}\le C(\log n/n)^{1/2}$ and $L_R\le Cr_n\le\frac12$ for all sufficiently large $n$, this proves \eqref{eq:onestage-centered}.
\end{proof}

\subsection{Proof of \texorpdfstring{\cref{thm:onestage}}{Theorem \ref{thm:onestage}}}\label{app:onestage-ub}
\begin{proof}
Let $\Omega_{\rm thr}$ be the event of \cref{lem:onestage-cal} with the given $\Omega_r$, so that $\bP(\Omega_{\rm thr})\ge\bP(\Omega_r)-Cn^{-c}$. We work on $\Omega_{\rm thr}$ and let $n$ be sufficiently large.

\emph{Step 1: the residual at the oracle set.} Since $r_n\to0$, the conditions of \cref{lem:hpd-sets} hold at every $x\in\cX_0$ with $\eps_{\rm sup}\le r_n$. Let $x\in\cX_0$ and $A_x:=A_x(\gamma_\star)$. By \cref{lem:hpd-sets}(ii), we have $R(\gamma_\star,x)=\int_{A_x}({\mathsf p}-\hp)(y\mid x)\,\d y$ and $|v_x(\gamma_\star)-t_\star(x)|\le Cr_n$. Since $A_x=\{y:\hp(y\mid x)\ge v_x(\gamma_\star)\}$ and $|\hp-{\mathsf p}|\le r_n$, we have $A_x\triangle\Cstar(x)\subset\mathcal J_x(Cr_n)$, and $\Leb(\mathcal J_x(Cr_n))\le Cr_n$ by \eqref{eq:margin-upper} with $\mathfrak L\equiv1$. Hence
	\begin{align*}
		\abs[1]{R(\gamma_\star,x)-R^\star(x)}
		&=\abs[3]{\int_{A_x}({\mathsf p}-\hp)(y\mid x)\,\d y-\int_{\Cstar(x)}({\mathsf p}-\hp)(y\mid x)\,\d y}\\
		&\le r_n\Leb\bigl(A_x\triangle\Cstar(x)\bigr)\le Cr_n^2 .
	\end{align*}
Averaging over $X$ given $X\in\cX_0$, we obtain $|\bar R(\gamma_\star)-\bar R^\star|\le Cr_n^2$, with $\bar R(\gamma_\star)$ as in \cref{lem:onestage-cal}(ii).

\emph{Step 2: the bounds.} By \eqref{eq:onestage-centered} and Step 1, for all $x\in\cX_0$,
	\begin{align*}
		\abs[1]{\cov^{\rm one}(x)-\gamma_\star-\{R^\star(x)-\bar R^\star\}}&\le\abs[1]{R(\gamma_\star,x)-R^\star(x)}+|\bar R(\gamma_\star)-\bar R^\star|\\
		&\quad+C\Bigl\{\Bigl(\frac{\log n}n\Bigr)^{1/2}+r_n|\bar R(\gamma_\star)|\Bigr\}.
	\end{align*}
Since $|R^\star|\le\Leb(\cY)r_n$, Step 1 gives $|\bar R(\gamma_\star)|\le Cr_n$, and the first bound follows. The first bound and $|R^\star|,|\bar R^\star|\le\Leb(\cY)r_n$ give the second. For all sufficiently large $n$, we have $Cr_n\le c_0$ and $r_n\le\varepsilon_0/8$. By \cref{lem:onestage-cal}(ii), $\hC^{\rm one}(x)=A_x(\hat\tau)$ is a super-level set of $\hp(\cdot\mid x)$, and its coverage lies in the window. Hence \cref{prop:supshape} with $r=r_n$ and $\mathfrak L\equiv1$ gives $\cS(\hC^{\rm one}(x))\le Cr_n^2$ and $\cE(\hC^{\rm one}(x))\le C\{r_n^2+|\cov^{\rm one}(x)-\gamma_\star|^2\}\le Cr_n^2$.
\end{proof}

\subsection{Proof of \texorpdfstring{\cref{thm:onestage-lower}}{Theorem \ref{thm:onestage-lower}}}\label{app:onestage-lower}
\begin{proof}
Let $\Omega$ be the event of \cref{thm:onestage} with $c=1$, so that $1-\bP(\Omega)\le\delta_n+Cn^{-1}$. Since $\hp(\cdot\mid x)$ and ${\mathsf p}(\cdot\mid x)$ are probability densities, we have $|R^\star|\le2$ and $|\bar R^\star|\le2$. Moreover, $|\cov^{\rm one}-\gamma_\star|\le1$. Hence, by the first bound of \cref{thm:onestage}, for every $x\in\cX_0$,
	\begin{align*}
		\abs[1]{\E\sbr{\cov^{\rm one}(x)-\gamma_\star}-\E\sbr[0]{R^\star(x)}+\E\sbr[0]{\bar R^\star}}\le C\Bigl\{r_n^2+\Bigl(\frac{\log n}n\Bigr)^{1/2}+\delta_n\Bigr\}=o(a_n).
	\end{align*}
By Fubini's theorem and the expansion of $\E\sbr[0]{R^\star(x)}$, which is uniform in $x$, we have $\E\sbr[0]{\bar R^\star}=a_n\bar{\mathfrak g}+o(a_n)$, where $\bar{\mathfrak g}:=\Pi(\cX_0)^{-1}\int_{\cX_0}\mathfrak g\,\d\Pi$. Therefore, uniformly in $x\in\cX_0$,
	\begin{align*}
		\E\sbr{\cov^{\rm one}(x)-\gamma_\star}=a_n\bigl\{\mathfrak g(x)-\bar{\mathfrak g}\bigr\}+o(a_n).
	\end{align*}
Since $|\E Z|\le\E|Z|$ for every integrable $Z$, Fubini's theorem gives
	\begin{align*}
		\E\sbr{\int_{\cX_0}|\cov^{\rm one}(x)-\gamma_\star|\,\d x}\ge\int_{\cX_0}\abs[1]{\E\sbr{\cov^{\rm one}(x)-\gamma_\star}}\,\d x\ge a_n\int_{\cX_0}|\mathfrak g-\bar{\mathfrak g}|\,\d x-o(a_n).
	\end{align*}
Since $\mathfrak g$ is not almost everywhere constant on $\cX_0$, we have $c_{\mathfrak g}:=\int_{\cX_0}|\mathfrak g-\bar{\mathfrak g}|\,\d x>0$, and the right-hand side is at least $c_{\mathfrak g}a_n/2$ for all sufficiently large $n$.
\end{proof}

\subsection{Proof of \texorpdfstring{\cref{cor:lp-onestage}}{Corollary \ref{cor:lp-onestage}}}\label{app:fixed-law}
The estimate of \cref{cor:lp-onestage} is $\hp^{\rm lp}_b:=\hat a^{\rm lp}_b/\widehat Z$ if $\widehat Z:=\int_\cY\hat a^{\rm lp}_b\,\d y\ge c_Z$, and $\hp^{\rm lp}_b:={\mathsf q}_0$ otherwise, with $c_Z$ as in \cref{app:lp-estimator}. The regularity conditions of \cref{cor:lp-onestage} are the following.
\begin{enumerate}[label=(R\arabic*)]
\item\label{it:R1} The univariate kernel $\mathcal K_1$ is a polynomial on $[-1,1]$, and the integer $J$ of the perturbation $\mathrm{pol}_J$ is even and larger than the constant $D_0$ of \cref{sec:cov-twostage-alg}.
\item\label{it:R2} The map $(x,y)\mapsto{\mathsf p}(y\mid x)$ is $\beta$ times continuously differentiable on the compact set $\cX_0^{2\delta_0}\times\cY$.
\item\label{it:R3} There exists a constant $C_{\rm low}>0$ such that $\Leb\{y\in\cY:{\mathsf p}(y\mid x)\le s\}\le C_{\rm low}s$ for all $x\in\cX_0$ and $s>0$.
\end{enumerate}
Condition \ref{it:R1} yields Lebesgue-null level sets through the perturbation, \ref{it:R2} yields the leading term of the bias, and \ref{it:R3} controls the clipping of the fit at low densities.

\begin{proof}
Let $\mathfrak g(x):=\int_{\Cstar(x)}B_\beta(x,y)\,\d y$. We verify the conditions of \cref{thm:onestage,thm:onestage-lower} with $a_n=b^\beta$ and $-\mathfrak g$ in place of $\mathfrak g$. The mean of $R^\star(x)$ is minus the bias integrated over $\Cstar(x)$, since the linear part of the noise has mean zero and the rest is of second order. We use the notation of \cref{app:lp-estimator} at the bandwidth $b$, write $\hp^{\rm lp}$ for $\hp^{\rm lp}_b$, and let $\Omega_{\rm vc}$ be the event of \cref{lem:vc} with $c:=d_X+2$.

\emph{Rates.} Since $b\ge h$ and $h\asymp(\log n/n)^{1/(2\beta+d)}$, we have $nb^d/\log n\to\infty$ and $\nu_b\le\nu_h\le C\zeta_n\le Cb^\beta$. Hence $\eps_b\le Cb^\beta$ and $\eps_b^2=o(b^\beta)$. Moreover, $(\log n/n)^{1/2}=o(\zeta_n)=o(b^\beta)$.

\emph{The conditions of \cref{thm:onestage}.} By \cref{lem:supnorm-b} and its proof, for all sufficiently large $n$, on $\Omega_{\rm vc}$ we have $|\widetilde Z-1|\le C\eps_b$ on $\cX_0^{\delta_0}$ and $|\hp-{\mathsf p}|\le C\eps_b$ on $\cX_0^{\delta_0}\times\cY$. As in the step \emph{\cref{cond:A}} of the proof of \cref{prop:lp-conditions}, we have $|\widehat Z-\widetilde Z|\le Cn^{-2}$ and $|\hat a/\widehat Z-\tilde a/\widetilde Z|\le Cn^{-2}$ on $\cX_0\times\cY$. Hence, for all sufficiently large $n$, we have $\widehat Z\ge c_Z$ and $\hp^{\rm lp}=\hat a/\widehat Z$ on $\Omega_{\rm vc}$, and
	\begin{align}\label{eq:lp-one-close}
		\sup_{\cX_0\times\cY}|\hp^{\rm lp}-\hp|\le Cn^{-2},\qquad\sup_{\cX_0\times\cY}|\hp^{\rm lp}-{\mathsf p}|\le C_{\rm sup}\eps_b=:r_n
	\end{align}
for a constant $C_{\rm sup}\ge1$. By \ref{it:R1}, the step \emph{Null level sets} of the same proof, which does not use the value of the bandwidth, applies. Hence the level sets of $\hat a(\cdot\mid x)$, and therefore those of $\hp^{\rm lp}(\cdot\mid x)$, are Lebesgue-null for every $x$ on $\Omega_{\rm vc}$. The sequence $r_n$ is deterministic and satisfies $r_n\to0$ and $r_n\ge\nu_b\ge(\log n/n)^{1/2}$. Thus the conditions of \cref{thm:onestage} hold with $\Omega_r:=\Omega_{\rm vc}$ and $\delta_n:=Cn^{-d_X-2}$, and, by the step \emph{Rates}, $r_n^2+(\log n/n)^{1/2}+\delta_n=o(b^\beta)$.

\emph{The leading bias of the fit.} Since $\beta$ is an integer, $\ell=\beta-1$. Let $D(s)$ be the vector of the monomials $s^{\boldsymbol k}$, $|\boldsymbol k|\le\ell$, with the multi-index notation of \cref{sec:setup-assumptions}, so that $D(z';z)=D((z'-z)/b)$. We call a point $z\in\cX_0^{\delta_0}\times\cY$ interior if its window is not truncated, that is, $z+bs\in\cX_0^{2\delta_0}\times\cY$ for all $s\in[-1,1]^d$.

Let $\Psi$ be as in \eqref{eq:Psi}, so that $\eta=\Psi({\mathsf p})$. Since ${\mathsf p}$ is $\beta$ times continuously differentiable on the compact set $\cX_0^{2\delta_0}\times\cY$ and $\Psi$ is smooth, $\eta$ is also $\beta$ times continuously differentiable there, with uniformly continuous derivatives of order $\beta$. Let $\mathcal T_z$ be the Taylor polynomial of degree $\ell$ of $\eta$ at $z$, as in the proof of \cref{lem:bias}. By Taylor's theorem with integral remainder and the uniform continuity of the derivatives of order $\beta$ of $\eta$,
	\begin{equation*}
		\eta(z')-\mathcal T_z(z')=\sum_{|\boldsymbol k|=\beta}\frac{\partial^{\boldsymbol k}\eta(z)}{\boldsymbol k!}(z'-z)^{\boldsymbol k}+o(b^\beta)\qquad\text{uniformly over }|z'-z|_\infty\le b .
	\end{equation*}
We insert this expansion into the representation
	\begin{align*}
		\bar\eta_b(z)-\eta(z)=\mathbf e_1^\top\Sigma_b(z)^{-1}\E\sbr[1]{(n_{\rm fit}b^d)^{-1}\sum_{\rm rows}W\,D\,\{\eta(z')-\mathcal T_z(z')\}}
	\end{align*}
of the proof of \cref{lem:bias}. As in that proof, the remainder $o(b^\beta)$ contributes $o(b^\beta)$. For the main terms, recall that a block has one row with density $f_X{\mathsf p}$ and $K$ rows with density $f_X{\mathsf q}_0$. Hence, for a function $v$ of a row and the rows $z_{ij}$ of block $i$ as in \cref{lem:linear}, we have $\E\sbr[0]{\sum_{j=0}^Kv(z_{ij})}=\int fv\,\d z'$. With $z'=z+bs$ and $(z'-z)^{\boldsymbol k}=b^\beta s^{\boldsymbol k}$, at an interior $z$,
	\begin{align*}
		\E\sbr[2]{(n_{\rm fit}b^d)^{-1}\sum_{\rm rows}W\,D\,(z'-z)^{\boldsymbol k}}&=b^\beta\!\int\mathcal K(s)D(s)s^{\boldsymbol k}f(z+bs)\,\d s\\
		&=b^\beta\Bigl\{f(z)\int\mathcal K(s)D(s)s^{\boldsymbol k}\,\d s+o(1)\Bigr\},\\
		\Sigma_b(z)&=\int\mathcal K(s)D(s)D(s)^\top f(z+bs)\,\d s=f(z)\Sigma_{\mathcal K}+o(1),
	\end{align*}
uniformly over interior $z$, because $f$ is uniformly continuous. Here
	\begin{align}\label{eq:Sigma-K}
		\Sigma_{\mathcal K}:=\int\mathcal KDD^\top\,\d s
	\end{align}
is invertible, and $f\ge f_{\min}K{\mathsf q}_0>0$. Hence the factor $f(z)$ cancels between $\Sigma_b(z)^{-1}$ and the expectation, and, uniformly over interior $z$,
	\begin{align*}
		\bar\eta_b-\eta=b^\beta \mathfrak b_\beta+o(b^\beta),\qquad \mathfrak b_\beta(z):=\sum_{|\boldsymbol k|=\beta}\frac{\partial^{\boldsymbol k}\eta(z)}{\boldsymbol k!}\,\mathbf e_1^\top\Sigma_{\mathcal K}^{-1}\!\int \mathcal K(s)D(s)s^{\boldsymbol k}\,\d s .
	\end{align*}

\emph{The bias after $\psi$ and the normalization.} By the choice of $h_0$ in \cref{app:lp-estimator}, for $b\le h_0$, which holds for all sufficiently large $n$, $\bar\eta_b$ takes values in $[-\frac12,\bar\eta'']$, on which $\psi$ is smooth. Let $\bar a$ be as in \eqref{eq:a-bar}, and recall that $K{\mathsf q}_0\psi(\eta)={\mathsf p}$. By a first-order expansion of $\psi$ and \cref{lem:bias}, we have
	\begin{align*}
		\bar a-{\mathsf p}=K{\mathsf q}_0\psi'(\eta)(\bar\eta_b-\eta)+O(b^{2\beta}).
	\end{align*}
Hence $\bar a-{\mathsf p}=b^\beta\mathfrak a_\beta+o(b^\beta)$ at interior points, where
	\begin{align}\label{eq:a-beta}
		\mathfrak a_\beta:=K{\mathsf q}_0\psi'(\eta)\mathfrak b_\beta,
	\end{align}
and $\bar a-{\mathsf p}=O(b^\beta)$ everywhere. For $x\in\cX_0$, the non-interior points $(x,y)$ form a layer of width $O(b)$ along $\partial\cY$, on which $\bar a-{\mathsf p}-b^\beta\mathfrak a_\beta=O(b^\beta)$. Hence $\int_\cY|\bar a-{\mathsf p}-b^\beta\mathfrak a_\beta|\,\d y=o(b^\beta)$. Since $\int{\mathsf p}(y\mid x)\,\d y=1$, the normalizing constant $\bar Z$ and $\bar{\mathsf p}$ satisfy
	\begin{align*}
		\bar Z&=\int\bar a\,\d y=1+b^\beta\int\mathfrak a_\beta\,\d y+o(b^\beta),\\
		\bar{\mathsf p}-{\mathsf p}&=\frac{\bar a}{\bar Z}-{\mathsf p}=\frac{(\bar a-{\mathsf p})-{\mathsf p}(\bar Z-1)}{\bar Z}.
	\end{align*}
Expanding $1/\bar Z=1+O(b^\beta)$, we obtain
	\begin{equation}\label{eq:leading-bias}
		\bar {\mathsf p}-{\mathsf p}=b^\beta B_\beta+o(b^\beta)\ \text{in }L^1(\d y),\qquad B_\beta:=\mathfrak a_\beta-{\mathsf p}\int\mathfrak a_\beta\,\d y ,
	\end{equation}
uniformly in $x\in\cX_0$.

\emph{The mean of the noise.} On $\Omega_{\rm vc}$, we have $\hp=\tilde a/\widetilde Z$, and we write $N=\hp-\bar{\mathsf p}$ as in \cref{app:lp-estimator}. Let $\cY^{\rm low}_x:=\{y:{\mathsf p}(y\mid x)\le C_{\rm clip}\eps_b\}$, with $C_{\rm clip}$ as in \cref{lem:supnorm-b}. With $\xi$ and $N_\Sigma$ as in \cref{lem:linear}, so that $\heta^{\rm raw}-\bar\eta_b=\xi+N_\Sigma$, we expand $\tilde a/\widetilde Z$ around $(\bar a,\bar Z)$ of \eqref{eq:a-bar} and $\psi$ around $\bar\eta_b$, and we write
	\begin{align*}
		N=\frac{\tilde a-\bar a}{\bar Z}-\frac{\bar a}{\bar Z^2}\int(\tilde a-\bar a)\,\d y'+N_2,\qquad\tilde a-\bar a=K{\mathsf q}_0\psi'(\bar\eta_b)(\xi+N_\Sigma)+N_3.
	\end{align*}
By the proof of \cref{lem:supnorm-b}, $|\tilde a-\bar a|$ and $|\widetilde Z-\bar Z|$ are at most $C\eps_b$, and $\widetilde Z,\bar Z\ge\frac12$. Hence $|N_2|\le C\eps_b^2$. On $\Omega_{\rm vc}$, by \cref{lem:linear}, we have $|N_\Sigma|\le C\nu_b^2$. Off $\cY^{\rm low}_x$, the clipping is inactive by \cref{lem:supnorm-b}(ii) and $\psi''$ is bounded on $[-\frac12,\bar\eta']$, and hence $|N_3|\le C\nu_b^2$. On $\cY^{\rm low}_x$, the term $N_3$ consists of the correction from clipping $\heta^{\rm raw}$ at zero and the Taylor remainder of $\psi$, and $|N_3|\le C\eps_b$. By \ref{it:R3}, we have $\Leb(\cY^{\rm low}_x)\le C_{\rm low}C_{\rm clip}\eps_b$, and therefore $\int|N_3|\,\d y\le C\eps_b^2$. Integrating the expansion of $N$ over $\Cstar(x)$ and exchanging the integrals in the normalization term, we obtain, on $\Omega_{\rm vc}$,
	\begin{align*}
		&\int_{\Cstar(x)}N\,\d y=L(x)+O(\eps_b^2),\\
		&L(x):=\frac1{n_{\rm fit}b^d}\sum_{i\le n_{\rm fit}}\bar\xi_i(x),\qquad\bar\xi_i(x):=\int w_{\rm int}(x,y)\,\xi_i(x,y)\,\d y,
	\end{align*}
uniformly in $x\in\cX_0$, where the $\xi_i$ are the block terms of \cref{lem:linear} and $w_{\rm int}$ is the deterministic, bounded weight
	\begin{align*}
		w_{\rm int}(x,y):=\frac{K{\mathsf q}_0\psi'(\bar\eta_b(x,y))}{\bar Z(x)}\bigl\{\ind_{\Cstar(x)}(y)-\bar\gamma(x)\bigr\},\qquad\bar\gamma(x):=\int_{\Cstar(x)}\frac{\bar a}{\bar Z}\,\d y.
	\end{align*}
By \cref{lem:linear}, we have $\E\bar\xi_i(x)=0$. A row of block $i$ contributes to $\xi_i(x,\cdot)$ only on a $y$-window of volume at most $Cb^{d_Y}$, with a bounded integrand, and hence $|\bar\xi_i(x)|\le Cb^{d_Y}$. Therefore $\E L(x)=0$ and $|L(x)|\le Cb^{-d_X}\le Cn^{d_X}$, since $b\ge h\ge n^{-1}$.

\emph{Conclusion.} Let $B=\bar{\mathsf p}-{\mathsf p}$ be the bias of \cref{app:lp-estimator}. On $\Omega_{\rm vc}$, by \eqref{eq:lp-one-close} and the last display,
	\begin{align*}
		R^\star(x)=\int_{\Cstar(x)}({\mathsf p}-\hp)(y\mid x)\,\d y+O(n^{-2})=-\int_{\Cstar(x)}B\,\d y-L(x)+O(\eps_b^2+n^{-2}).
	\end{align*}
Off $\Omega_{\rm vc}$, we have $|R^\star(x)|\le2$ and $|L(x)|\le Cn^{d_X}$, and $\bP(\Omega_{\rm vc}^c)\le Cn^{-d_X-2}$. Since $\E L(x)=0$ and $|B|\le C$, taking expectations gives $\E\sbr[0]{R^\star(x)}=-\int_{\Cstar(x)}B\,\d y+O(\eps_b^2+n^{-2})$. By the $L^1$ form of \eqref{eq:leading-bias} and the step \emph{Rates}, it follows that
	\begin{align*}
		\E\sbr[0]{R^\star(x)}=-b^\beta\mathfrak g(x)+o(b^\beta)\qquad\text{uniformly in }x\in\cX_0 .
	\end{align*}
The function $\mathfrak g$ is bounded, because $B_\beta$ is continuous on the compact set $\cX_0\times\cY$, and $\mathfrak g$ is not almost everywhere constant by assumption. Hence \cref{thm:onestage-lower} with $a_n=b^\beta$ and $-\mathfrak g$ in place of $\mathfrak g$ gives the bound.
\end{proof}

\subsection{A curvature condition for non-constant integrated bias}\label{app:gauss-bias}
Let $\beta=2$ (local linear fitting), and let the univariate kernel $\mathcal K_1$ of the product kernel $\mathcal K$ of \cref{app:tools} be symmetric, with $\int\mathcal K\,\d s=1$. We write $\mu_2:=\int s_1^2\mathcal K(s)\,\d s$. Then $\Sigma_{\mathcal K}$ of \eqref{eq:Sigma-K} is $\mathrm{diag}(1,\mu_2I_d)$. By the symmetry of $\mathcal K_1$, the terms of $\mathfrak b_\beta$ with mixed derivatives vanish. Hence, with $\operatorname{tr}\nabla^2_z$ the Laplacian in $z$ and $\operatorname{tr}\nabla^2_x$, $\operatorname{tr}\nabla^2_y$ those in $x$ and $y$,
	\begin{equation}\label{eq:gb-b}
		\mathfrak b_\beta=\frac{\mu_2}2\,\operatorname{tr}\nabla^2_z\eta=\frac{\mu_2}2\bigl(\operatorname{tr}\nabla^2_x\eta+\operatorname{tr}\nabla^2_y\eta\bigr).
	\end{equation}
Let $\mathfrak a_\beta$ be as in \eqref{eq:a-beta}. By the definition of $B_\beta$ in \eqref{eq:leading-bias} and ${\mathsf P}_x(\Cstar(x))=\gamma_\star$, the function $\mathfrak g$ of the proof of \cref{cor:lp-onestage} is
	\begin{align}\label{eq:gb-g}
		\mathfrak g(x)=\int_{\Cstar(x)}B_\beta(x,y)\,\d y=\int_\cY \mathfrak a_\beta(x,y)\bigl(\ind_{\Cstar(x)}(y)-\gamma_\star\bigr)\,\d y .
	\end{align}

\begin{lemma}[A curvature condition]\label{lem:gb-curv}
Let the distribution satisfy the conditions of \cref{cor:lp-onestage} on ${\mathsf P}$, other than the non-constancy of $\mathfrak g$, and let the kernel be as above. Let $x_a$ and $x_b$ be interior points of $\cX_0$ at which ${\mathsf p}(\cdot\mid x_a)={\mathsf p}(\cdot\mid x_b)=:{\mathsf p}_0$ and $\nabla_x{\mathsf p}(\cdot\mid x_a)=\nabla_x{\mathsf p}(\cdot\mid x_b)=0$ on $\cY$, and let ${\mathsf C}_0$ be the oracle set of ${\mathsf p}_0$. Then
	\begin{equation}\label{eq:gb-curv}
		\mathfrak g(x_a)-\mathfrak g(x_b)=\frac{\mu_2}2\int_{{\mathsf C}_0}\bigl\{\operatorname{tr}\nabla^2_x{\mathsf p}(y\mid x_a)-\operatorname{tr}\nabla^2_x{\mathsf p}(y\mid x_b)\bigr\}\,\d y .
	\end{equation}
If the right-hand side is not zero, then $\mathfrak g$ is not almost everywhere constant on $\cX_0$.
\end{lemma}

The integral in \eqref{eq:gb-curv} is the difference between the Laplacians of $x\mapsto{\mathsf P}_x({\mathsf C}_0)$, the conditional probability of the fixed set ${\mathsf C}_0$, at $x_a$ and at $x_b$.

\begin{proof}
Let $\Psi$ be as in \eqref{eq:Psi}. Since $\eta=\Psi({\mathsf p})$ and ${\mathsf p}=K{\mathsf q}_0\psi(\eta)$, the function $K{\mathsf q}_0\psi$ is the inverse of $\Psi$, and hence $K{\mathsf q}_0\psi'(\eta)\Psi'({\mathsf p})=1$. By the chain rule,
	\begin{align*}
		\operatorname{tr}\nabla^2_x\eta=\Psi'({\mathsf p})\operatorname{tr}\nabla^2_x{\mathsf p}+\Psi''({\mathsf p})|\nabla_x{\mathsf p}|^2 .
	\end{align*}
Let $\eta_0:=\Psi({\mathsf p}_0)$. Since $\nabla_x{\mathsf p}$ vanishes at $x_a$ and $x_b$, \eqref{eq:a-beta}, \eqref{eq:gb-b} and this identity give, for $x_j\in\{x_a,x_b\}$,
	\begin{align*}
		\mathfrak a_\beta(x_j,y)=\frac{\mu_2}2\bigl\{\operatorname{tr}\nabla^2_x{\mathsf p}(y\mid x_j)+K{\mathsf q}_0\psi'(\eta_0)\operatorname{tr}\nabla^2_y\eta_0(y)\bigr\}.
	\end{align*}
Since ${\mathsf p}(\cdot\mid x_a)={\mathsf p}(\cdot\mid x_b)$, the second term is the same at the two points, and both points have the oracle set ${\mathsf C}_0$. Since the second derivatives of ${\mathsf p}$ are continuous on the compact set $\cX_0^{\delta_0}\times\cY$, differentiating $\int_\cY{\mathsf p}(y\mid x)\,\d y=1$ twice in $x$ under the integral sign gives $\int_\cY\operatorname{tr}\nabla^2_x{\mathsf p}(y\mid x)\,\d y=0$. Subtracting \eqref{eq:gb-g} at $x_b$ from \eqref{eq:gb-g} at $x_a$, we therefore obtain \eqref{eq:gb-curv}.

For the last claim, the weight $\mathfrak a_\beta$ is continuous on $\cX_0^{\delta_0}\times\cY$, because ${\mathsf p}$ has continuous second derivatives and $\psi'$ is smooth on the range of $\eta$. Moreover, $\Leb(\Cstar(x)\triangle\Cstar(\tilde x))\le C|x-\tilde x|_\infty$ by \cref{lem:tstar-lip}(ii) with $\mathfrak L\equiv1$. Hence, by \eqref{eq:gb-g}, $\mathfrak g$ is continuous on $\cX_0$. If $\mathfrak g(x_a)\ne\mathfrak g(x_b)$, then, since $x_a$ and $x_b$ are interior points of $\cX_0$, $\mathfrak g$ takes different values on two neighborhoods of positive Lebesgue measure, and hence it is not almost everywhere constant on $\cX_0$.
\end{proof}

\emph{A heteroscedastic Gaussian distribution.} Let the constants and the design satisfy \cref{ass:constants}. We use the conventions and the notation of \cref{app:lower}, in which $\cY$ is centered at the origin. For $s>0$ let
	\begin{align*}
		f^{\rm sc}_s(y):=\frac{\e^{-|y|^2/(2s^2)}}{Z^{\rm sc}(s)},\qquad Z^{\rm sc}(s):=\int_\cY\e^{-|y|^2/(2s^2)}\,\d y,\qquad y\in\cY,
	\end{align*}
be the density of $\mathtt N(0,s^2I_{d_Y})$ restricted to $\cY$ and renormalized. The member $f^{\rm sc}_{\sigma_0}$ of this family is the base density $f_0$ of \eqref{eq:lb-f0}, whose level is $\bar t$ and whose oracle set is the ball ${\mathsf C}_0=\bar{\mathcal B}(0,r_\star)$, as shown in the proof of \cref{lem:constants-suffice}. We fix two points $x_a,x_b$ in the interior of $\cX_0$ and a smooth compactly supported function $\chi\ge0$ on $\R^{d_X}$ with $\chi(x)=|x-x_a|^2$ near $x_a$ and $\chi=0$ near $x_b$. For all sufficiently small $\theta>0$, let
	\begin{align*}
		{\mathsf p}(y\mid x):=f^{\rm sc}_{s_\theta(x)}(y),\qquad s_\theta(x):=\sigma_0+\theta\chi(x),\qquad e_\theta(x,y):={\mathsf p}(y\mid x)-f_0(y).
	\end{align*}

\emph{The distribution belongs to $\cP_2$ for all sufficiently small $\theta$.} The map $(s,y)\mapsto f^{\rm sc}_s(y)$ is smooth on $[\sigma_0,2\sigma_0]\times\cY$, and
	\begin{align*}
		e_\theta(x,y)=\theta\chi(x)\int_0^1\partial_s f^{\rm sc}_s(y)\big|_{s=\sigma_0+u\theta\chi(x)}\,\d u.
	\end{align*}
Hence the derivatives of $e_\theta$ of order at most two in $(x,y)$ are at most $C\theta$ on the convex set $[0,1]^{d_X}\times\cY$, and
	\begin{align*}
		\norm{e_\theta}_{\mathcal H^2}\le C\theta,\qquad\int_\cY e_\theta(x,y)\,\d y=0 .
	\end{align*}
Consequently, for all sufficiently small $\theta$, every $e_\theta(x,\cdot)$ is a small perturbation in the sense of \eqref{eq:lb-admissible}, and $\norm[0]{e_\theta}_{\mathcal H^2}\le L/2$. By \eqref{eq:lb-Delta-bound} and $\norm[0]{e_\theta(x,\cdot)}_1\le\Leb(\cY)\sup|e_\theta|$, we have $\Delta_{e_\theta(x,\cdot)}\le C\theta\le c_0/2$ for every $x$. Therefore \cref{lem:lb-perturb}(iii), applied to $x\mapsto e_\theta(x,\cdot)$, gives ${\mathsf P}\in\cP_2$.

Part (i) of \cref{lem:lb-perturb} gives \cref{ass:margin} with $\mathfrak L\equiv1$, and, since ${\mathsf p}\ge{\mathsf p}_\flat$, condition \ref{it:R3} holds with $C_{\rm low}:=\Leb(\cY)/{\mathsf p}_\flat$. Moreover, ${\mathsf p}$ is smooth, and hence \ref{it:R2} holds with $\beta=2$. The distribution thus satisfies the conditions of \cref{cor:lp-onestage} with $\beta=2$, other than the non-constancy of $\mathfrak g$.

\emph{The condition of \cref{lem:gb-curv}.} Since $\chi$ and $\nabla\chi$ vanish at $x_a$ and $x_b$, we have $s_\theta=\sigma_0$ and $\nabla s_\theta=0$ at both points. Hence ${\mathsf p}(\cdot\mid x_a)={\mathsf p}(\cdot\mid x_b)=f_0$, whose oracle set is ${\mathsf C}_0$, and $\nabla_x{\mathsf p}=\partial_s f^{\rm sc}_s|_{s=\sigma_0}\nabla s_\theta=0$ at both points. Moreover, by the chain rule and $\nabla s_\theta=0$, for $x\in\{x_a,x_b\}$,
	\begin{align*}
		\operatorname{tr}\nabla^2_x{\mathsf p}(y\mid x)=\partial_s f^{\rm sc}_s(y)\big|_{s=\sigma_0}\operatorname{tr}\nabla^2s_\theta(x),\qquad\operatorname{tr}\nabla^2s_\theta(x_a)=2\theta d_X,\qquad
		\operatorname{tr}\nabla^2s_\theta(x_b)=0 .
	\end{align*}
By \cref{lem:gb-curv},
	\begin{align*}
		\mathfrak g(x_a)-\mathfrak g(x_b)=\mu_2\theta d_XI_0,\qquad I_0:=\int_{{\mathsf C}_0}\partial_s f^{\rm sc}_s(y)\big|_{s=\sigma_0}\,\d y .
	\end{align*}
We determine the sign of $I_0$. Let $\E_s$ be the expectation under $Y\sim f^{\rm sc}_s$. Then
	\begin{align*}
		\partial_s\log f^{\rm sc}_s(y)=\frac{|y|^2}{s^3}-\frac{(Z^{\rm sc})'(s)}{Z^{\rm sc}(s)},\qquad\frac{(Z^{\rm sc})'(s)}{Z^{\rm sc}(s)}=\frac{\E_s|Y|^2}{s^3}.
	\end{align*}
Since $\partial_sf^{\rm sc}_s=f^{\rm sc}_s\,\partial_s\log f^{\rm sc}_s$ and $\int_{{\mathsf C}_0}f_0\,\d y=\gamma_\star$, it follows that
	\begin{align*}
		I_0=\frac{\gamma_\star}{\sigma_0^3}\Bigl\{\E_{\sigma_0}\bigl(|Y|^2\mid Y\in {\mathsf C}_0\bigr)-\E_{\sigma_0}|Y|^2\Bigr\}.
	\end{align*}
By the law of total expectation, we have
	\begin{align*}
		\E_{\sigma_0}|Y|^2=\gamma_\star\E_{\sigma_0}\bigl(|Y|^2\mid Y\in{\mathsf C}_0\bigr)+(1-\gamma_\star)\E_{\sigma_0}\bigl(|Y|^2\mid Y\notin{\mathsf C}_0\bigr).
	\end{align*}
Since $1-\gamma_\star>0$ and since $|y|\le r_\star$ on ${\mathsf C}_0$ and $|y|>r_\star$ off it, we obtain
	\begin{align*}
		I_0=-\frac{\gamma_\star(1-\gamma_\star)}{\sigma_0^3}\Bigl\{\E_{\sigma_0}\bigl(|Y|^2\mid Y\notin{\mathsf C}_0\bigr)-\E_{\sigma_0}\bigl(|Y|^2\mid Y\in{\mathsf C}_0\bigr)\Bigr\}<0 .
	\end{align*} Therefore $\mathfrak g(x_a)\ne \mathfrak g(x_b)$, and \cref{lem:gb-curv} shows that $\mathfrak g$ is not almost everywhere constant on $\cX_0$. The sign of $I_0$ reflects the fact that widening the conditional distribution moves mass out of the oracle set.

\subsection{Proof of \texorpdfstring{\cref{prop:partition-jump}}{Proposition \ref{prop:partition-jump}}}\label{app:partition-jump}
\begin{proof}
We show that, at a point of the common face of two adjacent cubes, $R^\star$ of \cref{sec:cov-onestage-thm} has one-sided limits that differ by order $h_X$, and we transfer this jump to the coverage by \cref{thm:onestage}. We use the notation of the proof of \cref{prop:partition}, and we work on its event $\Omega_{\rm p}$ with $c=1$, on which \eqref{eq:part-sup} holds and, by that proof, the level sets of $\hp$ are Lebesgue-null. Since $r_n\to0$ and $r_n\ge\nu_{\rm part}\ge(\log n/n)^{1/2}$, the conditions of \cref{thm:onestage} hold with the sequence $Cr_n$, where $C$ is the constant of \eqref{eq:part-sup}, and with $\Omega_r:=\Omega_{\rm p}$. Hence, on an event of probability tending to one,
	\begin{align}\label{eq:part-exp}
		\sup_{x\in\cX_0}\abs[1]{\cov^{\rm one}(x)-\gamma_\star-\{R^\star(x)-\bar R^\star\}}\le C\Bigl\{r_n^2+\Bigl(\frac{\log n}n\Bigr)^{1/2}\Bigr\}\le C(r_n'+h_X^2),
	\end{align}
where the last inequality uses $r_n=h_X+r_n'$, $r_n'\le2$ and $(\log n/n)^{1/2}\le\nu_{\rm part}\le r_n'$.

\emph{The two cubes.} For all sufficiently large $n$, $\mathcal B$ contains the closures of two cubes $\mathcal Q^-$ and $\mathcal Q^+$ of the grid, where $\mathcal Q^+$ is the translate of $\mathcal Q^-$ by $h_X$ in the first coordinate. Let $x_f\in\bar{\mathcal Q}^-\cap\bar{\mathcal Q}^+$ lie on their common face. We have $\sup_{x\in\bar{\mathcal Q}^-\cup\bar{\mathcal Q}^+}|x-x_f|_\infty\le h_X$. For $\mathcal Q\in\{\mathcal Q^-,\mathcal Q^+\}$, let $\hP_{\mathcal Q}$ and ${\mathsf P}_{\mathcal Q}$ be the distributions with densities $\hp_{\mathcal Q}$ and ${\mathsf p}_{\mathcal Q}$, and let $\bar V_{\mathcal Q}:=\Pi(\mathcal Q)^{-1}\int_{\mathcal Q}V_x(\gamma_\star)\,\Pi(\d x)$.

\emph{One-sided limits.} Let $\mathcal Q\in\{\mathcal Q^-,\mathcal Q^+\}$. For $x\in\mathcal Q$, we have $\hp(\cdot\mid x)=\hp_{\mathcal Q}$, and hence $R^\star(x)=\gamma_\star-\hP_{\mathcal Q}(\Cstar(x))$. By \cref{lem:tstar-lip}(ii) with $\mathfrak L\equiv1$, we have $\Leb(\Cstar(x)\triangle\Cstar(x_f))\le C|x-x_f|_\infty$, and $\hp_{\mathcal Q}$ is bounded. Hence $R^\star(x)\to R^\star_{\mathcal Q}:=\gamma_\star-\hP_{\mathcal Q}(\Cstar(x_f))$ as $x\to x_f$ within $\mathcal Q$. By \cref{lem:onestage-cal}(ii), on the event of \eqref{eq:part-exp}, the output at $x\in\mathcal Q$ is the set $A_x(\hat\tau)$, which does not depend on $x\in\mathcal Q$. Since ${\mathsf p}$ is continuous in $x$, the coverage $\cov^{\rm one}(x)$ has a limit as $x\to x_f$ within $\mathcal Q$. Let
	\begin{align*}
		\Delta^{\rm jump}:=\lim_{x\to x_f,\,x\in\mathcal Q^+}\cov^{\rm one}(x)-\lim_{x\to x_f,\,x\in\mathcal Q^-}\cov^{\rm one}(x).
	\end{align*}
Since $\bar R^\star$ does not depend on $x$, \eqref{eq:part-exp} gives
	\begin{align*}
		\Delta^{\rm jump}=R^\star_{\mathcal Q^+}-R^\star_{\mathcal Q^-}+O(r_n'+h_X^2)=\hP_{\mathcal Q^-}(\Cstar(x_f))-\hP_{\mathcal Q^+}(\Cstar(x_f))+O(r_n'+h_X^2).
	\end{align*}

\emph{The jump.} By the first bound of \eqref{eq:part-sup}, we have $|\hP_{\mathcal Q}(\Cstar(x_f))-{\mathsf P}_{\mathcal Q}(\Cstar(x_f))|\le Cr_n'$. Let $x\in\bar{\mathcal Q}$. Since ${\mathsf p}$ is Lipschitz, $\sup_y|{\mathsf p}(y\mid x_f)-{\mathsf p}(y\mid x)|\le Ch_X$, and hence $|{\mathsf P}_x(\Cstar(x_f))-\gamma_\star|\le Ch_X\le c_0$ for all sufficiently large $n$. Therefore \cref{prop:supshape} with $\hat a={\mathsf p}(\cdot\mid x_f)$ and $r=Ch_X$ bounds the excess mass of $\Cstar(x_f)$ at $x$ by $Ch_X^2$, and \cref{lem:identity} gives
	\begin{align*}
		{\mathsf P}_x(\Cstar(x_f))-\gamma_\star=t_\star(x)\bigl\{V_{x_f}(\gamma_\star)-V_x(\gamma_\star)\bigr\}+O(h_X^2).
	\end{align*}
By \cref{lem:tstar-lip}, $|t_\star(x)-t_\star(x_f)|$ and $|V_{x_f}(\gamma_\star)-V_x(\gamma_\star)|$ are $O(h_X)$, and hence replacing $t_\star(x)$ by $t_\star(x_f)$ changes the right-hand side by $O(h_X^2)$. Averaging over $x\in\mathcal Q$ under $\Pi$ and taking the difference between $\mathcal Q^-$ and $\mathcal Q^+$, in which $V_{x_f}(\gamma_\star)$ cancels, we obtain
	\begin{align*}
		\abs[1]{\Delta^{\rm jump}-t_\star(x_f)\bigl(\bar V_{\mathcal Q^+}-\bar V_{\mathcal Q^-}\bigr)}\le C(r_n'+h_X^2).
	\end{align*}

By \cref{ass:design}, $f_X$ is Lipschitz. Let $\tilde V_{\mathcal Q}:=\Leb(\mathcal Q)^{-1}\int_{\mathcal Q}V_x(\gamma_\star)\,\d x$ and $\bar f_{\mathcal Q}:=\Pi(\mathcal Q)/\Leb(\mathcal Q)$. Then
	\begin{align*}
		\abs[0]{\bar V_{\mathcal Q}-\tilde V_{\mathcal Q}}=\frac1{\Pi(\mathcal Q)}\abs[3]{\int_{\mathcal Q}\bigl\{V_x(\gamma_\star)-V_{x_f}(\gamma_\star)\bigr\}\bigl\{f_X(x)-\bar f_{\mathcal Q}\bigr\}\,\d x}\le Ch_X^2,
	\end{align*}
and $\varsigma(\tilde V_{\mathcal Q^+}-\tilde V_{\mathcal Q^-})\ge c_Vh_X$ by the monotonicity assumption on $V_x(\gamma_\star)$ in \cref{prop:partition-jump}. Since $r_n'/h_X\to0$ and $h_X\to0$, combining these bounds with the bound on $\Delta^{\rm jump}$ and $t_\star(x_f)\ge t_{\min}$ (\cref{lem:t-window}) gives, for all sufficiently large $n$,
	\begin{align*}
		\varsigma\Delta^{\rm jump}&\ge t_\star(x_f)\,\varsigma\bigl(\bar V_{\mathcal Q^+}-\bar V_{\mathcal Q^-}\bigr)-C(r_n'+h_X^2)\\
		&\ge t_\star(x_f)\,\varsigma\bigl(\tilde V_{\mathcal Q^+}-\tilde V_{\mathcal Q^-}\bigr)-C(r_n'+h_X^2)\\
		&\ge t_{\min}c_Vh_X-C(r_n'+h_X^2)\ge\frac{t_{\min}c_Vh_X}2.
	\end{align*}

Hence one of the two one-sided limits of $\cov^{\rm one}$ at $x_f$ differs from $\gamma_\star$ by at least $|\Delta^{\rm jump}|/2\ge t_{\min}c_Vh_X/4$. Since $\mathcal Q^\pm\subset\mathcal B\subset\cX_0$, this gives the bound with $c:=t_{\min}c_V/4$. Finally, $\nu_{\rm part}/h_X\to0$ and $\nu_{\rm part}\ge\{\log n/(nh_X^{d_X})\}^{1/2}$ imply $h_X/(\log n/n)^{1/(2+d_X)}\to\infty$, and $\rho_n\le(\log n/n)^{1/(2+d_X)}$ since $\beta\ge1$. Hence $h_X/\rho_n\to\infty$.
\end{proof}

\section{Proofs and additional results for \texorpdfstring{\cref{sec:shape}}{Section \ref{sec:shape}}}\label{app:shape}
Throughout this appendix, \cref{ass:margin} holds with $\mathfrak L\equiv1$, as in \cref{sec:shape}. For $q=\infty$, expressions in $q$ take their limits as $q\to\infty$, for instance $1/q=0$ and $2q/(q+1)=2$. We keep $x$ fixed, use the abbreviations of \cref{app:setup}, with $\cS$ and $\cE$ evaluated at the set at hand, and write $\hat a$ for $\hat a(\cdot\mid x)$.

\begin{lemma}[Misplaced points]\label{lem:swap}
Under \cref{ass:smooth,ass:margin}, let $A\subset\cY$ have coverage $\gamma:={\mathsf P}_x(A)$ in the window, and let $v\in\R$ and $\kappa>0$. If $\mathcal D^+\subset A\cap\{y:{\mathsf p}(y\mid x)\le v-\kappa\}$ and $\mathcal D^-\subset\{y:{\mathsf p}(y\mid x)\ge v+\kappa\}\setminus A$, then $\cS(A)\ge\kappa\min\{\Leb(\mathcal D^+),\Leb(\mathcal D^-)\}/{\mathsf p}_{\max}$.
\end{lemma}

\begin{proof}
By \cref{lem:t-window}, we have $t_\gamma>0$, and $t_\gamma\le\sup_\cY{\mathsf p}(\cdot\mid x)\le{\mathsf p}_{\max}$ by \cref{ass:smooth}. If $t_\gamma\le v$, then ${\mathsf p}(\cdot\mid x)\ge t_\gamma+\kappa$ on $\mathcal D^-$, and hence $\mathcal D^-\subset\Cstar_\gamma\setminus A$. If $t_\gamma>v$, then ${\mathsf p}(\cdot\mid x)<t_\gamma-\kappa$ on $\mathcal D^+$, and hence $\mathcal D^+\subset A\setminus\Cstar_\gamma$. Thus, in either case, $\Leb\{y\in A\triangle\Cstar_\gamma:|{\mathsf p}(y\mid x)-t_\gamma|\ge\kappa\}\ge\min\{\Leb(\mathcal D^+),\Leb(\mathcal D^-)\}$. Since $t_\gamma\le{\mathsf p}_{\max}$, the integral expression of $\cS(A)$ in \cref{sec:setup-shape} gives
	\begin{align*}
		\cS(A)=\frac1{t_\gamma}\int_{A\triangle\Cstar_\gamma}|{\mathsf p}(y\mid x)-t_\gamma|\,\d y\ge\frac\kappa{{\mathsf p}_{\max}}\min\{\Leb(\mathcal D^+),\Leb(\mathcal D^-)\}.
	\end{align*}
\end{proof}

\subsection{Proof of \texorpdfstring{\cref{thm:Lq}}{Theorem \ref{thm:Lq}}}
\begin{proof}
Let $\hat a_\phi:=\phi\circ\hat a$ and $\hat t_\phi:=\phi(\hat t)$. Since $\phi$ is strictly increasing, $\hC=\{y:\hat a(y)\ge\hat t\}=\{y:\hat a_\phi(y)\ge\hat t_\phi\}$, and $\sup_\cY|\hat a_\phi-{\mathsf p}(\cdot\mid x)|\le\varepsilon_0/8$ by assumption. Let $\mathfrak M:=\hC\triangle\Cstar_{\cov}$. By \cref{lem:shape-error}, applied to $\hat a_\phi$ and $\hat t_\phi$, we have $\mathfrak M\subset\mathcal J_x$, $\cS\le C\int_{\mathfrak M}|\hat a_\phi-{\mathsf p}(\cdot\mid x)|\,\d y$ and $\Leb(\mathfrak M)\le C\cS^{1/2}$. Since $\mathfrak M\subset\mathcal J_x$, H\"older's inequality gives
	\begin{align*}
		\cS\le Cr^{(q)}\Leb(\mathfrak M)^{1-1/q}\le C_2r^{(q)}\cS^{(q-1)/(2q)},
	\end{align*}
where $C_2$ does not depend on $q$, since $1-1/q\in[0,1]$. The claim is trivial if $\cS=0$. If $\cS>0$, dividing by $\cS^{(q-1)/(2q)}$ gives $\cS^{(q+1)/(2q)}\le C_2r^{(q)}$. Since $2q/(q+1)\in[1,2]$, it follows that $\cS\le\max\{C_2,C_2^2\}(r^{(q)})^{2q/(q+1)}$.
\end{proof}

\subsection{Proof of \texorpdfstring{\cref{thm:sep-smooth}}{Theorem \ref{thm:sep-smooth}}}
\begin{proof}
The estimate adds to ${\mathsf p}(\cdot\mid x)=f_0$ a smooth perturbation of height $\delta$, supported in a shell of half-width $w\gg\delta$ around the level set, nonnegative on the half-space $\{y:y_{d_Y}>0\}$ and nonpositive on the other half. The H\"older condition permits $w\asymp\delta^{1/\beta}$. On two antipodal caps of the shell, the perturbation equals $\delta$ and $-\delta$, respectively. Hence every super-level set of the estimate at a level within $\delta/4$ of $t_\star$ contains a layer of thickness of order $\delta$ in the first cap on which $f_0$ lies below the level by at least $\delta/4$, and misses a similar layer in the second cap on which $f_0$ lies above the level by at least $\delta/4$. By \cref{lem:swap}, it follows that $\cS\gtrsim\delta^2$.

We use the notation introduced in \cref{app:lower} before \cref{lem:lb-perturb}, except that $q$ remains the exponent of $L^q$. Here $t_\star=\bar t$, the $\gamma_\star$-level of $f_0$, and the level set $\{y:f_0(y)=\bar t\}$ is the sphere $\{y:|y|=r_\star\}$. By \cref{lem:lb-perturb}(iii) with $\varphi_x=0$, the distribution with $\Pi$ uniform and ${\mathsf p}(\cdot\mid x)=f_0$ belongs to $\cP_\beta$. Hence it satisfies \cref{ass:smooth,ass:margin}. By (L2), $\bar f(r_\star)=\bar t$ and the mean value theorem,
	\begin{equation}\label{eq:sep-radial}
		2\kappa_0\bigl||y|-r_\star\bigr|\le|f_0(y)-\bar t|\le C_{\bar f'}\bigl||y|-r_\star\bigr|\qquad(r_-\le|y|\le r_+).
	\end{equation}

\emph{Step 1: the perturbation.} Let $\chi\in C_c^\infty((-1,1))$ with $0\le\chi\le1$ and $\chi=1$ on $[-\frac12,\frac12]$, and let $\chi_{\rm odd}\in C^\infty(\R)$ be odd and nondecreasing, with $\chi_{\rm odd}=1$ on $[\frac12,\infty)$. Then $|\chi_{\rm odd}|\le1$ and $\chi_{\rm odd}=-1$ on $(-\infty,-\frac12]$. Let $y_{d_Y}$ denote the last coordinate of $y$ and, for $0<w\le\min\{1,r_\star/2\}$, let
	\begin{align*}
		\mathrm{dip}_w(y):=\chi\Bigl(\frac{|y|-r_\star}w\Bigr)\,\chi_{\rm odd}\Bigl(\frac{y_{d_Y}}{r_\star}\Bigr)\qquad(y\in\R^{d_Y}).
	\end{align*}
The function $\mathrm{dip}_w$ vanishes outside the shell $\mathcal A_w:=\{y:||y|-r_\star|<w\}$. Since $|y|>r_\star/2$ on $\mathcal A_w$, the map $y\mapsto|y|$ is smooth there. Therefore $\mathrm{dip}_w\in C^\infty(\R^{d_Y})$. Moreover, since $\chi_{\rm odd}$ is odd, $\mathrm{dip}_w(-y)=-\mathrm{dip}_w(y)$.

The derivatives of $y\mapsto|y|$ are bounded on $\{y:|y|\ge r_\star/2\}$, each derivative of $(|y|-r_\star)/w$ contains the factor $1/w$, and $w\le1$. Hence, by the chain and product rules, there exists a constant $C_{\rm dip}$, depending only on $\chi$, $\chi_{\rm odd}$, $r_\star$, $d_Y$ and $\beta$, such that $\sup|\partial^{\boldsymbol k}\mathrm{dip}_w|\le C_{\rm dip}w^{-|\boldsymbol k|}$ for $|\boldsymbol k|\le\ell+1$. For $|\boldsymbol k|=\ell$ and $z,z'\in\R^{d_Y}$, this bound and the mean value theorem give, since $\beta-\ell\in(0,1]$,
	\begin{align*}
		|\partial^{\boldsymbol k}\mathrm{dip}_w(z)-\partial^{\boldsymbol k}\mathrm{dip}_w(z')|\le2d_Y^{1/2}C_{\rm dip}w^{-\ell}\min\Bigl\{1,\frac{|z-z'|}w\Bigr\}\le C_{\rm sep}w^{-\beta}|z-z'|^{\beta-\ell},
	\end{align*}
with $C_{\rm sep}:=2d_Y^{1/2}C_{\rm dip}$. Since $w\le1$, the derivatives of order at most $\ell$ are also bounded by $C_{\rm sep}w^{-\beta}$, and hence $\norm[0]{\mathrm{dip}_w}_{\mathcal H^\beta_{\R^{d_Y}}}\le C_{\rm sep}w^{-\beta}$.
Let
	\begin{equation}\label{eq:sep-w}
		w:=\Bigl(\frac{2C_{\rm sep}\delta}L\Bigr)^{1/\beta}.
	\end{equation}
Then $w\asymp\delta^{1/\beta}$, and $\delta/w=\delta^{1-1/\beta}(L/(2C_{\rm sep}))^{1/\beta}$ tends to zero with $\delta$ because $\beta>1$. Hence we can take $\delta$ so small that
	\begin{equation}\label{eq:sep-small}
		w\le\min\Bigl\{1,\ \frac12\min\{r_\star-r_-,\,r_+-r_\star\}\Bigr\},\qquad4\delta\le \kappa_0w,\qquad C_{\bar f'}w+\delta\le\frac{\varepsilon_0}2,\qquad\delta\le\frac{\varepsilon_0}8,
	\end{equation}
together with $\frac54{\mathsf p}_{\max}C_m\delta\le c_0$, which is used in Step~2. Then, since $r_-=r_\star/2$, we have $w\le r_\star/4$, and $\mathcal A_w\subset\mathcal A\subset\cY$, with $\mathcal A$ as in \eqref{eq:lb-annulus}. We set
	\begin{align*}
		{\mathsf p}_\delta(\cdot\mid x):=f_0+\delta\mathrm{dip}_w .
	\end{align*}

\emph{Density.} Since $\mathrm{dip}_w$ is odd and $\cY=-\cY$, we have $\int_\cY\mathrm{dip}_w\,\d y=0$, and hence $\int_\cY{\mathsf p}_\delta(y\mid x)\,\d y=1$. Off $\mathcal A_w$, we have ${\mathsf p}_\delta(\cdot\mid x)=f_0\ge0$. On $\mathcal A_w$, we have $|f_0-\bar t|\le C_{\bar f'}w$ by \eqref{eq:sep-radial}. Therefore, by \eqref{eq:sep-small} and $\bar t\ge2\varepsilon_0$,
	\begin{align*}
		{\mathsf p}_\delta(\cdot\mid x)\ge\bar t-C_{\bar f'}w-\delta\ge\bar t-\frac{\varepsilon_0}2>0 .
	\end{align*}
Thus ${\mathsf p}_\delta(\cdot\mid x)$ is a probability density on $\cY$.

\emph{H\"older norm.} We regard ${\mathsf p}_\delta$ as a function of $(x,y)$ that does not depend on $x$. Since $|y-y'|\le|(x,y)-(x',y')|$, it suffices to bound its H\"older norm in $y$. By the bound on $\mathrm{dip}_w$ and \eqref{eq:sep-w}, we have $\norm[0]{\delta\mathrm{dip}_w}_{\mathcal H^\beta}\le C_{\rm sep}\delta w^{-\beta}=L/2$. Hence, by (L1) and subadditivity, $\norm[0]{{\mathsf p}_\delta}_{\mathcal H^\beta}\le L$, that is, ${\mathsf p}_\delta\in\mathcal H^\beta(L)$.

\emph{Norms.} Let
	\begin{align*}
		\Theta_+&:=\{\theta\in\mathbb S^{d_Y-1}:\theta_{d_Y}\ge\tfrac34\},\qquad\Theta_-:=-\Theta_+,\\
		\mathcal C_\pm&:=\{r\theta:\theta\in\Theta_\pm,\ |r-r_\star|\le w/2\},
	\end{align*}
and $c_\Theta:=S_{d_Y-1}(\Theta_+)>0$. For $r\theta\in\mathcal C_\pm$, we have $|r-r_\star|\le w/2\le r_\star/8$. Hence $\chi((r-r_\star)/w)=1$ and $\pm r\theta_{d_Y}\ge\frac{7r_\star}8\cdot\frac34>\frac{r_\star}2$. Therefore $\mathrm{dip}_w=\pm1$ on $\mathcal C_\pm$. Since $|\mathrm{dip}_w|\le\ind_{\mathcal A_w}$ and the disjoint sets $\mathcal C_\pm$ lie in $\mathcal A_w$, polar coordinates give, for $q\in\{1,2\}$,
	\begin{align*}
		2c_\Theta r_-^{d_Y-1}w\le\Leb(\mathcal C_+\cup\mathcal C_-)\le\int_\cY|\mathrm{dip}_w|^q\,\d y\le\Leb(\mathcal A_w)\le2|\mathbb S^{d_Y-1}|r_+^{d_Y-1}w .
	\end{align*}
Since $\norm[0]{{\mathsf p}_\delta(\cdot\mid x)-f_0}_{L^2}=\delta\norm[0]{\mathrm{dip}_w}_{L^2}$ and $\TV({\mathsf P}_{\delta,x},{\mathsf P}_x)=\frac\delta2\norm[0]{\mathrm{dip}_w}_{L^1}$, the last display and $w\asymp\delta^{1/\beta}$ give
	\begin{equation}\label{eq:sep-tv}
		\norm{{\mathsf p}_\delta(\cdot\mid x)-f_0}_{L^2}\asymp\delta w^{1/2}\asymp\delta^{1+1/(2\beta)},\qquad\TV({\mathsf P}_{\delta,x},{\mathsf P}_x)\asymp\delta w\asymp\delta^{1+1/\beta}.
	\end{equation}
Moreover, $\norm[0]{{\mathsf p}_\delta(\cdot\mid x)-f_0}_\infty=\delta$, since $|\mathrm{dip}_w|\le1$ with equality on $\mathcal C_\pm$. This proves the first two claims.

\emph{Step 2: the coverage of the super-level sets of ${\mathsf p}_\delta(\cdot\mid x)$.} Let $s$ satisfy $|s-\bar t|\le\delta/4$, and let $\gamma_s:={\mathsf P}_x\{{\mathsf p}_\delta(Y\mid x)\ge s\}$. Since $|{\mathsf p}_\delta(\cdot\mid x)-f_0|\le\delta$, \eqref{eq:sandwich} and \eqref{eq:margin-upper} give
	\begin{align*}
		\Leb\bigl(\{y:{\mathsf p}_\delta(y\mid x)\ge s\}\triangle\{y:f_0(y)\ge\bar t\}\bigr)\le\Leb\Bigl\{y:|f_0(y)-\bar t|\le\frac{5\delta}4\Bigr\}\le\frac54C_m\delta .
	\end{align*}
Since $f_0\le{\mathsf p}_{\max}$ and ${\mathsf P}_x\{f_0(Y)\ge\bar t\}=G(\bar t)=\gamma_\star$, it follows that
	\begin{equation}\label{eq:sep-gamma}
		|\gamma_s-\gamma_\star|\le\frac54{\mathsf p}_{\max}C_m\delta\le c_0
	\end{equation}
by the choice of $\delta$. Therefore the coverage of $\{y:{\mathsf p}_\delta(y\mid x)\ge s\}$ is in the window.

\emph{Step 3: the shape loss at levels near $\bar t$.} Let $s$ satisfy $|s-\bar t|\le\delta/4$, let $\hC_s:=\{y:{\mathsf p}_\delta(y\mid x)\ge s\}$ and $\gamma:=\gamma_s$, and let
	\begin{align*}
		I^+&:=\Bigl\{r\in[r_-,r_+]:s-\frac\delta2\le\bar f(r)\le s-\frac\delta4\Bigr\},&\mathcal D^+&:=\{r\theta:\theta\in\Theta_+,\ r\in I^+\},\\
		I^-&:=\Bigl\{r\in[r_-,r_+]:s+\frac\delta4\le\bar f(r)\le s+\frac\delta2\Bigr\},&\mathcal D^-&:=\{r\theta:\theta\in\Theta_-,\ r\in I^-\}.
	\end{align*}
By (L2), $\bar f$ maps $[r_-,r_+]$ onto an interval containing $[\bar t-4\varepsilon_0,\bar t+4\varepsilon_0]\supset[s-\delta/2,s+\delta/2]$. Since $\bar f$ is strictly decreasing with $|\bar f'|\le C_{\bar f'}$ on $[r_-,r_+]$, the sets $I^\pm$ are intervals of length at least $\delta/(4C_{\bar f'})$. For $r\in I^\pm$, we have $|\bar f(r)-\bar t|\le3\delta/4$, and hence, by \eqref{eq:sep-radial} and \eqref{eq:sep-small}, $|r-r_\star|\le3\delta/(8\kappa_0)\le w/2$. Consequently, $\mathcal D^\pm\subset\mathcal C_\pm$, where ${\mathsf p}_\delta(\cdot\mid x)=f_0\pm\delta$, and
	\begin{align*}
		{\mathsf p}_\delta(\cdot\mid x)&=f_0+\delta\ge s+\frac\delta2\quad\text{and}\quad f_0\le s-\frac\delta4\quad\text{on }\mathcal D^+,\\
		{\mathsf p}_\delta(\cdot\mid x)&=f_0-\delta\le s-\frac\delta2\quad\text{and}\quad f_0\ge s+\frac\delta4\quad\text{on }\mathcal D^- .
	\end{align*}
Thus $\mathcal D^+\subset\hC_s\cap\{f_0\le s-\delta/4\}$ and $\mathcal D^-\subset\{f_0\ge s+\delta/4\}\setminus\hC_s$. In polar coordinates, $\Leb(\mathcal D^\pm)\ge c_\Theta r_-^{d_Y-1}\delta/(4C_{\bar f'})$. By \eqref{eq:sep-gamma}, \cref{lem:swap} applies with $A=\hC_s$, $v=s$ and $\kappa=\delta/4$ and yields
	\begin{align*}
		\cS\ge\frac\delta4\cdot\frac{c_\Theta r_-^{d_Y-1}\delta}{4C_{\bar f'}{\mathsf p}_{\max}}=\frac{c_\Theta r_-^{d_Y-1}}{16C_{\bar f'}{\mathsf p}_{\max}}\,\delta^2,
	\end{align*}
and $\cE\ge\cS$ by \cref{lem:shape-split}. For the upper bounds, we apply \cref{prop:supshape} to $\hat a:={\mathsf p}_\delta(\cdot\mid x)$, with $\delta=\sup_\cY|{\mathsf p}_\delta(\cdot\mid x)-f_0|\le\varepsilon_0/8$ in place of $r$. By that proposition and \eqref{eq:sep-gamma}, we obtain
	\begin{align*}
		\cS\le C\delta^2,\qquad\cE\le C\bigl(\delta^2+|\gamma-\gamma_\star|^2\bigr)\le C\delta^2 .
	\end{align*}
This proves the third claim.

\emph{Step 4: the rates stated after the theorem.} Let $\eps:=\delta^{1+1/(2\beta)}$, as in \cref{sec:shape}, so that $\delta=\eps^{2\beta/(2\beta+1)}$. By \eqref{eq:sep-tv}, the $L^2$ error is of order $\eps$, and
	\begin{align*}
		\TV({\mathsf P}_{\delta,x},{\mathsf P}_x)\asymp\delta w\asymp\delta^{1+1/\beta}\le\delta^{1+1/(2\beta)}=\eps .
	\end{align*}
Hence the $L^2$ and total-variation errors are at most of order $\eps$, while, at every level within $\delta/4$ of $\bar t$,
	\begin{align*}
		\cS\asymp\delta^2=\eps^{4\beta/(2\beta+1)},
	\end{align*}
which is of larger order than $\eps^2$ because $4\beta/(2\beta+1)<2$.
\end{proof}

\subsection{Optimality of the exponent in \texorpdfstring{\cref{thm:Lq}}{Theorem \ref{thm:Lq}}}\label{app:Lq-lower}
\begin{theorem}[Optimality of the exponent]\label{thm:Lq-lower}
Suppose that \cref{ass:smooth} holds and that \cref{ass:margin} holds with $\mathfrak L\equiv1$. Then there exist constants $c,C>0$ such that, for every $\delta\in(0,c]$, there exists a conditional density ${\mathsf p}_\delta$, used in place of an estimate and equal to the true conditional density at covariates other than $x$, with the following properties.
\begin{enumerate}[label=(\roman*)]
\item For every $q\in[1,\infty]$ and every level $s$ with $|s-t_\star(x)|\le\delta$, with $\mathcal J_x$ as in \eqref{eq:Jx},
	\begin{align*}
		c\,\delta^{1+1/q}&\le\norm{{\mathsf p}_\delta(\cdot\mid x)-{\mathsf p}(\cdot\mid x)}_{L^q(\mathcal J_x)}\le C\delta^{1+1/q},\\
		c\,\delta^2&\le\cS\bigl(\{y:{\mathsf p}_\delta(y\mid x)\ge s\}\bigr)\le C\delta^2.
	\end{align*}
\item For every $\vartheta\in(0,1)$, if the calibration points are i.i.d.\ copies of $(X,Y)$, then, conditionally on $n_{\rm cal}^\circ\ge C\delta^{-2}\log(2/\vartheta)$, with probability at least $1-\vartheta$, the output of HPD-split applied to ${\mathsf p}_\delta$ is $\{y:{\mathsf p}_\delta(y\mid x)\ge s\}$ for some $s$ with $|s-t_\star(x)|\le\delta$.
\end{enumerate}
In particular, the exponent $2q/(q+1)$ in \cref{thm:Lq} cannot be improved for any $q\in[1,\infty]$.
\end{theorem}

\begin{proof}[Proof of \cref{thm:Lq-lower}]
The estimate raises the density by $2\delta$ on a set of measure of order $\delta$ on which ${\mathsf p}(\cdot\mid x)$ lies less than $\delta$ below the level, and lowers it by $2\delta$ on a set of the same measure on which ${\mathsf p}(\cdot\mid x)$ lies less than $\delta$ above the level. At every level within $\delta$ of $t_\star$, the first set therefore lies in the super-level set of the estimate and the second set outside it. The resulting estimate has $L^q$ error $r^{(q)}\asymp\delta^{1+1/q}$ and shape loss $\cS\asymp\delta^2$ at all these levels. Thus $\cS\asymp (r^{(q)})^{2q/(q+1)}$. The total-variation error of the estimate is $O(\delta^2)=o(\delta)$. Consequently the calibrated level stays within $\delta$ of $t_\star$.

\emph{Step 1: the construction.} Let $\delta>0$ be sufficiently small, with $2\delta\le\varepsilon_0/8$, and let $0<c_{\rm swap}\le c_m/(2{\mathsf p}_{\max})$. By \eqref{eq:margin-lower}, the bands
	\begin{align*}
		\{y:t_\star-\delta<{\mathsf p}(y\mid x)\le t_\star-a\delta\}\qquad\text{and}\qquad\{y:t_\star+a\delta\le{\mathsf p}(y\mid x)<t_\star+\delta\}
	\end{align*}
have ${\mathsf P}_x$-probability at least $c_m\delta$, and hence, since ${\mathsf p}(\cdot\mid x)\le{\mathsf p}_{\max}$, Lebesgue measure at least $c_m\delta/{\mathsf p}_{\max}$. We choose subsets $\mathcal D^+$ of the first and $\mathcal D^-$ of the second, each of Lebesgue measure $\mu_{\rm swap}:=c_{\rm swap}\delta$, and set
	\begin{align*}
		{\mathsf p}_\delta(\cdot\mid x):={\mathsf p}(\cdot\mid x)+2\delta\ind_{\mathcal D^+}-2\delta\ind_{\mathcal D^-} .
	\end{align*}
The function ${\mathsf p}_\delta(\cdot\mid x)$ integrates to one, because $\mathcal D^+$ and $\mathcal D^-$ have the same measure. It is nonnegative for all sufficiently small $\delta$, because ${\mathsf p}(\cdot\mid x)\ge t_\star+a\delta$ on $\mathcal D^-$. Since $\mathcal D^\pm\subset\mathcal J_x$, its $L^q$ error is
	\begin{align*}
		\norm{{\mathsf p}_\delta(\cdot\mid x)-{\mathsf p}(\cdot\mid x)}_{L^q(\mathcal J_x)}=2\delta(2\mu_{\rm swap})^{1/q}\asymp\delta^{1+1/q},
	\end{align*}
uniformly in $q$. Further,
	\begin{align*}
		{\mathsf p}_\delta(\cdot\mid x)&>t_\star-\delta+2\delta=t_\star+\delta\quad\text{on }\mathcal D^+,\\
		{\mathsf p}_\delta(\cdot\mid x)&<t_\star+\delta-2\delta=t_\star-\delta\quad\text{on }\mathcal D^-,
	\end{align*}
and hence, for every $s$ with $|s-t_\star|\le\delta$, ${\mathsf p}_\delta(\cdot\mid x)>s$ on $\mathcal D^+$ and ${\mathsf p}_\delta(\cdot\mid x)<s$ on $\mathcal D^-$.

\emph{Step 2: every level within $\delta$ of $t_\star$ has shape loss of order $\delta^2$.} Let $s$ satisfy $|s-t_\star|\le\delta$, and let $\hC_s:=\{y:{\mathsf p}_\delta(y\mid x)\ge s\}$ and $\gamma:={\mathsf P}_x(\hC_s)$. We have $\hC_s\triangle\{y:{\mathsf p}(y\mid x)\ge s\}\subset\mathcal D^+\cup\mathcal D^-$. Hence, since $G(t_\star)=\gamma_\star$ by \cref{lem:t-window}, \eqref{eq:margin-upper} gives
	\begin{align*}
		|\gamma-\gamma_\star|\le|G(s)-G(t_\star)|+{\mathsf P}_x(\mathcal D^+\cup\mathcal D^-)\le{\mathsf p}_{\max}C_m\delta+2{\mathsf p}_{\max}\mu_{\rm swap}\le c_0
	\end{align*}
for all sufficiently small $\delta$. Thus the coverage of $\hC_s$ is in the window. By Step 1, we have $\mathcal D^+\subset\hC_s\cap\{y:{\mathsf p}(y\mid x)\le t_\star-a\delta\}$ and $\mathcal D^-\subset\{y:{\mathsf p}(y\mid x)\ge t_\star+a\delta\}\setminus\hC_s$. Since $\Leb(\mathcal D^\pm)=\mu_{\rm swap}$, \cref{lem:swap} with $A=\hC_s$, $v=t_\star$ and $\kappa=a\delta$ gives $\cS\ge ac_{\rm swap}\delta^2/{\mathsf p}_{\max}$. The matching upper bound $\cS\lesssim\delta^2$ follows from \cref{thm:Lq} with $q=\infty$ and $\phi$ the identity, since $\sup_\cY|{\mathsf p}_\delta(\cdot\mid x)-{\mathsf p}(\cdot\mid x)|=r^{(\infty)}=2\delta\le\varepsilon_0/8$.

\emph{Step 3: the calibrated level lies within $\delta$ of $t_\star$.} By \eqref{eq:margin-lower},
	\begin{align*}
		G(t_\star-\delta)-G(t_\star)\ge{\mathsf P}_x\{y:t_\star-\delta<{\mathsf p}(y\mid x)\le t_\star-a\delta\}\ge c_m\delta .
	\end{align*}
Off $\mathcal D^-$ we have ${\mathsf p}_\delta(\cdot\mid x)\ge{\mathsf p}(\cdot\mid x)$. Hence, by the last display and the choice of $c_{\rm swap}$,
	\begin{align*}
		{\mathsf P}_x\{{\mathsf p}_\delta(Y\mid x)\ge t_\star-\delta\}\ge G(t_\star-\delta)-{\mathsf P}_x(\mathcal D^-)\ge\gamma_\star+c_m\delta-{\mathsf p}_{\max}\mu_{\rm swap}\ge\gamma_\star+\frac{c_m}2\delta .
	\end{align*}
Symmetrically, off $\mathcal D^+$ we have ${\mathsf p}_\delta(\cdot\mid x)\le{\mathsf p}(\cdot\mid x)$, and therefore
	\begin{align*}
		{\mathsf P}_x\{{\mathsf p}_\delta(Y\mid x)>t_\star+\delta\}\le{\mathsf P}_x\{{\mathsf p}(Y\mid x)>t_\star+\delta\}+{\mathsf P}_x(\mathcal D^+)\le\gamma_\star-\frac{c_m}2\delta .
	\end{align*}

We apply HPD-split with $\hp:={\mathsf p}_\delta$, which is the true conditional density at covariate values other than $x$. Then, with $\TV_x$, $\overline{\TV}$, $\eps_{\rm cal}$ and $\Delta$ as in \eqref{eq:tv-consts},
	\begin{align*}
		\TV_x&=\frac12\int|{\mathsf p}_\delta(y\mid x)-{\mathsf p}(y\mid x)|\,\d y=2\delta\mu_{\rm swap}\lesssim\delta^2,\\
		\overline{\TV}&=\Pi(\{x\}\mid X\in\cX_0)\,\TV_x\lesssim\delta^2.
	\end{align*}
If $n_{\rm cal}^\circ\ge C\delta^{-2}\log(2/\vartheta)$ with a sufficiently large constant $C$, then, since $\TV_x+\overline{\TV}\le2\TV_x=4c_{\rm swap}\delta^2$,
	\begin{align*}
		\TV_x+\Delta\le4c_{\rm swap}\delta^2+\eps_{\rm cal}+\frac2{n_{\rm cal}^\circ}<\frac{c_m}2\delta .
	\end{align*}

We work on the event of probability at least $1-\vartheta$ constructed in Step~1 of the proof of \cref{prop:tv}. As noted at the start of \cref{app:tv}, this construction applies to a general $\hp$. Since $1/(n_{\rm cal}^\circ+1)<\Delta<c_m\delta/2\le1-\gamma_\star$ for all sufficiently small $\delta$, the first case of Step~1 of the proof of \cref{prop:tv}, which requires $1-\gamma_\star<1/(n_{\rm cal}^\circ+1)$, is excluded. The second case of that step therefore gives $|\hat\tau-\gamma_\star|\le\Delta$, and hence $|\hat\tau-\gamma_\star|\le c_0$ for all sufficiently small $\delta$. We apply \cref{lem:hpd-sets} with $\hp:={\mathsf p}_\delta$. Its conditions hold for all sufficiently small $\delta$: we have
	\begin{align*}
		&\eps_{\rm sup}=\sup_y|{\mathsf p}_\delta(y\mid x)-{\mathsf p}(y\mid x)|=2\delta\le\frac{\varepsilon_0}8,\qquad\Leb(\cY)\eps_{\rm sup}<\min\cbr[2]{\frac{c_m\varepsilon_0}{32},\frac{c_0}2},\\
		&\{y:{\mathsf p}_\delta(y\mid x)=v\}\subset\bigcup_{v'\in\{v-2\delta,v,v+2\delta\}}\{y:{\mathsf p}(y\mid x)=v'\}\qquad(|v-t_\star|\le\varepsilon_0/2),
	\end{align*}
and the sets on the right are Lebesgue-null by \eqref{eq:margin-upper}. By \cref{lem:hpd-sets}(ii) with $u=\hat\tau$, for every value of the tie-breaking variable,
	\begin{align*}
		\hC=\{y:{\mathsf p}_\delta(y\mid x)\ge s\},\qquad\hP_x(\hC)=\hat\tau,\qquad s:=v_x(\hat\tau).
	\end{align*}
If $s>t_\star+\delta$, then $\hC\subset\{y:{\mathsf p}_\delta(y\mid x)>t_\star+\delta\}$, and if $s<t_\star-\delta$, then $\hC\supset\{y:{\mathsf p}_\delta(y\mid x)\ge t_\star-\delta\}$. Since $|\hP_x(A)-{\mathsf P}_x(A)|\le\TV_x$ for every set $A\subset\cY$, the bounds on ${\mathsf P}_x\{{\mathsf p}_\delta(Y\mid x)>t_\star+\delta\}$ and ${\mathsf P}_x\{{\mathsf p}_\delta(Y\mid x)\ge t_\star-\delta\}$, together with $\TV_x+\Delta<c_m\delta/2$, give
	\begin{align*}
		s>t_\star+\delta&\implies\hat\tau\le\gamma_\star-\frac{c_m}2\delta+\TV_x<\gamma_\star-\Delta,\\
		s<t_\star-\delta&\implies\hat\tau\ge\gamma_\star+\frac{c_m}2\delta-\TV_x>\gamma_\star+\Delta .
	\end{align*}
Both conclusions contradict $|\hat\tau-\gamma_\star|\le\Delta$. Hence $|s-t_\star|\le\delta$, and Step 2 applies to $\hC=\{y:{\mathsf p}_\delta(y\mid x)\ge s\}$.
\end{proof}

\section{Details of the numerical experiments}\label{app:numerics}
The code and the results of all runs are available at \url{https://github.com/ilsangohn/conformal-hpd}.

\subsection{Models and samples}
We write $\Phi$ for the standard normal distribution function and $\operatorname{expit}(s)=1/(1+\e^{-s})$ for the logistic function. The models of \cref{sec:simulations} are specified as follows.
\begin{itemize}
\item \emph{Gaussian}: $\mathrm{ctr}(x)=(0.6\,\Phi(\mathrm{idx}_1(x))-0.3,0,\dots,0)$ and $\mathrm{sd}(x)=0.30+0.30\,\Phi(\mathrm{idx}_2(x))$, where $\mathrm{idx}(x)=\sqrt{12}\,\mathrm{proj}^\top(x_1-\tfrac12,x_2-\tfrac12,x_3-\tfrac12)^\top$, and $\mathrm{proj}\in\R^{3\times2}$ has orthonormal columns and rows $(-0.0023,0.7628)$, $(0.5163,-0.5528)$, $(0.8564,0.3354)$.
\item \emph{multimodal}: the centers $\mathrm{ctr}_k\in[-1.6,1.6]^{d_Y}$ lie at pairwise distances of at least $1.2$, and the weights $(\pi_1(x),\dots,\pi_4(x))$ are the softmax of $3B(x_1-\tfrac12,x_2-\tfrac12,x_3-\tfrac12)^\top$ for a matrix $B\in\R^{4\times3}$. The centers, the assignment of the scales to them and $B$ are drawn once for each $d_Y$, in this order, from the NumPy generator \texttt{default\_rng}$(2026+d_Y)$: the centers uniformly on $[-1.6,1.6]^{d_Y}$, each redrawn until it lies at a distance of at least $1.2$ from the previous ones; the scales as a random permutation of the geometric sequence of four values from $0.02$ to $0.5$; and $B$ with independent standard normal entries. For $d_Y=2$, the centers $(-0.36,0.30)$, $(0.76,1.55)$, $(1.10,-0.32)$ and $(-0.35,-1.36)$ have scales $0.17$, $0.5$, $0.02$ and $0.059$.
\item \emph{manifold}: $\mathsf T$ is the distribution of $(0.8\cos\Theta,0.8\sin\Theta,0,\dots,0)+0.02\,Z$, with $\Theta\sim U[0,2\pi)$ and $Z\sim\mathtt N(0,1)$ independent, and $\mathrm{wt}(x)=\operatorname{expit}((x_1-0.5)/0.15)$.
\end{itemize}
The conditional distributions are truncated to $\cY$. The truncated density ${\mathsf p}(\cdot\mid x)/{\mathsf P}_x(\cY)$ is available in closed form, and all draws from it are obtained by rejection sampling.

Each replication draws $4000$ points for fitting and for the auxiliary step, and a calibration sample of $1000$ points. The last $750$ of the first $3000$ points are the tuning points of every density estimate. For every candidate value of its tuning parameters, an estimate is fitted on its fitting sample without the tuning points. The value with the smallest held-out loss on the tuning points is selected, and the estimate is refitted on its whole fitting sample, of $4000$ or $3000$ points, with this value. For the network and the LP, the held-out loss is the contrastive log-loss on $3000$ random contrastive rows of the tuning points.

\subsection{Evaluation}
Conditional coverage and shape loss are computed at twenty covariate points drawn once with a fixed seed, ten in the hard region and ten in the easy region. At each of them, we draw $L=20000$ responses from the truncated conditional distribution and pool them with $L$ uniform draws on $\cY$, shared by all covariates. The pooled draws form an importance sample from $q_x:=\frac12{\mathsf P}_x+\frac12\mathrm{Unif}(\cY)$. The coverage of a set is the fraction of the conditional draws that it contains, and the Lebesgue measure of a set $A$ is the average of $\ind_A/q_x$ over the pooled draws. The oracle level $t_{\cov}$ at coverage $\cov$ is the $(1-\cov)$-quantile of ${\mathsf p}(Y\mid x)$ over the conditional draws, and the shape loss is computed from its integral form,
	\begin{align*}
		\cS^\circ(\hC(x))=\frac1{t_{\cov}}\int_{\hC(x)\triangle\Cstar_{\cov}(x)}\abs[0]{{\mathsf p}(y\mid x)-t_{\cov}}\,\d y,
	\end{align*}
using the same importance sample, with $\cov$ the coverage that $\hC(x)$ attains, whatever its distance from $1-\alpha$. Shape loss is reported relative to $\Leb(\Cstar(x))$, averaged with equal weights over the twenty covariates. The ``region-averaged'' coverage error is the equal-weight average of the mean errors over the points in the two regions. Marginal coverage is the average conditional coverage over a further $100$ draws $X_1,\dots,X_{100}$ of $X$, each with $L=5000$. Over the same draws, the volume ratio of \cref{tab:main-gaussian,tab:main-multimodal,tab:main-manifold} is $\sum_{j\le100}\Leb(\hC(X_j))/\sum_{j\le100}\Leb(\Cstar(X_j))$. All methods use the same draws. The sets of one-stage calibration with the true density are evaluated in the same way, and their coverage errors measure the Monte Carlo and calibration noise. The coverages of \cref{fig:cov-curves} are computed in the same way at forty equally spaced covariates on the segment, with the density fits and calibrations of the main experiment. For \cref{fig:residual}, the HPD values are computed on a $60\times60$ grid of $\cY$, and the coverages at eighty covariates on the segment.

The HPD values of one-stage calibration are computed by self-normalized importance sampling. We draw $M$ proposal points $U_1,\dots,U_M$, shared by all covariates, from $q:=0.9\,q_{\rm s}+0.1\,\mathrm{Unif}(\cY)$, truncated to $\cY$. Here $q_{\rm s}$ is a Gaussian kernel smoothing of the responses of the first $3000$ points, in which each draw uses a bandwidth chosen uniformly from $\{0.05,0.15,0.35\}$, so that the proposal also covers the components of scale $0.02$. At a covariate $x$, the HPD value of $y$ is
	\begin{align*}
		\frac{\sum_{m\le M}w_m\ind\{\hp(U_m\mid x)>\hp(y\mid x)\}}{\sum_{m\le M}w_m},\qquad w_m:=\frac{\hp(U_m\mid x)}{q(U_m)}.
	\end{align*}
We use $M=16000$ for $d_Y=2$ and $M=40000$ for $d_Y=5$.

\subsection{Estimators}
\paragraph{Density estimators} All density estimators except the partition estimate use the contrastive device of \cref{sec:cov-nn} with $K=5$ references from ${\mathsf q}_0:=0.9\,\mathtt N_{\cY}(\mathrm m,\mathrm{diag}(\mathrm s^2))+0.1\,\mathrm{Unif}(\cY)$. Here $\mathrm m$ and $\mathrm s$ are the coordinatewise mean and $1.25$ times the coordinatewise standard deviation of the responses of the first $3000$ points, and $\mathtt N_{\cY}$ is the normal distribution truncated to $\cY$. The reference is thus part of the density fit and uses neither the auxiliary nor the calibration sample. The estimate of $\eta$ is clipped to $[0,0.9999]$ for the network and to $[0,0.995]$ for the LP before it is mapped to a density. Unlike the estimates of the theory, the network and the LP omit the perturbation $n^{-2}\mathrm{pol}_J$, and the LP uses the product Epanechnikov kernel, which is not differentiable at the edge of its support.
\begin{itemize}
\item \emph{Network (NN).} A multilayer perceptron with three hidden layers of width $128$ and SiLU activations, applied to standardized inputs $(x,y)$, trained on the cross-entropy loss with Adam (learning rate $2\times10^{-3}$, batches of $1024$ rows). The number of epochs, at most $200$, is the one with the smallest loss on the tuning rows under early stopping with patience $8$. The network is then retrained from the same initialization on all contrastive rows of its fitting sample for that many epochs.
\item \emph{Local polynomial (LP).} Local linear regression on the pooled contrastive rows with the product Epanechnikov kernel and common bandwidths $h^{(x)}$ in the coordinates of $x$ and $h^{(y)}$ in those of $y$. The pair is chosen by contrastive log-loss on the tuning rows from $h^{(x)}\in\{0.6,0.9,1.35,2.0\}$ and $h^{(y)}\in\{0.2,0.3,0.45,0.65\}$, with the further values $h^{(x)}=3.0$ and $h^{(y)}=0.95$ for $d_Y=5$, where the selections reach the upper ends of these grids. Where fewer than $3(d+1)$ rows fall in the kernel window, the local constant fit is used.
\item \emph{Partition estimate.} The estimate of \cref{prop:partition} on the cubes of side $1/m$ of a regular grid of $[0,1]^3$, with a Gaussian kernel in $y$ of bandwidth $h_Y$, normalized on $\cY$ without the division at each $y$ by the kernel mass of $\cY$. A cube without fitting points gets the uniform density. The pair $(m,h_Y)$ is chosen from $\{2,3,4,5\}\times\{0.1,0.15,0.2,0.3,0.45\}$, with $h_Y=0.65$ added for $d_Y=5$, by held-out conditional log-likelihood on the tuning points.
\end{itemize}

\paragraph{Two-stage calibration} The weights of \cref{alg:local} are the local linear weights with the product triweight kernel, and a window with fewer than $8$ auxiliary points gets zero weights. The bandwidth $h$ is chosen from $\{0.3,0.4,0.55,0.75,1.0,1.5\}$ by five-fold cross-validation over the auxiliary points. The selected weights therefore depend on the responses of the auxiliary points, a case that the theory of \cref{sec:cov-theory} does not cover. For a held-out point $(X_j,Y_j)$, let $\overline Q^{\,h}_j$ be the function $\overline Q_{X_j}$ of \cref{alg:local} computed with bandwidth $h$ from the auxiliary points of the other four folds, and let $t_j^{\rm lo}:=\inf\{t:\overline Q^{\,1.0}_j(t)\le0.97\}$ and $t_j^{\rm hi}:=\inf\{t:\overline Q^{\,1.0}_j(t)\le0.80\}$. The loss of $h$ at this point is
	\begin{align*}
		\frac1{t_j^{\rm hi}-t_j^{\rm lo}}\int_{t_j^{\rm lo}}^{t_j^{\rm hi}}\bigl(\ind\{\hat a(Y_j\mid X_j)\ge t\}-\overline Q^{\,h}_j(t)\bigr)^2\,\d t,
	\end{align*}
and the selected bandwidth minimizes the average loss over the held-out points with $t_j^{\rm lo}<t_j^{\rm hi}$. The pilot bandwidth $1.0$ fixes the range of levels, so that all bandwidths are compared on the same levels. Given the bandwidth, the score at a covariate is computed as in the paragraph \emph{Computation} of \cref{sec:cov-procedure}. Two-stage (own covariate) uses the same weights and grid, with the bandwidth chosen by five-fold cross-validation of the squared error of the local linear regression of $\ind\{\hat{\mathsf s}^{\rm hpd}(X_i,Y_i)\le\hat q\}$ on $X_i$ over the auxiliary points, where $\hat q$ is the $(1-\alpha)$-quantile of these HPD values. The predictions are clipped to $[0,1]$. The rectified scores of RCP and of two-stage (quantile) use the $k$ nearest auxiliary points, with $k\in\{25,50,100,200,400,600\}$ chosen by five-fold cross-validation of the pinball loss.

\paragraph{Baselines} CHCDS uses the network density. Let $\hat a_{\rm N}(\cdot\mid x)$ be the unnormalized density implied by the network, and let $\widehat Z(x)$ be its integral over $\cY$, computed with the importance sample of the HPD values. Let $\hp_{\rm N}=\hat a_{\rm N}/\widehat Z$ be the normalized estimate, and let $\hat t^{\rm raw}_{\rm N}(x)$ be the level of $\hat a_{\rm N}(\cdot\mid x)$ whose super-level set has coverage $1-\alpha$ under $\hp_{\rm N}(\cdot\mid x)$. The score of CHCDS is $\{\hat t^{\rm raw}_{\rm N}(x)-\hat a_{\rm N}(y\mid x)\}/\widehat Z(x)$, the difference between the level $\hat t^{\rm raw}_{\rm N}(x)/\widehat Z(x)$ of $\hp_{\rm N}(\cdot\mid x)$ and the density $\hp_{\rm N}(y\mid x)$. Probabilistic conformal prediction (PCP) draws $50$ points from the network density, by sampling-importance resampling from the proposal points, and uses the distance to the nearest draw as the score. The local ellipsoid \citep{johnstone2021conformal,messoudi2022ellipsoidal} uses a mean network (two hidden layers of width $64$, SiLU, Adam with learning rate $2\times10^{-3}$, $60$ epochs) and the covariance of the residuals of the $150$ nearest training points in $x$, with a ridge of $10^{-4}$, and its score is the Mahalanobis distance. The rectified versions (RCP) refit the base score without the auxiliary points and divide it by its local $0.9$-quantile, estimated from the scores of the $k$ nearest auxiliary points in $x$ as described in the paragraph \emph{Two-stage calibration}. OT-CP, used only for $d_Y=2$, fits the mean network without the auxiliary points and matches the $1000$ auxiliary residuals to $1000$ uniform draws on the unit disk by optimal assignment. The score of a residual is the norm of the average of the matched points of its five nearest auxiliary residuals. Of the two true-density references, the first is one-stage calibration with the true density and the second uses the score $-{\mathsf p}(y\mid x)$ with a single threshold calibrated on the calibration sample.

\paragraph{Computation} The HPD values of one-stage calibration require an integral over $\cY$ at each calibration and evaluation covariate, which here amounts to $M$ density evaluations at the proposal points. \Cref{alg:local} evaluates the density estimate only at the auxiliary points with nonzero weights and computes no integral. Its bandwidth rule evaluates $\hat a(Y_i\mid X_j)$ at all pairs of auxiliary points, since the largest window covers $[0,1]^3$.

\subsection{Results of all methods}\label{app:all-results}
\Cref{tab:main-gaussian,tab:main-multimodal,tab:main-manifold} report the results of all methods.

\begin{table}[!t]
\caption{Gaussian model: marginal coverage, region-averaged conditional coverage error ($\times100$), shape loss $\cS^\circ$ relative to $\Leb(\Cstar)$, averaged over the evaluation covariates, and the volume ratio $\sum_j\Leb(\hC(X_j))/\sum_j\Leb(\Cstar(X_j))$ over the random covariates $X_j$ of the marginal coverage; mean $\pm$ standard error over twenty replications.}\label{tab:main-gaussian}
\centering\footnotesize\setlength{\tabcolsep}{3pt}
\resizebox{\ifdim\width>\textwidth\textwidth\else\width\fi}{!}{%
\begin{tabular}{llcccc}
\toprule
\multicolumn{2}{l}{method} & marg.\ coverage & cov.\ error & shape loss & volume ratio\\
\midrule
\multicolumn{6}{l}{$d_Y=2$}\\
\addlinespace[2pt]
\textit{true density} & \textit{one-stage} & $0.896${\scriptsize${}\pm0.003$} & $1.05${\scriptsize${}\pm0.17$} & $2.52\times10^{-8}${\scriptsize${}\pm4.8\times10^{-9}$} & $0.986${\scriptsize${}\pm0.012$}\\
 & \textit{global threshold} & $0.892${\scriptsize${}\pm0.003$} & $3.40${\scriptsize${}\pm0.10$} & $3.41\times10^{-8}${\scriptsize${}\pm6.8\times10^{-9}$} & $0.938${\scriptsize${}\pm0.011$}\\
\addlinespace[2pt]
NN & one-stage & $0.897${\scriptsize${}\pm0.003$} & $1.62${\scriptsize${}\pm0.14$} & $0.00669${\scriptsize${}\pm4.8\times10^{-4}$} & $0.987${\scriptsize${}\pm0.011$}\\
 & \cref{alg:local} & $0.896${\scriptsize${}\pm0.003$} & $2.38${\scriptsize${}\pm0.13$} & $0.00697${\scriptsize${}\pm4.7\times10^{-4}$} & $0.993${\scriptsize${}\pm0.012$}\\
 & two-stage (own covariate) & $0.896${\scriptsize${}\pm0.002$} & $2.01${\scriptsize${}\pm0.14$} & $0.00690${\scriptsize${}\pm4.7\times10^{-4}$} & $0.998${\scriptsize${}\pm0.010$}\\
 & two-stage (quantile) & $0.896${\scriptsize${}\pm0.003$} & $1.87${\scriptsize${}\pm0.15$} & $0.00694${\scriptsize${}\pm4.8\times10^{-4}$} & $0.989${\scriptsize${}\pm0.011$}\\
\addlinespace[2pt]
LP & one-stage & $0.896${\scriptsize${}\pm0.002$} & $2.21${\scriptsize${}\pm0.10$} & $0.00980${\scriptsize${}\pm4.1\times10^{-4}$} & $0.974${\scriptsize${}\pm0.0100$}\\
 & \cref{alg:local} & $0.897${\scriptsize${}\pm0.002$} & $2.29${\scriptsize${}\pm0.11$} & $0.0117${\scriptsize${}\pm8.4\times10^{-4}$} & $1.01${\scriptsize${}\pm0.011$}\\
 & two-stage (own covariate) & $0.898${\scriptsize${}\pm0.002$} & $2.03${\scriptsize${}\pm0.11$} & $0.0118${\scriptsize${}\pm8.6\times10^{-4}$} & $1.02${\scriptsize${}\pm0.0086$}\\
 & two-stage (quantile) & $0.897${\scriptsize${}\pm0.002$} & $1.83${\scriptsize${}\pm0.13$} & $0.0113${\scriptsize${}\pm7.6\times10^{-4}$} & $0.997${\scriptsize${}\pm0.0094$}\\
\addlinespace[2pt]
partition & one-stage & $0.899${\scriptsize${}\pm0.002$} & $3.65${\scriptsize${}\pm0.08$} & $0.0297${\scriptsize${}\pm0.0015$} & $1.02${\scriptsize${}\pm0.0097$}\\
 & \cref{alg:local} & $0.897${\scriptsize${}\pm0.002$} & $2.29${\scriptsize${}\pm0.12$} & $0.0325${\scriptsize${}\pm0.0012$} & $1.03${\scriptsize${}\pm0.010$}\\
 & two-stage (own covariate) & $0.900${\scriptsize${}\pm0.002$} & $2.84${\scriptsize${}\pm0.10$} & $0.0340${\scriptsize${}\pm0.0013$} & $1.07${\scriptsize${}\pm0.011$}\\
 & two-stage (quantile) & $0.899${\scriptsize${}\pm0.002$} & $3.17${\scriptsize${}\pm0.12$} & $0.0332${\scriptsize${}\pm0.0013$} & $1.05${\scriptsize${}\pm0.015$}\\
\addlinespace[2pt]
baselines & CHCDS & $0.896${\scriptsize${}\pm0.003$} & $1.61${\scriptsize${}\pm0.13$} & $0.00665${\scriptsize${}\pm4.7\times10^{-4}$} & $0.987${\scriptsize${}\pm0.011$}\\
 & local ellipsoid & $0.897${\scriptsize${}\pm0.003$} & $1.94${\scriptsize${}\pm0.11$} & $0.0115${\scriptsize${}\pm6.4\times10^{-4}$} & $0.999${\scriptsize${}\pm0.014$}\\
 & RCP (ellipsoid score) & $0.899${\scriptsize${}\pm0.003$} & $1.97${\scriptsize${}\pm0.09$} & $0.0126${\scriptsize${}\pm8.1\times10^{-4}$} & $1.01${\scriptsize${}\pm0.012$}\\
 & PCP & $0.894${\scriptsize${}\pm0.002$} & $4.17${\scriptsize${}\pm0.15$} & $0.283${\scriptsize${}\pm0.0066$} & $1.23${\scriptsize${}\pm0.012$}\\
 & RCP (PCP score) & $0.898${\scriptsize${}\pm0.002$} & $2.88${\scriptsize${}\pm0.17$} & $0.280${\scriptsize${}\pm0.0053$} & $1.28${\scriptsize${}\pm0.014$}\\
 & OT-CP & $0.889${\scriptsize${}\pm0.002$} & $6.22${\scriptsize${}\pm0.09$} & $0.763${\scriptsize${}\pm0.14$} & $1.61${\scriptsize${}\pm0.12$}\\
\midrule
\multicolumn{6}{l}{$d_Y=5$}\\
\addlinespace[2pt]
\textit{true density} & \textit{one-stage} & $0.903${\scriptsize${}\pm0.002$} & $0.68${\scriptsize${}\pm0.10$} & $1.26\times10^{-8}${\scriptsize${}\pm3.1\times10^{-9}$} & $1.02${\scriptsize${}\pm0.014$}\\
 & \textit{global threshold} & $0.894${\scriptsize${}\pm0.002$} & $5.33${\scriptsize${}\pm0.10$} & $2.43\times10^{-7}${\scriptsize${}\pm10.0\times10^{-8}$} & $0.781${\scriptsize${}\pm0.015$}\\
\addlinespace[2pt]
NN & one-stage & $0.901${\scriptsize${}\pm0.002$} & $1.67${\scriptsize${}\pm0.08$} & $0.0177${\scriptsize${}\pm8.6\times10^{-4}$} & $0.981${\scriptsize${}\pm0.013$}\\
 & \cref{alg:local} & $0.898${\scriptsize${}\pm0.002$} & $2.95${\scriptsize${}\pm0.13$} & $0.0197${\scriptsize${}\pm0.0012$} & $0.950${\scriptsize${}\pm0.017$}\\
 & two-stage (own covariate) & $0.900${\scriptsize${}\pm0.002$} & $1.89${\scriptsize${}\pm0.10$} & $0.0197${\scriptsize${}\pm0.0011$} & $1.04${\scriptsize${}\pm0.022$}\\
 & two-stage (quantile) & $0.900${\scriptsize${}\pm0.002$} & $1.96${\scriptsize${}\pm0.13$} & $0.0196${\scriptsize${}\pm0.0012$} & $1.00${\scriptsize${}\pm0.025$}\\
\addlinespace[2pt]
LP & one-stage & $0.899${\scriptsize${}\pm0.003$} & $3.78${\scriptsize${}\pm0.25$} & $0.147${\scriptsize${}\pm0.026$} & $1.00${\scriptsize${}\pm0.025$}\\
 & \cref{alg:local} & $0.901${\scriptsize${}\pm0.002$} & $2.83${\scriptsize${}\pm0.12$} & $0.117${\scriptsize${}\pm0.0086$} & $1.23${\scriptsize${}\pm0.051$}\\
 & two-stage (own covariate) & $0.903${\scriptsize${}\pm0.002$} & $2.31${\scriptsize${}\pm0.10$} & $0.171${\scriptsize${}\pm0.040$} & $1.47${\scriptsize${}\pm0.16$}\\
 & two-stage (quantile) & $0.900${\scriptsize${}\pm0.003$} & $2.56${\scriptsize${}\pm0.23$} & $0.252${\scriptsize${}\pm0.13$} & $1.53${\scriptsize${}\pm0.34$}\\
\addlinespace[2pt]
partition & one-stage & $0.900${\scriptsize${}\pm0.002$} & $5.99${\scriptsize${}\pm0.09$} & $0.164${\scriptsize${}\pm0.0051$} & $0.991${\scriptsize${}\pm0.021$}\\
 & \cref{alg:local} & $0.900${\scriptsize${}\pm0.002$} & $2.77${\scriptsize${}\pm0.13$} & $0.181${\scriptsize${}\pm0.0057$} & $1.16${\scriptsize${}\pm0.024$}\\
 & two-stage (own covariate) & $0.905${\scriptsize${}\pm0.002$} & $3.08${\scriptsize${}\pm0.12$} & $0.208${\scriptsize${}\pm0.0067$} & $1.31${\scriptsize${}\pm0.031$}\\
 & two-stage (quantile) & $0.901${\scriptsize${}\pm0.002$} & $3.39${\scriptsize${}\pm0.14$} & $0.204${\scriptsize${}\pm0.0080$} & $1.36${\scriptsize${}\pm0.12$}\\
\addlinespace[2pt]
baselines & CHCDS & $0.901${\scriptsize${}\pm0.002$} & $1.60${\scriptsize${}\pm0.10$} & $0.0174${\scriptsize${}\pm8.2\times10^{-4}$} & $1.02${\scriptsize${}\pm0.026$}\\
 & local ellipsoid & $0.901${\scriptsize${}\pm0.002$} & $1.94${\scriptsize${}\pm0.10$} & $0.0758${\scriptsize${}\pm0.0026$} & $1.04${\scriptsize${}\pm0.022$}\\
 & RCP (ellipsoid score) & $0.901${\scriptsize${}\pm0.002$} & $2.08${\scriptsize${}\pm0.12$} & $0.0799${\scriptsize${}\pm0.0027$} & $1.07${\scriptsize${}\pm0.022$}\\
 & PCP & $0.891${\scriptsize${}\pm0.002$} & $7.20${\scriptsize${}\pm0.13$} & $1.23${\scriptsize${}\pm0.023$} & $1.51${\scriptsize${}\pm0.023$}\\
 & RCP (PCP score) & $0.900${\scriptsize${}\pm0.003$} & $3.25${\scriptsize${}\pm0.17$} & $0.904${\scriptsize${}\pm0.021$} & $1.77${\scriptsize${}\pm0.052$}\\
\bottomrule
\end{tabular}}
\end{table}

\begin{table}[!t]
\caption{Multimodal model: marginal coverage, region-averaged conditional coverage error ($\times100$), shape loss $\cS^\circ$ relative to $\Leb(\Cstar)$, averaged over the evaluation covariates, and the volume ratio $\sum_j\Leb(\hC(X_j))/\sum_j\Leb(\Cstar(X_j))$ over the random covariates $X_j$ of the marginal coverage; mean $\pm$ standard error over twenty replications.}\label{tab:main-multimodal}
\centering\footnotesize\setlength{\tabcolsep}{3pt}
\resizebox{\ifdim\width>\textwidth\textwidth\else\width\fi}{!}{%
\begin{tabular}{llcccc}
\toprule
\multicolumn{2}{l}{method} & marg.\ coverage & cov.\ error & shape loss & volume ratio\\
\midrule
\multicolumn{6}{l}{$d_Y=2$}\\
\addlinespace[2pt]
\textit{true density} & \textit{one-stage} & $0.899${\scriptsize${}\pm0.002$} & $0.77${\scriptsize${}\pm0.12$} & $4.14\times10^{-8}${\scriptsize${}\pm6.9\times10^{-9}$} & $0.994${\scriptsize${}\pm0.015$}\\
 & \textit{global threshold} & $0.896${\scriptsize${}\pm0.002$} & $4.79${\scriptsize${}\pm0.19$} & $4.20\times10^{-6}${\scriptsize${}\pm1.0\times10^{-6}$} & $0.813${\scriptsize${}\pm0.018$}\\
\addlinespace[2pt]
NN & one-stage & $0.898${\scriptsize${}\pm0.002$} & $2.76${\scriptsize${}\pm0.19$} & $0.180${\scriptsize${}\pm0.017$} & $1.15${\scriptsize${}\pm0.028$}\\
 & \cref{alg:local} & $0.903${\scriptsize${}\pm0.002$} & $3.18${\scriptsize${}\pm0.15$} & $0.246${\scriptsize${}\pm0.029$} & $1.21${\scriptsize${}\pm0.024$}\\
 & two-stage (own covariate) & $0.897${\scriptsize${}\pm0.002$} & $2.08${\scriptsize${}\pm0.10$} & $0.169${\scriptsize${}\pm0.010$} & $1.14${\scriptsize${}\pm0.021$}\\
 & two-stage (quantile) & $0.898${\scriptsize${}\pm0.002$} & $2.20${\scriptsize${}\pm0.13$} & $0.177${\scriptsize${}\pm0.011$} & $1.13${\scriptsize${}\pm0.020$}\\
\addlinespace[2pt]
LP & one-stage & $0.903${\scriptsize${}\pm0.003$} & $5.28${\scriptsize${}\pm0.11$} & $0.980${\scriptsize${}\pm0.035$} & $1.28${\scriptsize${}\pm0.030$}\\
 & \cref{alg:local} & $0.902${\scriptsize${}\pm0.002$} & $3.03${\scriptsize${}\pm0.15$} & $0.310${\scriptsize${}\pm0.023$} & $1.33${\scriptsize${}\pm0.029$}\\
 & two-stage (own covariate) & $0.904${\scriptsize${}\pm0.002$} & $3.25${\scriptsize${}\pm0.14$} & $0.521${\scriptsize${}\pm0.038$} & $1.27${\scriptsize${}\pm0.023$}\\
 & two-stage (quantile) & $0.906${\scriptsize${}\pm0.002$} & $4.40${\scriptsize${}\pm0.14$} & $0.623${\scriptsize${}\pm0.071$} & $1.34${\scriptsize${}\pm0.037$}\\
\addlinespace[2pt]
partition & one-stage & $0.903${\scriptsize${}\pm0.002$} & $5.86${\scriptsize${}\pm0.12$} & $1.17${\scriptsize${}\pm0.042$} & $1.34${\scriptsize${}\pm0.024$}\\
 & \cref{alg:local} & $0.903${\scriptsize${}\pm0.002$} & $2.93${\scriptsize${}\pm0.11$} & $0.772${\scriptsize${}\pm0.054$} & $1.55${\scriptsize${}\pm0.026$}\\
 & two-stage (own covariate) & $0.900${\scriptsize${}\pm0.002$} & $3.78${\scriptsize${}\pm0.11$} & $0.584${\scriptsize${}\pm0.024$} & $1.51${\scriptsize${}\pm0.027$}\\
 & two-stage (quantile) & $0.899${\scriptsize${}\pm0.002$} & $4.42${\scriptsize${}\pm0.19$} & $0.845${\scriptsize${}\pm0.13$} & $1.78${\scriptsize${}\pm0.15$}\\
\addlinespace[2pt]
baselines & CHCDS & $0.898${\scriptsize${}\pm0.002$} & $3.95${\scriptsize${}\pm0.39$} & $0.180${\scriptsize${}\pm0.020$} & $1.36${\scriptsize${}\pm0.11$}\\
 & local ellipsoid & $0.892${\scriptsize${}\pm0.002$} & $5.37${\scriptsize${}\pm0.14$} & $11.7${\scriptsize${}\pm0.28$} & $3.89${\scriptsize${}\pm0.056$}\\
 & RCP (ellipsoid score) & $0.894${\scriptsize${}\pm0.003$} & $5.58${\scriptsize${}\pm0.20$} & $11.0${\scriptsize${}\pm0.36$} & $4.01${\scriptsize${}\pm0.053$}\\
 & PCP & $0.897${\scriptsize${}\pm0.002$} & $5.43${\scriptsize${}\pm0.18$} & $2.90${\scriptsize${}\pm0.20$} & $1.76${\scriptsize${}\pm0.074$}\\
 & RCP (PCP score) & $0.896${\scriptsize${}\pm0.002$} & $3.37${\scriptsize${}\pm0.16$} & $1.58${\scriptsize${}\pm0.46$} & $1.75${\scriptsize${}\pm0.044$}\\
 & OT-CP & $0.901${\scriptsize${}\pm0.002$} & $6.71${\scriptsize${}\pm0.21$} & $27.5${\scriptsize${}\pm1.0$} & $5.51${\scriptsize${}\pm0.14$}\\
\midrule
\multicolumn{6}{l}{$d_Y=5$}\\
\addlinespace[2pt]
\textit{true density} & \textit{one-stage} & $0.899${\scriptsize${}\pm0.002$} & $0.97${\scriptsize${}\pm0.12$} & $8.72\times10^{-8}${\scriptsize${}\pm5.1\times10^{-8}$} & $0.994${\scriptsize${}\pm0.028$}\\
 & \textit{global threshold} & $0.897${\scriptsize${}\pm0.002$} & $4.86${\scriptsize${}\pm0.12$} & $2.80\times10^{-4}${\scriptsize${}\pm1.8\times10^{-4}$} & $0.636${\scriptsize${}\pm0.023$}\\
\addlinespace[2pt]
NN & one-stage & $0.892${\scriptsize${}\pm0.003$} & $5.23${\scriptsize${}\pm0.37$} & $4861${\scriptsize${}\pm722$} & $1.32${\scriptsize${}\pm0.079$}\\
 & \cref{alg:local} & $0.900${\scriptsize${}\pm0.002$} & $2.50${\scriptsize${}\pm0.09$} & $1248${\scriptsize${}\pm415$} & $1.45${\scriptsize${}\pm0.044$}\\
 & two-stage (own covariate) & $0.901${\scriptsize${}\pm0.002$} & $3.68${\scriptsize${}\pm0.27$} & $1708${\scriptsize${}\pm658$} & $1.52${\scriptsize${}\pm0.055$}\\
 & two-stage (quantile) & $0.901${\scriptsize${}\pm0.002$} & $3.93${\scriptsize${}\pm0.25$} & $2112${\scriptsize${}\pm517$} & $3.66${\scriptsize${}\pm0.86$}\\
\addlinespace[2pt]
LP & one-stage & $0.898${\scriptsize${}\pm0.002$} & $7.06${\scriptsize${}\pm0.13$} & $3.83\times10^{4}${\scriptsize${}\pm4779$} & $1.84${\scriptsize${}\pm0.072$}\\
 & \cref{alg:local} & $0.900${\scriptsize${}\pm0.002$} & $2.77${\scriptsize${}\pm0.11$} & $660${\scriptsize${}\pm190$} & $1.64${\scriptsize${}\pm0.054$}\\
 & two-stage (own covariate) & $0.903${\scriptsize${}\pm0.003$} & $3.50${\scriptsize${}\pm0.16$} & $665${\scriptsize${}\pm194$} & $1.79${\scriptsize${}\pm0.077$}\\
 & two-stage (quantile) & $0.900${\scriptsize${}\pm0.003$} & $4.03${\scriptsize${}\pm0.21$} & $572${\scriptsize${}\pm170$} & $4.25${\scriptsize${}\pm1.6$}\\
\addlinespace[2pt]
partition & one-stage & $0.893${\scriptsize${}\pm0.002$} & $5.99${\scriptsize${}\pm0.10$} & $1.42\times10^{4}${\scriptsize${}\pm1143$} & $1.53${\scriptsize${}\pm0.059$}\\
 & \cref{alg:local} & $0.898${\scriptsize${}\pm0.002$} & $2.64${\scriptsize${}\pm0.12$} & $36.3${\scriptsize${}\pm14$} & $2.11${\scriptsize${}\pm0.049$}\\
 & two-stage (own covariate) & $0.899${\scriptsize${}\pm0.002$} & $3.23${\scriptsize${}\pm0.10$} & $40.6${\scriptsize${}\pm17$} & $2.22${\scriptsize${}\pm0.053$}\\
 & two-stage (quantile) & $0.892${\scriptsize${}\pm0.003$} & $4.01${\scriptsize${}\pm0.21$} & $72.8${\scriptsize${}\pm39$} & $26.2${\scriptsize${}\pm5.2$}\\
\addlinespace[2pt]
baselines & CHCDS & $0.895${\scriptsize${}\pm0.003$} & $10.08${\scriptsize${}\pm0.61$} & $4.60\times10^{4}${\scriptsize${}\pm2.3\times10^{4}$} & $121${\scriptsize${}\pm11$}\\
 & local ellipsoid & $0.898${\scriptsize${}\pm0.002$} & $3.43${\scriptsize${}\pm0.12$} & $2202${\scriptsize${}\pm382$} & $6.57${\scriptsize${}\pm0.16$}\\
 & RCP (ellipsoid score) & $0.899${\scriptsize${}\pm0.002$} & $3.47${\scriptsize${}\pm0.14$} & $1012${\scriptsize${}\pm240$} & $7.12${\scriptsize${}\pm0.26$}\\
 & PCP & $0.896${\scriptsize${}\pm0.003$} & $8.85${\scriptsize${}\pm0.62$} & $1.28\times10^{6}${\scriptsize${}\pm2.8\times10^{5}$} & $41.8${\scriptsize${}\pm8.1$}\\
 & RCP (PCP score) & $0.893${\scriptsize${}\pm0.002$} & $8.46${\scriptsize${}\pm0.60$} & $5.84\times10^{5}${\scriptsize${}\pm2.2\times10^{5}$} & $50.8${\scriptsize${}\pm6.5$}\\
\bottomrule
\end{tabular}}
\end{table}

\begin{table}[!t]
\caption{Manifold model: marginal coverage, region-averaged conditional coverage error ($\times100$), shape loss $\cS^\circ$ relative to $\Leb(\Cstar)$, averaged over the evaluation covariates, and the volume ratio $\sum_j\Leb(\hC(X_j))/\sum_j\Leb(\Cstar(X_j))$ over the random covariates $X_j$ of the marginal coverage; mean $\pm$ standard error over twenty replications.}\label{tab:main-manifold}
\centering\footnotesize\setlength{\tabcolsep}{3pt}
\resizebox{\ifdim\width>\textwidth\textwidth\else\width\fi}{!}{%
\begin{tabular}{llcccc}
\toprule
\multicolumn{2}{l}{method} & marg.\ coverage & cov.\ error & shape loss & volume ratio\\
\midrule
\multicolumn{6}{l}{$d_Y=2$}\\
\addlinespace[2pt]
\textit{true density} & \textit{one-stage} & $0.899${\scriptsize${}\pm0.002$} & $0.78${\scriptsize${}\pm0.13$} & $5.40\times10^{-8}${\scriptsize${}\pm1.2\times10^{-8}$} & $1.00${\scriptsize${}\pm0.012$}\\
 & \textit{global threshold} & $0.899${\scriptsize${}\pm0.002$} & $2.10${\scriptsize${}\pm0.10$} & $1.10\times10^{-6}${\scriptsize${}\pm2.9\times10^{-7}$} & $0.969${\scriptsize${}\pm0.014$}\\
\addlinespace[2pt]
NN & one-stage & $0.899${\scriptsize${}\pm0.003$} & $1.66${\scriptsize${}\pm0.15$} & $0.0738${\scriptsize${}\pm0.0048$} & $1.04${\scriptsize${}\pm0.020$}\\
 & \cref{alg:local} & $0.901${\scriptsize${}\pm0.003$} & $2.59${\scriptsize${}\pm0.08$} & $0.0870${\scriptsize${}\pm0.0061$} & $1.09${\scriptsize${}\pm0.015$}\\
 & two-stage (own covariate) & $0.898${\scriptsize${}\pm0.003$} & $2.02${\scriptsize${}\pm0.16$} & $0.0785${\scriptsize${}\pm0.0054$} & $1.06${\scriptsize${}\pm0.018$}\\
 & two-stage (quantile) & $0.899${\scriptsize${}\pm0.002$} & $2.11${\scriptsize${}\pm0.17$} & $0.0811${\scriptsize${}\pm0.0055$} & $1.06${\scriptsize${}\pm0.017$}\\
\addlinespace[2pt]
LP & one-stage & $0.903${\scriptsize${}\pm0.002$} & $3.30${\scriptsize${}\pm0.06$} & $0.548${\scriptsize${}\pm0.0072$} & $1.15${\scriptsize${}\pm0.014$}\\
 & \cref{alg:local} & $0.902${\scriptsize${}\pm0.003$} & $2.47${\scriptsize${}\pm0.11$} & $0.385${\scriptsize${}\pm0.018$} & $1.18${\scriptsize${}\pm0.016$}\\
 & two-stage (own covariate) & $0.901${\scriptsize${}\pm0.002$} & $2.06${\scriptsize${}\pm0.13$} & $0.355${\scriptsize${}\pm0.017$} & $1.18${\scriptsize${}\pm0.017$}\\
 & two-stage (quantile) & $0.901${\scriptsize${}\pm0.003$} & $2.72${\scriptsize${}\pm0.14$} & $0.386${\scriptsize${}\pm0.019$} & $1.23${\scriptsize${}\pm0.045$}\\
\addlinespace[2pt]
partition & one-stage & $0.898${\scriptsize${}\pm0.003$} & $3.92${\scriptsize${}\pm0.04$} & $0.807${\scriptsize${}\pm0.015$} & $1.27${\scriptsize${}\pm0.017$}\\
 & \cref{alg:local} & $0.904${\scriptsize${}\pm0.002$} & $2.47${\scriptsize${}\pm0.09$} & $0.554${\scriptsize${}\pm0.018$} & $1.26${\scriptsize${}\pm0.014$}\\
 & two-stage (own covariate) & $0.901${\scriptsize${}\pm0.002$} & $2.77${\scriptsize${}\pm0.09$} & $0.598${\scriptsize${}\pm0.016$} & $1.29${\scriptsize${}\pm0.021$}\\
 & two-stage (quantile) & $0.902${\scriptsize${}\pm0.002$} & $3.95${\scriptsize${}\pm0.23$} & $0.662${\scriptsize${}\pm0.040$} & $1.73${\scriptsize${}\pm0.17$}\\
\addlinespace[2pt]
baselines & CHCDS & $0.899${\scriptsize${}\pm0.003$} & $1.83${\scriptsize${}\pm0.17$} & $0.0701${\scriptsize${}\pm0.0047$} & $1.05${\scriptsize${}\pm0.023$}\\
 & local ellipsoid & $0.905${\scriptsize${}\pm0.002$} & $7.22${\scriptsize${}\pm0.14$} & $0.0277${\scriptsize${}\pm0.0025$} & $1.28${\scriptsize${}\pm0.015$}\\
 & RCP (ellipsoid score) & $0.899${\scriptsize${}\pm0.003$} & $4.42${\scriptsize${}\pm0.16$} & $0.999${\scriptsize${}\pm0.072$} & $1.28${\scriptsize${}\pm0.010$}\\
 & PCP & $0.901${\scriptsize${}\pm0.002$} & $4.42${\scriptsize${}\pm0.15$} & $1.44${\scriptsize${}\pm0.031$} & $1.60${\scriptsize${}\pm0.017$}\\
 & RCP (PCP score) & $0.900${\scriptsize${}\pm0.002$} & $3.40${\scriptsize${}\pm0.14$} & $1.37${\scriptsize${}\pm0.022$} & $1.56${\scriptsize${}\pm0.016$}\\
 & OT-CP & $0.904${\scriptsize${}\pm0.002$} & $5.53${\scriptsize${}\pm0.10$} & $1.31${\scriptsize${}\pm0.17$} & $1.51${\scriptsize${}\pm0.073$}\\
\midrule
\multicolumn{6}{l}{$d_Y=5$}\\
\addlinespace[2pt]
\textit{true density} & \textit{one-stage} & $0.899${\scriptsize${}\pm0.003$} & $1.00${\scriptsize${}\pm0.15$} & $6.32\times10^{-8}${\scriptsize${}\pm3.7\times10^{-8}$} & $1.00${\scriptsize${}\pm0.026$}\\
 & \textit{global threshold} & $0.899${\scriptsize${}\pm0.003$} & $3.38${\scriptsize${}\pm0.09$} & $3.95\times10^{-4}${\scriptsize${}\pm1.8\times10^{-4}$} & $0.905${\scriptsize${}\pm0.031$}\\
\addlinespace[2pt]
NN & one-stage & $0.899${\scriptsize${}\pm0.002$} & $4.90${\scriptsize${}\pm0.45$} & $9.04${\scriptsize${}\pm3.6$} & $1.60${\scriptsize${}\pm0.086$}\\
 & \cref{alg:local} & $0.904${\scriptsize${}\pm0.003$} & $2.61${\scriptsize${}\pm0.12$} & $3.61${\scriptsize${}\pm1.2$} & $1.18${\scriptsize${}\pm0.028$}\\
 & two-stage (own covariate) & $0.902${\scriptsize${}\pm0.003$} & $2.79${\scriptsize${}\pm0.17$} & $34.5${\scriptsize${}\pm15$} & $1.21${\scriptsize${}\pm0.032$}\\
 & two-stage (quantile) & $0.902${\scriptsize${}\pm0.003$} & $3.34${\scriptsize${}\pm0.18$} & $3.40\times10^{4}${\scriptsize${}\pm2.3\times10^{4}$} & $3.55${\scriptsize${}\pm1.1$}\\
\addlinespace[2pt]
LP & one-stage & $0.902${\scriptsize${}\pm0.003$} & $5.49${\scriptsize${}\pm0.08$} & $1.18\times10^{4}${\scriptsize${}\pm882$} & $1.94${\scriptsize${}\pm0.077$}\\
 & \cref{alg:local} & $0.900${\scriptsize${}\pm0.003$} & $2.52${\scriptsize${}\pm0.09$} & $8.80${\scriptsize${}\pm2.3$} & $2.34${\scriptsize${}\pm0.094$}\\
 & two-stage (own covariate) & $0.899${\scriptsize${}\pm0.003$} & $2.21${\scriptsize${}\pm0.10$} & $18.7${\scriptsize${}\pm3.8$} & $2.18${\scriptsize${}\pm0.13$}\\
 & two-stage (quantile) & $0.897${\scriptsize${}\pm0.002$} & $3.83${\scriptsize${}\pm0.21$} & $1158${\scriptsize${}\pm420$} & $13.3${\scriptsize${}\pm2.3$}\\
\addlinespace[2pt]
partition & one-stage & $0.896${\scriptsize${}\pm0.003$} & $4.12${\scriptsize${}\pm0.07$} & $3256${\scriptsize${}\pm110$} & $1.40${\scriptsize${}\pm0.037$}\\
 & \cref{alg:local} & $0.904${\scriptsize${}\pm0.003$} & $2.53${\scriptsize${}\pm0.10$} & $9.61${\scriptsize${}\pm1.8$} & $1.69${\scriptsize${}\pm0.047$}\\
 & two-stage (own covariate) & $0.900${\scriptsize${}\pm0.003$} & $3.15${\scriptsize${}\pm0.11$} & $35.2${\scriptsize${}\pm28$} & $1.72${\scriptsize${}\pm0.045$}\\
 & two-stage (quantile) & $0.898${\scriptsize${}\pm0.002$} & $4.26${\scriptsize${}\pm0.26$} & $477${\scriptsize${}\pm174$} & $15.8${\scriptsize${}\pm3.4$}\\
\addlinespace[2pt]
baselines & CHCDS & $0.899${\scriptsize${}\pm0.003$} & $10.49${\scriptsize${}\pm1.02$} & $1.60\times10^{4}${\scriptsize${}\pm1.6\times10^{4}$} & $59.1${\scriptsize${}\pm10$}\\
 & local ellipsoid & $0.897${\scriptsize${}\pm0.003$} & $2.78${\scriptsize${}\pm0.09$} & $422${\scriptsize${}\pm26$} & $1.18${\scriptsize${}\pm0.047$}\\
 & RCP (ellipsoid score) & $0.900${\scriptsize${}\pm0.003$} & $3.09${\scriptsize${}\pm0.16$} & $421${\scriptsize${}\pm208$} & $1.35${\scriptsize${}\pm0.061$}\\
 & PCP & $0.899${\scriptsize${}\pm0.004$} & $4.87${\scriptsize${}\pm0.21$} & $5797${\scriptsize${}\pm600$} & $2.03${\scriptsize${}\pm0.074$}\\
 & RCP (PCP score) & $0.898${\scriptsize${}\pm0.003$} & $3.77${\scriptsize${}\pm0.31$} & $1880${\scriptsize${}\pm313$} & $2.04${\scriptsize${}\pm0.073$}\\
\bottomrule
\end{tabular}}
\end{table}

\FloatBarrier
\subsection{The same density fit}\label{app:samefit}
\Cref{tab:samefit} compares one-stage calibration with the density fitted on all $4000$ fitting points to one-stage calibration and \cref{alg:local} with the same density fitted on $3000$ points.

\begin{table}[!t]
\caption{The same density fit with one-stage calibration and with \cref{alg:local}: region-averaged conditional coverage error ($\times100$) and shape loss $\cS^\circ$ relative to $\Leb(\Cstar)$ of one-stage calibration with the density fitted on all $4000$ fitting points (all), and of one-stage calibration (one) and \cref{alg:local} (two) with the density fitted on the $3000$ points of two-stage calibration; mean $\pm$ standard error over twenty replications.}\label{tab:samefit}
\centering\small
\resizebox{\ifdim\width>\textwidth\textwidth\else\width\fi}{!}{%
\begin{tabular}{lllcccccc}
\toprule
 & & & \multicolumn{3}{c}{cond.\ cov.\ error ($\times100$)} & \multicolumn{3}{c}{shape loss $\cS^\circ$}\\
\cmidrule(lr){4-6}\cmidrule(lr){7-9}
model & $d_Y$ & estimate & all & one & two & all & one & two\\
\midrule
Gaussian & $2$ & NN & $1.62${\scriptsize${}\pm0.14$} & $1.81${\scriptsize${}\pm0.13$} & $2.38${\scriptsize${}\pm0.13$} & $0.00669${\scriptsize${}\pm4.8\times10^{-4}$} & $0.00698${\scriptsize${}\pm4.6\times10^{-4}$} & $0.00697${\scriptsize${}\pm4.7\times10^{-4}$}\\
 &  & LP & $2.21${\scriptsize${}\pm0.10$} & $2.25${\scriptsize${}\pm0.12$} & $2.29${\scriptsize${}\pm0.11$} & $0.00980${\scriptsize${}\pm4.1\times10^{-4}$} & $0.0109${\scriptsize${}\pm7.1\times10^{-4}$} & $0.0117${\scriptsize${}\pm8.4\times10^{-4}$}\\
 &  & partition & $3.65${\scriptsize${}\pm0.08$} & $3.65${\scriptsize${}\pm0.09$} & $2.29${\scriptsize${}\pm0.12$} & $0.0297${\scriptsize${}\pm0.0015$} & $0.0322${\scriptsize${}\pm0.0011$} & $0.0325${\scriptsize${}\pm0.0012$}\\
\addlinespace[2pt]
 & $5$ & NN & $1.67${\scriptsize${}\pm0.08$} & $1.89${\scriptsize${}\pm0.11$} & $2.95${\scriptsize${}\pm0.13$} & $0.0177${\scriptsize${}\pm8.6\times10^{-4}$} & $0.0204${\scriptsize${}\pm0.0012$} & $0.0197${\scriptsize${}\pm0.0012$}\\
 &  & LP & $3.78${\scriptsize${}\pm0.25$} & $3.77${\scriptsize${}\pm0.18$} & $2.83${\scriptsize${}\pm0.12$} & $0.147${\scriptsize${}\pm0.026$} & $0.0882${\scriptsize${}\pm0.0068$} & $0.117${\scriptsize${}\pm0.0086$}\\
 &  & partition & $5.99${\scriptsize${}\pm0.09$} & $6.16${\scriptsize${}\pm0.09$} & $2.77${\scriptsize${}\pm0.13$} & $0.164${\scriptsize${}\pm0.0051$} & $0.214${\scriptsize${}\pm0.0075$} & $0.181${\scriptsize${}\pm0.0057$}\\
\addlinespace[2pt]
multimodal & $2$ & NN & $2.76${\scriptsize${}\pm0.19$} & $2.48${\scriptsize${}\pm0.14$} & $3.18${\scriptsize${}\pm0.15$} & $0.180${\scriptsize${}\pm0.017$} & $0.200${\scriptsize${}\pm0.013$} & $0.246${\scriptsize${}\pm0.029$}\\
 &  & LP & $5.28${\scriptsize${}\pm0.11$} & $5.43${\scriptsize${}\pm0.17$} & $3.03${\scriptsize${}\pm0.15$} & $0.980${\scriptsize${}\pm0.035$} & $1.11${\scriptsize${}\pm0.033$} & $0.310${\scriptsize${}\pm0.023$}\\
 &  & partition & $5.86${\scriptsize${}\pm0.12$} & $5.97${\scriptsize${}\pm0.10$} & $2.93${\scriptsize${}\pm0.11$} & $1.17${\scriptsize${}\pm0.042$} & $1.31${\scriptsize${}\pm0.061$} & $0.772${\scriptsize${}\pm0.054$}\\
\addlinespace[2pt]
 & $5$ & NN & $5.23${\scriptsize${}\pm0.37$} & $5.19${\scriptsize${}\pm0.31$} & $2.50${\scriptsize${}\pm0.09$} & $4861${\scriptsize${}\pm722$} & $1.22\times10^{4}${\scriptsize${}\pm3113$} & $1248${\scriptsize${}\pm415$}\\
 &  & LP & $7.06${\scriptsize${}\pm0.13$} & $6.92${\scriptsize${}\pm0.16$} & $2.77${\scriptsize${}\pm0.11$} & $3.83\times10^{4}${\scriptsize${}\pm4779$} & $4.64\times10^{4}${\scriptsize${}\pm6767$} & $660${\scriptsize${}\pm190$}\\
 &  & partition & $5.99${\scriptsize${}\pm0.10$} & $6.06${\scriptsize${}\pm0.10$} & $2.64${\scriptsize${}\pm0.12$} & $1.42\times10^{4}${\scriptsize${}\pm1143$} & $1.67\times10^{4}${\scriptsize${}\pm1214$} & $36.3${\scriptsize${}\pm14$}\\
\addlinespace[2pt]
manifold & $2$ & NN & $1.66${\scriptsize${}\pm0.15$} & $1.83${\scriptsize${}\pm0.12$} & $2.59${\scriptsize${}\pm0.08$} & $0.0738${\scriptsize${}\pm0.0048$} & $0.0844${\scriptsize${}\pm0.0056$} & $0.0870${\scriptsize${}\pm0.0061$}\\
 &  & LP & $3.30${\scriptsize${}\pm0.06$} & $3.41${\scriptsize${}\pm0.07$} & $2.47${\scriptsize${}\pm0.11$} & $0.548${\scriptsize${}\pm0.0072$} & $0.581${\scriptsize${}\pm0.017$} & $0.385${\scriptsize${}\pm0.018$}\\
 &  & partition & $3.92${\scriptsize${}\pm0.04$} & $3.93${\scriptsize${}\pm0.04$} & $2.47${\scriptsize${}\pm0.09$} & $0.807${\scriptsize${}\pm0.015$} & $0.806${\scriptsize${}\pm0.013$} & $0.554${\scriptsize${}\pm0.018$}\\
\addlinespace[2pt]
 & $5$ & NN & $4.90${\scriptsize${}\pm0.45$} & $3.86${\scriptsize${}\pm0.30$} & $2.61${\scriptsize${}\pm0.12$} & $9.04${\scriptsize${}\pm3.6$} & $111${\scriptsize${}\pm38$} & $3.61${\scriptsize${}\pm1.2$}\\
 &  & LP & $5.49${\scriptsize${}\pm0.08$} & $5.44${\scriptsize${}\pm0.09$} & $2.52${\scriptsize${}\pm0.09$} & $1.18\times10^{4}${\scriptsize${}\pm882$} & $1.32\times10^{4}${\scriptsize${}\pm1189$} & $8.80${\scriptsize${}\pm2.3$}\\
 &  & partition & $4.12${\scriptsize${}\pm0.07$} & $4.21${\scriptsize${}\pm0.09$} & $2.53${\scriptsize${}\pm0.10$} & $3256${\scriptsize${}\pm110$} & $3572${\scriptsize${}\pm103$} & $9.61${\scriptsize${}\pm1.8$}\\
\bottomrule
\end{tabular}}
\end{table}

\FloatBarrier

\subsection{Bike-sharing data}\label{app:bike-mv}
We use the hourly records of the Capital Bikeshare system in Washington, D.C., for 2011--2012 \citep{fanaee2014event,fanaee2013bike}, which cover $17{,}379$ hours. The response is $Y=(\log(1+\mathrm{cnt}_{\rm casual}+\mathrm{jit}_1),\log(1+\mathrm{cnt}_{\rm registered}+\mathrm{jit}_2))$, where $\mathrm{cnt}_{\rm casual}$ and $\mathrm{cnt}_{\rm registered}$ are the numbers of rentals by casual and by registered users, and the independent $U[0,1]$ variables $\mathrm{jit}_1,\mathrm{jit}_2$ remove the discreteness of the counts. The covariates are the hour of the day, the normalized temperature and the humidity, so that $d_X=3$ and $d_Y=2$. We omit the working-day indicator, so that at a given hour the conditional distribution of $Y$ mixes commuting and leisure patterns and is often skewed or bimodal. The methods and tuning rules are those of the simulations, except for the bandwidth grids given in the paragraph \emph{Tuning}.

Each of ten random splits uses $11000$ hours for fitting, $1500$ for calibration and the remaining $4879$ for testing, and methods with an auxiliary step fit on $8000$ of the fitting hours and use the other $3000$ for that step. These random splits evaluate prediction for held-out hours from the same two-year period, rather than forecasting future hours. We estimate conditional coverage on fifteen groups of test points, formed by crossing five blocks of hours with the terciles of temperature. Sampling error alone gives a mean absolute deviation of the group coverages from $0.9$ of about $0.014$, that is, $0.14$ on the $\times10$ scale of \cref{tab:bike}. We call this value the sampling floor.

\begin{figure}[!t]
\centering
\includegraphics[width=0.95\textwidth]{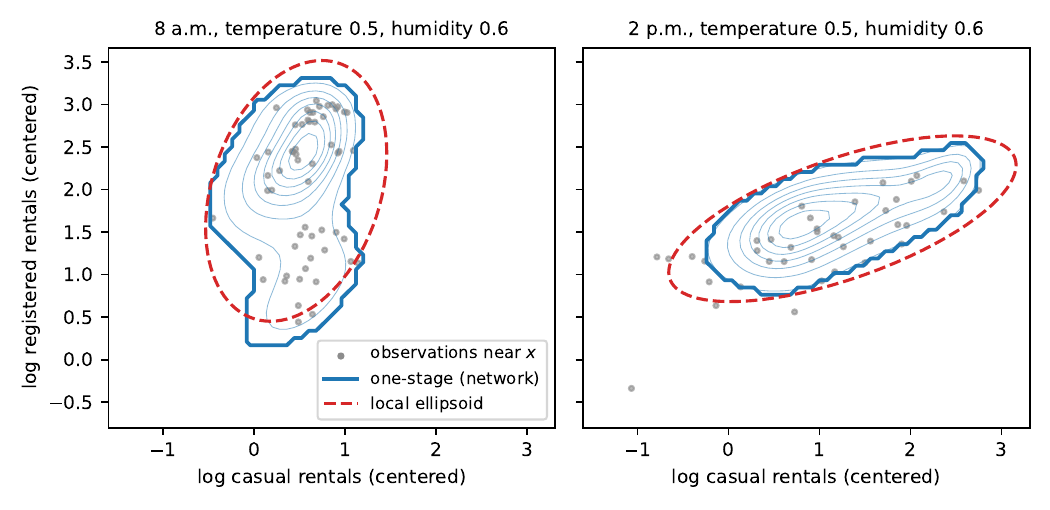}
\caption{Bike-sharing demand at 8 a.m.\ and 2 p.m., with normalized temperature $0.5$ and humidity $0.6$ (first split). Solid: one-stage calibration with the network density, as in \cref{tab:bike}; thin lines: contours of the estimated density; dashed: local ellipsoid; gray: all hours of the data set at the same hour of the day with temperature within $0.08$ and humidity within $0.12$ of these values. Both sets are calibrated to the nominal marginal coverage $0.9$.}\label{fig:bike}
\end{figure}

\paragraph{Preprocessing} The covariates are the hour divided by $23$ and the normalized temperature and humidity supplied with the data. In each split, fresh uniforms $\mathrm{jit}_1,\mathrm{jit}_2$ on $(0,1)$ are drawn for all $17{,}379$ hours, and each coordinate of $Y$ is mapped affinely so that its $0.1\%$ and $99.9\%$ quantiles over all hours are sent to $-2.25$ and $2.25$. Values outside $[-2.5,2.5]$, at most $0.1\%$ in each tail, are clipped to the boundary, and areas are converted back to log-count units. The axes of \cref{fig:bike} show the log counts minus the midpoint of these two quantiles. The map is common to all methods and is a symmetric function of all hours. Hence it does not affect exchangeability under the random split. The five blocks of hours are 0--5, 6--9, 10--15, 16--19 and 20--23, and the terciles of temperature are also computed over all hours and are used only for evaluation. A random permutation then assigns the hours to the parts described above. \Cref{fig:bike} uses the first split, with the network and the local ellipsoid fitted on the $11000$ fitting hours, as in the one-stage rows of \cref{tab:bike}, and both calibrated on the calibration hours. The sets are drawn on an $80\times80$ grid.

\paragraph{Tuning} On the grids of the simulations, the selected bandwidths of the LP and \cref{alg:local} lie at the lower ends. We therefore choose the LP bandwidths from $h^{(x)}\in\{0.08,0.12,0.18,0.27,0.4,0.6,0.9,1.35,2.0\}$ and $h^{(y)}\in\{0.14,0.2,0.3,0.45,0.65,1.0\}$, and the bandwidth of \cref{alg:local} from $\{0.1,0.15,0.22,0.3,0.4,0.55,0.75,1.0,1.5\}$. The other grids are those of the simulations, and the contrastive reference is fitted to the responses of the $8000$ fitting hours of the methods with an auxiliary step, among which lie the $1500$ tuning hours. The HPD values are computed on a $30\times30$ grid of $\cY$ instead of by importance sampling.

\begin{table}[!t]
\caption{Bike-sharing demand: mean $\pm$ standard error over 10 random splits; the sampling floor of the mean
group-coverage deviation is about $0.14$ on the $\times10$ scale. Columns marked ($\times10$) report ten times the value. Areas are in units of log-counts squared. The bandwidths of the LP and of \cref{alg:local} are chosen from grids extended below those of the simulations.}\label{tab:bike}
\centering\small
\resizebox{\ifdim\width>\textwidth\textwidth\else\width\fi}{!}{%
\begin{tabular}{llccccc}
\toprule
\multicolumn{2}{l}{method} & marginal coverage & \multicolumn{2}{c}{$|$group coverage$-0.9|$ ($\times10$)} & area & area relative\\
 & & & mean & max & & to NN, one-stage\\
\midrule
NN & one-stage & $0.893${\scriptsize${}\pm0.002$} & $0.211${\scriptsize${}\pm0.019$} & $0.631${\scriptsize${}\pm0.069$} & $4.280${\scriptsize${}\pm0.038$} & $1.000${\scriptsize${}\pm0.000$}\\
 & \cref{alg:local} & $0.896${\scriptsize${}\pm0.003$} & $0.339${\scriptsize${}\pm0.012$} & $0.823${\scriptsize${}\pm0.027$} & $5.210${\scriptsize${}\pm0.093$} & $1.218${\scriptsize${}\pm0.024$}\\
 & two-stage (own covariate) & $0.896${\scriptsize${}\pm0.002$} & $0.202${\scriptsize${}\pm0.013$} & $0.452${\scriptsize${}\pm0.038$} & $4.389${\scriptsize${}\pm0.055$} & $1.026${\scriptsize${}\pm0.014$}\\
 & two-stage (quantile) & $0.897${\scriptsize${}\pm0.002$} & $0.203${\scriptsize${}\pm0.014$} & $0.485${\scriptsize${}\pm0.042$} & $4.395${\scriptsize${}\pm0.054$} & $1.027${\scriptsize${}\pm0.013$}\\
\addlinespace[2pt]
LP & one-stage & $0.895${\scriptsize${}\pm0.003$} & $0.363${\scriptsize${}\pm0.023$} & $0.938${\scriptsize${}\pm0.105$} & $5.374${\scriptsize${}\pm0.100$} & $1.256${\scriptsize${}\pm0.027$}\\
 & \cref{alg:local} & $0.897${\scriptsize${}\pm0.003$} & $0.247${\scriptsize${}\pm0.014$} & $0.564${\scriptsize${}\pm0.028$} & $5.631${\scriptsize${}\pm0.086$} & $1.316${\scriptsize${}\pm0.022$}\\
 & two-stage (own covariate) & $0.897${\scriptsize${}\pm0.003$} & $0.216${\scriptsize${}\pm0.011$} & $0.554${\scriptsize${}\pm0.038$} & $5.613${\scriptsize${}\pm0.087$} & $1.312${\scriptsize${}\pm0.023$}\\
 & two-stage (quantile) & $0.899${\scriptsize${}\pm0.003$} & $0.186${\scriptsize${}\pm0.015$} & $0.447${\scriptsize${}\pm0.042$} & $5.886${\scriptsize${}\pm0.184$} & $1.375${\scriptsize${}\pm0.042$}\\
\addlinespace[2pt]
partition & one-stage & $0.895${\scriptsize${}\pm0.003$} & $0.304${\scriptsize${}\pm0.020$} & $0.684${\scriptsize${}\pm0.053$} & $6.473${\scriptsize${}\pm0.114$} & $1.513${\scriptsize${}\pm0.029$}\\
 & \cref{alg:local} & $0.896${\scriptsize${}\pm0.003$} & $0.249${\scriptsize${}\pm0.011$} & $0.640${\scriptsize${}\pm0.047$} & $7.171${\scriptsize${}\pm0.085$} & $1.677${\scriptsize${}\pm0.025$}\\
 & two-stage (own covariate) & $0.896${\scriptsize${}\pm0.003$} & $0.249${\scriptsize${}\pm0.015$} & $0.642${\scriptsize${}\pm0.061$} & $6.809${\scriptsize${}\pm0.100$} & $1.592${\scriptsize${}\pm0.030$}\\
 & two-stage (quantile) & $0.895${\scriptsize${}\pm0.004$} & $0.232${\scriptsize${}\pm0.010$} & $0.544${\scriptsize${}\pm0.050$} & $7.201${\scriptsize${}\pm0.366$} & $1.684${\scriptsize${}\pm0.087$}\\
\addlinespace[2pt]
baselines & CHCDS & $0.893${\scriptsize${}\pm0.002$} & $0.216${\scriptsize${}\pm0.020$} & $0.619${\scriptsize${}\pm0.069$} & $4.242${\scriptsize${}\pm0.045$} & $0.991${\scriptsize${}\pm0.005$}\\
 & local ellipsoid & $0.896${\scriptsize${}\pm0.002$} & $0.191${\scriptsize${}\pm0.012$} & $0.472${\scriptsize${}\pm0.036$} & $4.703${\scriptsize${}\pm0.031$} & $1.099${\scriptsize${}\pm0.010$}\\
 & RCP (ellipsoid score) & $0.894${\scriptsize${}\pm0.003$} & $0.191${\scriptsize${}\pm0.013$} & $0.482${\scriptsize${}\pm0.047$} & $4.799${\scriptsize${}\pm0.036$} & $1.122${\scriptsize${}\pm0.013$}\\
 & PCP & $0.899${\scriptsize${}\pm0.002$} & $0.346${\scriptsize${}\pm0.014$} & $0.881${\scriptsize${}\pm0.037$} & $5.514${\scriptsize${}\pm0.061$} & $1.289${\scriptsize${}\pm0.015$}\\
 & RCP (PCP score) & $0.897${\scriptsize${}\pm0.002$} & $0.198${\scriptsize${}\pm0.014$} & $0.540${\scriptsize${}\pm0.045$} & $5.621${\scriptsize${}\pm0.048$} & $1.314${\scriptsize${}\pm0.016$}\\
 & OT-CP & $0.897${\scriptsize${}\pm0.002$} & $0.504${\scriptsize${}\pm0.020$} & $1.799${\scriptsize${}\pm0.104$} & $5.958${\scriptsize${}\pm0.137$} & $1.392${\scriptsize${}\pm0.027$}\\
\bottomrule
\end{tabular}}
\end{table}

\paragraph{Results} Marginal coverage is close to $0.9$ for every method (\cref{tab:bike}), and the sets built from the network density by one-stage calibration and by CHCDS are the smallest. \Cref{fig:bike} illustrates why. At 8 a.m.\ the estimated conditional distribution has a sharp mode of commuting demand, with many registered users, and a long tail with few registered users, as on non-working days. The estimated HPD set follows this skewed shape, whereas the symmetric ellipsoid adds area beside the narrow tail, where there are few observations. At 2 p.m.\ the estimated high-density region is elongated, and the ellipse has nearly the same orientation but is larger.

At the group level, two-stage calibration reduces the coverage deviations of the LP and of the partition estimate and increases those of the network. Its regression pools hours within the bandwidth, which is $0.22$ in most splits, or about five hours of the day on either side. Over this range, the conditional distribution changes with the commuting pattern. Two-stage (own covariate) and two-stage (quantile), which regress HPD values, have smaller deviations than \cref{alg:local} for the network and the LP; in the simulations they do so where the estimate is accurate.

\FloatBarrier

\begingroup
\renewcommand{\doi}[1]{\href{https://doi.org/#1}{doi:\nolinkurl{#1}}}
\putbib
\endgroup
\end{bibunit}
\end{appendices}

\end{document}